\documentclass[12pt,leqno]{article}

\usepackage{amsfonts,amssymb,amsthm,amsmath}
\usepackage[margin=1in]{geometry}
\usepackage{hyperref}
\usepackage{setspace}
\usepackage[all]{xy}
\usepackage{longtable}
\usepackage{url}

\numberwithin{paragraph}{subsection}
\numberwithin{equation}{subparagraph}
\renewcommand{\theparagraph}{(\arabic{section}.\arabic{subsection}.\arabic{paragraph})}
\renewcommand{\thesubparagraph}{\arabic{section}.\arabic{subsection}.\arabic{paragraph}}

\allowdisplaybreaks

\begin{document}

\title{A Global Trace Formula for Reductive Lie Algebras\\ and the Harish-Chandra Transform\\ on the Space of Characteristic Polynomials}
\date{}
\author{Shuyang Cheng}
\maketitle

\begin{abstract}

In this paper an integral transform between spaces of nonstandard test functions on the affine space of dimension $n$ is constructed. The integral transform satisfies a summation formula of Poisson type, which is derived from an analogue of the Arthur-Selberg trace formula for the Lie algebra of $n\times n$ matrices.

\end{abstract}

\section*{Introduction}
\addcontentsline{toc}{section}{Introduction}

The classical summation formula of Poisson
\begin{eqnarray*}
\sum_{n\in\mathbb{Z}} f(n)
&=&
\sum_{n\in\mathbb{Z}} \hat{f}(n)
\end{eqnarray*}
has been related to the modularity of theta functions and the functional equation of zeta functions since the days of Jacobi and Riemann.

In his seminal paper \cite{1859} Riemann derived the functional equation of the zeta function which now bears his name from the modular identity of Jacobi's theta function, which is a consequence of the Poisson summation formula. 

In his thesis \cite{Tate} Tate recast Riemann's argument over the adeles and derived the functional equation of a large class of $L$-functions directly from his adelic Poisson summation formula. Tate's method was later generalized to the standard $L$-functions of higher rank general linear groups and central simple algebras by Godement and Jacquet in \cite{GoJa}. In their argument the group $\mathrm{GL}(n)$ is embedded into its Lie algebra $\mathrm{gl}(n)$, and the functional equation of standard $L$-functions of $\mathrm{GL}(n)$ is derived from the adelic Poisson summation formula for the vector space $\mathrm{gl}(n)$.

Recently in the work of Braverman and Kazhdan \cite{BrKa}, Ng\^{o} \cite{Ngo}, Sakellaridis \cite{Sake} and Lafforgue \cite{Laff}, the argument of Godement and Jacquet was partially generalized to more general automorphic $L$-functions:

Let $\mathrm{G}$ be a reductive group and let $\rho$ be a representation of the dual group $^L\mathrm{G}$ satisfying some mild assumptions, then the group $\mathrm{G}$ embeds into an algebraic monoid $\mathrm{M}_\rho$ which depends on $\rho$, and the functional equation of the automorphic $L$-functions $L(s,\pi,\rho)$ would follow from a conjectural summation formula of Poisson type for the monoid $\mathrm{M}_\rho$. In general $\mathrm{M}_\rho$ is nonlinear and possibly singular, and the conjectural Poisson summation formula would involve a nonstandard Fourier transform and spaces of nonstandard test functions.

\newpage~\\
The purpose of this paper is to construct a toy example of such a Poisson summation formula for the space $\mathcal{A}_n$ of characteristic polynomials of $n\times n$ matrices:

\paragraph*{Proposition \ref{poisson summation formula simplified statement}}
\emph{There exists a dense subspace $\mathcal{S}_0(\mathcal{A}_n(\mathbb{A}))$ of $L^2(\mathcal{A}_n(\mathbb{A}))$ of nonstandard Schwartz functions on $\mathcal{A}_n(\mathbb{A})$ and an endomorphism $\mathcal{H}_0$ of $\mathcal{S}_0(\mathcal{A}_n(\mathbb{A}))$ such that:}

\emph{For each point $X$ in $\mathcal{A}_n(\mathbb{Q})$ there exists a nonzero constant $a(X)$ such that for all $\varphi$ in $\mathcal{S}_0(\mathcal{A}_n(\mathbb{A}))$ which is cuspidal at two distinct place of $\mathbb{Q}$
\begin{eqnarray*}
\sum_{X\in\mathcal{A}_n(\mathbb{Q})} a(X)\cdot\varphi(X)
&=&
\sum_{X\in\mathcal{A}_n(\mathbb{Q})} a(X)\cdot\mathcal{H}_0(\varphi)(X),
\end{eqnarray*}
where the summand for a polynomial $X$ with repeated roots needs to be interpreted appropriately.}\\

The Poisson summation formula for $\mathcal{A}_n$ is deduced from an analogue of the Arthur-Selberg trace formula for the Lie algebra $\mathrm{gl}(n)$. In this paper such a trace formula is established for a restricted class of test functions subject to some mild local conditions on general reductive Lie algebra $\mathfrak{g}$. Some possible problems for future investigation include removing the restrictions on the test function $\varphi$, and establishing the Poisson summation formula for the affine space $\mathcal{A}_\mathrm{G}$ defined as the adjoint quotient of $\mathfrak{g}$ for a general reductive group $\mathrm{G}$. Such a Poisson summation formula will likely require a stable trace formula for $\mathfrak{g}$.

\section*{Acknowledgements}
\addcontentsline{toc}{section}{Acknowledgements}

I wish to express my gratitude  to my advisor Ng\^{o} Bao Ch\^{a}u, without whose continual support and encouragement this work would not have been possible. I wish to thank Pierre-Henri Chaudouard for spotting an error in an earlier version of this paper and providing very valuable suggestions. I also wish to thank Bill Casselman, Eddie Herman, Robert Kottwitz and Yiannis Sakellaridis for illuminating conversations. The proof of Corollary \ref{local harish chandra transform bijectivity} is suggested by Sakellaridis.

\newpage

\tableofcontents

\newpage

\section{Preliminaries}

\subsection{Notations and definitions}

\paragraph{}\setcounter{equation}{0}
Let $\mathrm{G}$\label{sym1} be a connected reductive group defined over the field of rational numbers $\mathbb{Q}$, let $\mathfrak{g}$ be its Lie algebra equipped with the adjoint action of $\mathrm{G}$ from the right. More generally an algebraic group will be denoted by a capital roman letter, its Lie algebra by the corresponding lowercase fraktur letter, except for $\mathfrak{a}$ which is reserved for a Euclidean vector space. 

Denote by $\mathbb{Q}_v$\label{sym3} the completion of $\mathbb{Q}$ at a place $v$\label{sym2}. If $S$ is a finite set of places, denote by $\mathbb{Q}_S$ the direct product of $\mathbb{Q}_v$ for all $v$ in $S$. Denote by $\mathbb{A}$ the ring of adeles of $\mathbb{Q}$. Define $S$-local and global norms\label{sym4} by
\begin{eqnarray}
\forall x\in\mathbb{Q}_S
\quad
|x|_S=\prod_{v\in S}|x|_v,
&&
\forall x\in\mathbb{A}
\quad
|x|_\mathbb{A}=\lim_S|x|_S.
\end{eqnarray}
The local norms are normalized in such a way that 
\begin{eqnarray}
\label{product formula}
\forall x\in\mathbb{Q}-\{0\}
&&
|x|_\mathbb{A}=1.
\label{product formula}
\end{eqnarray}

The groups $\mathrm{G}(\mathbb{Q}_v)$\label{sym5}, $\mathrm{G}(\mathbb{Q}_S)$ and $\mathrm{G}(\mathbb{A})$ of $\mathbb{Q}_v$, $\mathbb{Q}_S$ and $\mathbb{A}$-valued points of $\mathrm{G}$ are locally compact with respect to the analytic topology. The group $\mathrm{G}(\mathbb{Q})$ of $\mathbb{Q}$-valued points of $\mathrm{G}$ is a discrete subgroup of $\mathrm{G}(\mathbb{A})$ with respect to the analytic topology. 

\paragraph{}\setcounter{equation}{0}
Fix a minimal parabolic subgroup $\mathrm{P}_0$\label{sym6} of $\mathrm{G}$. Fix a Levi subgroup $\mathrm{M}_0$ of $\mathrm{P}_0$ with split component $\mathrm{A}_0$. A parabolic subgroup $\mathrm{P}$\label{sym7} of $\mathrm{G}$ is said to be \emph{standard} if $\mathrm{P}$ contains $\mathrm{P}_0$. Denote by $\mathrm{N}_\mathrm{P}$ the unipotent radical of $\mathrm{P}$, by $\mathrm{M}_\mathrm{P}$ the unique Levi subgroup of $\mathrm{P}$ containing $\mathrm{M}_0$, by $\mathrm{A}_\mathrm{P}$ the split component of $\mathrm{M}_\mathrm{P}$. Such a Levi subgroup $\mathrm{M}_\mathrm{P}$ is said to be \emph{standard}. Denote by $\overline{\mathrm{P}}$\label{sym9} the parabolic subgroup group opposite to $\mathrm{P}$, by $\overline{\mathrm{N}}_\mathrm{P}$ its unipotent radical. To simplify notations the standard Levi, split and unipotent components of $\mathrm{P}_i$\label{sym8} will be denoted by $\mathrm{M}_i$, $\mathrm{A}_i$ and $\mathrm{N}_i$ where $i$ is a natural number.

Let $\mathrm{M}$ and $\mathrm{L}$ be Levi subgroups of $\mathrm{G}$ such that $\mathrm{M}$ is contained in $\mathrm{L}$, denote by $\mathcal{F}^\mathrm{L}(\mathrm{M})$\label{sym10} the set of parabolic subgroups of $\mathrm{L}$ that contain $\mathrm{M}$, by $\mathcal{P}^\mathrm{L}(\mathrm{M})$ the set of parabolic subgroups of $\mathrm{L}$ whose Levi component is $\mathrm{M}$, by $\mathcal{L}^\mathrm{L}(\mathrm{M})$ the set of Levi subgroups of $\mathrm{L}$ that contain $\mathrm{M}$. To simplify notations denote by $\mathcal{F}(\mathrm{M})$ the set $\mathcal{F}^\mathrm{G}(\mathrm{M})$, by $\mathcal{F}^\mathrm{L}$ the set $\mathcal{F}^\mathrm{L}(\mathrm{M}_0)$, by $\mathcal{F}$ the set $\mathcal{F}^\mathrm{G}(\mathrm{M}_0)$. Similar notations apply to $\mathcal{P}$ and $\mathcal{L}$.

Let $\mathrm{P}$ be a parabolic subgroup of $\mathrm{G}$, denote by $X(\mathrm{M}_\mathrm{P})$\label{sym11} the group of rational characters of $\mathrm{M}_\mathrm{P}$
\begin{eqnarray}
X(\mathrm{M}_\mathrm{P}) 
&=&
\mathrm{Hom}_{\mathrm{Grp}/\mathbb{Q}}(\mathrm{M}_\mathrm{P},\mathrm{GL}(1,\mathbb{Q})).
\end{eqnarray}
Let $\mathfrak{a}_\mathrm{P}$\label{sym12} denote the real vector space 
\begin{eqnarray}
\mathfrak{a}_\mathrm{P}
&=&
\mathrm{Hom}_\mathbb{Z}(X(\mathrm{M}_\mathrm{P}),\mathbb{R}),
\end{eqnarray}
let $\mathfrak{a}_\mathrm{P}^*$ denote the dual space
\begin{eqnarray}
\mathfrak{a}_\mathrm{P}^*
&=&
X(\mathrm{M}_\mathrm{P})\otimes_\mathbb{Z}\mathbb{R}.
\end{eqnarray}
Let\label{sym13}
\begin{eqnarray}
\Phi_\mathrm{P},~\Delta_\mathrm{P}\subset\mathfrak{a}_\mathrm{P}^*,
&&
\Phi_\mathrm{P}^\vee,~\Delta_\mathrm{P}^\vee\subset\mathfrak{a}_\mathrm{P}
\end{eqnarray}
denote respectively the set of \emph{roots}, \emph{simple roots}, \emph{coroots}, \emph{simple coroots} of $\mathrm{A}_\mathrm{P}$ in $\mathfrak{g}$ and $\mathfrak{n}_\mathrm{P}$. The quadruple $(X(\mathrm{M}_0),\Phi_0,X(\mathrm{M}_0)^*,\Phi_0^\vee)$ is called the \emph{root datum} of $\mathrm{G}$.

Let $\mathrm{P}_1$ and $\mathrm{P}_2$ be parabolic subgroups of $\mathrm{G}$ with $\mathrm{P}_1$ contained in $\mathrm{P}_2$, denote
\begin{eqnarray}
\mathrm{N}_1^2=\mathrm{N}_1\cap\mathrm{M}_2, 
&&
\overline{\mathrm{N}}_1^2=\overline{\mathrm{N}}_1\cap\mathrm{M}_2.
\end{eqnarray}
Let $\Delta_1^2$\label{sym14} and $\Delta_1^{2,\vee}$ be the set of simple roots and coroots of $\mathrm{A}_1$ in $\mathfrak{n}_1^2$. There are canonical splittings
\begin{eqnarray}
\mathfrak{a}_1=\mathfrak{a}_1^2\oplus\mathfrak{a}_2,
&& 
\mathfrak{a}_1^*=\mathfrak{a}_1^{2,*}\oplus\mathfrak{a}_2^*.
\end{eqnarray}
The sets $\Delta_1^2$ and $\Delta_1^{2,\vee}$ form bases of $\mathfrak{a}_1^{2,*}$\label{sym15} and $\mathfrak{a}_1^2$. The respective dual bases are called the \emph{coweights} and \emph{weights} and denoted by $\hat{\Delta}_1^2$ and $\hat{\Delta}_1^{2,\vee}$.

Let $\mathrm{W}_0^\mathrm{G}$\label{sym16} be the \emph{Weyl group} of the pair $(\mathrm{G},\mathrm{A}_0)$. Let $\mathrm{M}_1$ and $\mathrm{M}_2$ be two Levi subgroups of $\mathrm{G}$, define the \emph{Weyl set} $\mathrm{W}(\mathfrak{a}_1,\mathfrak{a}_2)$ to be the set of linear isomorphisms from $\mathfrak{a}_1$ to $\mathfrak{a}_2$ obtained by restricting the action of elements of the Weyl group. The group $\mathrm{W}_0^\mathrm{G}$ operates on the root datum of $\mathrm{G}$, hence on $\mathfrak{a}_0$ and $\mathfrak{a}_0^*$. Fix Euclidean inner products on $\mathfrak{a}_0$ and $\mathfrak{a}_0^*$ which are compatible with each other and the underlying root datum, hence invariant under the Weyl group action. The inner product on $\mathfrak{a}_0$ induces an inner product on $\mathfrak{a}_1^2$.

Let $\tau_1^2$\label{sym17} denote the characteristic function on $\mathfrak{a}_0$ of the points that are positive with respect to every element of $\Delta_1^2$, let $\hat{\tau}_1^2$ denote the characteristic function on $\mathfrak{a}_0$ of the set of points that are positive with respect to every element of $\hat{\Delta}_1^{2,\vee}$:
\begin{eqnarray}
\tau_1^2=\mathbb{I}_{\big\{H\in\mathfrak{a}_0:~\alpha(H)>0,~\forall\alpha\in\Delta_1^2\big\}},
&&
\hat{\tau}_1^2=\mathbb{I}_{\big\{H\in\mathfrak{a}_0 :~\varpi(H)>0,~\forall\varpi\in\hat{\Delta}_1^{2,\vee}\big\}}.
\end{eqnarray}
To simplify notations denote by $\tau_1$ the set $\tau_1^\mathrm{G}$, by  $\hat{\tau}_1$ the set $\hat{\tau}_1^\mathrm{G}$.

\paragraph{}\setcounter{equation}{0}
Define $\mathrm{G}(\mathbb{A})^1$\label{sym18} to be the subgroup of $\mathrm{G}(\mathbb{A})$ consisting of  elements $g$ of $\mathrm{G}(\mathbb{Q})$ such that 
\begin{eqnarray}
\forall\chi\in X(\mathrm{G})
&&
|\chi(g)|_\mathbb{A}=1. 
\end{eqnarray}
Fix an admissible maximal compact subgroup\label{sym19}
\begin{eqnarray}
\mathrm{K} &=& \prod_v \mathrm{K}_v 
\end{eqnarray}
of $\mathrm{G}(\mathbb{A})$ such that the Iwasawa decomposition 
\begin{eqnarray}
\mathrm{G}(\mathbb{A})
&=&
\mathrm{P}(\mathbb{A})\mathrm{K}
\\\nonumber&=&
\mathrm{M}_\mathrm{P}(\mathbb{A})^1\exp(\mathfrak{a}_\mathrm{P})\mathrm{N}_\mathrm{P}(\mathbb{A})\mathrm{K}
\end{eqnarray} 
holds. Denote by\label{sym20}
\begin{eqnarray*}
H_\mathrm{P}~:~\mathrm{G}(\mathbb{A})\quad\longrightarrow\quad\mathfrak{a}_ \mathrm{P}
\end{eqnarray*} 
the natural projection. 

The \emph{Tamagawa measure} on $\mathrm{G}(\mathbb{A})$ is the measure induced from the choice of a basis of rational 1-forms on $\mathrm{G}$, which is well-defined by the product formula \eqref{product formula}. The Euclidean vector space $\mathfrak{a}_\mathrm{G}$ has a translation invariant measure, which without loss of generality assigns the coweight lattice 
\begin{eqnarray}
\mathrm{Hom}_\mathbb{Z}(X(\mathrm{G}),\mathbb{Z})
&\subset&
\mathfrak{a}_\mathrm{G}
\end{eqnarray}
covolume one.

The various measures are compatible under the Iwasawa decomposition in the sense that
\begin{eqnarray}
&&
\int_{\mathrm{G}(\mathbb{A})} f(g) ~\mathrm{d}g
\\\nonumber\\\nonumber
&=&
\int_{\mathrm{M}_\mathrm{P}(\mathbb{A})^1} \int_{\mathfrak{a}_\mathrm{P}} \int_{\mathrm{N}_\mathrm{P}(\mathbb{A})} \int_\mathrm{K} f(mank) ~\mathrm{d}k\mathrm{d}n\mathrm{d}a\mathrm{d}m
\\\nonumber\\\nonumber
&=&
\int_{\mathrm{M}_\mathrm{P}(\mathbb{A})^1} \int_{\mathfrak{a}_\mathrm{P}} \int_{\mathrm{N}_\mathrm{P}(\mathbb{A})} \int_\mathrm{K} f(nmak) e^{-2\rho_\mathrm{P}(H_0(a))} ~\mathrm{d}k\mathrm{d}n\mathrm{d}a\mathrm{d}m
\end{eqnarray}
where $\mathrm{d}n$ is the Tamagawa measure on $\mathrm{N}_\mathrm{P}(\mathbb{A})$ which could also be characterized by assigning $\mathrm{N}_\mathrm{P}(\mathbb{Q})$ covolume one in $\mathrm{N}_\mathrm{P}(\mathbb{A})$, the point $\rho_\mathrm{P}$\label{sym21} in $\mathfrak{a}_0^*$ is the \emph{Weyl vector} defined as the half sum of the roots of $\mathfrak{n}_\mathrm{P}$. The choices of measures on $\mathrm{G}(\mathbb{A})$ and $\mathfrak{a}_\mathrm{G}$ determine a Haar measure on $\mathrm{G}(\mathbb{A})^1$, hence a measure on the automorphic quotient $\mathbb{G}(\mathbb{Q})\backslash\mathrm{G}(\mathbb{A})^1$.

Let $T'$\label{sym22} be a point in $\mathfrak{a}_0$, let $\omega$\label{sym23} be a compact subset of $\mathrm{N}_0(\mathbb{A})\mathrm{M}_0(\mathbb{A})^1$. The \emph{Siegel set} $\mathfrak{S}(T',\omega)$ is the subset of $\mathrm{G}(\mathbb{A})^1$ defined as
\begin{eqnarray}
&&
\bigg\{ x=pak:~ p\in\omega,~a\in\exp(\mathfrak{a}_0),~k\in\mathrm{K},~
\beta(H_0(a)-T')\geq0 ~\forall\beta\in\Delta_0 \bigg\}. 
\end{eqnarray}
A Siegel set is said to be a \emph{Siegel domain} if it contains a fundamental domain for $\mathrm{G}(\mathbb{Q})\backslash\mathrm{G}(\mathbb{A})^1$:
\begin{eqnarray}
\mathrm{G}(\mathbb{A})^1 
&=& 
\mathrm{G}(\mathbb{Q})\mathfrak{S}(T',\omega).
\end{eqnarray}
Let $T$ be a sufficiently positive point in $\mathfrak{a}_0$, such a $T$ is said to be a \emph{truncation parameter}. Let $\mathfrak{S}(T',\omega)$\label{sym24} be a Siegel domain. The \emph{truncated Siegel domain} $\mathfrak{S}^T(T',\omega)$ is defined as
\begin{eqnarray}
\label{truncated siegel domain}
&&
\bigg\{ x\in\mathfrak{S}(T',\omega):~\varpi(H_0(x)-T)\leq0 ~\forall\varpi\in\hat{\Delta}_0 \bigg\}.
\end{eqnarray}

\paragraph*{Proposition}
(Borel, Harish-Chandra)\\
\emph{There exist a point $T'$ and a compact set $\omega$ such that $\mathfrak{S}(T',\omega)$ is a Siegel domain.}

\proof
See \S13 of \cite{Borel}.
\qed\\

The Siegel domain $\mathfrak{S}(T',\omega)$, therefore $\mathrm{G}(\mathbb{Q})\backslash\mathrm{G}(\mathbb{A})^1$, has finite volume. The truncated Siegel domain $\mathfrak{S}^T(T',\omega)$ is compact and exhausts $\mathfrak{S}(T',\omega)$ as $T$ approaches infinity. Fix a pair ($T',\omega$) such that $\mathfrak{S}(T',\omega)$ is a Siegel domain, denote by $F^\mathrm{G}(x,T)$\label{sym25} the characterstic function on $\mathrm{G}(\mathbb{Q})\backslash \mathrm{G}(\mathbb{A})^1$ of $\mathfrak{S}^T(T',\omega)$. There are analogues $F^\mathrm{P}(x,T)$ on $\mathrm{P}(\mathbb{Q})\backslash\mathrm{G}(\mathbb{A})^1$.

\paragraph{}\setcounter{equation}{0}
The space $\mathcal{S}(\mathfrak{g}(\mathbb{A}))$\label{sym26} of \emph{Schwartz functions} on $\mathfrak{g}(\mathbb{A})$ is defined as the tensor product 
\begin{eqnarray}
\bigotimes_p^\mathrm{res}C_c^\infty(\mathfrak{g}(\mathbb{Q}_p))\otimes\mathcal{S}(\mathfrak{g}(\mathbb{R}))
\end{eqnarray}
restricted at all but finitely many finite primes with respect to the unit vector $\mathbb{I}_{\mathfrak{g}(\mathbb{Z}_p)}$\label{sym27} in $C_c^\infty (\mathfrak{g}(\mathbb{Q}_p))$, equipped with the final topology with respect to 
\begin{eqnarray}
\label{global lie algebra schwartz space topology}
\mathcal{S}(\mathfrak{g}(\mathbb{A}))
&=&
\varinjlim_S\mathcal{S}(\mathfrak{g}(\mathbb{Q}_S)).
\end{eqnarray}

Fix a nondegenerate G-invariant rational bilinear form $\langle~,~\rangle$\label{sym28} on $\mathfrak{g}$, fix a global additive unitary character $\psi$ on $\mathbb{A}$
\begin{eqnarray}
\psi:\mathbb{A}\rightarrow\mathrm{U}(1)
\end{eqnarray}
such that 
\begin{eqnarray}
\forall x\in\mathbb{Q}
&&
\psi(x)=1,
\end{eqnarray} 
and $\psi(\langle~,~\rangle)$ identifies $\mathbb{Q}$ and $\mathbb{A}/\mathbb{Q}$ as Pontryagin duals of each other. The \emph{Fourier transform} on $\mathcal{S}(\mathfrak{g}(\mathbb{A}))$ is defined by\label{sym29}
\begin{eqnarray}
f\hat{~}(X) 
&=&
\int_{\mathfrak{g}(\mathbb{A})} f(Y)\psi(\langle X,Y \rangle) ~\mathrm{d}Y.
\end{eqnarray}
The global Fourier transform on $\mathfrak{g}(\mathbb{A})$ factorizes as the tensor product of the local Fourier transforms on $\mathfrak{g}(\mathbb{Q}_v)$ with respect to compatible choices of $\psi_v(\langle~,~\rangle)$. Denote by $\vee$ the inverse Fourier transform.

\paragraph*{Proposition}(Poisson summation formula, Tate)\\
\emph{For every Schwartz function $f$ on $\mathfrak{g}(\mathbb{A})$,
\begin{eqnarray}
\label{poisson summation formula, tate}
\sum_{X\in\mathfrak{g}(\mathbb{Q})}f(X) 
&=& 
\sum_{X\in\mathfrak{g}(\mathbb{Q})}f\hat{~}(X),
\end{eqnarray}
the sums are the absolutely convergent.}

\proof
See \S4.2 of \cite{Tate}.
\qed

\paragraph{}\setcounter{equation}{0}
Define an equivalence relation $\sim$\label{sym30} on $\mathfrak{g}(\mathbb{Q})$ by
\begin{eqnarray}
X\sim Y 
&\textrm{if}& 
\exists g\in\mathrm{G}(\mathbb{Q})
~Y_\mathrm{ss}=X_\mathrm{ss}\cdot\mathrm{ad}(g).
\end{eqnarray}\label{sym31}
In general $\sim$ is weaker than conjugacy by $\mathrm{G}(\mathbb{Q})$. A typical equivalence class will be denoted by $\mathfrak{o}$.

Let $D$\label{sym32} be the \emph{discriminant function} on $\mathfrak{g}(\mathbb{Q})$. Let $X$ be an element of $\mathfrak{g}(\mathbb{Q})$, define $D(X)$ to be the coefficient of the characteristic polynomial of $\mathrm{ad}(X)$, acting on $\mathfrak{g}(\mathbb{Q})$ as a linear endomorphism, in degree $r$, the absolute rank of $\mathrm{G}$:
\begin{eqnarray}
r
&=&
\mathrm{rank}_{\overline{\mathbb{Q}}}(\mathrm{G}\otimes_\mathbb{Q}\overline{\mathbb{Q}}).
\end{eqnarray}
Alternatively $D$ could be defined as the product of the roots of $\mathfrak{g}$ over an algebraically closed field. Denote by $D^\mathrm{M}$ the discriminant function on $\mathfrak{m}$ for a Levi subgroup $\mathrm{M}$ of $\mathrm{G}$.

Let $X$ be a semisimple element of $\mathfrak{g}(\mathbb{Q})$. Then $X$ is said to be \emph{regular} if $D(X)$ does not vanish. Denote by $\mathfrak{g}_\mathrm{reg.ss}$\label{sym33} the locus of regular semisimple points on $\mathfrak{g}$. If an equivalence class $\mathfrak{o}$ contains a regular semisimple element, the set $\mathfrak{o}$ is a $\mathrm{G}(\mathbb{Q})$-orbit consisting of regular semisimple elements. Such an $\mathfrak{o}$ is said to be \emph{regular}, otherwise $\mathfrak{o}$ is said to be \emph{singular}. A semisimple element $X$ is said to be \emph{$\mathbb{Q}$-elliptic} if it is stablized under the adjoint action by a maximal torus that is anisotropic over $\mathbb{Q}$ modulo the center of $\mathrm{G}$.

Denote by $\mathrm{G}_X$\label{sym34} the centralizer of $X$ in $\mathrm{G}$, let $\mathrm{G}_X^0$ be the connected component of the identity of $\mathrm{G}_\Sigma$, let $\pi_0(\mathrm{G}_X)$ be the group of connected components of $\mathrm{G}_X$. There is an exact sequence
\begin{eqnarray}
\xymatrix{1\ar[r]&\mathrm{G}_X^0\ar[r]
&\mathrm{G}_X\ar[r]&\pi_0(\mathrm{G}_X)\ar[r]&1}.
\end{eqnarray}

\paragraph{}\setcounter{equation}{0}
Let $S$ be a finite set of places of $\mathbb{Q}$, denote by $\mathrm{G}_S$\label{sym35} the base change $\mathrm{G}\otimes_\mathbb{Q}\mathbb{Q}_S$ of $\mathrm{G}$. Let $v$ be a place of $\mathbb{Q}$, denote by $\mathrm{G}_v$ the base change $\mathrm{G}\otimes_\mathbb{Q}\mathbb{Q}_v$ of $\mathrm{G}$. Similar notations apply to the Lie algebra $\mathfrak{g}$. The underlying topological groups of $\mathrm{G}(\mathbb{Q}_S)$ and $\mathrm{G}_S(\mathbb{Q}_S)$ are the same, however
\begin{eqnarray}
X(\mathrm{G})
&=&
\mathrm{Hom}_{\mathrm{Grp}/\mathbb{Q}}(\mathrm{G},\mathrm{GL}(1,\mathbb{Q})),
\\\nonumber
X(\mathrm{G}_S)
&=&
\mathrm{Hom}_{\mathrm{Grp}/\mathbb{Q}_S}(\mathrm{G}_S,\mathrm{GL}(1,\mathbb{Q}_S))
\end{eqnarray} 
are in general different.

All the constructions above generalize to the local and $S$-local settings with similar caveats:

\begin{itemize}
\item Fix a minimal parabolic subgroup $\mathrm{P}_{S,0}$\label{sym36} of $\mathrm{G}_S$ contained in $\mathrm{P}_{0,S}$, fix a Levi subgroup $\mathrm{M}_{S,0}$ contained in $\mathrm{M}_{0,S}$ with split component $\mathrm{A}_{S,0}$ containing $\mathrm{A}_{0,S}$. The choice of $\mathrm{P}_{S,0}$ is equivalent to a choice of a minimal parabolic $\mathrm{P}_{v,0}$ of $\mathrm{G}_v$ for each $v$ in $S$, similarly for $\mathrm{M}_{S,0}$ and $\mathrm{A}_{S,0}$.

\item Let $\mathrm{M}_S$ and $\mathrm{L}_S$ be Levi subgroups of $\mathrm{G}_S$ such that $\mathrm{M}_S$ is contained in $\mathrm{L}_S$, denote by $\mathcal{F}^{\mathrm{L}_S}(\mathrm{M}_S)$, $\mathcal{P}^{\mathrm{L}_S}(\mathrm{M}_S)$ and $\mathcal{L}^{\mathrm{L}_S}(\mathrm{M}_S)$ the analogous sets of parabolic and Levi subgroups. The Levi subgroups $\mathrm{M}_S$ and $\mathrm{L}_S$ determine local Levi subgroups $\mathrm{M}_v$ and $\mathrm{L}_v$ for each $v$ in $S$. There are bijections
\begin{eqnarray}
\mathcal{F}^{\mathrm{L}_S}(\mathrm{M}_S)
&=&
\prod_{v\in S}\mathcal{F}^{\mathrm{L}_v}(\mathrm{M}_v),
\\\nonumber
\mathcal{P}^{\mathrm{L}_S}(\mathrm{M}_S)
&=&
\prod_{v\in S}\mathcal{P}^{\mathrm{L}_v}(\mathrm{M}_v),
\\\nonumber
\mathcal{L}^{\mathrm{L}_S}(\mathrm{M}_S)
&=&
\prod_{v\in S}\mathcal{L}^{\mathrm{L}_v}(\mathrm{M}_v).
\end{eqnarray}
If $\mathrm{M}_S$ and $\mathrm{L}_S$ are the base change of Levi subgroups $\mathrm{M}$ and $\mathrm{L}$ of $\mathrm{G}$ from $\mathbb{Q}$ to $\mathbb{Q}_S$, there are diagonal inclusions 
\begin{eqnarray}
\mathcal{F}^\mathrm{L}(\mathrm{M})
&\subset&
\mathcal{F}^{\mathrm{L}_S}(\mathrm{M}_S),
\\\nonumber
\mathcal{P}^\mathrm{L}(\mathrm{M})
&\subset&
\mathcal{P}^{\mathrm{L}_S}(\mathrm{M}_S),
\\\nonumber
\mathcal{L}^\mathrm{L}(\mathrm{M})
&\subset&
\mathcal{L}^{\mathrm{L}_S}(\mathrm{M}_S).
\end{eqnarray}

\item Let $\mathrm{P}_S$ be a parabolic subgroup of $\mathrm{G}_S$, define real vector spaces
\begin{eqnarray}
\mathfrak{a}_{\mathrm{P}_S}=
\mathrm{Hom}_\mathbb{Z}(X(\mathrm{M}_{\mathrm{P}_S}),\mathbb{R}),
&&
\mathfrak{a}_{\mathrm{P}_S}^*=
X(\mathrm{M}_{\mathrm{P}_S})\otimes_\mathbb{Z}\mathbb{R}.
\end{eqnarray}
If $\mathrm{P}_S$ is the base change of a parabolic subgroup $\mathrm{P}$ of $\mathrm{G}$ from $\mathbb{Q}$ to $\mathbb{Q}_S$, there are diagonal inclusions
\begin{eqnarray}
\mathfrak{a}_\mathrm{P}\subset\mathfrak{a}_{\mathrm{P}_S}=
\bigoplus_{v\in S}\mathfrak{a}_{\mathrm{P}_v},
&&
\mathfrak{a}_\mathrm{P}^*\subset\mathfrak{a}_{\mathrm{P}_S}^*=
\bigoplus_{v\in S}\mathfrak{a}_{\mathrm{P}_v}^*.
\end{eqnarray}

\item A maximal torus $\mathrm{T}_S$\label{sym37} of $\mathrm{G}_S$ is equivalent to the choice of a maximal torus $\mathrm{T}_v$ of $\mathrm{G}_v$ for each $v$ in $S$. The associated Cartan subalgebra $\mathfrak{t}_S$, which is a free $\mathbb{Q}_S$-module, is equal as an abelian group to the direct sum of the $\mathbb{Q}_v$-vector spaces $\mathfrak{t}_v$ for all $v$ in $S$.

A maximal torus $\mathrm{T}_S$ is said to be \emph{elliptic} in $\mathrm{G}_S$ if it is anisotropic modulo $\mathrm{A}_{\mathrm{G}_S}$ over $\mathbb{Q}_S$. A maximal torus $\mathrm{T}_S$ is elliptic in $\mathrm{G}_S$ if and only if $\mathrm{T}_v$ is elliptic in $\mathrm{G}_v$ for each $v$ in $S$. If this is the case, $\mathrm{T}_S(\mathbb{Q}_S)$ is compact modulo $\mathrm{A}_{\mathrm{G}_S}(\mathbb{Q}_S)$ in the analytic topology.

Denote by $\mathcal{T}_\mathrm{ell}(\mathrm{G}_S)$\label{sym38} the set of conjugacy classes of elliptic maximal tori of $\mathrm{G}_S$.

\item Let $\mathrm{W}_{S,0}^{\mathrm{G}_S}$ be the Weyl group of the pair $(\mathrm{G}_S,\mathrm{A}_{S,0})$, there is a bijection 
\begin{eqnarray}
\mathrm{W}_{S,0}^{\mathrm{G}_S} 
&=&
\prod_{v\in S} \mathrm{W}_{v,0}^{\mathrm{G}_v},
\end{eqnarray}
the linear representation of $\mathrm{W}_{S,0}^{\mathrm{G}_S}$ on $\mathfrak{a}_{S,0}$ is the direct sum of the local representations.

Denote by $\mathrm{W}(\mathrm{G}_S,\mathrm{T}_S)$\label{sym39} the Weyl group of the pair $(\mathrm{G}_S(\mathbb{Q}_S),\mathrm{T}_S(\mathbb{Q}_S))$. There is a bijection
\begin{eqnarray}
\mathrm{W}(\mathrm{G}_S,\mathrm{T}_S) 
&=&
\prod_{v\in S} \mathrm{W}(\mathrm{G}_v,\mathrm{T}_v).
\end{eqnarray}

\item The Schwartz space $\mathcal{S}(\mathfrak{g}_S(\mathbb{Q}_S))$, the Fourier transform $\wedge$, the discriminant function $D^{\mathrm{G}_S}$, and the regular semisimple locus $\mathfrak{g}_{S,\mathrm{reg.ss}}(\mathbb{Q}_S)$ are unchanged under base change from $\mathbb{Q}$ to $\mathbb{Q}_S$.
\end{itemize}

\paragraph*{Proposition}(Weyl integration formula)\\
\emph{If $f_S$ is a Schwartz function on $\mathfrak{g}_S(\mathbb{Q}_S)$, then}
\begin{eqnarray}
\label{weyl integration formula}
&&
\int_{\mathfrak{g}_S(\mathbb{Q}_S)} f_S(X)~\mathrm{d}X
\\\nonumber\\\nonumber
&=&
\sum_{\mathrm{M}_S\in\mathcal{L}^{\mathrm{G}_S}} 
|\mathrm{W}_{S,0}^{\mathrm{M}_S}| 
|\mathrm{W}_{S,0}^{\mathrm{G}_S}|^{-1} 
\sum_{\mathrm{T}_S\in\mathcal{T}_\mathrm{ell}(\mathrm{M}_S)} 
|\mathrm{W}(\mathrm{M}_S,\mathrm{T}_S)|^{-1}\times 
\\\nonumber\\\nonumber
&& \quad \times \int_{\mathfrak{t}_S(\mathbb{Q}_S)} 
|D^{\mathrm{G}_S}(X)|_S
\int_{\mathrm{A}_{\mathrm{M}_S}(\mathbb{Q}_S)\backslash\mathrm{G}_S (\mathbb{Q}_S)} f_S(X\cdot\mathrm{ad}(x)) ~\mathrm{d}x\mathrm{d}X.
\end{eqnarray}

\proof
For the $p$-adic case see \S7.11 of \cite{Kott}. For the real case see Lemma 2 on page 35 of \cite{Var} and the references therein.
\qed

\paragraph*{Definition}
Let $X$ be an element of $\mathfrak{g}_{S,\mathrm{reg.ss}}(\mathbb{Q}_S)$. Define the $S$-local \emph{orbital integral} $I_\mathrm{G}^\mathrm{G}(X,~)$\label{sym40} to be the invariant distribution on $\mathfrak{g}_S(\mathbb{Q}_S)$ such that
\begin{eqnarray}
\label{standard orbital integral}
&&
\forall f_S\in\mathcal{S}(\mathfrak{g}_S(\mathbb{Q}_S))
\\\nonumber
&&
I_\mathrm{G}^\mathrm{G}(X,f_S)
=
|D^\mathrm{G}(X)|_S^{1/2}
\int_{\mathrm{G}_{S,X}^0(\mathbb{Q}_S)\backslash\mathrm{G}_S(\mathbb{Q}_S)} f_S(X\cdot\mathrm{ad}(x)) ~\mathrm{d}x
\end{eqnarray}
where $\mathrm{G}_{S,X}^0$ denotes the connected component of the identity of the stablizer subgroup of $X$ in $\mathrm{G}_S$.

\section{The non-invariant trace formula}

In this chapter a preliminary version of the trace formula for the reductive Lie algebra $\mathfrak{g}$ established in Chaudouard \cite{Chaud} is recalled.

\subsection{A motivating example}

\paragraph{Definition}\setcounter{equation}{0}
Let $f$ be a Schwartz function on $\mathfrak{g}(\mathbb{A})$, let $\mathfrak{o}$ be a $\sim$ equivalence class on $\mathfrak{g}(\mathbb{Q})$. Define the \emph{kernel functions} $K(x,f)$\label{sym41} and $K_\mathfrak{o}(x,f)$ by
\begin{eqnarray}
\quad
\forall x\in \mathrm{G}(\mathbb{Q})\backslash\mathrm{G}(\mathbb{A})
&&
K(x,f) = \sum_{X\in\mathfrak{g}(\mathbb{Q})}f(X\cdot\mathrm{ad}(x)),
\\\nonumber\\\nonumber
&&
K_\mathfrak{o}(x,f) = \sum_{X\in\mathfrak{o}}f(X\cdot\mathrm{ad}(x)).
\end{eqnarray}

\paragraph{Remark}\setcounter{equation}{0}
By the Poisson summation formula \eqref{poisson summation formula, tate}, the function $K(x,f)$ satisfies the functional equation 
\begin{eqnarray}
\label{kernel function functional equation}
K(x,f)
&=&
K(x,f\hat{~}).
\end{eqnarray}

\paragraph{Proposition}\setcounter{equation}{0}
\emph{Let $f$ be a Schwartz function on $\mathfrak{g}(\mathbb{A})$. If $\mathrm{G}$ is anisotropic over $\mathbb{Q}$, then
\begin{eqnarray}
&&
\lim_S
\sum_{\substack{\mathfrak{o}\in\mathfrak{g}(\mathbb{Q})/\sim\\\textrm{regular}}} 
\mathrm{Vol}(\mathrm{G}_{X_\mathfrak{o}}^0(\mathbb{Q})
\backslash\mathrm{G}_{X_\mathfrak{o}}^0(\mathbb{A})) 
\cdot I_\mathrm{G}^\mathrm{G}(X_\mathfrak{o},f_S)
\\\nonumber\\\nonumber
&&
\quad
+\sum_{\substack{\mathfrak{o}\in\mathfrak{g}(\mathbb{Q})/\sim\\\textrm{singular}}}
\mathrm{Vol}(\mathrm{G}_{X_\mathfrak{o}}^0(\mathbb{Q})
\backslash\mathrm{G}_{X_\mathfrak{o}}^0(\mathbb{A}))
\cdot I_\mathfrak{o}(f)
\\\nonumber\\\nonumber
&=&
\lim_S
\sum_{\substack{\mathfrak{o}\in\mathfrak{g}(\mathbb{Q})/\sim\\\textrm{regular}}}
\mathrm{Vol}(\mathrm{G}_{X_\mathfrak{o}}^0(\mathbb{Q})
\backslash \mathrm{G}_{X_\mathfrak{o}} ^0(\mathbb{A})) 
\cdot I_\mathrm{G}^\mathrm{G}(X_\mathfrak{o},f\hat{~}_S)
\\\nonumber\\\nonumber
&&
\quad
+\sum_{\substack{\mathfrak{o}\in\mathfrak{g}(\mathbb{Q})/\sim\\\textrm{singular}}}
\mathrm{Vol}(\mathrm{G}_{X_\mathfrak{o}}^0(\mathbb{Q})
\backslash\mathrm{G}_{X_\mathfrak{o}}^0(\mathbb{A}))
\cdot I_\mathfrak{o}(f\hat{~}),
\end{eqnarray}
where for each class $\mathfrak{o}$ choose an element $X_\mathfrak{o}$ of $\mathfrak{o}$, for each singular class $\mathfrak{o}$ define}\label{sym42}
\begin{eqnarray}
I_\mathfrak{o}(f)
&=&
\int_{\mathrm{G}_{X_\mathfrak{o}}^0(\mathbb{A})\backslash\mathrm{G}(\mathbb{A})} 
f(X_\mathfrak{o}\cdot\mathrm{ad}(x))~\mathrm{d}x.
\end{eqnarray}

\proof
The function $K(x,f)$ is continuous. By assumption $\mathrm{G}$ is anisotropic, so $\mathrm{G}(\mathbb{Q})\backslash\mathrm{G}(\mathbb{A})$ is compact. Therefore $K(x,f)$ is absolutely integrable, hence
\begin{eqnarray}
&&
\int_{\mathrm{G}(\mathbb{Q})\backslash\mathrm{G}(\mathbb{A})}
K(x,f)~\mathrm{d}x
\\\nonumber\\\nonumber
&=&
\int_{\mathrm{G}(\mathbb{Q})\backslash\mathrm{G}(\mathbb{A})} 
\sum_{\mathfrak{o}\in\mathfrak{g}(\mathbb{Q})/\sim} 
K_\mathfrak{o}(x,f)~\mathrm{d}x
\\\nonumber\\\nonumber 
&=& 
\sum_{\mathfrak{o}\in\mathfrak{g}(\mathbb{Q})/\sim} 
\int_{\mathrm{G}(\mathbb{Q})\backslash\mathrm{G}(\mathbb{A})}  
K_\mathfrak{o}(x,f)~\mathrm{d}x
\\\nonumber\\\nonumber 
&=&
\sum_{\mathfrak{o}\in\mathfrak{g}(\mathbb{Q})/\sim} 
\mathrm{Vol}\big(\mathrm{G}_{X_\mathfrak{o}}^0(\mathbb{Q})
\backslash \mathrm{G}_{X_\mathfrak{o}} ^0(\mathbb{A})\big) 
\int_{\mathrm{G}_{X_\mathfrak{o}}^0(\mathbb{A})\backslash\mathrm{G}(\mathbb{A})}  f(X_\mathfrak{o}\cdot\mathrm{ad}(x))~\mathrm{d}x
\\\nonumber\\\label{discriminant equals one}
&=&
\lim_S
\sum_{\substack{\mathfrak{o}\in\mathfrak{g}(\mathbb{Q})/\sim\\\mathrm{regular}}} 
\mathrm{Vol}\big(\mathrm{G}_{X_\mathfrak{o}}^0(\mathbb{Q})
\backslash \mathrm{G}_{X_\mathfrak{o}} ^0(\mathbb{A})\big) 
\cdot I_\mathrm{G}^\mathrm{G}(X_\mathfrak{o},f_S)
\\\nonumber\\\nonumber
&&
\quad
+\sum_{\substack{\mathfrak{o}\in\mathfrak{g}(\mathbb{Q})/\sim\\\mathrm{singular}}}
\mathrm{Vol}\big(\mathrm{G}_{X_\mathfrak{o}}^0(\mathbb{Q})
\backslash\mathrm{G}_{X_\mathfrak{o}}^0(\mathbb{A})\big)
\cdot I_\mathfrak{o}(f\hat{~}),
\end{eqnarray}
the equality \eqref{discriminant equals one} follows from 
\begin{eqnarray}
\forall X_\mathfrak{o}\in\mathfrak{g}(\mathbb{Q})
~\exists S
~\forall S'\supset S
&&
|D(X_\mathfrak{o})|_{S'}=1
\end{eqnarray}
by the product formula \eqref{product formula}. The proposition follows from the functional equation \eqref{kernel function functional equation}.
\qed

\subsection{The non-invariant trace formula of Chaudouard}

\paragraph{Definition}\setcounter{equation}{0}
\label{truncated kernel definition}
Let $f$ be a Schwartz function on $\mathfrak{g}(\mathbb{A})$, let $\mathfrak{o}$ be a $\sim$ equivalence class on $\mathfrak{g}(\mathbb{Q})$, let $T$ be a truncation parameter. Define the \emph{truncated kernel function} $k_\mathfrak{o}^T(~,f)$\label{sym43} on $\mathrm{G}(\mathbb{Q})\backslash\mathrm{G}(\mathbb{A})^1$ by
\begin{eqnarray}
&&
\forall x\in\mathrm{G}(\mathbb{Q})\backslash\mathrm{G}(\mathbb{A})^1
\\\nonumber
&&
k_\mathfrak{o}^T(x,f)
=
\sum_{\substack{\mathrm{P}\in\mathcal{F}\\\mathrm{standard}}}
(-1)^{\dim(\mathrm{A}_\mathrm{P}/\mathrm{A}_\mathrm{G})}
\sum_{\delta\in\mathrm{P}(\mathbb{Q})\backslash\mathrm{G}(\mathbb{Q})}
\hat{\tau}_\mathrm{P}(H_0(\delta x)-T) K_{\mathrm{P},\mathfrak{o}}(\delta x)
\end{eqnarray}
where
\begin{eqnarray}
K_{\mathrm{P},\mathfrak{o}}(x,f) &=& \sum_{X\in\mathfrak{m}_\mathrm{P}(\mathbb{Q})\cap\mathfrak{o}}
\int_{\mathfrak{n}_\mathrm{P}(\mathbb{A})}f((X+N)\cdot\mathrm{ad}(x))~\mathrm{d}N.
\end{eqnarray}
Define distributions $J_\mathfrak{o}^T$\label{sym44} and $J^T$ on $\mathfrak{g}(\mathbb{A})$ by
\begin{eqnarray}
\forall f\in\mathcal{S}(\mathfrak{g}(\mathbb{A}))
&&
J_\mathfrak{o}^T(f)=\int_{\mathrm{G}(\mathbb{Q})\backslash\mathrm{G}(\mathbb{A})^1} k_\mathfrak{o}^T(x,f)~\mathrm{d}x,
\\\nonumber\\\nonumber
&&
J^T(f)=\int_{\mathrm{G}(\mathbb{Q})\backslash\mathrm{G}(\mathbb{A})^1}
\sum_{\mathfrak{o}\in\mathfrak{g}(\mathbb{Q})/\sim}
k_\mathfrak{o}^T(x,f)~\mathrm{d}x.
\end{eqnarray}

\paragraph{Proposition}\setcounter{equation}{0}
\emph{Let $f$ be a Schwartz function on $\mathfrak{g}(\mathbb{A})$, let $T$ be a truncation parameter, then}
\begin{eqnarray}
\label{truncated kernel integrability}
\int_{\mathrm{G}(\mathbb{Q})\backslash\mathrm{G}(\mathbb{A})^1}
\sum_{\mathfrak{o}\in\mathfrak{g}(\mathbb{Q})/\sim}
\big|k_\mathfrak{o}^T(x,f)\big| ~\mathrm{d}x &<& \infty.
\end{eqnarray}

\proof
This is Th\'{e}or\`{e}me 3.1 of \cite{Chaud}. The following is a sketch of the argument, for the details please refer to loc. cit.

The proof depends on the following lemma.
\paragraph*{Lemma} (Combinatorial lemma of Langlands)\\
\emph{Let $\mathrm{P}_1$ and $\mathrm{P}_3$ be parabolic subgroups of $\mathrm{G}$ such that $\mathrm{P}_1$ is contained in $\mathrm{P}_3$, then}
\begin{eqnarray}
\label{langlands combinatorial lemma}
\sum_{\substack{\mathrm{P}_2\in\mathcal{F}\\\mathrm{P}_1\subset\mathrm{P}_2\subset\mathrm{P}_3}} 
(-1)^{\dim(\mathrm{A}_2 / \mathrm{A}_3)} 
\tau_1^2(H) \hat{\tau}_2^3(H) 
&=& 
\left\{ \begin{array}{ll} 1 & \textrm{\emph{if} } \mathrm{P}_1=\mathrm{P}_3, 
\\\\ 0 & \textrm{\emph{otherwise}}. \end{array} \right.  
\end{eqnarray}

\proof See \S6 of \cite{Arthur_CoarseI}.\qed\\

By \eqref{langlands combinatorial lemma}
\begin{eqnarray}
\label{truncated kernel integrability first upperbound}
&&
\int_{\mathrm{G}(\mathbb{Q})\backslash\mathrm{G}(\mathbb{A})^1}
\sum_{\mathfrak{o}\in\mathfrak{g}(\mathbb{Q})/\sim}
\big|k_\mathfrak{o}^T(x,f)\big| ~\mathrm{d}x
\\\nonumber\\\nonumber
&\leq&
\sum_{\substack{\mathrm{P}_1,\mathrm{P}_4\in\mathcal{F}\\
\mathrm{standard}\\\mathrm{P}_1\subset\mathrm{P}_4}} 
\int_{\mathrm{P}_1(\mathbb{Q})\backslash\mathrm{G}(\mathbb{A})^1}
\sum_{\mathfrak{o}\in\mathfrak{g}(\mathbb{Q})/\sim}
F^{\mathrm{P}_1}(x,T)
\sigma_1^4(H_0(x)-T) \times
\\\nonumber\\\nonumber
&&\qquad
\times\Bigg|
\sum_{\substack{\mathrm{P}_3\in\mathcal{F}\\\mathrm{P}_1\subset\mathrm{P}_3\subset\mathrm{P}_4}} 
(-1)^{\dim(\mathrm{A}_3 / \mathrm{A}_\mathrm{G})} 
K_{\mathrm{P}_3,\mathfrak{o}}(x,f) \Bigg| ~\mathrm{d}x
\end{eqnarray}
where
\begin{eqnarray}
\label{sym45}
\sigma_1^4(H) 
&=& 
\sum_{\substack{\mathrm{P}_5\in\mathcal{F}\\\mathrm{P}_4\subset\mathrm{P}_5}}
(-1)^{\dim(\mathrm{A}_4 / \mathrm{A}_5)} 
\tau_1^5(H) \hat{\tau}_5(H).
\end{eqnarray}

The second factor of the integrand of the right hand side of \eqref{truncated kernel integrability first upperbound} satisfies the inequality
\begin{eqnarray}
\label{truncated kernel integrability second inequality}
&&
\Bigg| 
\sum_{\substack{\mathrm{P}_3\in\mathcal{F}\\\mathrm{P}_1\subset\mathrm{P}_3\subset\mathrm{P}_4}}
(-1)^{\dim(\mathrm{A}_3 / \mathrm{A}_\mathrm{G})} 
K_{\mathrm{P}_3,\mathfrak{o}}(x,f) \Bigg|
\\\nonumber\\\nonumber
&\leq& 
\sum_{\substack{\mathrm{P}_2\in\mathcal{F}\\\mathrm{P}_1\subset\mathrm{P}_2\subset\mathrm{P}_4}} 
\Bigg| 
\sum_{\substack{\mathrm{P}_3\in\mathcal{F}\\\mathrm{P}_2\subset\mathrm{P}_3\subset\mathrm{P}_4}} (-1)^{\dim(\mathrm{A}_3 /\mathrm{A}_\mathrm{G})} 
\sum_{X\in\mathfrak{m}_1^2 (\mathbb{Q})'\cap\mathfrak{o}}
\\\nonumber\\\nonumber 
&& 
\quad 
\times \sum_{Y\in\mathfrak{n}_2^3(\mathbb{Q})}  
\int_{\mathfrak{n}_3(\mathbb{A})} 
f((X+Y+N)\cdot\mathrm{ad}(x)) ~\mathrm{d}N \Bigg|
\\\nonumber\\ 
\label{truncated kernel integrability second equality}
&=& 
\Bigg| \sum_{X\in\mathfrak{m}_1^4(\mathbb{Q})'\cap\mathfrak{o}}  
\int_{\mathfrak{n}_4(\mathbb{A})} 
f((X+N)\cdot\mathrm{ad}(x)) ~\mathrm{d}N \Bigg|+
\\\nonumber\\\nonumber
&& 
\quad 
+\sum_{\substack{\mathrm{P}_2\in\mathcal{F}\\\mathrm{P}_1\subset\mathrm{P}_2\subsetneq\mathrm{P}_4}} 
\Bigg| 
\sum_{\substack{\mathrm{P}_3\in\mathcal{F}\\\mathrm{P}_2\subset\mathrm{P}_3\subset\mathrm{P}_4}}
(-1)^{\dim(\mathrm{A}_3 /\mathrm{A}_\mathrm{G})}
\sum_{X\in\mathfrak{m}_1^2(\mathbb{Q})'\cap\mathfrak{o}}
\\\nonumber\\\nonumber
&& 
\quad\quad 
\times \sum_{\overline{Y}\in\overline{\mathfrak{n}}_2^3(\mathbb{Q})'} 
\int_{\mathfrak{n}_2(\mathbb{A})} 
f((X+N)\cdot\mathrm{ad}(x))
\psi(\langle N,\overline{Y}\rangle) ~\mathrm{d}N \Bigg|
\end{eqnarray}
where $\mathfrak{m}_\mathrm{P}^\mathrm{Q}(\mathbb{Q})'$\label{sym46} denotes the set of points of $\mathfrak{m}_\mathrm{P}^\mathrm{Q}(\mathbb{Q})$ not contained in any proper parabolic subalgebra of $\mathfrak{q}$. Similar notation applies to $\mathfrak{n}$. The equality \eqref{truncated kernel integrability second equality} follows from the Poisson summation formula applied to $\mathfrak{n}_2^3(\mathbb{Q})$ as a lattice in $\mathfrak{n}_2^3(\mathbb{A})$.

By the inclusion-exclusion principle applied to
\begin{eqnarray}
\sum_{\substack{\mathrm{P}_3\in\mathcal{F}\\\mathrm{P}_2\subset\mathrm{P}_3\subset\mathrm{P}_4}}
(-1)^{\mathrm{dim}(\mathrm{A}_3 /\mathrm{A}_\mathrm{G})}
\end{eqnarray}
the last expression in \eqref{truncated kernel integrability second inequality} reduces to a majorant of the form
\begin{eqnarray}
\prod_{\alpha\in\Delta_1^4} \int_0^\infty (1+t_\alpha)^{p_\alpha}e^{-q_\alpha t_\alpha} ~\mathrm{d}t_\alpha
\end{eqnarray}
for some natural numbers $p_\alpha$ and $q_\alpha$, which is finite.
\qed

\paragraph{Definition}\setcounter{equation}{0}
Let $\mathrm{P}_2$ be a standard parabolic subgroup of $\mathrm{G}$, let $T$ be a truncation parameter in $\mathfrak{a}_0$. Define the \emph{geometric gamma$'$ function} $\Gamma'_{\mathrm{P}_2}(~,T)$\label{sym47} on $\mathfrak{a}_0$ by
\begin{eqnarray}
\label{gamma' function definition}
\forall H\in\mathfrak{a}_0
&&
\Gamma'_{\mathrm{P}_2}(H,T) 
=
\sum_{\substack{\mathrm{P}_3\in\mathcal{F}\\\mathrm{P}_2\subset\mathrm{P}_3}} 
(-1)^{\dim(\mathrm{A}_3 /\mathrm{A}_\mathrm{G})} 
\tau_2^3(H) \hat{\tau}_3(H-T).
\end{eqnarray}

\paragraph{Remark}\setcounter{equation}{0}
For each parabolic subgroup $\mathrm{P}_1$ of $\mathrm{G}$ the geometric gamma$'$ functions satisfy the identity
\begin{eqnarray}
\label{gamma' function identity}
\hat{\tau}_1(H-T) 
&=& 
\sum_{\substack{\mathrm{P}_2\in\mathcal{F}\\\mathrm{P}_1\subset\mathrm{P}_2}} 
(-1)^{\dim(\mathrm{A} _2 /\mathrm{A}_\mathrm{G})} 
\hat{\tau}_1^2(H) \Gamma'_{\mathrm{P}_2}(H,T).
\end{eqnarray}
For a proof see page 13 of \cite{Arthur_Inv}.

\paragraph{Proposition}\setcounter{equation}{0}
\label{truncation polynomial}
\emph{Let $f$ be a Schwartz function on $\mathfrak{g}(\mathbb{A})$, let $\mathfrak{o}$ be a $\sim$ equivalence class, let $T$ be a truncation parameter. Then $J_\mathfrak{o}^T(f)$ and $J^T(f)$ are polynomials in $T$ of degree at most $\dim(\mathrm{A}_0/\mathrm{A}_\mathrm{G})$.}

\proof
This is Th\'{e}or\`{e}me 4.2 of \cite{Chaud}. The following is a sketch of the argument for $J_\mathfrak{o}^T(f)$, the argument for $J^T(f)$ is the same. For the details please refer to loc. cit.

Fix a point $T'$ in $\mathfrak{a}_0$ that is sufficiently positive and assume that $T$ dominates $T'$. By \eqref{gamma' function identity}
\begin{eqnarray}
&&
J_\mathfrak{o}^T(f)
\\\nonumber\\\nonumber
&=& 
\int_{\mathrm{G}(\mathbb{Q})\backslash\mathrm{G}(\mathbb{A})^1} 
\Bigg( \sum_{\substack{\mathrm{P}_1\in\mathcal{F}\\\mathrm{standard}}}
(-1)^{\dim(\mathrm{A}_1/\mathrm{A}_\mathrm{G})} \times 
\\\nonumber\\\nonumber 
&&
\quad 
\times \sum_{\delta\in\mathrm{P}_1(\mathbb{Q})\backslash\mathrm{G}(\mathbb{Q})}
\hat{\tau}_1(H_0(\delta x)-T)\cdot 
K_{\mathrm{P}_1,\mathfrak{o}}(\delta x)\Bigg) ~\mathrm{d}x 
\\\nonumber\\\nonumber 
&=&  
\int_{\mathrm{G}(\mathbb{Q})\backslash \mathrm{G}(\mathbb{A})^1} 
\Bigg( 
\sum_{\substack{\mathrm{P}_1,\mathrm{P}_2\in\mathcal{F}\\
\mathrm{standard}\\\mathrm{P}_1\subset\mathrm{P}_2}}
(-1)^{\dim(\mathrm{A}_1/\mathrm{A}_2)} 
\times 
\\\nonumber\\\nonumber 
&& 
\quad 
\times \sum_{\delta\in\mathrm{P}_1(\mathbb{Q})\backslash\mathrm{G}(\mathbb{Q})}
\hat{\tau}_1^2(H_{\mathrm{P}_1}(\delta x)-T') 
\Gamma_{\mathrm{P}_2}'(H_{\mathrm{P}_2}(\delta x)-T',T-T') \times
\\\nonumber\\\nonumber
&&
\quad\quad
\times \sum_{X\in\mathfrak{m}_1(\mathbb{Q})\cap\mathfrak{o}}
\int_{\mathfrak{n}_1(\mathbb{A})}
f((X+N)\cdot\mathrm{ad}(\delta x)) ~\mathrm{d}N \Bigg) ~\mathrm{d}x 
\\\nonumber\\\nonumber 
&=&  
\sum_{\substack{\mathrm{P}_2\in\mathcal{F}\\\mathrm{standard}}} 
\int_{\mathrm{M}_2(\mathbb{Q})\backslash \mathrm{M}_2(\mathbb{A})^1} 
\Bigg( \sum_{\substack{\mathrm{P}_1\in\mathcal{F}\\
\mathrm{standard}\\\mathrm{P}_1\subset\mathrm{P}_2}}
(-1)^{\dim(\mathrm{A}_1/\mathrm{A}_2)} \times
\\\nonumber\\\nonumber
&&
\quad
\times \sum_{\delta\in\mathrm{P}_1(\mathbb{Q})\cap\mathrm{M}_2(\mathbb{Q})
\backslash \mathrm{M}_2(\mathbb{Q})}
\hat{\tau}_1^2(H_{\mathrm{P}_1}(\delta x)-T') 
\sum_{X\in\mathfrak{m}_1(\mathbb{Q})\cap\mathfrak{o}} 
\int_{\mathfrak{n}_1^2(\mathbb{A})} 
\\\nonumber\\\nonumber 
&& 
\quad\quad
\times \bigg( \int_\mathrm{K} \int_{\mathfrak{n}_2(\mathbb{A})} 
f((X+N)\cdot\mathrm{ad}(\delta xk)+N'\cdot \mathrm{ad}(k)) 
~\mathrm{d}N'\mathrm{d}k \bigg) ~\mathrm{d}N \Bigg) ~\mathrm{d}x \times 
\\\nonumber\\\nonumber 
&& 
\quad\quad\quad
\times \int_{\mathfrak{a}_2^\mathrm{G}} 
\Gamma'_{\mathrm{P}_2}(H-T',T-T') ~\mathrm{d}H 
\\\nonumber\\\nonumber 
&=& 
\sum_{\substack{\mathrm{P}_2\in\mathcal{F}\\\mathrm{standard}}} 
J_\mathfrak{o}^{\mathrm{M}_2,T'}(f_{\mathrm{P}_2}) 
\int_{\mathfrak{a}_2^\mathrm{G}} 
\Gamma'_{\mathrm{P}_2}(H,T-T') ~\mathrm{d}H
\end{eqnarray}
where  $J_\mathfrak{o}^\mathrm{M}$\label{sym48} is defined to be the sum of $J_{\mathfrak{o}'}^\mathrm{M}$ over all the $\mathrm{M}(\mathbb{Q})\sim$ equivalence classes $\mathfrak{o}'$ contained in $\mathfrak{o}$ and $f_\mathrm{P}$\label{sym49} is defined by
\begin{eqnarray}
f_\mathrm{P}(X) &=& \int_\mathrm{K} \int_{\mathfrak{n}_\mathrm{P}(\mathbb{A})} f((X+N) \cdot \mathrm{ad}(k)) ~\mathrm{d}N\mathrm{d}k. 
\end{eqnarray} 

Because $\int \Gamma_\mathrm{P}'(H,T-T') ~\mathrm{d}H$ is a polynomial in $T$ which is homogeneous of degree $\dim(\mathrm{A}_\mathrm{P}/\mathrm{A}_\mathrm{G})$, Proposition \ref{truncation polynomial} follows by induction.
\qed

\paragraph{Proposition}\setcounter{equation}{0}
\emph{Let $f$ be a Schwartz function on $\mathfrak{g}(\mathbb{A})$, let $T$ be a truncation parameter. For every positive $\epsilon$
\begin{eqnarray}
\label{asymptotic expansion}
\Bigg|J^T(f)-\int_{\mathrm{G}(\mathbb{Q})\backslash\mathrm{G}(\mathbb{A})^1}
F^\mathrm{G}(x,T) 
\sum_{X\in\mathfrak{g}(\mathbb{Q})}
f(X\cdot\mathrm{ad}(x))~\mathrm{d}x \Bigg|
&=& 
O(e^{-\epsilon \|T\|}) 
\end{eqnarray}
where $\|~\|$\label{sym50} denotes the Euclidean norm on $\mathfrak{a}_0$, as $T$ approaches infinity such that $T$ is uniformly bounded away from the walls of the positive chamber.}

\proof
This is a corollary of the proof of Lemma \ref{refined asymptotic expansion}.
\qed

\paragraph{Proposition}\setcounter{equation}{0} (Non-invariant trace formula of Chaudouard)\\
\emph{Let $f$ be a Schwartz function on $\mathfrak{g}(\mathbb{A})$, then
\begin{eqnarray}
\label{chaudouard trace formula}
\sum_{\mathfrak{o}\in\mathfrak{g}(\mathbb{Q})/\sim}J_\mathfrak{o}^T(f) &=& \sum_{\mathfrak{o}\in\mathfrak{g}(\mathbb{Q})/\sim}J_\mathfrak{o}^T(f\hat{~})
\end{eqnarray}
holds as an equality between polynomials in $T$.}

\proof
By \eqref{kernel function functional equation}
\begin{eqnarray}
&&
\int_{\mathrm{G}(\mathbb{Q})\backslash\mathrm{G}(\mathbb{A})^1}
F^\mathrm{G}(x,T) 
\sum_{X\in\mathfrak{g}(\mathbb{Q})}
f(X\cdot\mathrm{ad}(x))~\mathrm{d}x
\\\nonumber\\\nonumber
&=&
\int_{\mathrm{G}(\mathbb{Q})\backslash\mathrm{G}(\mathbb{A})^1}
F^\mathrm{G}(x,T) 
K(x,f)~\mathrm{d}x
\\\nonumber\\\nonumber
&=&
\int_{\mathrm{G}(\mathbb{Q})\backslash\mathrm{G}(\mathbb{A})^1}
F^\mathrm{G}(x,T) 
K(x,f\hat{~})~\mathrm{d}x.
\end{eqnarray}
By \eqref{asymptotic expansion} the difference between $J^T(f)$ and $J^T(f\hat{~})$ converges to zero as $T$ approaches infinity. By Proposition \ref{truncation polynomial} $J^T(f)$ and $J^T(f\hat{~})$ are both polynomials in $T$, hence equal to each other.
\qed

\paragraph{Lemma}\setcounter{equation}{0}
\emph{There exists a unique point $T_0$\label{sym51} in $\mathfrak{a}_0^\mathrm{G}$ such that
\begin{eqnarray}
\label{point t0}
\forall s\in\mathrm{W}_0^\mathrm{G}
&&
H_0(w_s^{-1})=T_0-s^{-1}T_0
\end{eqnarray}
where $w_s$\label{sym52} denotes a representative of $s$ in $\mathrm{G}(\mathbb{Q})$.}

\proof
See Lemma 1.1 of \cite{Arthur_Inv}.
\qed

\paragraph{Definition}\setcounter{equation}{0}
Define distributions $J$\label{sym53} and $J_\mathfrak{o}$ on $\mathfrak{g}(\mathbb{A})$ by
\begin{eqnarray}
\forall f\in\mathcal{S}(\mathfrak{g}(\mathbb{A}))
&&
J(f)=J^{T_0}(f),~J_\mathfrak{o}(f)=J_\mathfrak{o}^{T_0}(f).
\end{eqnarray}

\paragraph{Remark}\setcounter{equation}{0}
The coefficients of the polynomials $J^T(f)$ and $J_\mathfrak{o}^T(f)$ in positive degrees depend only on the orbital integrals of $f_\mathrm{P}$ along proper Levi subalgebras of $\mathfrak{g}$. Therefore the constant terms $J(f)$ and $J_\mathfrak{o}(f)$ contain the essential information.
The choice of $T_0$ implies that the distributions $J$ and $J_\mathfrak{o}$ are independent of the choice of the minimal parabolic subgroup $\mathrm{P}_0$.

\paragraph{Proposition}\setcounter{equation}{0}
\emph{Let $f$ be a Schwartz function on $\mathfrak{g}(\mathbb{A})$. If $\mathfrak{o}$ is the conjugacy class of a regular semisimple element $X$ which is contained in a standard parabolic $\mathfrak{p}$ and elliptic in its Levi component $\mathfrak{m}$, then
\begin{eqnarray}
\label{elementary refinement}
&&
J_\mathfrak{o}(f) 
= 
|\pi_0(\mathrm{G}_X)|^{-1}
\mathrm{Vol}(\mathrm{M}_X^0(\mathbb{Q})\backslash\mathrm{M}_X^0(\mathbb{A})^1) \int_{\mathrm{G}_X^0(\mathbb{A})\backslash\mathrm{G}(\mathbb{A})} f(X\cdot\mathrm{ad}(x))~v_\mathrm{M}(x)\mathrm{d}x
\end{eqnarray}
where the weight factor $v_\mathrm{M}(x)$\label{sym54} is defined to be the volume of the convex hull of}
\begin{eqnarray}
\bigg\{ -H_\mathrm{P}(x):~ \mathrm{P}\in\mathcal{P}(\mathrm{M}) \bigg\} 
&\subset&
\mathfrak{a}_\mathrm{M}^\mathrm{G}.
\end{eqnarray}

\proof
This is a corollary of Theorem \ref{refined expansion}. See also (5.4) of \cite{Chaud}.
\qed

\paragraph{Proposition}\setcounter{equation}{0}
\emph{Let $f$ be a Schwartz function on $\mathfrak{g}(\mathbb{A})$, let $\mathfrak{o}$ be a $\sim$ equivalence class, then
\begin{eqnarray}
\label{elementary variance}
\forall x\in\mathrm{G}(\mathbb{A})^1
&&
J_\mathfrak{o}(f\circ\mathrm{ad}(x))=\sum_{\mathrm{P}\in\mathcal{F}} 
|\mathrm{W}_0^{\mathrm{M}_\mathrm{P}}| 
|\mathrm{W}_0^\mathrm{G}|^{-1}
 J_\mathfrak{o}^{\mathrm{M}_\mathrm{P}}(f_{\mathrm{P},x})
\end{eqnarray}
where $J_\mathfrak{o}^\mathrm{M}$ is defined to be the sum of $J_{\mathfrak{o}'}^\mathrm{M}$ over all the $\mathrm{M}(\mathbb{Q})\sim$ equivalence classes $\mathfrak{o}'$ contained in $\mathfrak{o}$, and the function $f_{\mathrm{P},x}$\label{sym55} on $\mathfrak{m}_\mathrm{P}(\mathbb{A})$ is defined as
\begin{eqnarray}
\label{elementary variance parabolic descent' definition}
f_{\mathrm{P},x}(X) &=& \int_\mathrm{K}\int_{\mathfrak{n}_\mathrm{P} (\mathbb{A})} f((X+N)\cdot\mathrm{ad}(k)) ~v_\mathrm{P}'(kx) \mathrm{d}N\mathrm{d}k
\end{eqnarray}
where the weight factor $v_\mathrm{P}'(x)$ is defined by}
\begin{eqnarray}
v_\mathrm{P}'(x) &=& \int_{\mathfrak{a}_\mathrm{P}^\mathrm{G}} \Gamma_\mathrm{P}'(H,-H_\mathrm{P}(x)) ~\mathrm{d}H.
\end{eqnarray}

\proof
By definition
\begin{eqnarray}
\label{elementary variance first equation}
J_\mathfrak{o}^T(f\circ\mathrm{ad}(x)) 
&=& 
\int_{\mathrm{G}(\mathbb{Q})\backslash \mathrm{G}(\mathbb{A})^1} 
\Bigg( \sum_{\substack{\mathrm{P}_1\in\mathcal{F}\\\mathrm{standard}}}
(-1)^{\dim(\mathrm{A}_1\mathrm{A}_\mathrm{G})} \times
\\\nonumber\\\nonumber 
&& 
\quad
\times \sum_{\delta\in\mathrm{P}_1(\mathbb{Q})\backslash\mathrm{G}(\mathbb{Q})}
\hat{\tau}_1(H_0(\delta yx)-T)
\cdot K_{\mathrm{P}_1,\mathfrak{o}}(\delta yx) ~\mathrm{d}y \Bigg).
\end{eqnarray}
By \eqref{gamma' function identity}
\begin{eqnarray}
\hat{\tau}_1(H_{\mathrm{P}_1}(\delta yx)-T) 
&=& 
\sum_{\substack{\mathrm{P}_2\in\mathcal{F}\\\mathrm{P}_1\subset\mathrm{P}_2}} 
(-1)^{\dim(\mathrm{A}_2/\mathrm{A}_\mathrm{G})} 
\hat{\tau}_1^2(H_{\mathrm{P}_1}(\delta y)-T)\times 
\\\nonumber && \quad \times 
\Gamma_{\mathrm{P}_2}'(H_{\mathrm{P}_1}(\delta y)-T,-H_{\mathrm{P}_1}(kx))
\end{eqnarray}
where $k$ is a $\mathrm{K}$ component of $\delta y$ under the decomposition of $\mathrm{G}(\mathbb{A})$ as $\mathrm{P}_1(\mathbb{A})\mathrm{K}$, the point $H_{\mathrm{P}_1}(kx)$ is independent of the choice of the element $k$. Hence the right hand side of \eqref{elementary variance first equation} is equal to
\begin{eqnarray}
\label{elementary variance second expression}
&&
\sum_{\substack{\mathrm{P}_2\in\mathcal{F}\\\mathrm{standard}}} 
\int_{\mathrm{M}_2(\mathbb{Q})\backslash \mathrm{M}_2(\mathbb{A})^1} 
\Bigg( \sum_{\substack{\mathrm{P}_1\in\mathcal{F}\\
\mathrm{standard}\\\mathrm{P}_1\subset\mathrm{P}_2}}
(-1)^{\dim(\mathrm{A}_1/\mathrm{A}_2)} \times
\\\nonumber\\\nonumber
&&\quad\times
\sum_{\delta\in\mathrm{P}_1(\mathbb{Q})\cap\mathrm{M}_2(\mathbb{Q})
\backslash \mathrm{M}_2(\mathbb{Q})}
\hat{\tau}_1^2(H_{\mathrm{P}_1}(\delta y)-T)
\sum_{X\in\mathfrak{m}_1(\mathbb{Q})\cap\mathfrak{o}}
\int_{\mathfrak{n}_1^2(\mathbb{A})} 
\\\nonumber\\\nonumber
&&
\qquad  
\times \bigg( \int_\mathrm{K} 
\int_{\mathfrak{n}_2(\mathbb{A})} 
f((X+N)\cdot\mathrm{ad}(\delta yx)+N'\cdot\mathrm{ad}(kx)) \times 
\\\nonumber\\\nonumber 
&& 
\qquad\quad
\times \Big(
\int_{\mathfrak{a}_2^\mathrm{G}} 
\Gamma_{\mathrm{P}_2}'(H-T,-H_{\mathrm{P}_2}(kx)) 
~\mathrm{d}H \Big) ~\mathrm{d}N'\mathrm{d}k~\bigg) ~\mathrm{d}N  \Bigg) ~\mathrm{d}y.
\end{eqnarray}
Substituting $T_0$ for $T$, the expression \eqref{elementary variance second expression} becomes
\begin{eqnarray}
\sum_{\substack{\mathrm{P}_2\in\mathcal{F}\\\mathrm{standard}}} 
J_\mathfrak{o}^{\mathrm{M}_2}(f_{{\mathrm{P}_2},x}) 
&=& 
\sum_{\mathrm{P}_2\in\mathcal{F}} 
|\mathrm{W}_0^{\mathrm{M}_2}| 
|\mathrm{W}_0^\mathrm{G}|^{-1} J_\mathfrak{o}^{\mathrm{M}_2}(f_{\mathrm{P}_2,x})
\end{eqnarray}
since the distribution $J_\mathfrak{o}^\mathrm{M}$ is independent of the choice of the minimal parabolic subgroup $\mathrm{P}_0$.
\qed

\section{Refined expansions}

In this chapter the distribution $J_\mathfrak{o}$ is decomposed as a linear combination of weighted orbital integrals following the methods of Arthur \cite{Arthur_Uni} \cite{Arthur_Orb}.

\subsection{Weighted orbital integrals}

\paragraph{Definition}\setcounter{equation}{0} 
Let $\mathrm{M}$ be a standard Levi subgroup of $\mathrm{G}$. A collection of complex-valued functions 
\begin{eqnarray}
\bigg\{c_\mathrm{P}\in C^\infty(i\mathfrak{a}_\mathrm{M}^*):~\mathrm{P}\in\mathcal{P}(\mathrm{M})\bigg\}
\end{eqnarray}
is said to be a \emph{$(\mathrm{G},\mathrm{M})$-family} if for each pair of adjacent parabolic subgroups $\mathrm{P}$ and $\mathrm{P}'$ in $\mathcal{P}(\mathrm{M})$, the functions $c_\mathrm{P}$ and $c_{\mathrm{P}'}$ agree on the hyperplane spanned by the common wall of the positive chambers of  $i\mathfrak{a}_\mathrm{M}^*$ defined by $\mathrm{P}$ and $\mathrm{P}'$.

Let $(c_\mathrm{P})$\label{sym57} be a $(\mathrm{G},\mathrm{M})$-family, define the function $c_\mathrm{M}$ on the complement of the coroot hyperplanes in $i\mathfrak{a}_\mathrm{M}^*$ by
\begin{eqnarray}
\label{weight constant}
\forall \lambda\in i\mathfrak{a}_\mathrm{M}^*
&&
c_\mathrm{M}(\lambda) 
=
\sum_{\mathrm{P}\in\mathcal{P}(\mathrm{M})} 
c_\mathrm{P}(\lambda)
\theta_\mathrm{P}(\lambda)^{-1}
\end{eqnarray}
where\label{sym56}
\begin{eqnarray}
\quad \theta_\mathrm{P}(\lambda) 
&=& 
\mathrm{Vol}(\mathfrak{a}_\mathrm{M}^\mathrm{G}\slash \mathbb{Z}(\Delta_\mathrm{P}^\vee)^{-1})
\prod_{\alpha\in\Delta_\mathrm{P}}
\langle\alpha\check{~},\lambda\rangle.
\end{eqnarray}
The function $c_\mathrm{M}$ extends smoothly over $i\mathfrak{a}_\mathrm{M}^*$. Denote by $c_\mathrm{M}$ its value at the origin of $i\mathfrak{a}_\mathrm{M}^*$.

\paragraph{Definition}\setcounter{equation}{0}
Let $\mathrm{Q}$ be a parabolic subgroup in $\mathcal{F}(\mathrm{M})$ with Levi component $\mathrm{L}$, let $(c_\mathrm{P})$\label{sym58} be a $(\mathrm{G},\mathrm{M})$-family. Let
\begin{eqnarray}
i_\mathrm{Q}^\mathrm{G}
&:&
\mathcal{P}^\mathrm{L}(\mathrm{M})\rightarrow\mathcal{P}^\mathrm{G}(\mathrm{M})
\end{eqnarray}
be the map that sends a parabolic subgroup $\mathrm{P}$ in $\mathcal{P}^\mathrm{L}(\mathrm{M})$ to the unique parabolic subgroup in $\mathcal{P}^\mathrm{G}(\mathrm{M})$ that is contained in $\mathrm{Q}$ whose intersection with $\mathrm{L}$ is $\mathrm{P}$. Denote by $(c_\mathrm{P}^\mathrm{Q})$ the $(\mathrm{L},\mathrm{M})$-family
\begin{eqnarray}
\bigg\{c_{i_\mathrm{Q}^\mathrm{G}(\mathrm{P})}:
~\mathrm{P}\in\mathcal{P}^\mathrm{L}(\mathrm{M})\bigg\}.
\end{eqnarray}
Let
\begin{eqnarray}
j_\mathrm{M}^\mathrm{L}
&:&
i\mathfrak{a}_\mathrm{L}^*\rightarrow i\mathfrak{a}_\mathrm{M}^*
\end{eqnarray}
be the natural inclusion map. Denote by the $(c_\mathrm{P})$ the $(\mathrm{G},\mathrm{L})$-family
\begin{eqnarray}
\bigg\{ j_\mathrm{M}^{\mathrm{L},*} (c_{\mathrm{P}'}):
~\mathrm{P}\in\mathcal{P}^\mathrm{G}(\mathrm{L}) \bigg\}
\end{eqnarray}
where $\mathrm{P}'$ is a parabolic subgroup in $\mathcal{P}^\mathrm{G}(\mathrm{M})$ contained in $\mathrm{P}$, and the function $j_\mathrm{M}^{\mathrm{L},*}(c_{\mathrm{P}'})$ is independent of the choice of $\mathrm{P}'$ in $\mathcal{P}^\mathrm{G}(\mathrm{M})$.

Let $(c_\mathrm{P})$ and $(d_\mathrm{P})$ be two $(\mathrm{G},\mathrm{M})$-families, denote by $((cd)_\mathrm{P})$ the product of $(c_\mathrm{P})$ and $(d_\mathrm{P})$, which is a $(\mathrm{G},\mathrm{M})$-family. For each parabolic subgroup $\mathrm{Q}$ in $\mathcal{F}(\mathrm{M})$ there exists a function $c_\mathrm{Q}'$\label{sym59} on $i\mathfrak{a}_\mathrm{M}^*$ such that
\begin{eqnarray}
\label{weight constant'}
\forall \lambda\in i\mathfrak{a}_\mathrm{M}^*
&&
(cd)_\mathrm{M}(\lambda)
=
\sum_{\mathrm{Q}\in\mathcal{F}(\mathrm{M})} 
c_\mathrm{Q}'(\lambda)d_\mathrm{M}^\mathrm{Q}(\lambda).
\end{eqnarray}
Denote by $c_\mathrm{Q}'$ the value $c_\mathrm{Q}'(0)$.

\paragraph{Definition}\setcounter{equation}{0}
\label{(g,m)-orthogonal set equation}
Let $\mathrm{M}$ be a standard Levi subgroup of $\mathrm{G}$. A collection of points\label{sym60}
\begin{eqnarray}
\mathcal{Y}_\mathrm{M}
&=&
\bigg\{
Y_\mathrm{P}\in\mathfrak{a}_\mathrm{M}:~\mathrm{P}\in\mathcal{P}(\mathrm{M})
\bigg\}
\end{eqnarray}
is said to be a \emph{$(\mathrm{G},\mathrm{M})$-orthogonal set} if for each pair of adjacent parabolic subgroups $\mathrm{P}$ and $\mathrm{P}'$ in $\mathcal{P}(\mathrm{M})$, the vector 
\begin{eqnarray}
Y_\mathrm{P}-Y_{\mathrm{P}'}
&\in& 
\mathfrak{a}_\mathrm{M}
\end{eqnarray}
is orthogonal to the hyperplane spanned by the common wall of the positive chambers defined by $\mathrm{P}$ and $\mathrm{P}'$.

Let $\alpha_{\mathrm{P}'}^\mathrm{P}\check{~}$ be the unique coroot that is positive for $\mathrm{P}$ and negative for $\mathrm{P}'$. A $(\mathrm{G},\mathrm{M})$-orthogonal set $\mathcal{Y}_\mathrm{M}$ is \emph{positive} if 
\begin{eqnarray}
\exists t>0
&&
Y_\mathrm{P}-Y_{\mathrm{P}'}=t\alpha_{\mathrm{P}'}^\mathrm{P}\check{~}.
\end{eqnarray}

\paragraph{Remark}\setcounter{equation}{0}
Let $\mathcal{Y}_\mathrm{M}$ be a positive $(\mathrm{G},\mathrm{M})$-orthogonal set, then the collection of functions\label{sym61}
\begin{eqnarray}
\label{(g,m)-orthogonal set to (g,m)-family}
\bigg\{
v_\mathrm{P}(\mathcal{Y}_\mathrm{M})(\lambda) = e^{\langle\lambda,Y_\mathrm{P}\rangle}
:~ \mathrm{P}\in\mathcal{P}(\mathrm{M})
\bigg\}
\end{eqnarray}
forms a $(\mathrm{G},\mathrm{M})$-family. The associated constant $v_\mathrm{M}(\mathcal{Y}_\mathrm{M})$ as in \eqref{weight constant} is equal to the volume of the convex hull of $\mathcal{Y}_\mathrm{M}$ in $\mathfrak{a}_\mathrm{M}$.

\paragraph{Definition}\setcounter{equation}{0}
Let $\mathrm{M}$ be a standard Levi subgroup of $\mathrm{G}$, let $x$ be an element of $\mathrm{G}(\mathbb{Q}_S)$. The collection of points\label{sym62}
\begin{eqnarray}
\label{weighted orbital integral (g,m)-orthogonal set}
\mathcal{Y}_\mathrm{M}(x)
&=&
\bigg\{-H_\mathrm{P}(x):~\mathrm{P}\in\mathcal{P}(\mathrm{M})\bigg\}
\end{eqnarray}
forms a positive $(\mathrm{G},\mathrm{M})$-orthogonal set. Define the \emph{weight factor} $v_\mathrm{M}(x)$ to be the associated constant 
\begin{eqnarray}
\label{weight factor definition}
v_\mathrm{M}(x)
&=&
v_\mathrm{M}(\mathcal{Y}_\mathrm{M}(x)).
\end{eqnarray}

Let $X$ be a point in $\mathfrak{m}(\mathbb{Q}_S)$ such that $\mathrm{G}_X^0(\mathbb{Q}_S)$ is contained in $\mathrm{M}(\mathbb{Q}_S)$, define the \emph{weighted orbital integral $J_\mathrm{M}^\mathrm{G}(X,~)$}\label{sym63} to be the distribution on $\mathfrak{g}(\mathbb{Q}_S)$ such that
\begin{eqnarray}
\label{regular weighted orbital integral}
&&
\forall f\in\mathcal{S}(\mathfrak{g}(\mathbb{Q}_S))
\\\nonumber
&&
J_\mathrm{M}^\mathrm{G}(X,f)=|D^\mathrm{G}(X)|_S^{1/2} \int_{\mathrm{G}_X^0(\mathbb{Q}_S) \backslash \mathrm{G}(\mathbb{Q}_S)} f(X\cdot\mathrm{ad}(x)) ~v_\mathrm{M}(x)\mathrm{d}x.
\end{eqnarray}

For a general element $X$ in $\mathfrak{m}(\mathbb{Q}_S)$ define the \emph{weighted orbital integral $J_\mathrm{M}^\mathrm{G}(X,~)$} to be the distribution on $\mathfrak{g}(\mathbb{Q}_S)$ such that
\begin{eqnarray}
\label{general weighted orbital integral}
&&
\forall f\in\mathcal{S}(\mathfrak{g}(\mathbb{Q}_S))
\\\nonumber
&&
J_\mathrm{M}^\mathrm{G}(X,f)=\lim_{A\rightarrow0} \sum_{\mathrm{L}\in\mathcal{L}(\mathrm{M})} r_\mathrm{M}^\mathrm{L}(\exp(X_\mathrm{nil}),\exp(A)) J_\mathrm{L}^\mathrm{G}(X+A,f)
\end{eqnarray}
where $A$ is a sequence of $\mathbb{Q}_S$-points of $\mathrm{Lie}(\mathrm{A}_\mathrm{M})$ such that
\begin{eqnarray}
\mathrm{G}_{X+A}^0(\mathbb{Q}_S)
&\subset&
\mathrm{M}(\mathbb{Q}_S)
\end{eqnarray}
and $r_\mathrm{P}^\mathrm{L}(x,a)$\label{sym64} is an auxiliary ($\mathrm{L},\mathrm{M}$)-family constructed by Arthur in \S5 of  \cite{Arthur_Loc}, see also \S2.4 of \cite{Hoff} and {\S}III.2 of \cite{Wal}.

\paragraph{Lemma}\setcounter{equation}{0}
\emph{Let $\mathrm{M}$ be a standard Levi subgroup of $\mathrm{G}$. Let $X$ be an element of $\mathfrak{m}(\mathbb{Q}_S)$, let $x$ be an element of $\mathrm{G}(\mathbb{Q}_S)$. Let $f$ be a Schwartz function on $\mathcal{S}(\mathfrak{g}(\mathbb{Q}_S))$. Then
\begin{eqnarray}
\label{weighted orbital integral variance}
J_\mathrm{M}^\mathrm{G}(X,f\circ\mathrm{ad}(x)) &=& \sum_{\mathrm{P}\in\mathcal{F}(\mathrm{M})} J_\mathrm{M}^{\mathrm{M}_\mathrm{P}}(X,f_{\mathrm{P},x})
\end{eqnarray}
where $f_{\mathrm{P},x}$\label{sym65} is the Schwartz function on $\mathfrak{m}_\mathrm{P}(\mathbb{Q}_S)$ defined as
\begin{eqnarray}
f_{\mathrm{P},x}(X) &=& \int_\mathrm{K}\int_{\mathfrak{n}_\mathrm{P} (\mathbb{A})} f((X+N)\cdot\mathrm{ad}(k)) ~v_\mathrm{P}'(kx) \mathrm{d}N\mathrm{d}k
\end{eqnarray}
where for each $y$ in $\mathrm{G}(\mathbb{Q}_S)$ and for each parabolic subgroup $\mathrm{Q}$ in $\mathcal{F}(\mathrm{M})$ the constant $v_\mathrm{Q}'(y)$ is the constant $v_\mathrm{Q}'$ intervening in \ref{weight constant'} associated to the $(\mathrm{G},\mathrm{M})$-orthogonal set $\mathcal{Y}_\mathrm{M}(x)$ defined in \ref{weighted orbital integral (g,m)-orthogonal set}.}

\proof
See III.3.(f) of \cite{Wal}.
\qed

\subsection{Nilpotent orbits}

\paragraph{Definition}\setcounter{equation}{0}
Let $\mathfrak{g}_\mathrm{nil}$\label{sym66} be the $\sim$ equivalence class of the origin in $\mathfrak{g}(\mathbb{Q})$, hence the nilpotent locus of $\mathfrak{g}$. Let $J_\mathrm{nil}$ and $J_\mathrm{nil}^T$ denote respectively $J_{\mathfrak{g}_\mathrm{nil}}$ and $J_{\mathfrak{g}_\mathrm{nil}}^T$. Let $J_\mathrm{nil}^\mathrm{G}$ denote $J_\mathrm{nil}$ for convenience in inductive arguments involving Levi subgroups.

\paragraph{Lemma}\setcounter{equation}{0}
\emph{There exists a continuous seminorm $\|~\|$\label{sym67} on $\mathcal{S}(\mathfrak{g}(\mathbb{A}))$ such that for every truncation parameter $T$ in $\mathfrak{a}_0$
\begin{eqnarray}
\label{refined asymptotic expansion}
&&
\forall f\in\mathcal{S}(\mathfrak{g}(\mathbb{A}))
\\\nonumber
&&
\Bigg|J_\mathrm{nil}^T(f)
-\int_{\mathrm{G}(\mathbb{Q})\backslash\mathrm{G}(\mathbb{A})^1}
F^\mathrm{G}(x,T) 
\sum_{X\in\mathfrak{g}_\mathrm{nil}(\mathbb{Q})}
f(X\cdot\mathrm{ad}(x)) ~\mathrm{d}x\Bigg| 
\leq 
\| f \| e^{- \frac{\mathrm{d}(T)}{2}}
\end{eqnarray}
where $\mathrm{d}(T)$ denotes the distance from $T$ to the root hyperplanes.}

\proof
The argument is valid for a general class $\mathfrak{o}$ and the sum of all the classes $\mathfrak{o}$, hence contains \ref{asymptotic expansion} as a special case. The necessary estimates are established by Chaudouard in the proof of Proposition 4.4 of \cite{Chaud}. The following is a sketch of the argument, for the details please refer to loc. cit. See also Theorem 3.1 of \cite{Arthur_Uni}.

Following the first part of the proof of \eqref{truncated kernel integrability}, by the combinatorial lemma of Langlands, expand $J_\mathrm{nil}^T(f)$ as a sum indexed by triples of nested standard parabolic subgroups whose leading term corresponding to $(\mathrm{G},\mathrm{G},\mathrm{G})$ is
\begin{eqnarray}
\int_{\mathrm{G}(\mathbb{Q})\backslash\mathrm{G}(\mathbb{A})^1}
F^\mathrm{G}(x,T) 
\sum_{X\in\mathfrak{g}_\mathrm{nil}(\mathbb{Q})}
f(X\cdot\mathrm{ad}(x)) ~\mathrm{d}x.
\end{eqnarray}
Hence the left hand side of the inequality in \eqref{refined asymptotic expansion} is bounded by
\begin{eqnarray}
\label{refined asymptotic expansion first majorant}
&&
\sum_{\substack{\mathrm{P}_1,\mathrm{P}_2,\mathrm{P}_3\in\mathcal{F}\\
\mathrm{standard}\\
\mathrm{P}_1\subset\mathrm{P}_2\subsetneq\mathrm{P}_3}}
\int_{\mathrm{P}_1(\mathbb{Q})\backslash\mathrm{G}(\mathbb{A})^1} 
F^{\mathrm{P}_1}(x,T) 
\sigma_1^3(H_0(x)-T) 
\sum_{X\in\mathfrak{m}_1^2(\mathbb{Q})'\cap\mathfrak{o}} 
\Bigg|\sum_{\overline{Y} \in\overline{\mathfrak{n}}_2^3(\mathbb{Q})'} 
\Phi_X(x,\overline{Y}) \Bigg| ~\mathrm{d}x
\end{eqnarray}
where $\Phi$\label{sym68} denotes the partial Fourier transform of $f$
\begin{eqnarray}
\Phi_X(x,\overline{Y}) 
&=& 
\int_{\mathfrak{n}_2(\mathbb{A})} 
f((X+N)\cdot\mathrm{ad}(x)) 
\cdot \psi(\langle N,\overline{Y}\rangle) ~\mathrm{d}N.
\end{eqnarray}
Changing variables each summand of \eqref{refined asymptotic expansion first majorant} is bounded by
\begin{eqnarray}
\label{refined asymptotic expansion second majorant}
&&
\sup_{y\in\Gamma} 
\Bigg( \int_{\mathrm{A}_1^\mathrm{G}(\mathbb{R})_0} 
\int_{\mathrm{A}_{0,T'}^{1,T}(\mathbb{R})}
\delta_0^2(a_1a)^{-1}
\sigma_1^3(H_{\mathrm{P}_1}(a_1)-T) \times 
\\\nonumber\\\nonumber 
&& \quad
\times \sum_{X\in\mathfrak{m}_1^2(\mathbb{Q})' \cap\mathfrak{o}}  
~\sum_{\overline{Y}\in\overline{ \mathfrak{n}}_2^3(\mathbb{Q})'} 
|\Phi_{X\cdot\mathrm{ad}(a_1a)} (y,\overline{Y}\cdot\mathrm{ad}(a_1a))| ~\mathrm{d}a\mathrm{d}a_1 \Bigg)
\end{eqnarray}
where $\Gamma$\label{sym69} is a fixed compact subset of $\mathrm{M}_0(\mathbb{A})^1$ which is independant of the truncation parameter $T$,  the subset $\mathrm{A}_{0,T'}^{1,T}(\mathbb{R})$ of $\mathrm{A}_0^1(\mathbb{R})$ is the $T',T$-truncated part as in \eqref{truncated siegel domain}, and $\delta_0^2$ is the modulus function of $\mathrm{P}_2$.

Let $n$ be a natural number, let $\mathcal{D}$\label{sym70} be an invariant differential operator on $\mathfrak{g}(\mathbb{R})$ of degree $n$, denote by $\Phi^\mathcal{D}$ the partial Fourier transform of $\mathcal{D}f$ where $\mathcal{D}$ operates on $f$ via its archimedean component $f_\infty$, then 
\begin{eqnarray}
&&
\big|\Phi_{X\cdot\mathrm{ad}(a_1a)} (y,\overline{Y}\cdot\mathrm{ad}(a_1a))\big| 
=  
C(y)\cdot \big|\overline{Y}\cdot\mathrm{ad}(a_1a)\big|^{-n} 
\cdot \big|\Phi_{X\cdot\mathrm{ad}(a_1a)}^\mathcal{D} (y,\overline{Y}\cdot\mathrm{ad}(a_1a))\big|
\end{eqnarray}
for some constant $C(y)$ that depends continuously on $y$. For a Schwartz function $f$ on $\mathfrak{g}(\mathbb{A})$ define $N(f)$\label{sym71} to be the smallest natural number $N$ such that $f$ is supported on $\frac{1}{N}\widehat{\mathbb{Z}}\times\mathbb{R}$ where $\widehat{\mathbb{Z}}$\label{sym72} denotes the profinite completion of $\mathbb{Z}$. Let the natural number $n$ be large enough so that
\begin{eqnarray}
\sum_{\overline{Y}\in\overline{\mathfrak{n}}_2^3(\frac{1}{N(f)}\mathbb{Z})'} 
\big|\overline{Y}\cdot\mathrm{ad}(a_1a)\big|^{-n} 
&\leq& 
C \prod_{\alpha\in\Delta_0^3} 
e^{-k_\alpha\alpha(H_0(a_1a))}
\end{eqnarray}
for some constant $C$ and natural numbers $k_\alpha$. Define the seminorm $\|~\|'$\label{sym73} by 
\begin{eqnarray}
\|f\|'&=& \sup_{X\in\mathfrak{g}(\mathbb{A})} |\mathcal{D}f(X)|.
\end{eqnarray}
A Schwartz function on the $\mathbb{A}$-valued points of a rational vector space is bounded by a product of Schwartz functions on each coordinate, hence for $Z$ in $\mathfrak{p}_2(\mathbb{A})$
\begin{eqnarray}
|\mathcal{D}f(Z\cdot\mathrm{ad}(y))| &\leq& \bigg( \prod_{\mu\in\Phi_0-\Phi_2} \phi_\mu(Z_\mu) \bigg) \phi_{\mathfrak{n}_2}(Z_{\mathfrak{n}_2})
\end{eqnarray}
where $Z_\mu$\label{sym74} and $Z_{\mathfrak{n}_2}$ are the components of $Z$ on the weight space of $\mu$ and $\mathfrak{n}_2$, and $\phi_\bullet$\label{sym75} are positive Schwartz functions. If $a_0$ is an element of $\mathrm{A}_0(\mathbb{R})$, denote by $\Psi(a_0)$ for the sum
\begin{eqnarray}
\Psi(a_0) 
&=& 
\sum_{X\in\mathfrak{m}_1^2(\frac{1}{N(f)}\mathbb{Z})' \cap \mathfrak{o}} 
\bigg( \prod_{\mu\in\Phi_0-\Phi_2} \phi_\mu(\mu(a_0)^{-1}X_\mu) \bigg).
\end{eqnarray}
Then \eqref{refined asymptotic expansion second majorant} is bounded by
\begin{eqnarray}
\label{refined asymptotic expansion third majorant}
&&
\sup_{y\in\Gamma}~ C(y)  
\int_{\mathrm{A}_1^\mathrm{G}(\mathbb{R})} 
\int_{\mathrm{A}_{0,T'}^{1, T} (\mathbb{R})} 
\bigg( \delta_0^2(a_1a)^{-1}
\sigma_1^2(H_{\mathrm{P}_1}(a_1)-T)\times
\\\nonumber\\\nonumber 
&& \quad
\times \Psi(a_1a) \cdot C 
\prod_{\alpha\in\Delta_0^2} 
e^{-k_\alpha\alpha(H_0(a_1a))} \bigg) ~\mathrm{d}a\mathrm{d}a_1,
\end{eqnarray}
which only depends on the Schwartz function $f$ via its componentwise bounds $\phi_\bullet$ and the lattice $\overline{\mathfrak{n}}_2^3(\frac{1}{N(f)}\mathbb{Z})$, hence is proportional to the seminorm $\|~\|$ defined by
\begin{eqnarray}
\forall f\in\mathcal{S}(\mathfrak{g}(\mathbb{A}))
&&
\|f\|=\sup_{y\in\Gamma}~ C(y)\cdot C\cdot\|f\|'\cdot N(f)^n.
\end{eqnarray}
The constant $\sup_{y\in\Gamma} C(y)\cdot C$ is independant of both the Schwartz function $f$ and the truncation parameter $T$.

Substituting the definition of the seminorm $\|~\|$, the majorant \eqref{refined asymptotic expansion third majorant} reduces to
\begin{eqnarray}
\label{refined asymptotic expansion fourth majorant}
\|f\| \cdot \mathrm{Vol}(\mathrm{A}_{0,T'}^{1,T}(\mathbb{R})) 
\prod_{\alpha\in \Delta_1^3} 
\bigg( e^{-\alpha(T)}\int_0^\infty p(t_\alpha) e^{-t_\alpha} ~\mathrm{d}t_\alpha\bigg)
\end{eqnarray}
where $p(t)$ is a polynomial. Because 
\begin{eqnarray}
\prod_{\alpha\in \Delta_1^3} e^{-\alpha(T)}\int_0^\infty p(t)e^{-t} ~\mathrm{d}t
&\leq&
e^{-\mathrm{d}(T)},
\end{eqnarray}
and $\mathrm{Vol}(\mathrm{A}_{0,T'}^{1,T}(\mathbb{R}))$ is of polynomial growth in $T$, the majorant \eqref{refined asymptotic expansion fourth majorant} reduces to 
\begin{eqnarray}
\|f\|e^{-\frac{\mathrm{d}(T)}{2}}.
\end{eqnarray}
\qed

\paragraph{Definition}\setcounter{equation}{0}
\label{definition f epsilon}
Let $v$ be a place of $\mathbb{Q}$, let $\beta_v$\label{sym76} be a bump function on $\mathbb{Q}_v$, let $\nu$\label{sym77} be a $\mathrm{G}$-orbit contained in $\mathfrak{g}_\mathrm{nil}$, let $\{p_1,\dots,p_l\}$ be a collection of polynomials rational with rational coefficients cutting out the Zariski closure $\overline{\nu}$. Let $f$ be a Schwartz function on $\mathfrak{g}(\mathbb{A})$, let $\epsilon$ be a positive real number. Define the \emph{truncated function} $f_{\nu,v}^\epsilon$\label{sym78} on $\mathfrak{g}(\mathbb{A})$ as in \cite{Arthur_Uni} by
\begin{eqnarray}
\label{f epsilon}
\forall X\in\mathfrak{g}(\mathbb{A})
&&
f_{\nu,v}^\epsilon (X) 
=
f(X)\beta_v(\epsilon^{-1}|p_1(X)|_v)\dots\beta_v(\epsilon^{-1}|p_l(X)|_v ).
\end{eqnarray}

\paragraph{Lemma}\setcounter{equation}{0}
\label{f epsilon estimate}
\emph{Let the place $v$, the bump function $\beta_v$, the orbit $\nu$ and the polynomials $\{p_1,\dots,p_l\}$ be as in \ref{definition f epsilon}, then there exists a natural number $m$ and another seminorm $\|~\|_1$\label{sym79} for which the inequality (\ref{refined asymptotic expansion}) holds such that}
\begin{eqnarray}
\label{f epsilon inequality}
\forall \epsilon \text{ \emph{such that} }0<\epsilon<1~ 
\forall f\in\mathcal{S}(\mathfrak{g}(\mathbb{A}))
&&
\|f_{\nu,v}^\epsilon\| \leq \epsilon^{-ml}\|f\|_1.
\end{eqnarray}

\proof
It is enough to consider the case when there is a single polynomial $p$. There are two cases:
\begin{itemize}
\item If $v$ is the archimedean place then 
\begin{eqnarray}
N(f_{\nu,v}^\epsilon) &=& N(f), 
\end{eqnarray}
so by properties of the derivative
\begin{eqnarray}
\|f_{\nu,v}^\epsilon\|' &\leq& \epsilon^{-m}C\|f\|'.
\end{eqnarray}
Define $\|f\|_1'$ to be $C\|f\|$.
\item If $v$ is finite then 
\begin{eqnarray}
\|f_{\nu,v}^\epsilon\|'=\|f\|',
\end{eqnarray}
so 
\begin{eqnarray}
N(f_{\nu,v}^\epsilon) &\leq& \epsilon^{-m} N(f)\end{eqnarray}
since $N(f)$ only depends on the support of $f_\mathrm{finite}$ and the assignment 
\begin{eqnarray}
f
&\mapsto&
f_{\nu,v}^\epsilon=f \prod\beta_v(\epsilon^{-1}|p|_v)
\end{eqnarray}
shrinks the support of $f$ by a factor of $\epsilon^m$, upto a multiplicative constant. Define $\|f\|'$ to be $\|f\|$.
\end{itemize} 
Hence \eqref{f epsilon inequality} follows since
\begin{eqnarray}
\|f\| &=& \|f\|'N(f).
\end{eqnarray}
\qed

\paragraph{Lemma}\setcounter{equation}{0}
\emph{Let the place $v$, the bump function $\beta_v$, the orbit $\nu$ and the polynomials $\{p_1,\dots,p_l\}$ be as in \ref{definition f epsilon}, then there exists a positve real number $r$ such that
\begin{eqnarray}
\label{f epsilon asymptotic expansion}
&&
\forall \epsilon>0~
\forall f\in\mathcal{S}(\mathfrak{g}(\mathbb{A}))
\\\nonumber
&&\quad
\int_{\mathrm{G}(\mathbb{Q})\backslash \mathrm{G}(\mathbb{A})^1}
F^\mathrm{G}(x,T) 
\sum_{X\in\mathfrak{g}(\mathbb{Q})-\overline{\nu}(\mathbb{Q})}
\big|f_{\nu,v}^\epsilon(X\cdot\mathrm{ad}(x)) \big| ~\mathrm{d}x 
\leq
\|f\| \epsilon^r (1+|T|)^{d_0}
\end{eqnarray}
where $d_0$ is $\dim(\mathrm{A}_0/\mathrm{A}_\mathrm{G})$, the split rank of $\mathrm{G}$.}

\proof
This is Lemma 4.1 of \cite{Arthur_Uni}.
\qed

\paragraph{Proposition}\setcounter{equation}{0}
\emph{Let $T$ be a truncation parameter. For each nilpotent orbit $\nu$ there exists a distribution $J_\nu^T$\label{sym80} on $\mathfrak{g}(\mathbb{A})$ such that for each Schwartz function $f$ on $\mathfrak{g}(\mathbb{A})$, the expression $J_\nu^T(f)$ is a polynomial in $T$ of degree at most $d_0$, and
\begin{eqnarray}
J_\mathrm{nil}^T(f) &=& 
\sum_{\substack{\nu\subset\mathfrak{g}_\mathrm{nil}\\\textrm{orbit}}} J_\nu^T(f).
\end{eqnarray}
There exists a continuous seminorm $\|~\|$\label{sym81} on $\mathcal{S}(\mathfrak{g}(\mathbb{A}))$ and a positive real number $\epsilon$ such that}
\begin{eqnarray}
&&
\forall f\in\mathcal{S}(\mathfrak{g}(\mathbb{A}))
\\\nonumber
&&
\bigg|J_\nu^T(f)-\int_{\mathrm{G}(\mathbb{Q})\backslash\mathrm{G}(\mathbb{A})^1}F^\mathrm{G}(x,T) \sum_{X\in\nu(\mathbb{Q})}f(X\cdot\mathrm{ad}(x)) ~\mathrm{d}x\bigg|
\leq
\| f \| e^{- \epsilon\mathrm{d}(T)}.
\end{eqnarray}

\proof
Define the polynomials recursively by the formula
\begin{eqnarray}
J_{\overline{\nu}}^T(f) &=& \lim_{\epsilon\rightarrow0}J_\mathrm{nil}^T(f_{\nu,v}^\epsilon)
\end{eqnarray}
where $J_{\overline{\nu}}^T$\label{sym82} denotes the sum of $J_{\nu'}^T$ for those orbits $\nu'$ contained in $\overline{\nu}$.

The limit of $J_\mathrm{nil}^T(f_{\nu,v}^\epsilon)$ as $\epsilon$ approaches 0 exists because
\begin{eqnarray}
\sum_{X\in\overline{\nu}(\mathbb{Q})} f(X\cdot\mathrm{ad}(x)) &=& \sum_{X\in\overline{\nu}(\mathbb{Q})} f_{\nu,v}^\epsilon(X\cdot\mathrm{ad}(x))
\end{eqnarray}
implies that
\begin{eqnarray}
\label{nilpotent asymptotic polynomial proof first inequality}
&&
\bigg|J_\mathrm{nil}^T(f_{\nu,v}^\epsilon)-\int_{\mathrm{G}(\mathbb{Q})\backslash\mathrm{G}(\mathbb{A})^1}F^\mathrm{G}(x,T) \sum_{X\in{\overline{\nu}}(\mathbb{Q})}f(X\cdot\mathrm{ad}(x)) ~\mathrm{d}x\bigg|
\\\nonumber\\\nonumber
&\leq&
\bigg|J_\mathrm{nil}^T(f_{\nu,v}^\epsilon)-\int_{\mathrm{G}(\mathbb{Q}) \backslash \mathrm{G}(\mathrm{A})^1}F^\mathrm{G}(x,T)\sum_{X\in\mathfrak{g}_\mathrm{nil}(\mathbb{Q})} f_{\nu,v}^\epsilon (X\cdot\mathrm{ad}(x)) ~\mathrm{d}x \bigg| +  
\\\nonumber\\\nonumber
&&
\quad
 + \int_{\mathrm{G}(\mathbb{Q}) \backslash \mathrm{G}(\mathrm{A})^1}F^\mathrm{G}(x,T)\sum_{X\in\mathfrak{g}_\mathrm{nil}(\mathbb{Q})-\overline{\nu}(\mathbb{Q})} \bigg|f_{\nu,v}^\epsilon (X\cdot\mathrm{ad}(x)) \bigg| ~\mathrm{d}x.
\end{eqnarray}
By \eqref{refined asymptotic expansion} the first summand of the right hand side of \eqref{nilpotent asymptotic polynomial proof first inequality} is bounded by 
\begin{eqnarray}
\|f_{\nu,v}^\epsilon\|e^{-\frac{\mathrm{d}(T)}{2}}
&\leq&
\epsilon^{-lm}\cdot\|f\|_1\cdot e^{-\frac{\mathrm{d}(T)}{2}}
\end{eqnarray}
which follows from \eqref{f epsilon inequality}.

The second summand of the right hand side of \eqref{nilpotent asymptotic polynomial proof first inequality} is bounded by summing over the complement of $\overline{\nu}(\mathbb{Q})$ in $\mathfrak{g}(\mathbb{Q})$ instead of in $\mathfrak{g}_\mathrm{nil}(\mathbb{Q})$, hence by \eqref{f epsilon asymptotic expansion} is bounded by
\begin{eqnarray}
\|f\|_1\cdot \epsilon^r\cdot(1+|T|)^{d_0}.
\end{eqnarray}
Therefore the right hand side of \eqref{nilpotent asymptotic polynomial proof first inequality} is bounded by
\begin{eqnarray}
\label{nilpotent asymptotic polynomial proof second majorant}
\|f\|_1 \bigg( \epsilon^{-lm}e^{-\frac{\mathrm{d}(T)}{2}}+\epsilon^r(1+|T|)^{d_0} \bigg).
\end{eqnarray}
It suffices to take $\delta^n$ as $\epsilon$ with $\delta$ bounded strictly between 0 and 1 and $n$ a sequence natural numbers approaching infinity. Then the majorant \eqref{nilpotent asymptotic polynomial proof second majorant} satisfies the inequality
\begin{eqnarray}
\|f\|_1 \bigg( e^{|\log(\delta)|lmn-\frac{\mathrm{d}(T)}{2}}+\delta^{rn}(1+|T|)^{d_0} \bigg) &\leq& \|f\|\cdot\delta^{rn}\cdot(1+|T|)^{d_0}
\end{eqnarray}
provided $\mathrm{d}(T)$ is bounded below by $C|\log(\delta)|n$ for some constant $C$ and $\|~\|$ is another continuous seminorm with the same properties.

Therefore for every natural number $n$ and every $T$ in the open subcone 
\begin{eqnarray}
\bigg\{T:\mathrm{d}(T)>C|\log(\delta)|(n+1)\bigg\}
\end{eqnarray}
of the positive chamber, the following inequality holds
\begin{eqnarray}
|J_\mathrm{nil}^T(f_{\nu,v}^{\delta^n})-J_\mathrm{nil}^T(f_{\nu,v}^{\delta^{n+1}})| &\leq& 2\|f\|(1+|T|)^{d_0}\delta^{rn}.
\end{eqnarray}
The left hand side is a polynomial in $T$ of degree at most $d_0$, hence the polynomial extrapolation lemma in Lemma 5.2 of \cite{Arthur_Eis} applies:

There exists a constant $A$ such that for all $T$,
\begin{eqnarray}
\qquad
|J_\mathrm{nil}^T(f_{\nu,v}^{\delta^n})-J_\mathrm{nil}^T(f_{\nu,v}^{\delta^{n+1}})| &\leq& A \cdot \|f\| \cdot (1+|T|)^{d_0} \cdot \big( |\log(\delta)|(n+1) \big)^{d_0}  \cdot \delta^{rn}.
\end{eqnarray}
Since
\begin{eqnarray}
\sum_{n=0}^\infty \big( |\log(\delta)|(n+1) \big)^{d_0} \delta^{rn}&<& \infty,
\end{eqnarray}
by telescoping the series 
\begin{eqnarray}
\sum_{n=0}^\infty
J_\mathrm{nil}^T(f_{\nu,v}^{\delta^n})-J_\mathrm{nil}^T(f_{\nu,v}^{\delta^{n+1}})
\end{eqnarray}
the sequence $J_\mathrm{nil}^T(f_{\nu,v}^{\delta^n})$ has a limit as $n$ approaches infinity. The limit is a polynomial in $T$ of degree at most $d_0$, denoted by $J_{\overline{\nu}}^T(f)$.

By construction
\begin{eqnarray}
\qquad\qquad
&& 
\bigg|J_{\overline{\nu}}^T(f)-\int_{\mathrm{G}(\mathbb{Q})\backslash\mathrm{G}(\mathbb{A})^1}F^\mathrm{G}(x,T) \sum_{X\in{\overline{\nu}}(\mathbb{Q})}f(X\cdot\mathrm{ad}(x)) ~\mathrm{d}x\bigg| 
\\\nonumber\\\nonumber 
&\leq& 
\bigg| J_\mathrm{nil}^T(f_{\nu,v}^{\delta^n})-J_{\overline{\nu}}^T(f) \bigg| + \bigg|J_\mathrm{nil}^T(f_{\nu,v}^{\delta^n})-\int_{\mathrm{G}(\mathbb{Q})\backslash\mathrm{G}(\mathbb{A})^1}F^\mathrm{G}(x,T) \sum_{X\in{\overline{\nu}}(\mathbb{Q})}f(X\cdot\mathrm{ad}(x)) ~\mathrm{d}x\bigg| 
\\\nonumber\\
\label{nilpotent asymptotic polynomial proof third inequality} 
&\leq& 
\sum_{n=0}^\infty A\|f\|(1+|T|)^{d_0} \big( |\log(\delta)|(n+1) \big)^{d_0} \delta^{rn}+\|f\|\delta^{rn}(1+|T|)^{d_0}
\end{eqnarray}
where \eqref{nilpotent asymptotic polynomial proof third inequality} holds for a fixed $\delta$, a sufficiently positive $T$ and $n$ the largest natural number such that 
\begin{eqnarray}
\mathrm{d}(T) &\geq& C|\log(\delta)|n.
\end{eqnarray}

Therefore it is possible choose a new seminorm $\|~\|$ and a new constant $\epsilon>0$ such that
\begin{eqnarray}
\Bigg|J_{\overline{\nu}}^T(f)-\int_{\mathrm{G}(\mathbb{Q})\backslash\mathrm{G}(\mathbb{A})^1}F^\mathrm{G}(x,T) \sum_{X\in{\overline{\nu}}(\mathbb{Q})}f(X\cdot\mathrm{ad}(x)) ~\mathrm{d}x\Bigg| &\leq& \|f\|e^{-\epsilon\mathrm{d}(T)}.
\end{eqnarray}
The corresponding statements for $\nu$ instead of $\overline{\nu}$ follow recursively by setting
\begin{eqnarray}
J_{\overline{\nu}}^T(f) &=& \sum_{\substack{\nu'\subset\overline{\nu}\\\mathrm{orbit}}} J_{\nu'}^T(f).
\end{eqnarray} 
\qed

\paragraph{Proposition}\setcounter{equation}{0} 
\emph{Let $S$ be a finite set of places of $\mathbb{Q}$ containing the archimedean place, let $f$ be a Schwartz function on $\mathfrak{g}(\mathbb{A})$ such that $f_p$ is the characteristic function of the standard lattice $\mathfrak{g}(\mathbb{Z}_p)$ whenever $p$ is not contained in $S$. Let $\mathrm{M}$ be a standard Levi subgroup of $\mathrm{G}$. Denote by $\mathfrak{m}_\mathrm{nil}(\mathbb{Q})_{\mathrm{M},S}$\label{sym83} the set of $\mathrm{M}(\mathbb{Q}_S)$-conjugacy classes in $\mathfrak{m}_\mathrm{nil}(\mathbb{Q})$. For each nilpotent conjugacy class $\nu$ in $\mathfrak{m}_\mathrm{nil}(\mathbb{Q})_{\mathrm{M},S}$ here exists a constant $a^\mathrm{M}(S,\nu)$\label{sym84} such that}
\begin{eqnarray}
\label{nilpotent refined expansion} 
J_\mathrm{nil}^\mathrm{G}(f) &=& \sum_{\mathrm{M}\in\mathcal{L}} |\mathrm{W}_0^\mathrm{M}| 
|\mathrm{W}_0^\mathrm{G}|^{-1}
\sum_{\nu\in\mathfrak{m}_\mathrm{nil}(\mathbb{Q})_{\mathrm{M},S}} a^\mathrm{M}(S,\nu)J_\mathrm{M}^\mathrm{G} (\nu,f_S). 
\end{eqnarray}

\paragraph*{Remark}
By Lemma 7.1 of \cite{Arthur_Uni}, the sets $\mathfrak{m}_\mathrm{nil}(\mathbb{Q})_{\mathrm{M},S}$ and $\mathfrak{m}_\mathrm{nil}(\mathbb{Q}_S)/\mathrm{M}(\mathbb{Q}_S)^1$ are in natural bijection, hence the expression $J_\mathrm{M}^\mathrm{G}(\nu,f_S)$ makes sense.
 
\proof
The argument is based on the following two lemmas:

\paragraph*{Lemma}
\emph{For every $x$ in $\mathrm{G}(\mathbb{A})^1$
\begin{eqnarray}
\label{nilpotent truncated kernel variance}
J_\mathrm{nil}^\mathrm{G}(f\circ\mathrm{ad}(x)) 
&=& 
\sum_{\mathrm{Q}\in\mathcal{F}} 
|\mathrm{W}_0^{\mathrm{M}_\mathrm{Q}}| 
|\mathrm{W}_0^\mathrm{G}|^{-1} 
J_\mathrm{nil}^{\mathrm{M}_\mathrm{Q}}(f_{\mathrm{Q},x}), 
\end{eqnarray}
where the function $f_{\mathrm{Q},x}$ on $\mathfrak{m}_\mathrm{Q}(\mathbb{A})$ is defined as in (\ref{elementary variance parabolic descent' definition}).}

\proof
This is a special case of \eqref{elementary variance}.
\qed

\paragraph*{Lemma}
\emph{Let $\mathrm{M}$ be a standard Levi subgroup. For every point $x$ in $\mathrm{G}(\mathbb{Q}_S)^1$ and $\nu$ in $\mathfrak{m}_\mathrm{nil}(\mathbb{Q}_S)\slash\mathrm{M}(\mathbb{Q}_S)^1$}
\begin{eqnarray}
\label{nilpotent weighted orbital integral variance}
J_\mathrm{M}^\mathrm{G}(\nu,f_S\circ\mathrm{ad}(x)) 
&=& 
\sum_{\mathrm{Q}\in\mathcal{F}(\mathrm{M})}
 J_\mathrm{M}^{\mathrm{M}_\mathrm{Q}}(\nu,f_{S,\mathrm{Q},x}). 
\end{eqnarray}

\proof
This is a special case of \eqref{weighted orbital integral variance}.
\qed\\

Argue by induction. Assume that the constants $a^\mathrm{L}(S,\nu)$ and the identities
\begin{eqnarray}
J_\mathrm{nil}^\mathrm{L}(f) &=& \sum_{\mathrm{M}\in\mathcal{L}^\mathrm{L}} |\mathrm{W}_0^\mathrm{M}| 
|\mathrm{W}_0^\mathrm{L}|^{-1}
\sum_{\nu\in\mathfrak{m}_\mathrm{nil}(\mathbb{Q})_{\mathrm{M},S}} a^\mathrm{M}(S,\nu)J_\mathrm{M}^\mathrm{L} (\nu,f_S)
\end{eqnarray}
are known for every proper Levi subgroup $\mathrm{L}$ of $\mathrm{G}$.

Define the distribution $T^\mathrm{G}$\label{sym85} on $\mathfrak{g}(\mathbb{Q}_S)$ by
\begin{eqnarray}
\forall f_S\in\mathcal{S}(\mathfrak{g}(\mathbb{Q}_S))
\\\nonumber
T^\mathrm{G}(f_S)
&=&
J_\mathrm{nil}^\mathrm{G}\Big(f_S\otimes\bigotimes_{p\notin S}\mathbb{I}_{\mathfrak{g}(\mathbb{Z}_p)}\Big) - \sum_{\substack{\mathrm{M}\in\mathcal{L}\\\mathrm{M} \neq \mathrm{G}}} |\mathrm{W}_0^\mathrm{M}||\mathrm{W}_0^\mathrm{G}|^{-1} \times
\\\nonumber\\\nonumber
&&\quad
\times \sum_{\nu \in\mathfrak{m}_\mathrm{nil} (\mathbb{Q})_{\mathrm{M},S}} a^\mathrm{M}(S,\nu) J_\mathrm{M}^\mathrm{G}(\nu,f_S).
\end{eqnarray}
The distribution $T^\mathrm{G}$ is supported on $\mathfrak{g}_\mathrm{nil}(\mathbb{Q}_S)$. By \eqref{nilpotent truncated kernel variance} and \eqref{nilpotent weighted orbital integral variance}
\begin{eqnarray}
&&
T^\mathrm{G}(f_S\circ\mathrm{ad}(x)) - T^\mathrm{G}(f_S) 
\\\nonumber\\\nonumber
&=& 
\Bigg(\sum_{\mathrm{Q}\in\mathcal{F}} 
|\mathrm{W}_0^{\mathrm{M}_\mathrm{Q}}| 
|\mathrm{W}_0^\mathrm{G}|^{-1} 
J_\mathrm{nil}^{\mathrm{M}_\mathrm{Q}} 
\Big(( f_S\otimes\bigotimes_{p\notin S}\mathbb{I}_{\mathfrak{g}(\mathbb{Z}_p)})_{\mathrm{Q},x} \Big) 
\\\nonumber\\\nonumber
&& 
\quad
-\sum_{\substack{\mathrm{M}\in\mathcal{L}\\\mathrm{M}\neq\mathrm{G}}} 
\sum_{\mathrm{Q}\in\mathcal{F}(\mathrm{M})} 
|\mathrm{W}_0^\mathrm{M}|
|\mathrm{W}_0^\mathrm{G}|^{-1} \sum_{\nu\in\mathfrak{m}_\mathrm{nil} (\mathbb{Q})_{\mathrm{M},S}} 
a^\mathrm{M}(S,\nu) 
J_\mathrm{M}^{\mathrm{M}_\mathrm{Q}} (\nu,f_{S,\mathrm{Q},x}) \Bigg) 
\\\nonumber\\\nonumber
&& 
\qquad
-\Bigg(J_\mathrm{nil}^\mathrm{G}\Big(f_S\otimes\bigotimes_{p\notin S}\mathbb{I}_{\mathfrak{g}(\mathbb{Z}_p)}\Big)
- \sum_{\substack{\mathrm{M}\in\mathcal{L}\\\mathrm{M}\neq\mathrm{G}}} 
|\mathrm{W}_0^\mathrm{M}|
|\mathrm{W}_0^\mathrm{G}|^{-1} \times \\\nonumber\\\nonumber 
&&
\qquad\quad
\times \sum_{\nu\in\mathfrak{m}_\mathrm{nil} (\mathbb{Q})_{\mathrm{M},S}} a^\mathrm{M}(S,\nu) 
J_\mathrm{M}^\mathrm{G} (\nu,f_S) \Bigg) 
\\\nonumber\\\nonumber 
&=& 
\sum_{\substack{Q\in\mathcal{F}\\\mathrm{Q}\neq\mathrm{G}}} 
|\mathrm{W}_0^{\mathrm{M}_\mathrm{Q}}|
|\mathrm{W}_0^\mathrm{G}|^{-1} 
\Bigg( J_\mathrm{nil}^{\mathrm{M}_\mathrm{Q}} \Big(( f_S\otimes\bigotimes_{p\notin S}\mathbb{I}_{\mathfrak{g}(\mathbb{Z}_p)})_{\mathrm{Q},x}\Big) 
\\\nonumber\\\nonumber 
&& 
\quad
- \sum_{\mathrm{M} \in \mathcal{L}^{\mathrm{M}_\mathrm{Q}}} 
|\mathrm{W}_0^\mathrm{M}|
|\mathrm{W}_0^{\mathrm{M}_\mathrm{Q}}|^{-1}  \sum_{\nu\in\mathfrak{m}_\mathrm{nil} (\mathbb{Q})_{\mathrm{M},S}} 
a^\mathrm{M}(S,\nu) 
J_\mathrm{M}^{\mathrm{M}_\mathrm{Q}}(\nu,f_{S,\mathrm{Q},x}) \Bigg)
\end{eqnarray}
which vanishes by the inductive hypothesis that $T^\mathrm{L}$ is invariant applied to the Levi subgroup $\mathrm{M}_\mathrm{Q}$ with $\mathrm{Q}$ a proper parabolic subgroup of $\mathrm{G}$. Therefore $T^\mathrm{G}$ is invariant under the action of $\mathrm{G}(\mathbb{Q}_S)^1$.

Construct the constants $a^\mathrm{G}(S,\nu)$ subject to the new identity
\begin{eqnarray}
\label{t distribution identity}
T^\mathrm{G}(f_S) &=& \sum_{\nu\in\mathfrak{g}_\mathrm{nil}(\mathbb{Q})_{\mathrm{G},S}} a^\mathrm{G}(S,\nu)J_\mathrm{G}^\mathrm{G}(\nu,f_S).
\end{eqnarray}
Stratify the nilpotent locus $\mathfrak{g}_\mathrm{nil}(\mathbb{Q}_S)$ equivariantly by codimension. More precisely define for each natural number $d$ the open set\label{sym86}
\begin{eqnarray}
\mathfrak{g}_{\mathrm{nil},d}(\mathbb{Q}_S) &=& \bigcup_{\substack{\nu\in\mathfrak{g}_\mathrm{nil}(\mathbb{Q})_{\mathrm{G},S}\\\mathrm{codim}(\nu)\leq d}} \nu(\mathbb{Q}_S).
\end{eqnarray}
Denote by $T_d^\mathrm{G}$\label{sym87} the distribution obtained by restricting $T^\mathrm{G}$ to $\mathfrak{g}_{\mathrm{nil},d}(\mathbb{Q}_S)$.

The open set $\mathfrak{g}_{\mathrm{nil},0}(\mathbb{Q}_S)$ is the regular nilpotent orbit. Since $T_0^\mathrm{G}$ is invariant it is equal to a multiple of $J_\mathrm{G}^\mathrm{G}(\nu_\mathrm{reg},~)$, where $\nu_\mathrm{reg}$ denotes the regular nilpotent orbit. Define $a^\mathrm{G}(S,\nu_\mathrm{reg})$ to be the constant of proportionality.

The other constants $a^\mathrm{G}(S,\nu)$ are constructed by induction on the codimension $d$ which ranges among $0,1,2,\dots,\dim(\mathfrak{g}_\mathrm{nil})$. Let $T^{\mathrm{G},d}$ be the distribution on the complement of $\mathfrak{g} _{\mathrm{nil},d-1}(\mathbb{Q}_S)$ in $\mathfrak{g}_\mathrm{nil}(\mathbb{Q}_S)$ defined by
\begin{eqnarray}
T^{\mathrm{G},d}(f_S) &=& T^\mathrm{G} (f_S) -\sum_{\substack{\nu\in\mathfrak{g}_\mathrm{nil}(\mathbb{Q})_{\mathrm{G},S}\\\mathrm{codim}(\nu)<d}} a^\mathrm{G}(S,\nu) J_\mathrm{G}^\mathrm{G}(\nu,f_S).
\end{eqnarray}
Denote by $T_d^{\mathrm{G},d}$ its restriction to the complement of $\mathfrak{g}_{\mathrm{nil},d-1}(\mathbb{Q}_S)$ in $\mathfrak{g}_{\mathrm{nil},d}(\mathbb{Q}_S)$.

Since the complement of $\mathfrak{g}_{\mathrm{nil},d-1}(\mathbb{Q}_S)$ in $\mathfrak{g}_{\mathrm{nil},d}(\mathbb{Q}_S)$ has an open partition by the nilpotent orbits of codimension $d$, and the distribution $T_d^{\mathrm{G},d}$ is invariant, there exist constants $a^\mathrm{G}(S,\nu)$, one for each $\nu$ of codimension $d$, such that
\begin{eqnarray}
T_d^{\mathrm{G},d}(f_S) &=& \sum_{\substack{\nu\in\mathfrak{g}_\mathrm{nil}(\mathbb{Q})_{\mathrm{G},S}\\\mathrm{codim}(\nu)=d}} a^\mathrm{G}(S,\nu)J_\mathrm{G}^\mathrm{G}(\nu,f_S).
\end{eqnarray}

The constants $a^\mathrm{G}(S,\nu)$ are required to satisfy \eqref{t distribution identity}, which is equivalent to
\begin{eqnarray}
T^\mathrm{G}(f_S) &=& T_0^{\mathrm{G},0}(f_S)+ T_1^{\mathrm{G},1}(f_S)+ T_2^{\mathrm{G},2}(f_S)+ \dots+ T_{\dim(\mathfrak{g}_\mathrm{nil})}^{\mathrm{G},\dim(\mathfrak{g}_\mathrm{nil})}(f_S).
\end{eqnarray}
Let $v\in S$ be a place of $\mathbb{Q}_S$, let $\nu$ be a nilpotent orbit, define the function $f_{S,\nu,v}^\epsilon$ by the same formula \eqref{f epsilon} as for $f_{\nu,v}^\epsilon$. Let $\nu_d$ denote the complement of $\mathfrak{g}_{\mathrm{nil},d-1}$ in $\mathfrak{g}_{\mathrm{nil},d}$, the union of the nilpotent orbits $\nu$ of codimension $d$. Then the expression $f_{S,\nu_d,v}^\epsilon$ makes sense, and
\begin{eqnarray}
T_d^{\mathrm{G},d}(f_S) &=& \lim_{\epsilon\rightarrow0} T^\mathrm{G}(f_{S,\nu_d,v}^\epsilon) - \lim_{\epsilon\rightarrow0} T^\mathrm{G}(f_{S,\nu_{d+1},v}^\epsilon).
\end{eqnarray}
Therefore
\begin{eqnarray}
&& T_0^{\mathrm{G},0}(f_S)+ T_1^{\mathrm{G},1}(f_S)+ T_2^{\mathrm{G},2}(f_S)+ \dots+ T_{\dim(\mathfrak{g}_\mathrm{nil})}^{\mathrm{G},\dim(\mathfrak{g}_\mathrm{nil})}(f_S) 
\\\nonumber\\\nonumber
&=& 
\lim_{\epsilon\rightarrow0} \bigg(  \underbrace{T^\mathrm{G}(f_{S,\nu_0,v}^\epsilon}_{\textrm{this is } T^\mathrm{G}(f_{S})} ) -  T^\mathrm{G}(f_{S,\nu_1,v}^\epsilon) + T^\mathrm{G}(f_{S,\nu_1,v}^\epsilon) -  T^\mathrm{G}(f_{S,\nu_2,v}^\epsilon) + T^\mathrm{G}(f_{S,\nu_2,v}^\epsilon) - \dots 
\\\nonumber\\\nonumber 
&&  
\quad
\dots -  T^\mathrm{G}(f_{S,\nu_{\dim(\mathfrak{g}_\mathrm{nil})},v}^\epsilon) + T^\mathrm{G}(f_{S,\nu_{\dim(\mathfrak{g}_\mathrm{nil})},v}^\epsilon) -  T^\mathrm{G}( \underbrace{f_{S,\nu_{\dim(\mathfrak{g}_\mathrm{nil})+1},v}^\epsilon}_{\textrm{this is void}} ) \bigg) 
\\\nonumber\\\nonumber
&=& T^\mathrm{G}(f_S).
\end{eqnarray}
\qed

\paragraph{Remark}\setcounter{equation}{0}
If $\nu$ is the orbit consisting of the origin, the coefficient $a^\mathrm{M}(S,\nu)$ is independent of $S$ and equal to the corresponding Tamagawa number
\begin{eqnarray}
a^\mathrm{M}(S,0) &=& \mathrm{Vol}(\mathrm{M}(\mathbb{Q})\backslash \mathrm{M}(\mathbb{A})^1).
\end{eqnarray}

\subsection{General orbits}

\paragraph{Definition}\setcounter{equation}{0}
\label{semisimple descent first notations}
Let $\mathfrak{o}$ be a $\sim$ equivalence class on $\mathfrak{g}(\mathbb{Q})$. Let $\mathrm{P}_1$ be a parabolic subgroup of $\mathrm{G}$, let $\mathrm{M}_1$ be the Levi component of $\mathrm{P}_1$ which is standard. Fix a semisimple element $\Sigma$\label{sym88} in $\mathfrak{o}$ such that $\Sigma$ is contained in the Levi subalgebra $\mathfrak{m}_1$, but not in any proper parabolic subalgebra of $\mathfrak{p}_1$. The group $\mathrm{P}_{1,\Sigma}^0$\label{sym89} is a minimal parabolic subgroup of $\mathrm{G}_\Sigma^0$ with minimal Levi component $\mathrm{M}_{1,\Sigma}^0$. Denote by $\mathcal{F}^\Sigma$ the set of parabolic subgroups of $\mathrm{G}_\Sigma^0$ containing $\mathrm{M}_{1,\Sigma}^0$. A parabolic subgroup $\mathrm{Q}$ in $\mathcal{F}^\Sigma$ is said to be standard if $\mathrm{Q}$ contains $\mathrm{P}_{1,\Sigma}^0$. Fix a maximal compact subgroup $\mathrm{K}_\Sigma$ of $\mathrm{G}_\Sigma^0(\mathbb{A})$ that is admissible with respect to $\mathrm{M}_{1,\Sigma}^0$ such that for each parabolic subgroup $\mathrm{Q}$ in $\mathcal{F}^\Sigma$\label{sym91} there is the associated function 
\begin{eqnarray}
H_\mathrm{Q}:\quad \mathrm{G}_\Sigma^0(\mathbb{A}) 
&\rightarrow& \mathfrak{a}_\mathrm{Q}.
\end{eqnarray}
There is a unique point $T_{\Sigma,1}$\label{sym90} in $\mathfrak{a}_1$ modulo $\mathfrak{a}_{\mathrm{G}_\Sigma^0}$ defined in the same manner as $T_0$ in $\mathfrak{a}_0^\mathrm{G}$ that satisfies an identity analogous to \eqref{point t0}. Let $\mathrm{L}$ be a Levi subgroup of $\mathrm{G}_\Sigma^0$ containing $\mathrm{M}_{1,\Sigma}^0$, denote by $\mathrm{W}_1^\mathrm{L}$ the Weyl group of $\mathrm{L}$ with respect to the split torus $\mathrm{A}_1$.

Let $\pi_\Sigma$\label{sym92} be the surjection from $\mathcal{F}(\mathrm{M}_1)$ onto $\mathcal{F}^\Sigma$ defined by
\begin{eqnarray}
\forall \mathrm{P}\in\mathcal{F}(\mathrm{M}_1)
&&
\pi_\Sigma(\mathrm{P})=\mathrm{P}_\Sigma^0.
\end{eqnarray}
Let $\mathrm{Q}$ be a parabolic subgroup in $\mathcal{F}^\Sigma$. Let $\mathcal{F}_\mathrm{Q}(\mathrm{M}_1)$ be the inverse image $\pi_\Sigma^{-1}(\mathrm{Q})$. Define subsets $\mathring{\mathcal{F}}_\mathrm{Q}(\mathrm{M}_1)$ and $\bar{\mathcal{F}}_\mathrm{Q}(\mathrm{M}_1)$ of $\mathcal{F}(\mathrm{M}_1)$ by
\begin{eqnarray}
\mathring{\mathcal{F}}_\mathrm{Q}(\mathrm{M}_1)
&=& 
\bigg\{\mathrm{P}\in\mathcal{F}(\mathrm{M}_1):~\pi_\Sigma(\mathrm{P})=\mathrm{Q},
~\mathfrak{a}_\mathrm{P}=\mathfrak{a}_\mathrm{Q}\bigg\}
\\\nonumber
\bar{\mathcal{F}}_\mathrm{Q}(\mathrm{M}_1)
&=&
\bigg\{\mathrm{P}\in\mathcal{F}(\mathrm{M}_1):~\pi_\Sigma(\mathrm{P})\supset\mathrm{Q}\bigg\}.
\end{eqnarray}
There are inclusions
\begin{eqnarray}
\forall \mathrm{Q}\in\mathcal{F}^\Sigma
&&
\mathring{\mathcal{F}}_\mathrm{Q}(\mathrm{M}_1)
\subset\mathcal{F}_\mathrm{Q}(\mathrm{M}_1)
\subset\bar{\mathcal{F}}_\mathrm{Q}(\mathrm{M}_1)
\subset\mathcal{F}(\mathrm{M}_1).
\end{eqnarray}

\paragraph{Definition}\setcounter{equation}{0}
Let $\mathrm{P}_1$, $\mathrm{M}_1$, $\mathfrak{o}$, $\Sigma$ be as in \ref{semisimple descent first notations}. Let $\mathrm{Q}$ be a parabolic subgroup in $\mathcal{F}^\Sigma$. Let $\mathcal{Y}$\label{sym93} be a collection of points
\begin{eqnarray}
\mathcal{Y}
&=&
\bigg\{
Y_\mathrm{P}\in\mathfrak{a}_0:~
\mathrm{P}\in\mathring{\mathcal{F}}_\mathrm{Q}(\mathrm{M}_1)
\bigg\}
\end{eqnarray}
satisfying the compatibility conditions defining a $(\mathrm{G},\mathrm{M}_1)$-orthogonal set in \ref{(g,m)-orthogonal set equation}, namely that for each pair of adjacent parabolic subgroups $\mathrm{P}$ and $\mathrm{P}'$ in $\mathring{\mathcal{F}}_\mathrm{Q}(\mathrm{M}_1)$ the difference between $Y_\mathrm{P}$ and $Y_{\mathrm{P}'}$ is orthogonal to the common wall of the positive chambers defined by $\mathrm{P}$ and $\mathrm{P}'$. 

The collection $\mathcal{Y}$ has a unique extension from $\mathring{\mathcal{F}}_\mathrm{Q}(\mathrm{M}_1)$ to $\bar{\mathcal{F}}_\mathrm{Q}(\mathrm{M}_1)$. Let $\mathrm{Q}'$ be a parabolic subgroup in $\mathcal{F}^\Sigma$ containing $\mathrm{Q}$, denote by $\mathcal{Y}_{\mathrm{Q}'}$ the collection of points
\begin{eqnarray}
\label{(g,m)-orthogonal set semisimple descent generalization}
\mathcal{Y}_{\mathrm{Q}'}
&=&
\bigg\{Y_\mathrm{P}:~ \mathrm{P}\in\mathcal{F}_{\mathrm{Q}'}(\mathrm{M}_1)\bigg\}.
\end{eqnarray}

Let $\mathrm{Q}_3$ be a parabolic subgroup in $\mathcal{F}^\Sigma$. Define the \emph{gamma$'$ function} $\Gamma'_{\mathrm{Q}_3}(~,\mathcal{Y}_{\mathrm{Q}_3})$\label{sym94} on $\mathfrak{a}_1$ by
\begin{eqnarray}
\label{gamma' function definition var}
&&
\forall H\in\mathfrak{a}_1
\\\nonumber
&&
\Gamma_{\mathrm{Q}_3}'(H,\mathcal{Y}_{\mathrm{Q}_3}) 
=
\sum_{\substack{\mathrm{Q}_4\in\mathcal{F}^\Sigma\\\mathrm{Q}_4\supset\mathrm{Q}_3}} 
\tau_3^4(H)
\Bigg(\sum_{\mathrm{P}\in\mathcal{F}_{\mathrm{Q}_4}(\mathrm{M}_1)} 
(-1)^{\dim(\mathrm{A}_\mathrm{P} / \mathrm{A}_\mathrm{G})} 
\hat{\tau}_\mathrm{P}(H-Y_\mathrm{P})  \Bigg). 
\end{eqnarray}

\paragraph{Remark}\setcounter{equation}{0}
The function $\Gamma'_{\mathrm{Q}_3}(~,\mathcal{Y}_{\mathrm{Q}_3})$ factorizes through the projection from $\mathfrak{a}_1$ onto $\mathfrak{a}_3^\mathrm{G}$ and depends continuously on $\mathcal{Y}$. For each parabolic subgroup $\mathrm{Q}_2$ in $\mathcal{F}^\Sigma$
\begin{eqnarray}
\label{gamma'function identity var}
&&
\sum_{\mathrm{P}\in\mathcal{F}_{\mathrm{Q}_2}(\mathrm{M}_1)} 
(-1)^{\dim(\mathrm{A}_\mathrm{P}/\mathrm{A}_\mathrm{G})} 
\hat{\tau}_\mathrm{P}(H-Y_\mathrm{P}) 
= 
\sum_{\substack{\mathrm{Q}_3\in\mathcal{F}^\Sigma\\\mathrm{Q_3}\supset\mathrm{Q}_2}} (-1)^{\dim(\mathrm{A}_2/\mathrm{A}_3)} 
\hat{\tau}_2^3(H) 
\Gamma_{\mathrm{Q}_3}'(H,\mathcal{Y}_{\mathrm{Q}_3}).
\end{eqnarray}
See \S4 of \cite{Arthur_Orb}.

\paragraph{Lemma}\setcounter{equation}{0}
\emph{Let $\mathrm{Q}$ be a parabolic subgroup in $\mathcal{F}^\Sigma$, let $\mathcal{Y}$ be a collection of points as in (\ref{(g,m)-orthogonal set semisimple descent generalization}). The function $\Gamma'_\mathrm{Q}(~,\mathcal{Y}_\mathrm{Q})$ is compactly supported as a function on $\mathfrak{a}_\mathrm{Q}^\mathrm{G}$.}

\emph{Let $(c_\mathrm{P})$ be the ($\mathrm{G},\mathrm{M}_1$)-family associated with a ($\mathrm{G},\mathrm{M}_1$)-orthogonal set that contains $\mathcal{Y}_\mathrm{Q}$ as in (\ref{(g,m)-orthogonal set to (g,m)-family}). Let $c_\mathrm{P}'$ be the functions intervening in (\ref{weight constant'}). Let $\Gamma'_\mathrm{Q}\hat{~}(~,\mathcal{Y}_\mathrm{Q})$ denote the Fourier transform of $\Gamma'_\mathrm{Q}(~,\mathcal{Y}_\mathrm{Q})$. Then}
\begin{eqnarray}
\label{gamma' function var fourier transform}
\forall \lambda\in i\mathfrak{a}_\mathrm{Q}^{\mathrm{G},*}
&&
\Gamma'_\mathrm{Q}\hat{~}(\lambda,\mathcal{Y}_\mathrm{Q})
=
\sum_{\mathrm{P}\in\mathring{\mathcal{F}}_\mathrm{Q}(\mathrm{M}_1)} c_\mathrm{P}'(\lambda).
\end{eqnarray}

\proof
See Lemma 4.1 of \cite{Arthur_Orb}.
\qed

\paragraph{Remark}
The main ingredient of the proof is the following generalization of \eqref{langlands combinatorial lemma}:\\
Let $\mathrm{Q}$ be a parabolic subgroup in $\mathcal{F}^\Sigma$, let $\mathrm{P}$ be a parabolic subgroup in $\bar{\mathcal{F}}_\mathrm{Q}(\mathrm{M}_1)$, then
\begin{eqnarray}
&&
\sum_{\substack{\mathrm{P}'\in\bar{\mathcal{F}}_\mathrm{Q}(\mathrm{M}_1)\\
\mathrm{P}'\subset\mathrm{P}}} 
(-1)^{\dim(\mathrm{A}_{\mathrm{P}'}/\mathrm{A}_\mathrm{P})} 
\tau_\mathrm{Q}^{\pi_\Sigma(\mathrm{P}')} (H) \hat{\tau}_{\mathrm{P}'}^\mathrm{P} (H) 
= 
\left\{ \begin{array}{ll} 1 & \textrm{if } \mathrm{P}\in\mathcal{F}_\mathrm{Q}(\mathrm{M}_1)
\textrm{ and } H\in\mathfrak{a}_\mathrm{P},
\\\\ 0 & \textrm{otherwise}. 
\end{array}\right.
\end{eqnarray}
This is Lemma 4.2 of \cite{Arthur_Orb}.

\paragraph{Lemma}\setcounter{equation}{0} (Global semisimple descent)\\
\emph{Let $\mathrm{P}_1$, $\mathrm{M}_1$, $\mathfrak{o}$, $\Sigma$ be as in \ref{semisimple descent first notations}, let $f$ be a Schwartz function on $\mathfrak{g}(\mathbb{A})$, then
\begin{eqnarray}
\label{global semisimple descent}
&&
J_\mathfrak{o}(f) 
= 
|\pi_0(\mathrm{G}_\Sigma)|^{-1}
\int_{\mathrm{G}_\Sigma^0(\mathbb{A})\backslash\mathrm{G}(\mathbb{A})} 
\Bigg( \sum_{\mathrm{Q}\in\mathcal{F}^\Sigma} 
|\mathrm{W}_1^{\mathrm{M}_\mathrm{Q}}| 
|\mathrm{W}_1^{\mathrm{G}_\Sigma^0}|^{-1}  
J_\mathrm{nil}^{\mathrm{M}_\mathrm{Q}}(\Phi_{\mathrm{Q},x}^{T_0-T_{\Sigma,1}}) 
\Bigg)~\mathrm{d}x
\end{eqnarray}
where for a truncation parameter $T$ in $\mathfrak{a}_0$ the function $\Phi_{\mathrm{Q},x}^T$\label{sym95} on $\mathfrak{m}_\mathrm{Q}(\mathbb{A})$ is defined by
\begin{eqnarray}
&&
\forall X\in\mathfrak{m}_\mathbb{Q}(\mathbb{A})
\\\nonumber
&&
\Phi_{\mathrm{Q},x}^T(X) 
=\int_{\mathrm{K}_\Sigma}
\int_{\mathfrak{n}_\mathrm{Q}(\mathbb{A})}
f\Big( \big(\Sigma+(X+N) \cdot \mathrm{ad}(k)\big)
\cdot\mathrm{ad}(x) \Big)
~v_\mathrm{Q}'(kx,T) \mathrm{d}N \mathrm{d}k
\end{eqnarray}
where the weight factor $v_\mathrm{Q}'$ is defined as
\begin{eqnarray}
v_\mathrm{Q}'(kx,T)
 &=& 
\int_{\mathfrak{a}_\mathrm{Q}^\mathrm{G}} \Gamma_\mathrm{Q}'(H,\mathcal{Y}_\mathrm{Q}^T(k,x))
 ~\mathrm{d}H
\end{eqnarray}
where $\mathcal{Y}^T(k,x)$ is the collection of points defined by
\begin{eqnarray}
\forall \mathrm{P}\in\mathcal{F}_\mathrm{Q}(\mathrm{M}_1)
&&
Y_\mathrm{P}^T(k,x) = -H_\mathrm{P}(kx)-T_\Sigma+T
\label{(g,m) orthogonal set var}
\end{eqnarray} 
where $T_\Sigma$\label{sym96} is a truncation parameter in $\mathfrak{a}_1$ such that $T_\Sigma-T_{\Sigma,1}$ is the projection of $T-T_0$.}

\proof
This argument follows the proof of Lemma 6.2 of \cite{Arthur_Orb}.

Let $T$ be a truncation parameter in $\mathfrak{a}_0$. Define a second truncated kernel function $j_\mathfrak{o}^T(~,f)$\label{sym97} on $\mathrm{G}(\mathbb{Q})\backslash\mathrm{G}(\mathbb{A})^1$ by
\begin{eqnarray}
&&
\forall x\in\mathrm{G}(\mathbb{Q})\backslash\mathrm{G}(\mathbb{A})^1
\\\nonumber
&&
j_\mathfrak{o}^T(x,f)=
\sum_{\substack{\mathrm{P}\in\mathcal{F}\\\mathrm{standard}}}
(-1)^{\dim(\mathrm{A}_\mathrm{P}/\mathrm{A}_\mathrm{G})}
\sum_{\delta\in\mathrm{P}(\mathbb{Q})\backslash\mathrm{G}(\mathbb{Q})}
\hat{\tau}_\mathrm{P}(H_0(\delta x)-T) J_{\mathrm{P},\mathfrak{o}}(\delta x)
\end{eqnarray}
where
\begin{eqnarray}
&&
J_{\mathrm{P},\mathfrak{o}}(x,f)=\sum_{X\in\mathfrak{m}_\mathrm{P}(\mathbb{Q})\cap\mathfrak{o}}
~ \sum_{\eta\in\mathrm{N}_{\mathrm{P},X_{\mathrm{ss}}} (\mathbb{Q})\backslash \mathrm{N}_\mathrm{P}(\mathbb{Q})} \int_{\mathfrak{n}_{\mathrm{P},X_\mathrm{ss}}(\mathbb{A})} f((X+N)\cdot\mathrm{ad}(\eta x))~\mathrm{d}N.
\end{eqnarray}

\paragraph*{Lemma}
\emph{The function $j_\mathfrak{o}^T(~,f)$ is integrable on $\mathrm{G}(\mathbb{Q})\backslash\mathrm{G}(\mathbb{A})^1$, and}
\begin{eqnarray}
\label{global semisimple descent proof first lemma}
J_\mathfrak{o}^T(f) &=& \int_{\mathrm{G}(\mathbb{Q})\backslash \mathrm{G}(\mathbb{A})^1} j_\mathfrak{o}^T(x,f) ~\mathrm{d}x.
\end{eqnarray}

\proof
Following the proof of \eqref{truncated kernel integrability}, by the combinatorial lemma of Langlands \eqref{langlands combinatorial lemma}, it suffices to prove that
\begin{eqnarray}
&&
\sum_{\substack{\mathrm{P}_2,\mathrm{P}_5\in\mathcal{F}\\
\mathrm{standard}\\\mathrm{P}_2\subset\mathrm{P}_5}}
 \int_{\mathrm{P}_2(\mathbb{Q})\backslash\mathrm{G}(\mathbb{A})^1} 
F^{\mathrm{P}_2}(x,T)
\sigma_2^5(H_0(x)-T) \times
\\\nonumber\\\nonumber
&&
\times \bigg| \sum_{\substack{\mathrm{P}_4\in\mathcal{F}\\
\mathrm{P}_2\subset\mathrm{P}_4\subset\mathrm{P}_5}} 
(-1)^{\dim(\mathrm{A}_4 / \mathrm{A}_\mathrm{G})} 
J_{\mathrm{P}_4,\mathfrak{o}}(x,f) \bigg| ~\mathrm{d}x
\end{eqnarray}
is finite. Decompose the sum defining $J_{\mathrm{P}_4,\mathfrak{o}}(x,f)$ over the set of parabolic subgroups contained in $\mathrm{P}_4$,
\begin{eqnarray}
&&
J_{\mathrm{P}_4,\mathfrak{o}}(x,f) 
\\\nonumber\\\nonumber
&=& 
\sum_{\substack{\mathrm{P}_3\in\mathcal{F}\\\mathrm{P}_2\subset\mathrm{P}_3\subset\mathrm{P}_4}}  
\sum_{X\in\mathfrak{m}_2^3 (\mathbb{Q})'\cap\mathfrak{o}}  
~\sum_{Y\in\mathfrak{n}_3^4(\mathbb{Q})} 
\\\nonumber\\\nonumber
&&\quad
\times\sum_{\eta\in\mathrm{N}_{4,(X+Y)_{\mathrm{ss}}} (\mathbb{Q})
\backslash \mathrm{N}_4(\mathbb{Q})} 
~\int_{\mathfrak{n}_{4,(X+Y)_\mathrm{ss}}(\mathbb{A})} 
f\big((X+Y+N)\cdot\mathrm{ad}(\eta x)\big)~\mathrm{d}N
\\\nonumber\\
\label{global semisimple descent proof first equality}
&=& 
\sum_{\substack{\mathrm{P}_3\in\mathcal{F}\\\mathrm{P}_2\subset\mathrm{P}_3\subset\mathrm{P}_4}}  
\sum_{X\in\mathfrak{m}_2^3 (\mathbb{Q})'\cap\mathfrak{o}}  
~\sum_{Z\in\mathfrak{n}_{3,X_\mathrm{ss}}^4(\mathbb{Q})} 
~\sum_{\delta\in\mathrm{N}_{3,X_\mathrm{ss}}^4(\mathbb{Q})
\backslash\mathrm{N}_3^4(\mathbb{Q})} 
\\\nonumber\\\nonumber
&&  \quad 
\times \sum_{\eta\in\mathrm{N}_{4,(X+Z)_{\mathrm{ss}}\cdot\mathrm{ad}(\delta)} (\mathbb{Q})
\backslash \mathrm{N}_4(\mathbb{Q})} 
\\\nonumber\\\nonumber 
&&  \qquad 
\times \int_{\mathfrak{n}_{4,(X+Z)_\mathrm{ss} \cdot \mathrm{ad}(\delta)} (\mathbb{A})} 
f\Big(\big(X+Z+N\cdot\mathrm{ad}(\delta^{-1})\big)\cdot\mathrm{ad}(\delta \eta x)\Big)~\mathrm{d}N
\\\nonumber\\\nonumber 
&=& 
\sum_{\substack{\mathrm{P}_3\in\mathcal{F}\\\mathrm{P}_2\subset\mathrm{P}_3\subset\mathrm{P}_4}}  
\sum_{X\in\mathfrak{m}_2^3 (\mathbb{Q})'\cap\mathfrak{o}}  
~\sum_{Z\in\mathfrak{n}_{3,X_\mathrm{ss}}^4(\mathbb{Q})} 
~\sum_{\eta\in\mathrm{N}_{3,X_{\mathrm{ss}}} (\mathbb{Q})\backslash \mathrm{N}_3(\mathbb{Q})} 
\\\nonumber\\\nonumber
&& \quad 
\times \int_{\mathfrak{n}_{4,(X+Z)_\mathrm{ss}}(\mathbb{A})} 
f\big((X+Z+N)\cdot\mathrm{ad}(\eta x)\big)~\mathrm{d}N
\\\nonumber\\
\label{global semisimple descent proof second equality} 
&=& 
\sum_{\substack{\mathrm{P}_3\in\mathcal{F}\\\mathrm{P}_2\subset\mathrm{P}_3\subset\mathrm{P}_4}}  
\sum_{X\in\mathfrak{m}_2^3 (\mathbb{Q})'\cap\mathfrak{o}}  
~ \sum_{\eta\in\mathrm{N}_{3,X_{\mathrm{ss}}} (\mathbb{Q})\backslash \mathrm{N}_3(\mathbb{Q})} 
\\\nonumber\\\nonumber 
&& \quad  
\times \sum_{\overline{Z} \in\overline{\mathfrak{n}}_{3,X_\mathrm{ss}} ^4(\mathbb{Q})} 
\int_{\mathfrak{n}_{3,X_\mathrm{ss}}(\mathbb{A})} 
f\big((X+N)\cdot\mathrm{ad}(\eta x)\big)
\cdot\psi(\langle N,\overline{Z}\rangle) ~\mathrm{d}N
\end{eqnarray}
where \eqref{global semisimple descent proof first equality} follows from Corollary 2.4 of \cite{Chaud} which states that
\begin{eqnarray}
 \sum_{Z\in\mathfrak{n}_{X_\mathrm{ss}}(\mathbb{Q})} ~ \sum_{\delta\in\mathrm{N}_{X_\mathrm{ss}}(\mathbb{Q}) \backslash \mathrm{N}(\mathbb{Q})} f((X+Z)\cdot\mathrm{ad}(\delta))
 &=& 
\sum_{Y\in\mathfrak{n}(\mathbb{Q})} f(X+Y)
\label{Chaudouard 2.4}
\end{eqnarray}
and \eqref{global semisimple descent proof second equality} follows from the Poisson summation formula. The right hand side of \eqref{global semisimple descent proof second equality} is independant of $\mathrm{P}_4$, hence the alternating sum in $\mathrm{P}_4$ cancels by the inclusion-exclusion principle, therefore the same estimates in the proof of \eqref{truncated kernel integrability} implies the integrability of $j_\mathfrak{o}^T(~,f)$.

For the integral representation of $J_\mathfrak{o}^T(f)$, consider the $(\mathrm{P}_2,\mathrm{P}_5)$ summand of the integral of $j_\mathrm{o}^T(~,f)$ over $\mathrm{G}(\mathbb{Q})\backslash\mathrm{G}(\mathbb{A})^1$:
\begin{eqnarray}
\qquad\qquad
&& 
\int_{\mathrm{P}_2(\mathbb{Q})\backslash\mathrm{G}(\mathbb{A})^1} 
F^{\mathrm{P}_2}(x,T)
\sigma_2^5(H_0(x)-T)  \times
\\\nonumber\\\nonumber
&&\quad
\times \sum_{\substack{\mathrm{P}_4\in\mathcal{F}\\\mathrm{P}_2\subset\mathrm{P}_4\subset\mathrm{P}_5}}
(-1)^{\dim(\mathrm{A}_4 / \mathrm{A}_\mathrm{G})}
 J_{\mathrm{P}_4,\mathfrak{o}}(x,f)  ~\mathrm{d}x
\\\nonumber\\\nonumber
&=& \int_{\mathrm{M}_2(\mathbb{Q})\mathrm{N}_2(\mathbb{A})
\backslash\mathrm{G}(\mathbb{A})^1}
 F^{\mathrm{P}_2}(x,T)
\sigma_2^5(H_0(x)-T)  \times
\\\nonumber\\\nonumber
&&\quad
\times \sum_{\substack{\mathrm{P}_4\in\mathcal{F}\\\mathrm{P}_2\subset\mathrm{P}_4\subset\mathrm{P}_5}}
(-1)^{\dim(\mathrm{A}_4 / \mathrm{A}_\mathrm{G})} 
\sum_{X\in\mathfrak{m}_4(\mathbb{Q})\cap\mathfrak{o}} 
\Bigg(\int_{\mathrm{N}_2(\mathbb{Q}) 
\backslash \mathrm{N}_2(\mathbb{A})}  
\sum_{\eta\in\mathrm{N}_{4,X_{\mathrm{ss}}} (\mathbb{Q})
\backslash \mathrm{N}_4(\mathbb{Q})} 
\\\nonumber\\\nonumber 
&& \qquad 
\times \bigg( 
\int_{\mathfrak{n}_{4,X_\mathrm{ss}}(\mathbb{A})} 
f((X+N)\cdot\mathrm{ad}(\eta n_2x))
~\mathrm{d}N\bigg)\mathrm{d}n_2\Bigg)\mathrm{d}x
\\\nonumber\\
\label{global semisimple descent proof third equality}
&=& \int_{\mathrm{M}_2(\mathbb{Q})\mathrm{N}_2(\mathbb{A})
\backslash\mathrm{G}(\mathbb{A})^1}
 F^{\mathrm{P}_2}(x,T)
\sigma_2^5(H_0(x)-T)  \times
\\\nonumber\\\nonumber
&&\quad
\times
\sum_{\substack{\mathrm{P}_4\in\mathcal{F}\\\mathrm{P}_2\subset\mathrm{P}_4\subset\mathrm{P}_5}}
(-1)^{\dim(\mathrm{A}_4 / \mathrm{A}_\mathrm{G})} 
\sum_{X\in\mathfrak{m}_4(\mathbb{Q})\cap\mathfrak{o}} 
\Bigg( 
\int_{\mathrm{N}_2(\mathbb{Q}) \backslash \mathrm{N}_2(\mathbb{A})}  
\int_{\mathrm{N}_{4,X_{\mathrm{ss}}} (\mathbb{Q})\backslash \mathrm{N}_4(\mathbb{A})} 
\\\nonumber\\\nonumber 
&& \qquad 
\times \bigg(
\int_{\mathfrak{n}_{4,X_\mathrm{ss}}(\mathbb{A})} 
f((X+N)\cdot\mathrm{ad}(nn_2x))~\mathrm{d}N\bigg)\mathrm{d}n\mathrm{d}n_2
\Bigg)\mathrm{d}x
\\\nonumber\\
\label{global semisimple descent proof fourth equality}
&=& \int_{\mathrm{M}_2(\mathbb{Q})\mathrm{N}_2(\mathbb{A})
\backslash\mathrm{G}(\mathbb{A})^1}
 F^{\mathrm{P}_2}(x,T)
\sigma_2^5(H_0(x)-T)  
\sum_{\substack{\mathrm{P}_4\in\mathcal{F}\\\mathrm{P}_2\subset\mathrm{P}_4\subset\mathrm{P}_5}}
(-1)^{\dim(\mathrm{A}_4 / \mathrm{A}_\mathrm{G})} \times
\\\nonumber\\\nonumber 
&& \quad 
\times \Bigg( 
\sum_{X\in\mathfrak{m}_4(\mathbb{Q})\cap\mathfrak{o}} 
\int_{\mathrm{N}_2(\mathbb{Q}) \backslash \mathrm{N}_2(\mathbb{A})}   
\int_{\mathfrak{n}_4(\mathbb{A})} 
f((X+N)\cdot\mathrm{ad}( n_2x)) ~\mathrm{d}N\mathrm{d}n_2
\Bigg)\mathrm{d}x.
\end{eqnarray}
The equality \eqref{global semisimple descent proof third equality} holds since $\mathrm{N}_4(\mathbb{Q})\backslash\mathrm{N}_4(\mathbb{A})$ has volume 1, and the equality \eqref{global semisimple descent proof fourth equality} follows from Corollary 2.5 of \cite{Chaud} which is the integral analogue of \eqref{Chaudouard 2.4}. Reversing the combinatorial manipulations to the right hand side of \eqref{global semisimple descent proof fourth equality} skipping the step \eqref{global semisimple descent proof third equality},
\begin{eqnarray}
&& \int_{\mathrm{M}_2(\mathbb{Q})\mathrm{N}_2(\mathbb{A})\backslash\mathrm{G}(\mathbb{A})^1} 
F^{\mathrm{P}_2}(x,T)
\sigma_2^5(H_0(x)-T)  
\sum_{\substack{\mathrm{P}_4\in\mathcal{F}\\\mathrm{P}_2\subset\mathrm{P}_4\subset\mathrm{P}_5}}
(-1)^{\dim(\mathrm{A}_4 / \mathrm{A}_\mathrm{G})} \times
\\\nonumber\\\nonumber 
&& \quad
\times \Bigg( \sum_{X\in\mathfrak{m}_4(\mathbb{Q})\cap\mathfrak{o}} 
\int_{\mathrm{N}_2(\mathbb{Q}) \backslash \mathrm{N}_2(\mathbb{A})}   
\int_{\mathfrak{n}_4(\mathbb{A})} 
f((X+N)\cdot\mathrm{ad}( n_2x)) ~\mathrm{d}N\mathrm{d}n_2 \Bigg)\mathrm{d}x \
\\\nonumber\\\nonumber
&=& \int_{\mathrm{P}_2(\mathbb{Q})\backslash\mathrm{G}(\mathbb{A})^1} 
F^{\mathrm{P}_2}(x,T)
\sigma_2^5(H_0(x)-T)  \times
\\\nonumber\\\nonumber
&&\quad
\times\sum_{\substack{\mathrm{P}_4\in\mathcal{F}\\\mathrm{P}_2\subset\mathrm{P}_4\subset\mathrm{P}_5}}
(-1)^{\dim(\mathrm{A}_4 / \mathrm{A}_\mathrm{G})} 
K_{\mathrm{P}_4,\mathfrak{o}}(x,f)  ~\mathrm{d}x,
\end{eqnarray}
which is the $(\mathrm{P}_2,\mathrm{P}_5)$ summand of the integral of $k_\mathrm{o}^T(~,f)$ over $\mathrm{G}(\mathbb{Q})\backslash\mathrm{G}(\mathbb{A})^1$, hence $J_\mathfrak{o}^T(f)$.
\qed\\

By \eqref{global semisimple descent proof first lemma} it is enough to consider the function $J_{\mathrm{P},\mathfrak{o}}(~,f)$ instead of $K_{\mathrm{P},\mathfrak{o}}(~,f)$ for a parabolic subgroup $\mathrm{P}$ of $\mathrm{G}$ containing $\mathrm{P}_1$ for the proof of \eqref{global semisimple descent}. Every element $X$ in $\mathfrak{m}_\mathrm{P}(\mathbb{Q})\cap\mathfrak{o}$ is conjugate under the adjoint action of $\mathrm{G}(\mathbb{Q})$ to the sum of $\Sigma$ and $N_\Sigma$ for some $N_\Sigma$ in $\mathfrak{g}_{\Sigma,\mathrm{nil}}(\mathbb{Q})$. More precisely
\begin{eqnarray}
&&
\exists \mathrm{P}_1'\in\mathcal{F},~\mathrm{P}_1'~\textrm{standard},~\mathrm{P}_1'\subset\mathrm{P}
\quad\exists s\in\mathrm{W}(\mathfrak{a}_{\mathrm{P}_1},\mathfrak{a}_{\mathrm{P}_1'})
\quad\exists \mu\in\mathrm{M}_\mathrm{P}^0(\mathbb{Q})
\\\nonumber
&& \quad
\quad\exists N_\Sigma\in\mathfrak{m}_\mathrm{P}(\mathbb{Q})\cdot\mathrm{ad}(w_s)
\cap\mathfrak{g}_{\Sigma,\mathrm{nil}}(\mathbb{Q})
\quad X=(\Sigma+N_\Sigma)\cdot\mathrm{ad}(w_s^{-1}\mu)
\end{eqnarray} 
where 
\begin{itemize}
\item the element $w_s$ in $\mathrm{G}(\mathbb{Q})$ is a representative of $s$;
\item the double coset
\begin{eqnarray}
[s]
&\in&
\mathrm{W}_0^{\mathrm{M}_\mathrm{P}}\Big\backslash
\mathrm{W}(\mathfrak{a}_{\mathrm{P}_1},\mathfrak{a}_{\mathrm{P}_1'})\Big\slash
\mathrm{W}_1^{\mathrm{G}_\Sigma^0}
\end{eqnarray}
is uniquely determined;
\item for a fixed choice of the element $s$, the coset
\begin{eqnarray}
[\mu]
&\in&
\mathrm{M}_\mathrm{P}(\mathbb{Q})\cap w_s\mathrm{G}_\Sigma(\mathbb{Q}) w_s^{-1}
\Big\backslash\mathrm{M}_\mathrm{P}(\mathbb{Q})
\end{eqnarray}
is uniquely determined;
\item for a fixed choice of the representative $w_s$ and the element $\mu$, the element $N_\Sigma$ is uniquely determined.
\end{itemize}
Let $\mathrm{W}(\mathfrak{a}_1;\mathrm{M}_\mathrm{P}^+,\mathrm{G}_\Sigma^{0,+})$ be the subset of the Weyl group of $\mathrm{G}$ defined by
\begin{eqnarray}
&&
\mathrm{W}(\mathfrak{a}_1;\mathrm{M}_\mathrm{P}^+,\mathrm{G}_\Sigma^{0,+})
\\\nonumber
&=&
\bigcup_{\substack{\mathrm{P}_1'\in\mathcal{F}\\\mathrm{standard}\\\mathrm{P}_1'\subset\mathrm{P}}}
\bigg\{
s\in\mathrm{W}(\mathfrak{a}_{\mathrm{P}_1},\mathfrak{a}_{\mathrm{P}_1'})
:~
\forall \alpha\in\Delta_{\mathrm{P}_1'}^\mathrm{P}
~s^{-1}\alpha>0,
~\forall \beta\in\Delta_1^{\mathrm{G}_\Sigma^0}
~s\beta>0
\bigg\}.
\end{eqnarray}
The map
\begin{eqnarray}
&&
\Big(\mathrm{W}(\mathfrak{a}_1;\mathrm{M}_\mathrm{P}^+,\mathrm{G}_\Sigma^{0,+})\Big)\times
\Big(\mathrm{M}_\mathrm{P}(\mathbb{Q})\cap w_s \mathrm{G}_\Sigma^0(\mathbb{Q}) w_s^{-1} 
\Big\backslash \mathrm{M}_\mathrm{P}(\mathbb{Q})\Big)\times
\\\nonumber
&&\quad
\times\Big(\mathfrak{m}_\mathrm{P}(\mathbb{Q}) \cdot \mathrm{ad}(w_s)\cap
\mathfrak{g}_{\Sigma,\mathrm{nil}}(\mathbb{Q})\Big)
\\\nonumber
&\longrightarrow&
\mathfrak{m}_\mathrm{P}(\mathbb{Q})\cap\mathfrak{o}
\end{eqnarray}
defined by
\begin{eqnarray}
(s,\mu,N_\Sigma)
&\mapsto&
(\Sigma+N_\Sigma)\cdot\mathrm{ad}(w_s^{-1}\mu)
\end{eqnarray}
is surjective and the group $\pi_0(\mathrm{G}_\Sigma)$ operates simply transitively on each fiber. Hence
\begin{eqnarray}
\label{global semisimple descent proof fifth equality}
&&
\forall x\in\mathrm{G}(\mathbb{Q})\backslash\mathrm{G}(\mathbb{A})^1
\\\nonumber
&&
J_{\mathrm{P},\mathfrak{o}}(x,f)
\\\nonumber\\\nonumber
&=&
\sum_{s\in \mathrm{W}(\mathfrak{a}_1;\mathrm{M}_\mathrm{P}^+,\mathrm{G}_\Sigma^{0,+})}
\sum_{\mu\in \mathrm{M}_\mathrm{P}(\mathbb{Q})\cap w_s \mathrm{G}_\Sigma^0(\mathbb{Q}) w_s^{-1} 
\big\backslash \mathrm{M}_\mathrm{P}(\mathbb{Q})}
\\\nonumber\\\nonumber
&&\quad
\times \sum_{N_\Sigma\in\mathfrak{m}_\mathrm{P}(\mathbb{Q}) \cdot \mathrm{ad}(w_s)\cap
\mathfrak{g}_{\Sigma,\mathrm{nil}}(\mathbb{Q})}
\\\nonumber\\\nonumber
&&\qquad\times
|\pi_0(\mathrm{G}_\Sigma)|^{-1}
\sum_{\eta\in \mathrm{N}_{\mathrm{P},\Sigma\cdot\mathrm{ad}(w_s^{-1}\mu)}(\mathbb{Q}) 
\backslash \mathrm{N}_\mathrm{P}(\mathbb{Q})}
\int_{\mathfrak{n}_{\mathrm{P},\Sigma\cdot\mathrm{ad}(w_s^{-1}\mu)}(\mathbb{A})} 
\\\nonumber\\\nonumber
&&\qquad\quad\times
f\Big(\big((\Sigma+N_\Sigma)\cdot\mathrm{ad}(w_s^{-1}\mu)+N\big)
\cdot\mathrm{ad}(\eta x)\Big)~\mathrm{d}N
\\\nonumber\\
\label{global semisimple descent proof sixth equality}
&=&
|\pi_0(\mathrm{G}_\Sigma)|^{-1}
\sum_{s\in \mathrm{W}(\mathfrak{a}_1;\mathrm{M}_\mathrm{P}^+,\mathrm{G}_\Sigma^{0,+})}
~\sum_{N_\Sigma\in\mathfrak{m}_\mathrm{P}(\mathbb{Q}) \cdot \mathrm{ad}(w_s)\cap
\mathfrak{g}_{\Sigma,\mathrm{nil}}(\mathbb{Q})}
\\\nonumber\\\nonumber
&&\quad\times
\sum_{\pi\in \mathrm{P}(\mathbb{Q}) \cap w_s \mathrm{G}_\Sigma^0(\mathbb{Q}) w_s^{-1}\Big\backslash \mathrm{P}(\mathbb{Q})}
\int_{\mathfrak{n}_{\mathrm{P},\Sigma}(\mathbb{A})\cdot\mathrm{ad}(w_s)} 
\\\nonumber\\\nonumber
&&\qquad\times
f\Big(\big(\Sigma+N_\Sigma+N\big)
\cdot\mathrm{ad}(w_s^{-1}\pi x)\Big)~\mathrm{d}N
\end{eqnarray}
where \eqref{global semisimple descent proof sixth equality} follows from the change of variables
\begin{eqnarray}
&&
\mathrm{M}_\mathrm{P}(\mathbb{Q})\cap w_s \mathrm{G}_\Sigma^0(\mathbb{Q}) w_s^{-1} 
\Big\backslash \mathrm{M}_\mathrm{P}(\mathbb{Q})
\times
\mathrm{N}_{\mathrm{P},\Sigma\cdot\mathrm{ad}(w_s^{-1}\mu)}(\mathbb{Q}) 
\backslash \mathrm{N}_\mathrm{P}(\mathbb{Q})
\\\nonumber
&&\quad\stackrel{\sim}{\longrightarrow}\quad
\mathrm{P}(\mathbb{Q}) \cap w_s \mathrm{G}_\Sigma^0(\mathbb{Q}) w_s^{-1}\Big\backslash \mathrm{P}(\mathbb{Q}).
\end{eqnarray}
Substitute \eqref{global semisimple descent proof fifth equality} into the formula \eqref{global semisimple descent proof first lemma},
\begin{eqnarray}
&&
J_\mathfrak{o}^T(f) 
\\\nonumber\\\nonumber
&=& 
|\pi_0(\mathrm{G}_\Sigma)|^{-1}
\int_{\mathrm{G}(\mathbb{Q})\backslash\mathrm{G}(\mathbb{A})^1}
\bigg(\sum_{\substack{\mathrm{P}\in\mathcal{F}\\\mathrm{standard}}}
(-1)^{\mathrm{dim}(\mathrm{A}_\mathrm{P}/\mathrm{A}_\mathrm{G})}
\sum_{\delta\in\mathrm{P}(\mathbb{Q})\backslash\mathrm{G}(\mathbb{Q})}
\\\nonumber\\\nonumber
&&\quad
\times\sum_{s\in \mathrm{W}(\mathfrak{a}_1;\mathrm{M}_\mathrm{P}^+,\mathrm{G}_\Sigma^{0,+})}
~\sum_{N_\Sigma\in\mathfrak{m}_\mathrm{P}(\mathbb{Q}) \cdot \mathrm{ad}(w_s)\cap
\mathfrak{g}_{\Sigma,\mathrm{nil}}(\mathbb{Q})}
~\sum_{\pi\in \mathrm{P}(\mathbb{Q}) \cap w_s \mathrm{G}_\Sigma^0(\mathbb{Q}) w_s^{-1}\big\backslash \mathrm{P}(\mathbb{Q})}
\\\nonumber\\\nonumber
&&\qquad\times
\int_{\mathfrak{n}_{\mathrm{P},\Sigma}(\mathbb{A})\cdot\mathrm{ad}(w_s)} 
f\Big(\big(\Sigma+N_\Sigma+N\big)
\cdot\mathrm{ad}(w_s^{-1}\pi\delta x)\Big)~\mathrm{d}N\times
\\\nonumber\\\nonumber
&&\qquad\quad
\times\hat{\tau}_\mathrm{P}(H_\mathrm{p}(\delta x)-T)
\bigg)~\mathrm{d}x
\\\nonumber\\
\label{global semisimple descent proof seventh equality}
&=& 
|\pi_0(\mathrm{G}_\Sigma)|^{-1}
\int_{\mathrm{G}(\mathbb{Q})\backslash\mathrm{G}(\mathbb{A})^1}
\bigg(\sum_{\substack{\mathrm{P}\in\mathcal{F}\\\mathrm{standard}}}
(-1)^{\mathrm{dim}(\mathrm{A}_\mathrm{P}/\mathrm{A}_\mathrm{G})}
\sum_{s\in \mathrm{W}(\mathfrak{a}_1;\mathrm{M}_\mathrm{P}^+,\mathrm{G}_\Sigma^{0,+})}
\\\nonumber\\\nonumber
&&\quad
\times
\sum_{\xi\in \mathrm{Q}(\mathbb{Q}) \backslash \mathrm{G}(\mathbb{Q})}
~\sum_{N_\Sigma\in\mathfrak{m}_{\mathrm{Q},\mathrm{nil}}(\mathbb{Q})} 
\int_{\mathfrak{n}_\mathrm{Q}(\mathbb{A})} f((\Sigma+N_\Sigma+N)\cdot\mathrm{ad}(\xi x)) ~\mathrm{d}N\times
\\\nonumber\\\nonumber
&&\qquad
\times\hat{\tau}_\mathrm{P}(H_\mathrm{p}(w_s\xi x)-T)
\bigg)~\mathrm{d}x
\\\nonumber\\\nonumber
&=&
|\pi_0(\mathrm{G}_\Sigma)|^{-1}
\int_{\mathrm{G}(\mathbb{Q})\backslash\mathrm{G}(\mathbb{A})^1}
\bigg(\sum_{\substack{\mathrm{Q}\in\mathcal{F}^\Sigma\\\mathrm{standard}}}
~\sum_{\xi\in \mathrm{Q}(\mathbb{Q}) \backslash \mathrm{G}(\mathbb{Q})}
\\\nonumber\\\nonumber
&&\quad
\times
~\sum_{N_\Sigma\in\mathfrak{m}_{\mathrm{Q},\mathrm{nil}}(\mathbb{Q})} 
\int_{\mathfrak{n}_\mathrm{Q} (\mathbb{A})} f((\Sigma+N_\Sigma+N)\cdot\mathrm{ad}(\xi x)) ~\mathrm{d}N\times
\\\nonumber\\\nonumber
&&\qquad
\times\sum_{\substack{\mathrm{P}\in\mathcal{F},~\mathrm{standard}\\
s\in \mathrm{W}(\mathfrak{a}_1;\mathrm{M}_\mathrm{P}^+,\mathrm{G}_\Sigma^{0,+})\\
w_s^{-1}\mathrm{P}w_s\cap\mathrm{G}_\Sigma^0=\mathrm{Q}}}
(-1)^{\mathrm{dim}(\mathrm{A}_\mathrm{P}/\mathrm{A}_\mathrm{G})}
\hat{\tau}_\mathrm{P}(H_\mathrm{p}(w_s\xi x)-T)
\bigg)~\mathrm{d}x.
\end{eqnarray}
where in \eqref{global semisimple descent proof seventh equality} $\xi$ denotes the product $w_s^{-1}\pi\delta$ and $\mathrm{Q}$ denotes the standard parabolic subgroup of $\mathrm{G}_\Sigma^0$ defined by
\begin{eqnarray}
\mathrm{Q}
&=&
w_s^{-1}\mathrm{P}w_s\cap\mathrm{G}_\Sigma^0
\end{eqnarray}
with Levi decomposition
\begin{eqnarray}
\mathfrak{m}_\mathrm{Q}=\mathfrak{m}_\mathrm{P}\cdot\mathrm{ad}(w_s)\cap\mathfrak{g}_\Sigma,
&&
\mathfrak{n}_\mathrm{Q}=\mathfrak{n}_\mathrm{P}\cdot\mathrm{ad}(w_s)\cap\mathfrak{g}_\Sigma.
\end{eqnarray}
The assignment
\begin{eqnarray}
\label{global semisimple descent proof eighth map}
(\mathrm{P},s)
&\mapsto&
\mathrm{P}'=w_s^{-1}\mathrm{P}w_s
\end{eqnarray}
defines a bijection
\begin{eqnarray}
&&
\bigg\{\mathrm{P}\in\mathcal{F},~
s\in \mathrm{W}(\mathfrak{a}_1;\mathrm{M}_\mathrm{P}^+,\mathrm{G}_\Sigma^{0,+}):~
\mathrm{P}~\mathrm{standard},~
w_s^{-1}\mathrm{P}w_s\cap\mathrm{G}_\Sigma^0=\mathrm{Q}\bigg\}
\\\nonumber\\\nonumber
&&\quad
\stackrel{\sim}{\longrightarrow}
\mathcal{F}_\mathrm{Q}(\mathrm{M}_1),
\end{eqnarray}
hence
\begin{eqnarray}
\label{global semisimple descent proof ninth equality}
&&
J_\mathfrak{o}^T(f)
\\\nonumber\\\nonumber
&=&
|\pi_0(\mathrm{G}_\Sigma)|^{-1}
\int_{\mathrm{G}(\mathbb{Q})\backslash\mathrm{G}(\mathbb{A})^1}
\bigg(\sum_{\substack{\mathrm{Q}\in\mathcal{F}^\Sigma\\\mathrm{standard}}}
~\sum_{\xi\in \mathrm{Q}(\mathbb{Q}) \backslash \mathrm{G}(\mathbb{Q})}
\\\nonumber\\\nonumber
&&\quad
\times
~\sum_{N_\Sigma\in\mathfrak{m}_{\mathrm{Q},\mathrm{nil}}(\mathbb{Q})} 
\int_{\mathfrak{n}_\mathrm{Q} (\mathbb{A})} f((\Sigma+N_\Sigma+N)\cdot\mathrm{ad}(\xi x)) ~\mathrm{d}N\times
\\\nonumber\\\nonumber 
&& \qquad
\times \sum_{\mathrm{P}'\in\mathcal{F}_\mathrm{Q}(\mathrm{M}_1)}
(-1)^{\dim(\mathrm{A}_{\mathrm{P}'}/\mathrm{A}_\mathrm{G})} 
\hat{\tau}_{\mathrm{P}'}(H_{\mathrm{P}'}(\xi x)-s^{-1}(T-T_0)-T_0) \bigg) ~\mathrm{d}x
\end{eqnarray}
where $s$ denotes the second component of the inverse image of $\mathrm{P}'$ under the map defined in \eqref{global semisimple descent proof eighth map}.

Let $T_\Sigma$ be a truncation parameter for the triple $(\mathrm{G}_\Sigma(\mathbb{A})^0,\mathrm{M}_{1,\Sigma}(\mathbb{A}),\mathrm{K}_\Sigma)$ in $\mathfrak{a}_1$ such that $T_\Sigma-T_{\Sigma,1}$ is the projection of $T-T_0$. By \eqref{gamma'function identity var}
\begin{eqnarray}
&& 
\sum_{\mathrm{P}\in\mathcal{F}_\mathrm{Q}(\mathrm{M}_1)} 
(-1)^{\dim(\mathrm{A}_\mathrm{P}/\mathrm{A}_\mathrm{G})} 
\hat{\tau}_\mathrm{P}\big(H_\mathrm{P}(\delta xy)-s^{-1}(T-T_0)-T_0\big) 
\\\nonumber\\\nonumber 
&=& 
\sum_{\mathrm{P}\in\mathcal{F}_\mathrm{Q}(\mathrm{M}_1)} 
(-1)^{\dim(\mathrm{A}_\mathrm{P}/\mathrm{A}_\mathrm{G})} 
\hat{\tau}_\mathrm{P}\big((H_\mathrm{Q}(\delta x)-T_\Sigma)-Y_\mathrm{P}^T(\delta x,y)\big) 
\\\nonumber\\\nonumber 
&=& 
\sum_{\substack{\mathrm{Q}'\in\mathcal{F}^\Sigma\\\mathrm{Q}'\supset\mathrm{Q}}} 
(-1)^{\dim(\mathrm{A}_\mathrm{Q}/\mathrm{A}_{\mathrm{Q}'})} 
\hat{\tau}_\mathrm{Q}^{\mathrm{Q}'} \big(H_\mathrm{Q}(\delta x)-T_\Sigma \big)\times
\\\nonumber\\\nonumber
&&
\quad\times\Gamma_{\mathrm{Q}'}'\big(H_{\mathrm{Q}'}(\delta x)-T_\Sigma,
\mathcal{Y}_{\mathrm{Q}'}^T(\delta x,y) \big)
\end{eqnarray}
where the family $\mathcal{Y}_{\mathrm{Q}'}^T(\delta x,y)$ is defined by
\begin{eqnarray}
\forall \mathrm{P}\in\mathcal{F}_\mathrm{Q}(\mathrm{M}_1)
&&
Y_\mathrm{P}^T(\delta x,y)=-H_\mathrm{P}(ky)+s^{-1}(T-T_0)-T_\Sigma+T_0,
\end{eqnarray} 
where $k$ is the $\mathrm{K}_\Sigma$ component of $\delta x$ under the Iwasawa decomposition with respect to $\mathrm{P}_\Sigma$ and $s$ is the second component of the inverse image of $\mathrm{P}$ under the map defined in \eqref{global semisimple descent proof eighth map}. On the right hand side of \eqref{global semisimple descent proof ninth equality} make the change of variables
\begin{eqnarray}
&&
\mathrm{Q}(\mathbb{Q})\backslash\mathrm{G}(\mathbb{Q}) 
\times \mathrm{G}(\mathbb{Q})\backslash\mathrm{G}(\mathbb{A})^1 
\\\nonumber
&&\quad
\stackrel{\sim}{\longrightarrow}\quad
\mathrm{Q}(\mathbb{Q})\backslash\mathrm{G}_\Sigma^0(\mathbb{Q}) 
\times \mathrm{G}_\Sigma^0(\mathbb{Q}) 
\backslash \mathrm{G}_\Sigma^0(\mathbb{A})\cap\mathrm{G}(\mathbb{A})^1 \times \mathrm{G}_\Sigma^0(\mathbb{A}) \backslash \mathrm{G}(\mathbb{A}).
\end{eqnarray}
Then
\begin{eqnarray}
&&
J_\mathfrak{o}^T(f)
\\\nonumber\\\nonumber
&=&
|\pi_0(\mathrm{G}_\Sigma)|^{-1}
\int_{\mathrm{G}_\Sigma^0(\mathbb{A})\backslash\mathrm{G}(\mathbb{A})} 
\int_{\mathrm{G}_\Sigma^0(\mathbb{Q})\backslash\mathrm{G}_\Sigma^0(\mathbb{A})
\cap\mathrm{G}(\mathbb{A})^1} 
\bigg( \sum_{\substack{\mathrm{Q}\in\mathcal{F}^\Sigma\\
\mathrm{standard}}} 
~\sum_{\substack{\mathrm{Q}'\in\mathcal{F}^\Sigma\\
\mathrm{Q}'\supset\mathrm{Q}}} 
\\\nonumber\\\nonumber
&& \quad
\times\sum_{\delta\in\mathrm{Q}(\mathbb{Q})\backslash\mathrm{G}_\Sigma^0(\mathbb{Q})} 
~\sum_{N_\Sigma\in\mathfrak{m}_{\mathrm{Q},\mathrm{nil}}(\mathbb{Q})} 
(-1)^{\dim(\mathrm{A}_\mathrm{Q}/\mathrm{A}_{\mathrm{Q}'})} \times
\\\nonumber\\\nonumber 
&& \qquad
\times \int_{\mathfrak{n}_\mathrm{Q} (\mathbb{A})} 
f((\Sigma+N_\Sigma+N)\cdot\mathrm{ad}(\delta xy)) ~\mathrm{d}N \times
\\\nonumber\\\nonumber 
&& \qquad\quad
\times  \hat{\tau}_\mathrm{Q}^{\mathrm{Q}'} \big(H_\mathrm{Q}(\delta x)-T_\Sigma \big) 
\Gamma_{\mathrm{Q}'}'\big(H_{\mathrm{Q}'}(\delta x)-T_\Sigma,
\mathcal{Y}_{\mathrm{Q}'}^T(\delta x,y)\big) \bigg) \mathrm{d}x\mathrm{d}y
\\\nonumber\\
\label{global semisimple descent proof tenth equality}
&=&
|\pi_0(\mathrm{G}_\Sigma)|^{-1}
\int_{\mathrm{G}_\Sigma^0(\mathbb{A})\backslash\mathrm{G}(\mathbb{A})}  
\Bigg( \sum_{\substack{\mathrm{Q}'\in\mathcal{F}^\Sigma\\
\mathrm{standard}}} 
\int_{\mathrm{K}_\Sigma} 
\int_{\mathrm{A}_{\mathrm{Q}'}(\mathbb{R})\cap \mathrm{G}(\mathbb{A})^1} 
\int_{\mathrm{M}_{\mathrm{Q}'}(\mathbb{Q}) \backslash \mathrm{M}_{\mathrm{Q}'}(\mathrm{A})^1} 
\\\nonumber\\\nonumber
&&\quad 
\times \bigg( \sum_{\substack{\mathrm{Q}\in\mathcal{F}^\Sigma\\
\mathrm{standard}\\\mathrm{Q}\subset\mathrm{Q}'}} 
~\sum_{\mu\in\mathrm{Q}(\mathbb{Q}) 
\cap \mathrm{M}_{\mathrm{Q}'}(\mathbb{Q}) \backslash \mathrm{M}_{\mathrm{Q}'}(\mathbb{Q})}   
~\sum_{N_\Sigma\in\mathfrak{m}_{\mathrm{Q},\mathrm{nil}}(\mathbb{Q})} 
(-1)^{\dim(\mathrm{A}_\mathrm{Q}/\mathrm{A}_{\mathrm{Q}'})} \times
\\\nonumber\\\nonumber
&&\qquad 
\times  \int_{\mathfrak{n}_\mathrm{Q} (\mathbb{A}) \cap \mathfrak{m}_{\mathrm{Q}'}(\mathbb{A})} 
\Phi_{\mathrm{Q}',a,k,y}^T\big((N_\Sigma+N) \cdot\mathrm{ad}(\mu m)\big) ~\mathrm{d}N \times 
\\\nonumber\\\nonumber
&&\qquad\quad
\times \hat{\tau}_\mathrm{Q}^{\mathrm{Q}'} \big(H_\mathrm{Q}(\mu m)-T_\Sigma \big) \bigg) 
~\mathrm{d}m\mathrm{d}a\mathrm{d}k \Bigg) ~\mathrm{d}y
\\\nonumber\\\nonumber
&=&
|\pi_0(\mathrm{G}_\Sigma)|^{-1}
\int_{\mathrm{G}_\Sigma^0(\mathbb{A})\backslash\mathrm{G}(\mathbb{A})} 
\sum_{\substack{\mathrm{Q}'\in\mathcal{F}^\Sigma\\\mathrm{standard}}} 
\\\nonumber\\\nonumber
&&\quad
\times \bigg(\int_{\mathrm{K}_\Sigma} 
\int_{\mathrm{A}_{\mathrm{Q}'}(\mathbb{R})\cap\mathrm{G}(\mathbb{A})^1} 
J_\mathrm{nil} ^{\mathrm{M}_{\mathrm{Q}'},T_\Sigma}( \Phi_{\mathrm{Q}',a,k,y}^T ) ~\mathrm{d}a\mathrm{d}k \bigg)~\mathrm{d}y.
\end{eqnarray}
where on the right hand side of \eqref{global semisimple descent proof tenth equality} $\Phi_{\mathrm{Q}',a,k,y}^T$ denotes the function on $\mathfrak{m}_{\mathrm{Q}'}(\mathbb{A})$ defined by
\begin{eqnarray}
\forall X\in\mathfrak{m}_{\mathrm{Q}'}(\mathbb{A})
\\\nonumber
\Phi_{\mathrm{Q}',a,k,y}^T(X) 
&=&
\int_{\mathfrak{m}_{\mathrm{Q}'}(\mathbb{A})} 
f((\Sigma+X+N)\cdot\mathrm{ad}(ky)) ~\mathrm{d}N\times
\\\nonumber\\\nonumber
&&\quad
\times \Gamma_{\mathrm{Q}'}' \big(H_{\mathrm{Q}'}(a)-T_\Sigma, \mathcal{Y}_{\mathrm{Q}'}^T(k,y) \big),
\end{eqnarray}
and the equality \eqref{global semisimple descent proof tenth equality} follows from the changes of variables
\begin{eqnarray}
&&
\mathrm{Q}(\mathbb{Q})\backslash \mathrm{G}_\Sigma^0(\mathbb{Q}) 
\times \mathrm{G}_\Sigma^0(\mathbb{Q}) \backslash\mathrm{G}_\Sigma^0(\mathbb{A})
\cap\mathrm{G}(\mathbb{A})^1 
\\\nonumber 
&&\quad\stackrel{\sim}{\longrightarrow}\quad
\mathrm{Q}(\mathbb{Q}) 
\cap\mathrm{M}_{\mathrm{Q}'}(\mathbb{Q}) \backslash \mathrm{M}_{\mathrm{Q}'}(\mathbb{Q}) 
\times \mathrm{Q}'(\mathbb{Q})\backslash 
\mathrm{G}_\Sigma^0(\mathbb{A})\cap\mathrm{G}(\mathbb{A})^1
\end{eqnarray}
and
\begin{eqnarray}
\mathfrak{n}_\mathrm{Q}(\mathbb{A})
&\stackrel{\sim}{\longrightarrow}&
\mathfrak{n}_\mathrm{Q}(\mathbb{A}) \cap \mathfrak{m}_{\mathrm{Q}'}(\mathbb{A}) \times \mathfrak{n}_{\mathrm{Q}'}(\mathbb{A})
\end{eqnarray}
and the Iwasawa decomposition for $\mathrm{G}_\Sigma^0(\mathbb{A})\cap\mathrm{G}(\mathbb{A})^1$.

Evaluate $T$ at $T_0$ and $T_\Sigma$ at $T_{\Sigma,1}$, and apply \eqref{gamma' function var fourier transform}, then
\begin{eqnarray}
J_\mathfrak{o}(f) 
&=&
|\pi_0(\mathrm{G}_\Sigma)|^{-1}
\int_{\mathrm{G}_\Sigma^0(\mathbb{A})\backslash \mathrm{G}(\mathbb{A})} 
\sum_{\substack{\mathrm{Q}'\in\mathcal{F}^\Sigma\\\mathrm{standard}}} 
J_\mathrm{nil}^{\mathrm{M}_{\mathrm{Q}'}}(\Phi_{\mathrm{Q}',y}^{T_0-T_{\Sigma,1}}) ~\mathrm{d}y.
\end{eqnarray}
By the defining property \eqref{point t0} of the points $T_0$ and $T_{\Sigma,1}$, 
\begin{eqnarray}
J_\mathrm{nil}^{\mathrm{M}_{\mathrm{Q}'}}(\Phi_{\mathrm{Q}',y}^{T_0-T_{\Sigma,1}}) 
&=& 
J_\mathrm{nil}^{\mathrm{M}_{\mathrm{Q}''}}(\Phi_{\mathrm{Q}'',y}^{T_0-T_{\Sigma,1}})
\end{eqnarray}
whenever $\mathrm{Q}'$ and $\mathrm{Q}''$ are parabolic subgroups in $\mathcal{F}^\Sigma$ conjugate under the Weyl group of $\mathrm{G}_\Sigma^0$, hence 
\begin{eqnarray}
&&
J_\mathfrak{o}(f)
\\\nonumber
&=& 
|\pi_0(\mathrm{G}_\Sigma)|^{-1}
\int_{\mathrm{G}_\Sigma^0(\mathbb{A})\backslash\mathrm{G}(\mathbb{A})} 
\Bigg( \sum_{\mathrm{Q}\in\mathcal{F}^{\mathrm{G}_\Sigma}} 
|\mathrm{W}_1^{\mathrm{M}_\mathrm{Q}}| 
|\mathrm{W}_1^{\mathrm{G}_\Sigma^0}|^{-1}  
J_\mathrm{nil}^{\mathrm{M}_\mathrm{Q}}(\Phi_{\mathrm{Q},x}^{T_0-T_{\Sigma,1}}) \Bigg)~\mathrm{d}x.
\end{eqnarray}
\qed

\paragraph{Definition}\setcounter{equation}{0}
Fix a finite set of places $S$. Let $\mathrm{M}$ be a standard Levi subgroup of $\mathrm{G}$, let $\Xi$\label{sym98} be a semisimple element of $\mathfrak{m}(\mathbb{Q})$, let $\nu$ be a nilpotent element of $\mathfrak{m}_\Xi(\mathbb{Q}_S)$ defined upto the adjoint action of $\mathrm{M}_\Xi(\mathbb{Q}_S)$. The element 
\begin{eqnarray}
\label{local semisimple descent point definition}
X=\Xi+\nu
&\in&
\mathfrak{m}(\mathbb{Q}_S)
\end{eqnarray}
is well-defined modulo $\mathrm{M}(\mathbb{Q}_S)$. Denote by $D^\mathrm{G}$ the discriminant function on $\mathfrak{g}$.

Let $T$ be a point in $\mathfrak{a}_\mathrm{M}$, let $x$ be an element of $\mathrm{G}(\mathbb{Q}_S)$. Let $(v_\mathrm{P}(x,T))$\label{sym99} be the $(\mathrm{G},\mathrm{M})$-family defined by
\begin{eqnarray}
\forall \mathrm{P}\in\mathcal{F}(\mathrm{M})
~\forall \lambda\in i\mathfrak{a}_\mathrm{M}^*
&&
v_\mathrm{P}(x,T)(\lambda)
=
v_\mathrm{P}(x)(\lambda)\times e^{\langle\lambda,T\rangle}.
\end{eqnarray}
Let $\mathrm{Q}$ be a parabolic subgroup of $\mathrm{G}_\Xi^0$ containing $\mathrm{M}_\Xi^0$. Define the \emph{weight factor $v_\mathrm{Q}'$} by
\begin{eqnarray}
\label{local semisimple descent weight factor' definition}
v_\mathrm{Q}'(x,T) 
&=& 
\sum_{\mathrm{P}\in\mathring{\mathcal{F}}_\mathrm{Q}(\mathrm{M})} v_\mathrm{P}'(x,T)
\end{eqnarray}
where the set $\mathring{\mathcal{F}}_\mathrm{Q}(\mathrm{M})$ is defined with respect to the semisimple element $\Xi$. Let $T_{\Xi,\mathrm{M}}$\label{sym100} be the point in $\mathfrak{a}_{\mathrm{M}_\Xi^0}$ modulo $\mathfrak{a}_{\mathrm{G}_\Xi^0}$ defined in the same manner as $T_0$ in $\mathfrak{a}_0^\mathrm{G}$ that satisfies an identity analogous to \eqref{point t0}.

\paragraph{Lemma}\setcounter{equation}{0}
\emph{Let $\mathrm{L}$ be a Levi subgroup in $\mathcal{L}(\mathrm{M})$, let $T$ be a point in $\mathfrak{a}_\mathrm{M}$, then
\begin{eqnarray}
\label{local semisimple descent weight factor' identity}
\forall x\in\mathrm{G}(\mathbb{Q}_S)
~\forall y\in\mathrm{G}_\Xi^0(\mathbb{Q}_S)
&&
v_\mathrm{L}(yx) 
= 
\sum_{\mathrm{Q}\in\mathcal{F}^\Xi(\mathrm{M}_\Xi^0)} 
v_{\mathrm{L}_\Xi^0}^\mathrm{Q}(y) v_\mathrm{Q}'(kx,T)
\end{eqnarray}
where $v_\mathrm{L}(x)$ is the weight factor defined in \eqref{weight factor definition} and $k$ is the $\mathrm{K}_\Xi$-component of $y$ under the Iwasawa decomposition of $\mathrm{G}_\Xi^0(\mathbb{Q}_S)$ with respect to the parabolic subgroup $\mathrm{Q}_\Xi^0$. }

\proof
This is Corollary 8.4 of \cite{Arthur_Loc}.
\qed

\paragraph{Lemma}\setcounter{equation}{0} (Local semisimple descent)\\
\emph{Let $f_S$ be a Schwartz function on $\mathfrak{g}(\mathbb{Q}_S)$, let $X$ and $\nu$ be as in \eqref{local semisimple descent point definition}, then
\begin{eqnarray}
\label{local semisimple descent}
\qquad\qquad
J_\mathrm{M}^\mathrm{G}(X,f_S) 
&=& 
|D^\mathrm{G}(\Xi)|_S^{1/2} 
\int_{\mathrm{G}_\Xi^0(\mathbb{Q}_S)\backslash\mathrm{G}(\mathbb{Q}_S)}
\bigg( \sum_{\mathrm{Q}\in\mathcal{F}^\Xi(\mathrm{M}_\Xi^0)} 
J_{\mathrm{M}_\Xi^0}^{\mathrm{M}_\mathrm{Q}} 
(\nu,\Phi_{S,\mathrm{Q},x}^{T_0-T_{\Xi,\mathrm{M}}}) \bigg)~\mathrm{d}x
\end{eqnarray}
where $\Phi_{S,\mathrm{Q},x}^T$\label{sym101} is the function on $\mathfrak{m}_\mathrm{Q}(\mathbb{Q}_S)$ defined by
\begin{eqnarray}
&&
\forall Y\in\mathfrak{m}_\mathbb{Q}(\mathbb{Q}_S)
\\\nonumber
&&\quad
\Phi_{S,\mathrm{Q},x}^T(Y) 
\\\nonumber\\\nonumber
&&\quad
=\int_{\mathrm{K}_\Xi}
\int_{\mathfrak{n}_\mathrm{Q}(\mathbb{Q}_S)}f_S\Big( \big(\Xi+(Y+N) \cdot \mathrm{ad}(k)\big) \cdot\mathrm{ad}(x) \Big)~v_\mathrm{Q}'(kx,T) \mathrm{d}N \mathrm{d}k
\end{eqnarray}
where the weight factor $v_\mathrm{Q}'(x,T)$ is defined as in \eqref{local semisimple descent weight factor' definition}.}

\proof
This argument follows the proof of Corollary 8.7 of \cite{Arthur_Loc}, see also the discussion in \S2.6 of \cite{Hoff}.

Since the point $T_0-T_{\Xi,\mathrm{M}}$ lies in $\mathfrak{a}_\mathrm{M}$, the identity \eqref{local semisimple descent weight factor' identity} is valid for the weight factor $v_\mathrm{Q}'(x,T_0-T_{\Xi,\mathrm{M}})$, hence
\begin{eqnarray}
&&
J_\mathrm{M}^\mathrm{G}(X,f_S) 
\\\nonumber\\\nonumber
&=& 
\lim_{A\rightarrow0} \Bigg( |D^\mathrm{G}(X+A)|_S^{1/2} \int_{\mathrm{G}_{X_\mathrm{ss}}^0(\mathbb{Q}_S) \backslash \mathrm{G}(\mathbb{Q}_S)} 
\int_{\mathrm{G}_{X+A}^0(\mathbb{Q}_S)\backslash\mathrm{G}_{X_\mathrm{ss}}^0(\mathbb{Q}_S)}
\\\nonumber\\\nonumber 
&&  
\quad \times f_S((X+A)\cdot\mathrm{ad}(yx)) \times
\\\nonumber\\\nonumber
&& \qquad
\times \bigg( \sum_{\mathrm{L}\in\mathcal{L}(\mathrm{M})} 
r_\mathrm{M}^\mathrm{L}\big(\exp(X_\mathrm{nil}),\exp(A)\big) 
v_\mathrm{L}(yx) \bigg) \mathrm{d}y\mathrm{d}x \Bigg) 
\\\nonumber\\\nonumber
&=& 
\lim_{A\rightarrow0} 
\Bigg( |D^\mathrm{G}(X+A)|_S^{1/2} \int_{\mathrm{G}_\Xi^0(\mathbb{Q}_S)\backslash \mathrm{G}(\mathbb{Q}_S)} \int_{\mathrm{M}_X^0(\mathbb{Q}_S)\backslash \mathrm{G}_\Xi^0(\mathbb{Q}_S)}
\\\nonumber\\\nonumber
&& \quad
\times f_S((X+A)\cdot\mathrm{ad}(yx)) \times 
\\\nonumber\\\nonumber 
&& 
\qquad 
\times \bigg( \sum_{\mathrm{Q}\in\mathcal{F}^\Xi(\mathrm{M}_\Xi^0)} 
~\sum_{\mathrm{L}\in\mathcal{L} ^{\mathrm{M}_\mathrm{Q}} (\mathrm{M}_\Xi^0)} r_{\mathrm{M}_\Xi^0}^\mathrm{L}\big(\exp(X_\mathrm{nil}),\exp(A)\big) \times
\\\nonumber\\\nonumber
&& \qquad\quad
\times v_\mathrm{L}^\mathrm{Q}(y) v_\mathrm{Q}'(kx,T_0-T_{\Xi,\mathrm{M}}) \bigg)\mathrm{d}y\mathrm{d}x \Bigg)
\\\nonumber\\
\label{local semisimple descent proof equality}
&=&
\lim_{A\rightarrow0} 
\Bigg(  |D^\mathrm{G}(X+A)|_S^{1/2} 
~|D^{\mathrm{G}_\Xi}(X+A)|_S^{-1/2} 
\int_{\mathrm{G}_\Xi^0(\mathbb{Q}_S)\backslash \mathrm{G}(\mathbb{Q}_S)} 
\\\nonumber\\\nonumber
&&\quad
\times\bigg( \sum_{\mathrm{Q}\in\mathcal{F}^\Xi(\mathrm{M}_\Xi^0)}  \sum_{\mathrm{L}\in\mathcal{L} ^{\mathrm{M}_\mathrm{Q}} (\mathrm{M}_\Xi)} r_{\mathrm{M}_\Xi^0}^\mathrm{L}\big(\exp(X_\mathrm{nil}),\exp(A)\big) \times
\\\nonumber\\\nonumber  
&&\qquad
\times J_\mathrm{L}^{\mathrm{M}_\mathrm{Q}} (X_\mathrm{nil}+A,\Phi_{S,\mathrm{Q},y}^{T_0-T_{\Xi,\mathrm{M}}}) \bigg)\mathrm{d}y \Bigg)
\\\nonumber\\\nonumber
&=&
|D^\mathrm{G}(\Xi)|_S^{1/2} \int_{\mathrm{G}_\Xi^0(\mathbb{Q}_S)\backslash\mathrm{G}(\mathbb{Q}_S)} \bigg( \sum_{\mathrm{Q}\in\mathcal{F}^\Xi(\mathrm{M}_\Xi^0)} J_{\mathrm{M}_\Xi}^{\mathrm{M}_\mathrm{Q}} (X_\mathrm{nil},\Phi_{S,\mathrm{Q},y}^{T_0-T_{\Xi,\mathrm{M}}}) \bigg) \mathrm{d}y
\end{eqnarray}
where the equality \eqref{local semisimple descent proof equality} follows from the change of variables
\begin{eqnarray}
\mathrm{M}_X^0(\mathbb{Q}_S)\backslash\mathrm{G}_\Xi^0(\mathbb{Q}_S) &\stackrel{\sim}{\longrightarrow}& \mathrm{M}_X(\mathbb{Q}_S)\backslash\mathrm{M}_{\mathrm{Q}_\Xi^0}(\mathbb{Q}_S) \times \mathrm{N}_\mathrm{Q}(\mathbb{Q}_S) \times \mathrm{K}_\Xi 
\end{eqnarray}
for each parabolic subgroup $\mathrm{Q}$ in $\mathcal{F}^\Xi(\mathrm{M}_\Xi^0)$.
\qed

\paragraph{Definition}\setcounter{equation}{0}
Let $\equiv$\label{sym102} be the equivalence relation on $\mathfrak{m}(\mathbb{Q})\cap\mathfrak{o}$ defined by
\begin{eqnarray}
&&
\forall X\in\mathfrak{m}(\mathbb{Q})\cap\mathfrak{o}\quad
\forall Y\in\mathfrak{m}(\mathbb{Q})\cap\mathfrak{o}\quad
X\equiv Y\quad
\mathrm{if}
\\\nonumber
&&\quad
\exists\delta\in\mathrm{M}(\mathbb{Q})\quad
\exists\eta\in \mathrm{M}_{X_\mathrm{ss}}^0(\mathbb{Q}_S)\quad
X_\mathrm{ss}=Y_\mathrm{ss}\cdot\mathrm{ad}(\delta)\quad
X_\mathrm{nil}=Y_\mathrm{nil}\cdot \mathrm{ad}(\delta\eta).
\end{eqnarray}
Let $(\mathfrak{m}(\mathbb{Q})\cap\mathfrak{o})_{\mathrm{M},S}$\label{sym103} denote the collection of $\equiv$ equivalence classes in $\mathfrak{m}(\mathbb{Q})\cap\mathfrak{o}$. By abuse of notation denote the $\equiv$ equivalence class of an element $X$ in $\mathfrak{m}(\mathbb{Q})\cap\mathfrak{o}$ by the same symbol $X$.

\paragraph{Proposition}\setcounter{equation}{0}
\emph{Let $\mathfrak{o}$ be a $\sim$ equivalence class on $\mathfrak{g}(\mathbb{Q})$. Let $S_\mathfrak{o}$\label{sym104} be a sufficiently large finite set of places of $\mathbb{Q}$. Then for each standard Levi subgroup $\mathrm{M}$ of $\mathrm{G}$ and for each $\equiv$ equivalence class $X$ in $\mathfrak{m}(\mathbb{Q})\cap\mathfrak{o}$ there exists a constant $a^\mathrm{M}(S_\mathfrak{o},X)$\label{sym105} such that}
\begin{eqnarray}
\label{refined expansion}
&&
\forall f\in\mathcal{S}(\mathfrak{g}(\mathbb{A}))
\\\nonumber\\\nonumber
&&\quad
J_\mathfrak{o}^\mathrm{G}(f) = \sum_{\mathrm{M}\in\mathcal{L}}|\mathrm{W}_0^\mathrm{M}||\mathrm{W}_0^\mathrm{G}|^{-1}
\sum_{X\in(\mathfrak{m}(\mathbb{Q})\cap\mathfrak{o})_{\mathrm{M},S_\mathfrak{o}}} a^\mathrm{M}(S_\mathfrak{o},X)J_\mathrm{M}^\mathrm{G}(X,f_{S_\mathfrak{o}}).
\end{eqnarray}

\proof
Retain the notations of \eqref{global semisimple descent} and \eqref{local semisimple descent}. Let $\Xi$ be equal to $\Sigma$, hence 
\begin{eqnarray}
X=\Sigma+\nu
&\in&
\mathfrak{m}(\mathbb{Q}_S)\slash\mathrm{ad}(\mathrm{M}(\mathbb{Q}_S)),
\end{eqnarray} 
and 
\begin{eqnarray}
T_0-T_{\Xi,\mathrm{M}}
&=&
T_0-T_{\Sigma,1}.
\end{eqnarray}
Without loss of generality assume that the finite set $S_\mathfrak{o}$ is large enough such that
\begin{eqnarray}
\label{refined expansion proof first condition}
\forall p\notin S_\mathfrak{o}
&& \bullet
\Sigma\in\mathfrak{g}(\mathbb{Z}_p);
\\\nonumber
&&\bullet
|D(\Sigma)|_p=1;
\\\nonumber
&&\bullet
\mathfrak{g}(\mathbb{Z}_p)\cdot\mathrm{ad}(\mathrm{K}_p)=\mathfrak{g}(\mathbb{Z}_p);
\\\nonumber
&&\bullet
\mathrm{K}_p\cap\mathrm{G}_\Sigma^0(\mathbb{Q}_p)=\mathrm{K}_{\Sigma,p} \textrm{ is hyperspecial in } \mathrm{G}_\Sigma^0(\mathbb{Q}_p);
\\\nonumber
&&\bullet
\forall x_p\in\mathrm{G}(\mathbb{Q}_p)~
(\Sigma+\mathfrak{g}_{\Sigma,\mathrm{nil}}(\mathbb{Q}_p)) \cdot \mathrm{ad}(x_p)\cap \mathfrak{g}(\mathbb{Z}_p)\neq\emptyset
\\\nonumber
&&\quad
\Rightarrow x_p\in\mathrm{G}_\Sigma^0(\mathbb{Q}_p)\mathrm{K}_p.
\end{eqnarray} 
The last condition in \eqref{refined expansion proof first condition} is satisfied for large enough $S_\mathfrak{o}$ by Lemma 6.1 of \cite{Arthur_Orb}. Since 
\begin{eqnarray}
\bigg(H\mapsto\Gamma_\mathrm{Q}'\big( H-T_{\Sigma,1},\mathcal{Y}_\mathrm{Q}^{T_0}(k,x)\big)\bigg)^\wedge
&=&
v_\mathrm{Q}'(kx,T_0-T_{\Sigma,1})(\lambda)
\end{eqnarray}the conditions \eqref{refined expansion proof first condition} imply that 
\begin{eqnarray}
\Phi_{\mathrm{Q},x}^{T_0-T_{\Sigma,1}} &=& \Phi_{S_\mathfrak{o},\mathrm{Q},x_{S_\mathfrak{o}}}^{T_0-T_{\Sigma,1}} \otimes \bigotimes_{p\notin S_\mathfrak{o}}\mathbb{I}_{\mathrm{G}_\Sigma^0 (\mathbb{Q}_p)\mathrm{K}_p}(x)\cdot\mathbb{I}_ {\mathfrak{m}_\mathrm{Q}(\mathbb{Z}_p)} .
\end{eqnarray}
By\eqref{global semisimple descent}) and \eqref{nilpotent refined expansion},
\begin{eqnarray}
\label{refined expansion proof second equality}
&&
J_\mathfrak{o}^\mathrm{G}(f) 
\\\nonumber\\\nonumber
&=& 
|\pi_0(\mathrm{G}_\Sigma)|^{-1} \int_{\mathrm{G}_\Sigma^0(\mathbb{Q}_{S_\mathfrak{o}})\backslash\mathrm{G}(\mathbb{Q}_{S_\mathfrak{o}})} \bigg( \sum_{\mathrm{Q}\in\mathcal{F}^\Sigma} ~\sum_{\mathrm{L}\in\mathcal{L}^{\mathrm{M}_\mathrm{Q}}(\mathrm{M}_{1,\Sigma}^0)} |\mathrm{W}_1^\mathrm{L}| 
|\mathrm{W}_1^{\mathrm{G}_\Sigma^0}|^{-1} \times 
\\\nonumber\\\nonumber 
&&  \quad
\times 
\sum_{\nu\in(\mathfrak{l}_{\mathrm{nil}}(\mathbb{Q}))_{\mathrm{L},S_\mathfrak{o}}} 
 a^\mathrm{L}(S_\mathfrak{o},\nu)
 J_\mathrm{L}^{\mathrm{M}_\mathrm{Q}} 
(\nu,\Phi_{S_\mathfrak{o},\mathrm{Q},x}^{T_0-T_{\Sigma,1}})  \bigg) \mathrm{d}x  
\\\nonumber\\
\label{refined expansion proof third equality} 
&=& 
|\pi_0(\mathrm{G}_\Sigma)|^{-1}
\sum_{\mathrm{M} \in\mathring{\mathcal{L}}_\Sigma(\mathrm{M}_1)}
|\mathrm{W}_1^{\mathrm{M}_\Sigma^0}| 
|\mathrm{W}_1^{\mathrm{G}_\Sigma^0}|^{-1}  \times
\\\nonumber\\\nonumber
&&\quad \times
\sum_{\nu\in(\mathfrak{m}_{\Sigma,\mathrm{nil}}(\mathbb{Q}))_{\mathrm{M}_\Sigma^0,S_\mathfrak{o}}}  
a^{\mathrm{M}_\Sigma^0}(S_\mathfrak{o},\nu) J_\mathrm{M}^\mathrm{G} (\Sigma+\nu,f_{S_\mathfrak{o}})
\end{eqnarray}
where on the right hand side of \eqref{refined expansion proof third equality} the symbol $\mathring{\mathcal{L}}_\Sigma(\mathrm{M}_1)$\label{sym106} denotes the set
\begin{eqnarray}
\mathring{\mathcal{L}}_\Sigma(\mathrm{M}_1)
&=&
\bigg\{
\mathrm{M}\in\mathcal{L}(\mathrm{M}_1):~
\mathrm{A}_\mathrm{M}=\mathrm{A}_{\mathrm{M}_\Sigma^0}
\bigg\}.
\end{eqnarray}
The equality \eqref{refined expansion proof third equality} follows from \eqref{local semisimple descent} and the bijection
\begin{eqnarray}
\pi_\Sigma:\quad
\mathring{\mathcal{L}}_\Sigma(\mathrm{M}_1)
&\stackrel{\sim}{\longrightarrow}&
\mathcal{L}^{\mathrm{G}_\Sigma^0}(\mathrm{M}_{1,\Sigma}^0)
\end{eqnarray}
defined by
\begin{eqnarray}
\forall\mathrm{M}\in\mathring{\mathcal{L}}_\Sigma(\mathrm{M}_1)
&&
\pi_\Sigma(\mathrm{M})=\mathrm{M}_\Sigma^0=\mathrm{L}
\end{eqnarray}
and the equation
\begin{eqnarray}
|D(\Sigma)|_{S_\mathfrak{o}}&=&1
\end{eqnarray}
which follows from \eqref{refined expansion proof first condition}.

Define the constant $a^\mathrm{M}(S_\mathfrak{o},X)$ by
\begin{eqnarray}
\label{global coefficient}
a^\mathrm{M}(S_\mathfrak{o},X) 
&=& 
|\pi_0(\mathrm{M}_\Sigma)|^{-1}
\sum_{\substack{\nu\in(\mathfrak{m}_{\Sigma,\mathrm{nil}})_{\mathrm{M}_\Sigma^0,S_\mathfrak{o}}\\
\Sigma+\nu \equiv X \mathrm{mod}(\mathrm{M},S_\mathfrak{o}) }} a^{\mathrm{M}_\Sigma^0}(S_\mathfrak{o},\nu)
\end{eqnarray}
if $X_\mathrm{ss}$ is $\mathbb{Q}$-elliptic in $\mathfrak{m}(\mathbb{Q})$ and zero otherwise. 
For a fixed $\Sigma$ and let $\mathrm{M}$ vary in $\mathcal{L}(\mathrm{M}_1)$, then $\Sigma$ is $\mathbb{Q}$-elliptic in $\mathfrak{m}(\mathbb{Q})$ if and only if $A_\mathrm{M}$ and $A_{\mathrm{M}_\Sigma}$ are equal, hence by \eqref{refined expansion proof second equality}
\begin{eqnarray}
J_\mathfrak{o}^\mathrm{G}(f) 
&=& 
\sum_{\mathrm{M}\in\mathcal{L}(\mathrm{M}_1)} 
|\mathrm{W}_1^{\mathrm{M}_\Sigma^0}| |\mathrm{W}_1^{\mathrm{G}_\Sigma^0}|^{-1}
|\pi_0(\mathrm{M}_\Sigma)||\pi_0(\mathrm{G}_\Sigma)|^{-1} \times
\\\nonumber\\\nonumber
&&\quad
\times  \sum_{\substack{X\in(\mathfrak{m} (\mathbb{Q})\cap\mathfrak{o})_{\mathrm{M},S_\mathfrak{o}}
\\X_\mathrm{ss}=\Sigma}} 
a^\mathrm{M}(S_\mathfrak{o},X) 
J_\mathrm{M}^\mathrm{G}(X,f_{S_\mathfrak{o}}).
\end{eqnarray}
The set
\begin{eqnarray}
\label{refined expansion proof fourth set}
\Bigg\{~ (\mathrm{M},\sigma):~ \begin{array}{l} \mathrm{M}\in\mathcal{L},\sigma \textrm{ is a semisimple $\mathrm{M}(\mathbb{Q})$-}\\ \textrm{-orbit in } \mathfrak{m}(\mathbb{Q}), \textrm{ $\Sigma$ is conjugate}\\ \textrm{to a point in $\sigma$ under $\mathrm{G}(\mathbb{Q})$. }  \end{array} ~\Bigg\}
\end{eqnarray}
admits a natural action by $\mathrm{W}_0^\mathrm{G}$ which preserves $a^\mathrm{M}(S_\mathfrak{o},X)$ and $J_\mathrm{M}^\mathrm{G}(X,f_{S_\mathfrak{o}})$, where each pair $(\mathrm{M},\sigma)$ has stablizer $\mathrm{W}_0^\mathrm{M}$. Denote by $\mathrm{W}_1^{\mathrm{G}_\Sigma}$ the quotient of the normalizer of $\mathrm{A}_1$ in $\mathrm{G}_\Sigma$ by $\mathrm{M}_{1,\Sigma}$.  There is an exact sequence
\begin{eqnarray}
\xymatrix{1\ar[r]&\mathrm{W}_1^{\mathrm{G}_\Sigma^0}\ar[r]
&\mathrm{W}_1^{\mathrm{G}_\Sigma}\ar[r]&\pi_0(\mathrm{G}_\Sigma)\ar[r]
&1}.
\end{eqnarray}
Similar notations apply to $\mathrm{M}_\Sigma$. The subset of \eqref{refined expansion proof fourth set}
\begin{eqnarray}
\Bigg\{~ (\mathrm{M},\sigma):~ \mathrm{M}\in\mathcal{L}(\mathrm{M}_1) \textrm{ and } \sigma= \Sigma\cdot\mathrm{ad}(\mathrm{M}(\mathbb{Q})) ~\Bigg\} 
\end{eqnarray}
admits a natural action of $\mathrm{W}_1^{\mathrm{G}_\Sigma}$, where each pair $(\mathrm{M},\sigma)$ has stablizer $\mathrm{W}_1^{\mathrm{M}_\Sigma}$. The quotients
\begin{eqnarray}
\Bigg\{~ (\mathrm{M},\sigma):~ \begin{array}{l} \mathrm{M}\in\mathcal{L},\sigma \textrm{ is a semisimple $\mathrm{M}(\mathbb{Q})$-}\\ \textrm{-orbit in } \mathfrak{m}(\mathbb{Q}), \textrm{ $\Sigma$ is conjugate}\\ \textrm{to a point in $\sigma$ under $\mathrm{G}(\mathbb{Q})$. }  \end{array} ~\Bigg\}
\Bigg\slash\mathrm{W}_0^\mathrm{G}
\end{eqnarray}
and
\begin{eqnarray}
\Bigg\{~ (\mathrm{M},\sigma):~ \mathrm{M}\in\mathcal{L}(\mathrm{M}_1) \textrm{ and } \sigma= \Sigma\cdot\mathrm{ad}(\mathrm{M}(\mathbb{Q})) ~\Bigg\} 
\Bigg\slash\mathrm{W}_1^{\mathrm{G}_\Sigma}
\end{eqnarray}
are in bijection, therefore
\begin{eqnarray}
&&
J_\mathfrak{o}^\mathrm{G}(f) 
\\\nonumber\\\nonumber
&=& 
\sum_{\mathrm{M}\in\mathcal{L}(\mathrm{M}_1)} 
|\mathrm{W}_1^{\mathrm{M}_\Sigma}| |\mathrm{W}_1^{\mathrm{G}_\Sigma}|^{-1}
\sum_{\substack{X\in(\mathfrak{m} (\mathbb{Q})\cap\mathfrak{o})_{\mathrm{M},S_\mathfrak{o}}
\\X_\mathrm{ss}=\Sigma}} 
a^\mathrm{M}(S_\mathfrak{o},X) 
J_\mathrm{M}^\mathrm{G}(X,f_{S_\mathfrak{o}})
\\\nonumber\\\nonumber
&=&
\sum_{\mathrm{M}\in\mathcal{L}}  \frac{|\mathrm{W}_0^\mathrm{M}|}{|\mathrm{W}_0^\mathrm{G}|} \bigg(\frac{|\mathrm{W}_1^{\mathrm{M}_\Sigma}|}{|\mathrm{W}_1^{\mathrm{G}_\Sigma}|}\bigg)^{-1} |\mathrm{W}_1^{\mathrm{M}_\Sigma}| |\mathrm{W}_1^{\mathrm{G}_\Sigma}|^{-1} \times
\\\nonumber\\\nonumber
&&\quad\times
\sum_{X\in(\mathfrak{m} (\mathbb{Q}) \cap\mathfrak{o})_{\mathrm{M},S_\mathfrak{o}}} a^\mathrm{M}(S_\mathfrak{o},X) J_\mathrm{M}^\mathrm{G}(X,f_{S_\mathfrak{o}})   
\\\nonumber\\\nonumber 
 &=& 
\sum_{\mathrm{M}\in\mathcal{L}}|\mathrm{W}_0^\mathrm{M}||\mathrm{W}_0^\mathrm{G}|^{-1}
\sum_{X\in(\mathfrak{m}(\mathbb{Q})\cap\mathfrak{o})_{\mathrm{M},S_\mathfrak{o}}} a^\mathrm{M}(S_\mathfrak{o},X)J_\mathrm{M}^\mathrm{G}(X,f_{S_\mathfrak{o}}).
\end{eqnarray}
\qed

\paragraph{Remark}\setcounter{equation}{0}
For a semisimple element $X$ in $\mathfrak{o}$,
\begin{eqnarray}
\label{tamagawa number}
a^\mathrm{M}(S_\mathfrak{o},X) &=& |\pi_0(\mathrm{M}_X)|^{-1}\mathrm{Vol}\big(\mathrm{M}_X^0(\mathbb{Q})\backslash\mathrm{M}_X^0(\mathbb{A})^1\big)
\end{eqnarray}
if $X$ is $\mathbb{Q}$-elliptic in $\mathfrak{m}(\mathbb{Q})$, and vanishes otherwise. If $X$ is in addition assumed to be regular, the identity \eqref{refined expansion} reduces to \eqref{elementary refinement}.

\subsection{The refined trace formula}

\paragraph{Proposition}\setcounter{equation}{0}(Refined trace formula)\\
\label{refined trace formula statement}
\emph{For each sufficiently large finite set $S$ of places of $\mathbb{Q}$, for each $\sim$ equivalence class $\mathfrak{o}$ in $\mathfrak{g}(\mathbb{Q})$, for each $\equiv$ equivalence class $X$ in $(\mathfrak{m}(\mathbb{Q})\cap\mathfrak{o})_{\mathrm{M},S}$, there exists a constant $a^\mathrm{M}(S,X)$\label{sym107} such that}
\begin{eqnarray}
\label{refined trace formula}
&&
\forall f\in\mathcal{S}(\mathfrak{g}(\mathbb{A}))
\\\nonumber\\\nonumber
&&
\lim_S
\sum_{\mathfrak{o}\in\mathfrak{g}(\mathbb{Q})/\sim}
~\sum_{\mathrm{M}\in\mathcal{L}}
|\mathrm{W}_0^\mathrm{M}|
|\mathrm{W}_0^\mathrm{G}|^{-1}
\sum_{X\in(\mathfrak{m}(\mathbb{Q})\cap\mathfrak{o})_{\mathrm{M},S}}
a^\mathrm{M}(S,X)
J_\mathrm{M}^\mathrm{G}(X,f_S)
\\\nonumber\\\nonumber
&&\quad=
\lim_S
\sum_{\mathfrak{o}\in\mathfrak{g}(\mathbb{Q})/\sim}
~\sum_{\mathrm{M}\in\mathcal{L}}
|\mathrm{W}_0^\mathrm{M}|
|\mathrm{W}_0^\mathrm{G}|^{-1}
\sum_{X\in(\mathfrak{m}(\mathbb{Q})\cap\mathfrak{o})_{\mathrm{M},S}}
a^\mathrm{M}(S,X)
J_\mathrm{M}^\mathrm{G}(X,f\hat{~}_S).
\end{eqnarray}

\proof
By \eqref{chaudouard trace formula} it is enough to show that the left hand side of \eqref{refined trace formula} is equal to
\begin{eqnarray}
J(f)
&=&
\sum_{\mathfrak{o}\in\mathfrak{g}(\mathbb{Q})/\sim} J_\mathfrak{o}(f).
\end{eqnarray}
For each individual class $\mathfrak{o}$ by \eqref{refined expansion}
\begin{eqnarray}
J_\mathfrak{o}(f) &=& \lim_S \sum_\mathrm{M} |\mathrm{W}_0^\mathrm{M}||\mathrm{W}_0^\mathrm{G}|^{-1}
\sum_X a^\mathrm{M}(S,X)J_\mathrm{M}^\mathrm{G}(X,f_S)
\end{eqnarray}
where the limit stablizes as $S$ grows large enough.

\paragraph*{Lemma} 
\emph{Let $\Gamma$\label{sym108} be a compact subset of $\mathfrak{g}(\mathbb{A})$. Then there are only finitely many classes $\mathfrak{o}$ such that $\mathfrak{o}\cdot\mathrm{ad}(\mathrm{G}(\mathbb{A}))$ intersects $\Gamma$.}

\proof
The lemma follows from Corollary A.2 of \cite{Arthur_Orb} which states that there exists a compact set $\mathrm{G}_\Gamma$ contained in $\mathrm{G}(\mathbb{A})^1$ such that 
\begin{eqnarray}
\forall x\in\mathrm{G}(\mathbb{A})^1-\mathrm{G}(\mathbb{Q})\mathrm{G}_\Gamma
&&
\mathfrak{g}(\mathbb{Q})'\cdot\mathrm{ad}(x)\cap\Gamma=\emptyset  \end{eqnarray} 
where $\mathfrak{g}(\mathbb{Q})'$ denotes the set of points of $\mathfrak{g}(\mathbb{Q})$ not contained in any proper parabolic subalgebra. This is established by a reduction theory argument.

To prove the lemma first consider a class $\mathfrak{o}$ contained in $\mathfrak{g}(\mathbb{Q})'$. There exists a point $X$ in $\mathfrak{o}$ such that $X\cdot\mathrm{ad}(\mathrm{G}_\Gamma)$ is contained in $\Gamma$, hence there are only finitely many such $X$, hence finitely many such $\mathfrak{o}$. 

Next for an arbitrary class $\mathfrak{o}$ there is some Levi subalgebra $\mathfrak{m}$ such that $\mathfrak{o}$ intersects $\mathfrak{m}(\mathbb{Q})'$. Then by the Iwasawa decomposition
\begin{eqnarray}
&&
\mathfrak{o}\cdot\mathrm{ad}(\mathrm{G}(\mathbb{A})^1)
\\\nonumber
&=&
\mathfrak{o}\cdot\mathrm{ad}(\mathrm{M}(\mathbb{A})^1) \cdot\mathrm{ad}(\mathrm{N}(\mathbb{Q}) \backslash \mathrm{N}(\mathbb{A})) \cdot\mathrm{ad}(\mathrm{K})  
\\\nonumber 
 &=& 
 (\mathfrak{o}\cap\mathfrak{m}(\mathbb{Q})') \cdot\mathrm{ad}(\mathrm{M}_{\Gamma\cap\mathfrak{m}(\mathbb{A})}) \cdot\mathrm{ad}(\mathrm{N}(\mathbb{Q}) \backslash \mathrm{N}(\mathbb{A})) \cdot\mathrm{ad}(\mathrm{K}),
\end{eqnarray}
the lemma follows by induction.
\qed\\

The Schwartz function $f$ is compactly supported if its component at infinity $f_\infty$ is compactly supported. For such an $f$ the limit on the left hand side of \eqref{refined trace formula} stablizes for every $S$ such that
\begin{eqnarray}
S&\supset&\bigcup_{\mathfrak{o}\in\mathfrak{g}(\mathbb{Q})/\sim} S_\mathfrak{o}
\end{eqnarray}
where the union is supported on a finite set of $\mathfrak{o}$. Hence for a compactly supported $f$ the left hand side of \eqref{refined trace formula} is equal to $J(f)$. 

However the function spaces $C_c^\infty(\mathfrak{g}(\mathbb{R}))$ and hence $C_c^\infty(\mathfrak{g}(\mathbb{Q}_S))$ are dense in $\mathcal{S}(\mathfrak{g}(\mathbb{R}))$ and $\mathcal{S}(\mathfrak{g}(\mathbb{Q}_S))$ respectively, therefore $C_c^\infty(\mathfrak{g}(\mathbb{A}))$ is dense in $\mathcal{S}(\mathfrak{g}(\mathbb{A}))$ which by \eqref{global lie algebra schwartz space topology} is defined as
\begin{eqnarray}
\mathcal{S}(\mathfrak{g}(\mathbb{A}))
&=&
\varinjlim \mathcal{S}(\mathfrak{g}(\mathbb{Q}_S))
\end{eqnarray}
equipped with the final topology. 

Since $J(f)$ extends continuously to all Schwartz functions, the limit on the left hand side of \eqref{refined trace formula} exists,  and the interchangeability of the operations of taking the limit as $S$ approaches infinity and taking the sum over all classes $\mathfrak{o}$ extends from $C_c^\infty(\mathfrak{g}(\mathbb{A}))$ to $\mathcal{S}(\mathfrak{g}(\mathbb{A}))$. This establishes the identity \eqref{refined trace formula} for all Schwartz functions $f$ on $\mathfrak{g}(\mathbb{A})$.
\qed

\section{A simple invariant trace formula}

In this chapter following the arguments in \S7 of \cite{Arthur_Inv2} the refined trace formula \eqref{refined trace formula} is reduced to an identity between invariant distributions on $\mathfrak{g}(\mathbb{A})$ for a suitably restricted class of test functions.

\subsection{Parabolic descent and parabolic induction}

\paragraph{Definition}\setcounter{equation}{0}
Let $\mathrm{P}$ be a parabolic subgroup in $\mathcal{F}$ with Levi component $\mathrm{M}$ and unipotent radical $\mathrm{N}$. Let $f_S$ be a Schwartz function on $\mathfrak{g}(\mathbb{Q}_S)$. Define $f_{S,\mathrm{P}}$\label{sym114} to be the Schwartz function on $\mathfrak{m}(\mathbb{Q}_S)$ by \emph{parabolic descent} along $\mathrm{P}$
\begin{eqnarray}
\label{parabolic descent}
&&
\forall X\in\mathfrak{m}(\mathbb{Q}_S)
\\\nonumber\\\nonumber
&&
f_{S,\mathrm{P}}(X) 
=
\int_{\mathrm{K}_S}\int_{\mathfrak{n}(\mathbb{Q}_S)} f_S((X+N)\cdot\mathrm{ad}(k)) ~\mathrm{d}N\mathrm{d}k.
\end{eqnarray}
For each parabolic subgroup $\mathrm{P}_S$ in $\mathcal{F}^{\mathrm{G}_S}$ define the Schwartz function $f_{S,\mathrm{P}_S}$ on $\mathfrak{m}_S(\mathbb{Q}_S)$ analogously.

\paragraph{Remark}\setcounter{equation}{0}
Parabolic descent preserves orbital integrals, i.e.
\begin{eqnarray}
\label{parabolic descent orbital integral}
\forall X\in\mathfrak{m}_\mathrm{reg,ss}(\mathbb{Q}_S)
&&
I_\mathrm{G}^\mathrm{G}(X,f_S) 
=
I_\mathrm{M}^\mathrm{M}(X,f_{S,\mathrm{P}}).
\end{eqnarray}
See \S13.12 of \cite{Kott}. Hence the orbital integrals of $f_{S,\mathrm{P}}$ on $\mathfrak{m}_\mathrm{reg,ss}(\mathbb{Q}_S)$ is independent of the choice of the parabolic subgroup $\mathrm{P}$ in $\mathcal{P}(\mathrm{M})$. In this case denote $f_{S,\mathrm{P}}$ by the alternative notation $f_{S,\mathrm{M}}$ such that
\begin{eqnarray}
\forall X\in\mathfrak{m}_\mathrm{reg,ss}(\mathbb{Q}_S)
&&
I_\mathrm{M}^\mathrm{M}(X,f_{S,\mathrm{M}})
=
I_\mathrm{M}^\mathrm{M}(X,f_{S,\mathrm{P}})
\end{eqnarray}
is well-defined.

\paragraph{Remark}\setcounter{equation}{0}
With compatible choices of Fourier transforms on $\mathfrak{g}(\mathbb{Q}_S)$ and $\mathfrak{m}(\mathbb{Q}_S)$, parabolic descent intertwines the Fourier transforms, i.e. 
\begin{eqnarray}
\label{parabolic descent commutes with fourier transform}
(f_S\hat{~})_\mathrm{P}
&=&
(f_{S,\mathrm{P}})\hat{~}.
\end{eqnarray}
See \S13.13 of \cite{Kott}.

\paragraph{Definition}\setcounter{equation}{0}
Let $\mathrm{P}$ be a parabolic subgroup in $\mathcal{F}$ with Levi component $\mathrm{M}$ and unipotent radical $\mathrm{N}$. Let $X$ be an element in $\mathfrak{m}(\mathbb{Q}_S)$. Define $\mathrm{Ind}_{\mathrm{M},\mathrm{P}}^\mathrm{G}(X)$\label{sym115} to be the $\mathrm{ad}(\mathrm{G}(\mathbb{Q}_S))$-invariant subset of $\mathfrak{g}(\mathbb{Q}_S)$ by \emph{parabolic induction} along $\mathrm{P}$\begin{eqnarray}
\mathrm{Ind}_{\mathrm{M},\mathrm{P}}^\mathrm{G}(X)
&=&
\Big(\big(X\cdot\mathrm{ad}(\mathrm{M})+\mathfrak{n}\big)\cdot\mathrm{ad}(\mathrm{G})\Big)
_\mathrm{reg}(\mathbb{Q}_S)
\end{eqnarray}
where the subscript $\mathrm{reg}$ denotes the regular locus, the set of points with minimal dimensional isotropy group in $\mathrm{G}$.

\paragraph{Remark}\setcounter{equation}{0}
With a fixed choice of $\mathrm{M}$, parabolic induction along $\mathrm{P}$ is independent of the choice of the parabolic subgroup $\mathrm{P}$ in $\mathcal{P}(\mathrm{M})$, hence adopt the alternative notation $\mathrm{Ind}_\mathrm{M}^\mathrm{G}(X)$. For a proof see Satz 2.6 of \cite{Borho}. If the inducing Levi subgroup $\mathrm{M}$ is understood tacitly, denote $\mathrm{Ind}_\mathrm{M}^\mathrm{G}(X)$ by the alternative notation $X^\mathrm{G}$.

\paragraph{Remark}\setcounter{equation}{0}
If $X$ is regular semisimple then $\mathrm{Ind}_\mathrm{M}^\mathrm{G}(X)$ is equal to the $\mathrm{ad}(\mathrm{G}(\mathbb{Q}_S))$-orbit of $X$.
In general $\mathrm{Ind}_\mathrm{M}^\mathrm{G}(X)$ is a finite union of $\mathrm{ad}(\mathrm{G}(\mathbb{Q}_S))$-orbits which are geometrically conjugate to each other, but $\mathrm{Ind}_\mathrm{M}^\mathrm{G}(X)$ is still called the \emph{induced orbit} of $X$. See page 255 of \cite{Arthur_Loc} and \S2.1 of \cite{Borho}. If $\mathrm{L}$ is a Levi subgroup in $\mathcal{L}(\mathrm{M})$, define the invariant weighted orbital integral along the induced orbit $X^\mathrm{L}$ by
\begin{eqnarray}
\forall f_S\in\mathfrak{g}(\mathbb{Q}_S)
&&
J_\mathrm{L}^\mathrm{G}(X^\mathrm{L},f_S)
=
\sum_{\substack{\Omega\subset X^\mathrm{L}\\\mathrm{ad}(\mathrm{L}(\mathbb{Q}_S))\textrm{-orbit}}}
J_\mathrm{L}^\mathrm{G}(\Omega,f_S).
\end{eqnarray}

\paragraph{Remark}\setcounter{equation}{0}
Parabolic induction is transitive in nested chains of Levi subgroups, i.e.
\begin{eqnarray}
\forall \mathrm{L}\in\mathcal{L}^\mathrm{G}(\mathrm{M})
&&
\mathrm{Ind}_\mathrm{M}^\mathrm{G}(X)=\mathrm{Ind}_\mathrm{L}^\mathrm{G}
(\mathrm{Ind}_\mathrm{M}^\mathrm{L}(X)).
\end{eqnarray}
See \S2.3 of \cite{Borho}.

\paragraph{Remark}
Parabolic induction is compatible with Jordan decomposition, i.e.
\begin{eqnarray}
\label{parabolic induction jordan decomposition}
\mathrm{Ind}_\mathrm{M}^\mathrm{G}(X)
&=&
\Big(X_\mathrm{ss}+\mathrm{Ind}_{\mathrm{M}_{X_\mathrm{ss}}^0}^{\mathrm{G}_{X_\mathrm{ss}}^0}(X_\mathrm{nil})\Big)\cdot\mathrm{ad}(\mathrm{G}(\mathbb{Q}_S)).
\end{eqnarray}
See \S2.4 of \cite{Borho}. Also see Lemma 2 of \cite{Hoff_Ind}.

\paragraph{Lemma}\setcounter{equation}{0}
(Descent and splitting of ($\mathrm{G},\mathrm{M}$)-families)\\
\emph{Let $\mathrm{M}$ be a Levi subgroup in $\mathcal{L}$. There exist a function $d_\mathrm{M}^\mathrm{G}(~,~)$\label{sym116}
\begin{eqnarray}
d_\mathrm{M}^\mathrm{G}:\quad
\mathcal{L}(\mathrm{M})\times\mathcal{L}(\mathrm{M})
&\longrightarrow&
\mathbb{R},
\end{eqnarray}
and a partially defined map $s(~,~)$
\begin{eqnarray}
s:\quad
\mathcal{L}(\mathrm{M})\times\mathcal{L}(\mathrm{M})
&\longrightarrow&
\mathcal{F}(\mathrm{M})\times\mathcal{F}(\mathrm{M})
\end{eqnarray}
whose domain contains the pairs $(\mathrm{L}_1,\mathrm{L}_2)$ for which $d_\mathrm{M}^\mathrm{G}(\mathrm{L}_1,\mathrm{L}_2)$ is nonzero,  such that
\begin{itemize}
\item if $(\mathrm{L}_1,\mathrm{L}_2)$ is contained in the domain of $s$, then
\begin{eqnarray}
s(\mathrm{L}_1,\mathrm{L}_2)
&\in&
\mathcal{P}(\mathrm{L}_1)\times\mathcal{P}(\mathrm{L}_2);
\end{eqnarray}
\item if $(c_\mathrm{P})$ is a $(\mathrm{G},\mathrm{M})$-family and $\mathrm{L}$ is a Levi subgroup in $\mathcal{L}(\mathrm{M})$, then
\begin{eqnarray}
\label{(g,m)-family descent}
c_\mathrm{L} &=& \sum_{\mathrm{L}'\in\mathcal{L}(\mathrm{M})} d_\mathrm{M}^\mathrm{G}(\mathrm{L},\mathrm{L}')c_\mathrm{M}^{\mathrm{Q}'}
\end{eqnarray}
where $\mathrm{Q}'$ denotes the second component of $s(\mathrm{L},\mathrm{L}')$;
\item if $(c_\mathrm{P})$ and $(d_\mathrm{P})$ are $(\mathrm{G},\mathrm{M})$-families, then
\begin{eqnarray}
\label{(g,m)-family splitting}
(cd)_\mathrm{M} &=& \sum_{\mathrm{L}_1,\mathrm{L}_2\in\mathcal{L}(\mathrm{M})} d_\mathrm{M}^\mathrm{G}(\mathrm{L}_1,\mathrm{L}_2)c_\mathrm{M}^{\mathrm{Q}_1}d_\mathrm{M}^{\mathrm{Q}_2}
\end{eqnarray}
where 
\begin{eqnarray}
(\mathrm{Q}_1,\mathrm{Q}_2) &=& s(\mathrm{L}_1,\mathrm{L}_2).
\end{eqnarray}
\end{itemize}
Analogous results hold for the groups $\mathrm{G}_v$ and $\mathrm{G}_S$.}

\proof
See \S7 of \cite{Arthur_Inv1}.
\qed

\paragraph{Remark}\setcounter{equation}{0}
\label{angle jacobian}
The constant $d_\mathrm{M}^\mathrm{G}(\mathrm{L}_1,\mathrm{L}_2)$ is defined to be the volume in $\mathfrak{a}_\mathrm{M}^\mathrm{G}$ of the image of a fundamental parallelotope in the direct sum of $\mathfrak{a}_\mathrm{M}^{\mathrm{L}_1}$ and $\mathfrak{a}_\mathrm{M}^{\mathrm{L}_2}$ under the natural map
\begin{eqnarray}
\label{angle jacobian natural map}
\mathfrak{a}_\mathrm{M}^{\mathrm{L}_1} \oplus \mathfrak{a}_\mathrm{M}^{\mathrm{L}_2}
&\longrightarrow&
\mathfrak{a}_\mathrm{M}^\mathrm{G} 
\end{eqnarray}
if \eqref{angle jacobian natural map} is an isomorphism, and zero otherwise. In the former case $d_\mathrm{M}^\mathrm{G}(\mathrm{L}_1,\mathrm{L}_2)$ is equal to the volume in $\mathfrak{a}_{\mathrm{L}_2}^\mathrm{G}$ of the image of a fundamental parallelotope in $\mathfrak{a}_\mathrm{M}^{\mathrm{L}_1}$ under the natural isomorphism
\begin{eqnarray}
\mathfrak{a}_\mathrm{M}^{\mathrm{L}_1} &\longrightarrow& \mathfrak{a}_{\mathrm{L}_2}^\mathrm{G}.
\end{eqnarray}

\paragraph{Remark}\setcounter{equation}{0}
\label{section s}
The map $s$ depends on the choice of a vector $\xi$ in general position in $\mathfrak{a}_\mathrm{M}^\mathrm{G}$. Let $\mathrm{L}_1$ and $\mathrm{L}_2$ be Levi subgroups in $\mathcal{L}(\mathrm{M})$ such that \eqref{angle jacobian natural map} is an isomorphism, then
\begin{eqnarray}
\exists \xi_1\in\mathfrak{a}_{\mathrm{L}_1}^\mathrm{G}~
\exists \xi_2\in\mathfrak{a}_{\mathrm{L}_2}^\mathrm{G}
&&
\xi=\frac{\xi_1}{2}-\frac{\xi_2}{2}.
\end{eqnarray}
Define $s(\mathrm{L}_1,\mathrm{L}_2)$ to be $(\mathrm{Q}_1,\mathrm{Q}_2)$ where $\mathrm{Q}_i$ is the parabolic subgroup in $\mathcal{P}(\mathrm{L}_i)$ whose corresponding positive chamber in $\mathfrak{a}_{\mathrm{L}_i}^\mathrm{G}$ contains the vector $\xi_i$ where the index $i$ is 1 or 2.

\paragraph{Lemma}\setcounter{equation}{0}
(Descent and splitting of weighted orbital integrals)\\
\emph{Let $\mathrm{M}$ be a Levi subgroup in $\mathcal{L}$, let $X$ be an element of $\mathfrak{m}(\mathbb{Q}_S)$. Let $f_S$ be a Schwartz function on $\mathfrak{g}(\mathbb{Q}_S)$. Let $\xi$ be a vector in general position in $\mathfrak{a}_\mathrm{M}^\mathrm{G}$, let
\begin{eqnarray}
d_\mathrm{M}^\mathrm{G}
&:&
\mathcal{L}^\mathrm{G}(\mathrm{M})\times\mathcal{L}^\mathrm{G}(\mathrm{M})
\longrightarrow
\mathbb{R} 
\\\nonumber
s
&:&
\mathcal{L}^\mathrm{G}(\mathrm{M})\times\mathcal{L}^\mathrm{G}(\mathrm{M})
\longrightarrow
\mathcal{F}^\mathrm{G}(\mathrm{M})\times\mathcal{F}^\mathrm{G}(\mathrm{M}) 
\end{eqnarray}
be defined as in Remark \ref{angle jacobian} and Remark \ref{section s} with respect to $\xi$.
\begin{itemize}
\item Let $\mathrm{L}$ be a Levi subgroup in $\mathcal{L}^\mathrm{G}(\mathrm{M})$, then
\begin{eqnarray}
\label{weighted orbital integral descent}
J_\mathrm{L}^\mathrm{G}(X^\mathrm{L},f_S) 
&=& 
\sum_{\mathrm{L}'\in\mathcal{L}(\mathrm{M})} d_\mathrm{M}^\mathrm{G}(\mathrm{L},\mathrm{L}') J_\mathrm{M}^{\mathrm{L}'}(X,f_{S,\mathrm{Q}'})
\end{eqnarray}
where $\mathrm{Q}'$ denotes the second component of $s(\mathrm{L},\mathrm{L}')$.
\item Let $S$ be the set $\{v_1,v_2\}$. Let $f_S$ be of the form $f_{v_1}\otimes f_{v_2}$ where $f_{v_i}$ is a Schwartz function on $\mathfrak{g}(\mathbb{Q}_{v_i})$ where the index $i$ is 1 or 2. Then
\begin{eqnarray}
\label{weighted orbital integral splitting}
&&
J_\mathrm{M}^\mathrm{G}(X,f_S)
= 
\sum_{\mathrm{L}_1,\mathrm{L}_2\in\mathcal{L}^(\mathrm{M})} 
d_\mathrm{M}^\mathrm{G}(\mathrm{L}_1,\mathrm{L}_2) 
J_\mathrm{M}^{\mathrm{L}_1}(X,f_{v_1,\mathrm{Q}_1})
J_\mathrm{M}^{\mathrm{L}_2}(X,f_{v_2,\mathrm{Q}_2})
\end{eqnarray}
where 
\begin{eqnarray}
(\mathrm{Q}_1,\mathrm{Q}_2)&=&s(\mathrm{L}_1,\mathrm{L}_2).
\end{eqnarray}
\end{itemize}
Local identities analogous to (\ref{weighted orbital integral descent}) also hold for $\mathrm{G}_v$.}

\proof
For a regular semisimple $X$ this follows from \eqref{(g,m)-family descent} and \eqref{(g,m)-family splitting}. For a general $X$ this follows from the limit formula \eqref{general weighted orbital integral}.
\qed

\paragraph{Corollary}\setcounter{equation}{0}
(Descent and splitting of orbital integrals)\\
\emph{Let $X$ be an element of $\mathfrak{g}(\mathbb{Q}_S)$. Let $f_S$ be a Schwartz function on $\mathfrak{g}(\mathbb{Q}_S)$.
\begin{itemize}
\item Let $\mathrm{L}$ be a Levi subgroup in $\mathcal{L}$ such that $X$ is contained in $\mathfrak{l}(\mathbb{Q}_S)$, then
\begin{eqnarray}
\label{orbital integral descent}
I_\mathrm{G}^\mathrm{G}(X^\mathrm{G},f_S)
&=&
I_\mathrm{L}^\mathrm{L}(X,f_{S,\mathrm{L}}).
\end{eqnarray}
\item Let $S$ be the set $\{v_1,v_2\}$. Let $f_S$ be of the form $f_{v_1}\otimes f_{v_2}$ where $f_{v_i}$ is a Schwartz function on $\mathfrak{g}(\mathbb{Q}_{v_i})$ where the index $i$ is 1 or 2. Then
\begin{eqnarray}
\label{orbital integral splitting}
I_\mathrm{G}^\mathrm{G}(X,f_{v_1}\otimes f_{v_2})
&=&
I_\mathrm{G}^\mathrm{G}(X,f_{v_1})I_\mathrm{G}^\mathrm{G}(X,f_{v_2}).
\end{eqnarray}
\end{itemize}
Local identities analogous to (\ref{orbital integral descent}) also hold for $\mathrm{G}_v$.}

\proof
The identity \eqref{orbital integral descent} follows from \eqref{weighted orbital integral descent} where the only nonzero summand on the right hand side corresponds to the Levi subgroup $\mathrm{L}$ in $\mathcal{L}(\mathrm{L})$.

The identity \eqref{orbital integral splitting} follows from \eqref{weighted orbital integral splitting} where the only summand on the right hand side corresponds to the pair $(\mathrm{G},\mathrm{G})$ in $\mathcal{L}(\mathrm{G})\times\mathcal{L}(\mathrm{G})$.
\qed

\paragraph{Remark}\setcounter{equation}{0}
If $X$ is a regular semisimple element of $\mathfrak{m}(\mathbb{Q})$ then the local analogue of \eqref{weighted orbital integral descent} becomes
\begin{eqnarray}
\label{weighted orbital integral regular semisimple descent}
J_{\mathrm{L}_v}^{\mathrm{G}_v}(X,f_v)
&=& 
\sum_{\mathrm{L}'_v\in\mathcal{L}^{\mathrm{G}_v}(\mathrm{M}_v)} d_{\mathrm{M}_v}^{\mathrm{G}_v}(\mathrm{L}_v,\mathrm{L}'_v) J_{\mathrm{M}_v}^{\mathrm{L}'_v}(X,f_{v,\mathrm{Q}'_v}).
\end{eqnarray}
Repeatedly applying the identities \eqref{weighted orbital integral splitting} and \eqref{weighted orbital integral regular semisimple descent} reduces a regular semisimple weighted orbital integral $J_\mathrm{M}^\mathrm{G}(X,~)$ that appears generically in the refined trace formula \eqref{refined trace formula} to a linear combination of products of local weighted orbital integrals $J_{\mathrm{M}'_v}^{\mathrm{G}_v}(X,~)$ where $\mathrm{M}'_v$ is contained in $\mathrm{M}_v$ and $X$ is elliptic in $\mathfrak{m}'_v(\mathbb{Q}_v)$.

\subsection{A simple form of the trace formula}

\paragraph{Definition}\setcounter{equation}{0}
Let $v$ be a place of $\mathbb{Q}$, let $f_v$ be a Schwartz function on $\mathfrak{g}(\mathbb{Q}_v)$. The function $f_v$ is said to be \emph{cuspidal} if it is a finite linear combination of Schwartz functions $f_{v,i}$ on $\mathfrak{g}(\mathbb{Q}_v)$ such that
\begin{eqnarray}
\forall\mathrm{index}~i\quad\forall\mathrm{P}_v\in\mathcal{F}^{\mathrm{G}_v}-\{\mathrm{G}_v\}
&&
f_{v,i,\mathrm{P}_v}=0.
\end{eqnarray}
Denote the space of cuspidal Schwartz functions on $\mathfrak{g}(\mathbb{Q}_v)$ by $\mathcal{S}_\mathrm{cusp}(\mathfrak{g}(\mathbb{Q}_v))$.

\paragraph{Remark}\setcounter{equation}{0}
\label{cuspidality implies cuspidality}
A Schwartz function $f_v$ on $\mathfrak{g}(\mathbb{Q}_v)$ is cuspidal if
\begin{eqnarray}
I_\mathrm{G}^\mathrm{G}(X,f_v)
&=&
0
\end{eqnarray}
whenever $X$ is a regular semisimple element of $\mathfrak{g}(\mathbb{Q}_v)$ which is not $\mathbb{Q}_v$-elliptic.

\paragraph{Remark}\setcounter{equation}{0}
\label{fourier transform preserves cuspidality}
The space $\mathcal{S}_\mathrm{cusp}(\mathfrak{g}(\mathbb{Q}_v))$ is stable under the Fourier transform on $\mathfrak{g}(\mathbb{Q}_v)$ by \eqref{parabolic descent commutes with fourier transform}.

\paragraph{Proposition}\setcounter{equation}{0}
(Simple invariant trace formula)\\
\emph{For each sufficiently large finite set $S$ of places of $\mathbb{Q}$, for each $\sim$ equivalence class $\mathfrak{o}$ in $\mathfrak{g}(\mathbb{Q})$, for each $\equiv$ equivalence class $X$ in $(\mathfrak{m}(\mathbb{Q})\cap\mathfrak{o})_{\mathrm{M},S}$, there exists a constant $a^\mathrm{M}(S,X)$ such that for each Schwartz function $f$ on $\mathfrak{g}(\mathbb{A})$ which is cuspidal at two distinct places of $\mathbb{Q}$}
\begin{eqnarray}
\label{simple invariant trace formula}
&&
\lim_S
\sum_{\mathfrak{o}\in\mathfrak{g}(\mathbb{Q})/\sim}
~\sum_{X\in(\mathfrak{g}(\mathbb{Q})\cap\mathfrak{o})_{\mathrm{G},S}}
a^\mathrm{G}(S,X)
I_\mathrm{G}^\mathrm{G}(X,f_S)
\\\nonumber\\\nonumber
&=&
\lim_S
\sum_{\mathfrak{o}\in\mathfrak{g}(\mathbb{Q})/\sim}
~\sum_{X\in(\mathfrak{g}(\mathbb{Q})\cap\mathfrak{o})_{\mathrm{G},S}}
a^\mathrm{G}(S,X)
I_\mathrm{G}^\mathrm{G}(X,f\hat{~}_S).
\end{eqnarray}

\proof
This argument follows the proof of Theorem 7.1(b) of \cite{Arthur_Inv2}.

Let $v_1$ and $v_2$ denote two places at which $f$ is cuspidal. Without loss of generality the set $S$ contains both $v_1$ and $v_2$, and there exist two subsets $S_1$ and $S_2$ which partition $S$ such that $S_i$ contains $v_i$ where the index $i$ is 1 or 2.

With a choice of $\xi$ as in Remark \ref{section s}, it follows from \eqref{weighted orbital integral splitting} that
\begin{eqnarray}
\label{simple invariant trace formula proof first equality}
&&
J_\mathrm{M}^\mathrm{G}(X,f_S)
=
\sum_{\mathrm{L}_1,\mathrm{L}_2\in\mathcal{L}^(\mathrm{M})} 
d_\mathrm{M}^\mathrm{G}(\mathrm{L}_1,\mathrm{L}_2) 
J_\mathrm{M}^{\mathrm{L}_1}(X,f_{S_1,\mathrm{Q}_1})
J_\mathrm{M}^{\mathrm{L}_2}(X,f_{S_2,\mathrm{Q}_2})
\end{eqnarray}
where 
\begin{eqnarray}
(\mathrm{Q}_1,\mathrm{Q}_2)&=&s(\mathrm{L}_1,\mathrm{L}_2).
\end{eqnarray}
Since $f$ is cuspidal at $v_1$ and $v_2$, each summand on the right hand side of \eqref{simple invariant trace formula proof first equality} vanishes unless the parabolic subgroups $\mathrm{Q}_1$ and $\mathrm{Q}_2$ are both equal to $\mathrm{G}$. Hence the Levi subgroups $\mathrm{M}_1$ and $\mathrm{M}_2$ are both equal to $\mathrm{G}$, which implies that
\begin{eqnarray}
d_\mathrm{M}^\mathrm{G}(\mathrm{L}_1,\mathrm{L}_2)
&=&
d_\mathrm{M}^\mathrm{G}(\mathrm{G},\mathrm{G})
\end{eqnarray}
vanishes unless the Levi subgroup $\mathrm{M}$ is equal to $\mathrm{G}$ by Remark \ref{angle jacobian}.

Then \eqref{simple invariant trace formula} follows from the refined trace formula \eqref{refined trace formula} and Remark \ref{fourier transform preserves cuspidality}.
\qed

\paragraph{Remark}\setcounter{equation}{0}
Each distribution that appears in the summands of the simple invariant trace formula \eqref{simple invariant trace formula} factorizes into local distributions that are invariant under the adjoint action of $\mathrm{G}$ on $\mathfrak{g}$.

\section{The Harish-Chandra transform on the space\\of characteristic polynomials}

In this chapter an integral transform on the space of characteristic polynomials satisfying a summation formula of Poisson type is constructed.

\subsection{Preliminaries}

\paragraph{}\setcounter{equation}{0}
In this chapter $\mathrm{G}$\label{sym140} denotes the general linear group $\mathrm{GL}(n,\mathbb{Q})$ for some natural number $n$, with the standard choice of the minimal Levi subgroup $\mathrm{M}_0$ and the Borel subgroup $\mathrm{B}$ to be the subgroup consisting of the diagonal matrices and the upper triangular matrices respectively.

\paragraph{}\setcounter{equation}{0}
Let $\mathcal{A}_\mathrm{G}$\label{sym141} denote the affine space of characterstic polynomials of $n\times n$ matrices over $\mathbb{Q}$, i.e. 
\begin{eqnarray}
\mathcal{A}_\mathrm{G}
&=&
\mathrm{gl}(n,\mathbb{Q})/\!\!/\mathrm{GL}(n,\mathbb{Q})
\end{eqnarray}
where $/\!\!/$ denotes the affine quotient and $\mathrm{GL}(n,\mathbb{Q})$ acts on $\mathrm{gl}(n,\mathbb{Q})$ from the right by conjugation. The discriminant function $D$ on $\mathrm{gl}(n,\mathbb{Q})$ descends to a polynomial on $\mathcal{A}_\mathrm{G}$, denote by $\mathcal{A}_{\mathrm{G},\mathrm{reg}}$ the open subset where $D$ does not vanish.

Let $\mathrm{M}$ be a standard Levi subgroup of $\mathrm{G}$. Denote by $\mathcal{A}_\mathrm{M}$ the affine quotient of $\mathfrak{m}$ by the adjoint action of $\mathrm{M}$, then $\mathcal{A}_\mathrm{M}$ is an affine space and there exists a partition $(n_1,n_2,\dots,n_r)$ of $n$ such that
\begin{eqnarray}
\label{gln levi partition decomposition}
\mathcal{A}_\mathrm{M}
&=&
\mathcal{A}_{\mathrm{G}_1}\times\mathcal{A}_{\mathrm{G}_2}\times\dots\times\mathcal{A}_{\mathrm{G}_r}.
\end{eqnarray}
The embedding of $\mathfrak{m}$ into $\mathfrak{g}$ induces a map
\begin{eqnarray}
\pi_\mathrm{M}:\quad
\mathcal{A}_\mathrm{M}
&\longrightarrow&
\mathcal{A}_\mathrm{G}
\end{eqnarray}
which is finite of degree $|\mathrm{W}_0^\mathrm{M}|^{-1}|\mathrm{W}_0^\mathrm{G}|$ and \'{e}tale over $\mathcal{A}_{\mathrm{G},\mathrm{reg}}$.

\paragraph{}\setcounter{equation}{0}
Let $v$ be a place of $\mathbb{Q}$. Equip $\mathcal{A}_\mathrm{G}(\mathbb{Q}_v)$ with the measure\label{sym142}
\begin{eqnarray}
|D(X_v)|_v^{-1/2}\mathrm{d}X_v
\end{eqnarray}
where $\mathrm{d}X_v$ denotes the standard translation invariant measure on the affine space $\mathcal{A}_\mathrm{G}(\mathbb{Q}_v)$. This is equal to the pushforward of the translation invariant measure on $\mathfrak{m}_0(\mathbb{Q}_v)$ along the Chevalley morphism
\begin{eqnarray}
\pi_{\mathrm{M}_0}:\quad
\mathfrak{m}_0=\mathcal{A}_{\mathrm{M}_0}
&\longrightarrow&
\mathfrak{m}_0/\!\!/\mathrm{W}_0^\mathrm{G}=\mathcal{A}_\mathrm{G}.
\end{eqnarray}
The complement of $\mathcal{A}_{\mathrm{G},\mathrm{reg}}(\mathbb{Q}_v)$ in $\mathcal{A}_\mathrm{G}(\mathbb{Q}_v)$ is a null set.

For each standard Levi subgroup $\mathrm{M}$ of $\mathrm{G}$ equip $\mathcal{A}_\mathrm{M}(\mathbb{Q}_v)$ with the product measure 
\begin{eqnarray}
\prod_{i=1}^r
|D^{\mathrm{G}_i}(X_{i,v})|_v^{-1/2}\mathrm{d}X_{i,v}
\end{eqnarray}
where $\mathrm{G}_1,\mathrm{G}_2,\dots,\mathrm{G}_r$ are related to $\mathrm{M}$ as in \eqref{gln levi partition decomposition}.

Denote by $\mathcal{A}_{\mathrm{M},\mathrm{reg}}(\mathbb{Q}_v)_\mathrm{ell}$\label{sym143} the subset of $\mathcal{A}_{\mathrm{M},\mathrm{reg}}(\mathbb{Q}_v)$ consisting of the images of the $\mathbb{Q}_v$-elliptic elements in $\mathfrak{m}(\mathbb{Q}_v)$. Then
\begin{eqnarray}
\mathcal{A}_{\mathrm{G},\mathrm{reg}}(\mathbb{Q}_v) 
&=& 
\coprod_{\mathrm{M}\in\mathcal{L}}
\pi_\mathrm{M}(\mathcal{A}_{\mathrm{M},\mathrm{reg}}(\mathbb{Q}_v)_\mathrm{ell}).
\end{eqnarray}
This decomposition is compatible with the measures on $\mathcal{A}_\mathrm{G}(\mathbb{Q}_v)$ and $\mathcal{A}_\mathrm{M}(\mathbb{Q}_v)$.

The local measures induce $S$-local and global measures on $\mathcal{A}_\mathrm{G}(\mathbb{Q}_S)$ and $\mathcal{A}_\mathrm{G}(\mathbb{A})$.

\paragraph{Definition}\setcounter{equation}{0}
Let $S$ be a finite set of places of $\mathbb{Q}$, let $\mathrm{M}$ be a standard Levi subgroup of $\mathrm{G}$, let $X$ be an element of $\mathfrak{m}(\mathbb{Q}_S)$. Define the \emph{invariant weighted orbital integral $I_\mathrm{M}^\mathrm{G}(X,~)$}\label{sym109} to be the distribution on $\mathfrak{g}(\mathbb{Q}_S)$ such that
\begin{eqnarray}
\label{invariant weighted orbital integral definition}
&&
\forall f_S\in\mathcal{S}(\mathfrak{g}(\mathbb{Q}_S))
\\\nonumber\\\nonumber
&&
I_\mathrm{M}^\mathrm{G}(X,f_S)
\\\nonumber
&&\quad
= 
J_\mathrm{M}^\mathrm{G}(X,f_S)
-\sum_{\substack{ \mathrm{L}\in\mathcal{L}(\mathrm{M})\\\mathrm{L}\neq\mathrm{G}}} 
\Bigg( \sum_{\mathrm{L}'_S\in\mathcal{L}^{\mathrm{L}_S}} 
|\mathrm{W}_{S,0}^{\mathrm{L}'_S}| 
|\mathrm{W}_{S,0}^{\mathrm{L}_S}|^{-1} \times
\\\nonumber\\\nonumber
&&\qquad\times
\sum_{\mathrm{T}_S\in\mathcal{T}_\mathrm{ell}(\mathrm{L}'_S)} 
|\mathrm{W}(\mathrm{L}'_S,\mathrm{T}_S)|^{-1} \times \\\nonumber\\\nonumber
&&\qquad\quad \times
\int_{\mathfrak{t}_S(\mathbb{Q}_S)} 
J_\mathrm{L}^\mathrm{G}(Y,f_S\check{~})
\Big(I_\mathrm{M}^\mathrm{L}(X,~)\hat{~}(Y)\Big)
~|D^\mathrm{L}(Y)|_S^{1/2}\mathrm{d}Y \Bigg)
\end{eqnarray}
where $I_\mathrm{M}^\mathrm{M}(X,~)$ denotes the standard orbital integral on $\mathfrak{m}(\mathbb{Q}_S)$ defined by \eqref{standard orbital integral}.

\paragraph{Remark}\setcounter{equation}{0}
The invariant weighted orbital integral $I_\mathrm{M}^\mathrm{G}(X,~)$ is well-defined.
\begin{itemize}
\item The distribution $I_\mathrm{M}^\mathrm{G}(X,~)$ does not depend on the choice of the Fourier transforms on $\mathfrak{g}(\mathbb{Q}_S)$ and its Levi subalgebras $\mathfrak{l}(\mathbb{Q}_S)$ as long as these are compatible. See \cite{Wal} Lemme VI.5 for the $p$-adic case, the same argument also works for the real and $S$-local cases.
\item 
\label{invariant weighted orbital integral growth}
The integral in \eqref{invariant weighted orbital integral definition} converges since the distribution $I_\mathrm{M}^\mathrm{L}(X,~)\hat{~}$ is represented by a smooth function supported on $\mathfrak{g}_\mathrm{reg.ss}(\mathbb{Q}_S)$ such that the function
\begin{eqnarray}
Y &\mapsto& \Big(I_\mathrm{M}^\mathrm{L}(X,~)\hat{~}(Y)\Big)|D^\mathrm{L}(Y)|_S^{1/2}
\end{eqnarray}
is locally 
\begin{eqnarray}
O\Big(\max\big\{1,-\log(|D^\mathrm{G}(Y)|_S)^N\big\}\Big)
\end{eqnarray}
on $\mathfrak{t}(\mathbb{Q}_S)$ for some natural number $N$ and tempered at infinity, and the function
\begin{eqnarray}
Y \mapsto J_\mathrm{L}^\mathrm{G}(Y,f\check{~}_S)
\end{eqnarray}
is locally 
\begin{eqnarray}
O\Big(\max\big\{1,-\log(|D^\mathrm{G}(Y)|_S)^M\big\}\Big)
\end{eqnarray}
on $\mathfrak{t}(\mathbb{Q}_S)$ for some natural number $M$ and rapidly decreasing at infinity. For the $p$-adic case see Lemme VI.3(iv) and Corollaire III.6 of \cite{Wal}. For the real case see Proposition 9 on page 108 of \cite{Var} and Corollary 7.4 of \cite{Arthur_Discrete}.
\end{itemize}

\paragraph{Proposition}\setcounter{equation}{0}
\emph{Let $\mathrm{M}$ be a standard Levi subgroup of $\mathrm{G}$, let $X$ be an element of $\mathfrak{m}(\mathbb{Q})$. Then $I_\mathrm{M}^\mathrm{G}(X,~)$ is an invariant distribution on $\mathfrak{g}(\mathbb{Q}_S)$, i.e.}
\begin{eqnarray}
&&
\forall x\in\mathrm{G}(\mathbb{Q}_S)~
\forall f_S\in\mathcal{S}(\mathfrak{g}(\mathbb{Q}_S))
\\\nonumber
&&\quad
I_\mathrm{M}^\mathrm{G}(X,f_S\circ\mathrm{ad}(x))=I_\mathrm{M}^\mathrm{G}(X,f_S).
\end{eqnarray}

\proof
For the $p$-adic case see Proposition VI.1 of \cite{Wal}, the same argument also works for the $S$-local case.
\qed

\paragraph{Remark}\setcounter{equation}{0}
\label{duality remark}
To quote Waldspurger from {\S}Introduction of \cite{Wal}:\\
\begin{doublespace}
`` \emph{ On dispose de deux ensembles de distributions invariantes sur $\mathrm{G}$:}

\emph{(1) les int\'{e}grales orbitales associ\'{e}es aux \'{e}l\'{e}ments semi-simples de $\mathrm{G}$;}

\emph{(2) les caract\`{e}res de repr\'{e}sentations temp\'{e}r\'{e}es irr\'{e}ductibles de $\mathrm{G}$.}\\
( \dots \dots )
\end{doublespace}

\emph{Rempla\c{c}ons $G$ par son alg\`{e}bre de Lie $\mathfrak{g}$ et consid\'{e}rons l'espace de distributions invariantes par l'action adjointe de $G$. L'ensemble (1) a un analogue \'{e}vident: les int\'{e}grales orbitales associ\'{e}es aux \'{e}l\'{e}ments semi-simples de $\mathfrak{g}$. Le seul but de cet article est de fournir un support un peu consistant \`{a} l'id\'{e}e, d'ailleurs banale, que l'analogue de (2) est l'ensemble des transform\'{e}es de Fourier des int\'{e}grales orbitales pr\'{e}c\'{e}dents (invariantes). } ''\\

The duality between the orbital integrals and the Fourier transforms of the invariant weighted orbital integrals is embodied in the local trace formulae on the Lie algebra $\mathfrak{g}$.

\paragraph{Proposition}\setcounter{equation}{0}
(Local invariant trace formula of Waldspurger)\\
\emph{Let $v$ be a place of $\mathbb{Q}$, let $f_v$ and $g_v$ be Schwartz functions on $\mathfrak{g}_v(\mathbb{Q}_v)$, then}
\begin{eqnarray}
\label{local invariant trace formula}
&&
\sum_{\mathrm{M}_v\in\mathcal{L}^{\mathrm{G}_v}} 
|\mathrm{W}_{v,0}^{\mathrm{M}_v}| 
|\mathrm{W}_{v,0}^{\mathrm{G}_v}|^{-1} 
\sum_{\mathrm{T}_v\in\mathcal{T}_\mathrm{ell}(\mathrm{M}_v)} 
|\mathrm{W}(\mathrm{M}_v,\mathrm{T}_v)|^{-1} \times
\\\nonumber\\\nonumber
&&\quad 
\times \int_{\mathfrak{t}_v(\mathbb{Q}_v)} 
(-1)^{\dim(\mathrm{A}_{\mathrm{M}_v}/\mathrm{A}_{\mathrm{G}_v})} 
I_{\mathrm{M}_v}^{\mathrm{G}_v}(X_v,f_v\hat{~})
I_{\mathrm{G}_v}^{\mathrm{G}_v}(X_v,g_v)  ~\mathrm{d}X_v
\\\nonumber\\\nonumber
&=&
\sum_{\mathrm{M}_v\in\mathcal{L}^{\mathrm{G}_v}} 
|\mathrm{W}_{v,0}^{\mathrm{M}_v}| 
|\mathrm{W}_{v,0}^{\mathrm{G}_v}|^{-1} 
\sum_{\mathrm{T}_v\in\mathcal{T}_\mathrm{ell}(\mathrm{M}_v)} 
|\mathrm{W}(\mathrm{M}_v,\mathrm{T}_v)|^{-1} \times
\\\nonumber\\\nonumber
&&\quad 
\times \int_{\mathfrak{t}_v(\mathbb{Q}_v)} 
(-1)^{\dim(\mathrm{A}_{\mathrm{M}_v}/\mathrm{A}_{\mathrm{G}_v})} 
I_{\mathrm{M}_v}^{\mathrm{G}_v}(X_v,g_v\hat{~})
I_{\mathrm{G}_v}^{\mathrm{G}_v}(X_v,f_v)  ~\mathrm{d}X_v.
\end{eqnarray}

\proof
For the $p$-adic case see Th\'{e}or\`{e}me VII.1 of \cite{Wal}. The argument only uses combinatorial properties of $(\mathrm{G}_v,\mathrm{M}_v)$-families and standard results from harmonic analysis on a $\mathbb{Q}_v$-vector space, hence works equally well in the real case.
\qed

\paragraph{Definition}\setcounter{equation}{0}
Let $S$ be a finite set of places of $\mathbb{Q}$, let $\mathrm{M}$ be a Levi subgroup in $\mathcal{L}$, let $X$ be an element of $\mathfrak{m}(\mathbb{Q}_S)$. Define the \emph{vector-valued orbital integral} $\mathcal{I}_\mathrm{M}^\mathrm{G}(X,~)$\label{sym144} to be the vector-valued distribution on $\mathfrak{g}(\mathbb{Q}_S)$ taking values in $\bigoplus_{\mathrm{Q}\in\mathcal{F}(\mathrm{M})}\mathbb{C}$ such that
\begin{eqnarray}
&&
\forall f_S\in\mathcal{S}(\mathfrak{g}(\mathbb{Q}_S)
\\\nonumber\\\nonumber
&&
\mathcal{I}_\mathrm{M}^\mathrm{G}(X,f_S)=\bigg(
(-1)^{\mathrm{dim}(\mathrm{A}_\mathrm{M}/\mathrm{A}_{\mathrm{M}_\mathrm{Q}})}
I_\mathrm{M}^{\mathrm{M}_\mathrm{Q}}(X,f_{S,\mathrm{Q}})
\bigg)_{\mathrm{Q}\in\mathcal{F}(\mathrm{M}).}
\end{eqnarray}

\paragraph{Lemma}\setcounter{equation}{0}
(Induction and splitting)\\
\emph{Let $S$ be a finite set of places of $\mathbb{Q}$, let $\mathrm{M}$ be a Levi subgroup in $\mathcal{L}$, let $X$ be an element of $\mathfrak{m}(\mathbb{Q}_S)$, let $f_S$ be a Schwartz function on $\mathfrak{g}(\mathbb{Q}_S)$.
\begin{itemize}
\item If $\mathrm{L}$ is a Levi subgroup in $\mathcal{L}(\mathrm{M})$, then $\mathcal{I}_\mathrm{L}^\mathrm{G}(X,f_S)$ is completely determined by $\mathcal{I}_\mathrm{M}^\mathrm{G}(X,f_S)$.
\item If $S$ is the set $\{v_1,v_2,\dots,v_r\}$ for some natural number $r$ and
\begin{eqnarray}
&&
\forall i=1,2,\dots,r~
\exists f_{v_i}\in\mathcal{S}(\mathfrak{g}(\mathbb{Q}_{v_i}))
\\\nonumber
&&\quad
f_S=f_{v_1}\otimes f_{v_2}\otimes\dots\otimes f_{v_r},
\end{eqnarray}
then $\mathcal{I}_\mathrm{M}^\mathrm{G}(X,f_S)$ is completely determined by $\mathcal{I}_\mathrm{M}^\mathrm{G}(X_{v_i},f_{v_i})$ where $i$ ranges among $1,2,\dots,r$.
\end{itemize}}

\proof
This follows from analogues of the descent and splitting identities \eqref{weighted orbital integral descent} and \eqref{weighted orbital integral splitting} for the invariant weighted orbital integrals $I_\mathrm{M}^\mathrm{G}(X,~)$, which could be deduced by the same arguments as in \S7 of \cite{Arthur_Inv1}.
\qed

\paragraph{Lemma}\setcounter{equation}{0}
(Local trace formula)\\
\emph{Let $v$ be a place of $\mathbb{Q}$, let $f_v$ and $g_v$ be Schwartz functions on $\mathfrak{g}(\mathbb{Q}_v)$, then}
\begin{eqnarray}
\label{harish chandra local trace formula}
&&
\sum_{\mathrm{M}\in\mathcal{L}}
|\mathrm{W}_0^\mathrm{M}||\mathrm{W}_0^\mathrm{G}|^{-1} 
\sum_{\mathrm{T}_v\in\mathcal{T}_\mathrm{ell}(\mathrm{M})} 
|\mathrm{W}(\mathrm{M},\mathrm{T}_v)|^{-1}\times
\\\nonumber\\\nonumber
&&\quad
\times\int_{\mathfrak{t}_v(\mathbb{Q}_v)} 
\mathcal{I}_\mathrm{M}^\mathrm{G}(X_v,f_v\hat{~})
I_\mathrm{G}^\mathrm{G}(X_v,g_v)  ~\mathrm{d}X_v
\\\nonumber\\\nonumber
&=&
\sum_{\mathrm{M}\in\mathcal{L}}
|\mathrm{W}_0^\mathrm{M}||\mathrm{W}_0^\mathrm{G}|^{-1} 
\sum_{\mathrm{T}_v\in\mathcal{T}_\mathrm{ell}(\mathrm{M})} 
|\mathrm{W}(\mathrm{M},\mathrm{T}_v)|^{-1}\times
\\\nonumber\\\nonumber
&&\quad
\times
\int_{\mathfrak{t}_v(\mathbb{Q}_v)} 
\mathcal{I}_\mathrm{M}^\mathrm{G}(X_v,g_v\hat{~})
I_\mathrm{G}^\mathrm{G}(X_v,f_v)  ~\mathrm{d}X_v.
\end{eqnarray}

\proof
This follows from \eqref{local invariant trace formula}.

\paragraph{Lemma}\setcounter{equation}{0}
(Global trace formula)\\
\emph{Let $f$ be a Schwartz function on $\mathfrak{g}(\mathbb{A})$ which is cuspidal at two distinct places of $\mathbb{Q}$, then
\begin{eqnarray}
\label{harish chandra global trace formula}
&&
\lim_S
\sum_{\mathfrak{o}\in\mathfrak{g}(\mathbb{Q})/\sim}
~\sum_{X\in(\mathfrak{g}(\mathbb{Q})\cap\mathfrak{o})_{\mathrm{G},S}}
a^\mathrm{G}(S,X)
I_\mathrm{G}^\mathrm{G}(X,f_S)
\\\nonumber\\\nonumber
&=&
\lim_S
\sum_{\mathfrak{o}\in\mathfrak{g}(\mathbb{Q})/\sim}
~\sum_{X\in(\mathfrak{g}(\mathbb{Q})\cap\mathfrak{o})_{\mathrm{G},S}}
a^\mathrm{G}(S,X)
\mathcal{I}_\mathrm{G}^\mathrm{G}(X,f\hat{~}_S)
\end{eqnarray}
where the limit is taken over finite sets $S$ of places of $\mathbb{Q}$, and the coefficients $a^\mathrm{M}(S,X)$ are defined as in \eqref{global coefficient}.}

\proof
This follows from \eqref{simple invariant trace formula}.
\qed

\subsection{The local and global Schwartz spaces}

\paragraph{Definition}\setcounter{equation}{0}
Let $v$ be a place of $\mathbb{Q}$. Let $X_v$ be an element of $\mathcal{A}_{\mathrm{G},\mathrm{reg}}(\mathbb{Q}_v)$. Define the \emph{maximal orbital integral} $\mathcal{I}_\mathrm{max}^\mathrm{G}(X_v,~)$\label{sym146} to be the vector-valued distribution on $\mathfrak{g}(\mathbb{Q}_v)$ such that
\begin{eqnarray}
\forall f_v\in\mathcal{S}(\mathfrak{g}(\mathbb{Q}_v))
&&
\mathcal{I}_\mathrm{max}^\mathrm{G}(X_v,f_v) 
=
\mathcal{I}_{\mathrm{M}[\widetilde{X}_v\textrm{-ell}]}^\mathrm{G}(\widetilde{X}_v,f_v) 
\end{eqnarray}
where $\widetilde{X}_v$\label{sym147} is a regular semisimple element in $\mathfrak{g}(\mathbb{Q}_v)$ lifting $X_v$, and $\mathrm{M}[\widetilde{X}_v\textrm{-ell}]$ is a standard Levi subgroup of $\mathrm{G}$ such that $\mathfrak{m}[\widetilde{X}_v\textrm{-ell}](\mathbb{Q}_v)$ contains $\widetilde{X}_v$ as a $\mathbb{Q}_v$-elliptic element.

\paragraph{Definition}\setcounter{equation}{0}
Let $v$ be a place of $\mathbb{Q}$. Define the \emph{local Schwartz space} $\mathcal{S}_0(\mathcal{A}_\mathrm{G}(\mathbb{Q}_v))$\label{sym148} to be the space of complex-valued functions on $\mathcal{A}_{\mathrm{G},\mathrm{reg}}(\mathbb{Q}_v)$ such that
\begin{eqnarray}
\label{local schwartz space 0}
&&
\varphi_v\in\mathcal{S}_0(\mathcal{A}_\mathrm{G}(\mathbb{Q}_v))
\\\nonumber
&\Leftrightarrow&
\exists f_v\in\mathcal{S}(\mathfrak{g}(\mathbb{Q}_v))~
\forall X_v\in\mathcal{A}_{\mathrm{G},\mathrm{reg}}(\mathbb{Q}_v)
\\\nonumber
&&\quad
\varphi_v(X_v) = I_\mathrm{G}^\mathrm{G}(X_v,f_v).
\end{eqnarray}
Define the \emph{local Schwartz space} $\mathcal{S}_1(\mathcal{A}_\mathrm{G}(\mathbb{Q}_v))$ to be the space of vector-valued functions on $\mathcal{A}_{\mathrm{G},\mathrm{reg}}(\mathbb{Q}_v)$ such that
\begin{eqnarray}
\label{local schwartz space 1}
&&
\varphi_v\in\mathcal{S}_1(\mathcal{A}_\mathrm{G}(\mathbb{Q}_v))
\\\nonumber
&\Leftrightarrow&
\exists f_v\in\mathcal{S}(\mathfrak{g}(\mathbb{Q}_v))~
\forall X_v\in\mathcal{A}_{\mathrm{G},\mathrm{reg}}(\mathbb{Q}_v)
\\\nonumber
&&\quad
\varphi_v(X_v) = \mathcal{I}_\mathrm{max}^\mathrm{G}(X_v,f_v).
\end{eqnarray}

\paragraph{Definition}\setcounter{equation}{0}
Let $v$ be a place of $\mathbb{Q}$.  If $v$ is $p$-adic, denote by $\Lambda_0$\label{sym149} the standard lattice 
\begin{eqnarray}
\Lambda_0=\mathrm{gl}(n,\mathbb{Z})(\mathbb{Z}_p)
&\subset&
\mathrm{gl}(n,\mathbb{Q})(\mathbb{Q}_p)
\end{eqnarray}
and denote by $\mathbb{I}_{\Lambda_0}$ its characteristic function. If $v$ is archimedean, denote by $\mathcal{E}$ the Gaussian function on $\mathfrak{g}(\mathbb{R})$ which is self-dual with respect to the unitary Fourier transform, i.e.
\begin{eqnarray*}
\forall M\in\mathfrak{g}(\mathbb{R})
&&
\mathcal{E}(M)=e^{-\pi \mathrm{Tr}(M^\mathrm{T}M)} 
\end{eqnarray*}
where $M^\mathrm{T}$ denotes the transpose of $M$ and $\mathrm{Tr}$ denotes the trace of an $n\times n$ matrix.

Define the \emph{basic function} $\phi_{0,v}$\label{sym150} to be the element of $\mathcal{S}_0(\mathcal{A}_\mathrm{G}(\mathbb{Q}_v))$ such that
\begin{eqnarray}
\label{local basic function 0}
&&
\forall X_v\in\mathcal{A}_{\mathrm{G},\mathrm{reg}}(\mathbb{Q}_v)
\\\nonumber\\\nonumber
&&
\phi_{0,v}(X_v)=\left\{\begin{array}{ll}
I_\mathrm{G}^\mathrm{G}(X_v,\mathbb{I}_{\Lambda_0})
&\textrm{if $v$ is $p$-adic},\\\\
I_\mathrm{G}^\mathrm{G}(X_v,\mathcal{E})
&\textrm{if $v$ is archimedean}.
\end{array}\right.
\end{eqnarray}
Define the \emph{basic function} $\phi_{1,v}$ to be the element of $\mathcal{S}_1(\mathcal{A}_\mathrm{G}(\mathbb{Q}_v))$ such that
\begin{eqnarray}
\label{local basic function 1}
&&
\forall X_v\in\mathcal{A}_{\mathrm{G},\mathrm{reg}}(\mathbb{Q}_v)
\\\nonumber\\\nonumber
&&
\phi_{1,v}(X_v)=\left\{\begin{array}{ll}
\mathcal{I}_\mathrm{max}^\mathrm{G}(X_v,\mathbb{I}_{\Lambda_0})
&\textrm{if $v$ is $p$-adic},\\\\
\mathcal{I}_\mathrm{max}^\mathrm{G}(X_v,\mathcal{E})
&\textrm{if $v$ is archimedean}.
\end{array}\right.
\end{eqnarray}

\paragraph{Lemma}\setcounter{equation}{0}
If $v$ is a $p$-adic place, then a complex-valued function $\varphi_v$ on $\mathcal{A}_{\mathrm{G},\mathrm{reg}}(\mathbb{Q}_v)$ belongs to $\mathcal{S}_0(\mathcal{A}_\mathrm{G}(\mathbb{Q}_v))$ if and only if the following conditions are satisfied:
\begin{itemize}
\item $\varphi_v$ is a locally constant function on $\mathcal{A}_{\mathrm{G},\mathrm{reg}}(\mathbb{Q}_v)$;
\item after extending by zero to $\mathcal{A}_\mathrm{G}(\mathbb{Q}_v)$, $\varphi_v$ is compactly supported;
\item for each singular point $Z_v$\label{sym151} in the complement of $\mathcal{A}_{\mathrm{G},\mathrm{reg}}(\mathbb{Q}_S)$ in $\mathcal{A}_\mathrm{G}(\mathbb{Q}_v)$, there exists an open neighborhood $U_{\varphi_v,Z_v}$ of $Z_v$ in $\mathcal{A}_\mathrm{G}(\mathbb{Q}_v)$ such that 
\begin{eqnarray}
&&
\varphi_v(~)\big|_{U_{\varphi_v,Z_v}}
\in
\mathrm{span}\Bigg\{
\Gamma_\mathrm{G}^\mathrm{G}(~,\nu)\big|_{U_{\varphi_v,Z_v}}:~ 
\begin{array}{l}\nu\subset\mathfrak{g}_{\widetilde{Z}_v}(\mathbb{Q}_v)\\
\textrm{nilpotent orbit}\end{array} \Bigg\}
\end{eqnarray}
where $\varphi\big|_U$ denotes the restriction of the function $\varphi$ to the open subset $U$, the point $\widetilde{Z}_v$ is a regular semisimple element in $\mathfrak{g}(\mathbb{Q}_v)$ lifting $Z_v$, and $\Gamma_\mathrm{G}^\mathrm{G}(~,\nu)$ denotes the Shalika germ at the nilpotent orbit $\nu$, see \S17 of \cite{Kott}.
\end{itemize}

\proof
Away from the singular locus the function $\varphi_v$ is a linear combination of the characteristic functions of small compact open sets, and these all have lifts in $\mathcal{S}(\mathfrak{g}(\mathbb{Q}_v))$. 

Near the singular locus the theory of Shalika germs implies that the asymptotic behavior of $\varphi_v$ near the singular point $Z_v$ is necessary. It remains to show that all the Shalika germs appear in $\mathcal{S}_0(\mathcal{A}_\mathrm{G}(\mathbb{Q}_v))$, in other words there is no nontrivial linear relation among the possible asymptotes of $I_\mathrm{G}^\mathrm{G}(X_v,f_v)$ as $f_v$ ranges over $\mathcal{S}(\mathfrak{g}(\mathbb{Q}_v))$. The Shalika germ associated to the nilpotent orbit $\nu$ is homogeneous of degree equal to the codimension of $\nu$ in the nilpotent locus, hence the only possible linear relations are among the Shalika germs associated to the nilpotent orbits of the same dimension. But such nilpotent orbits are separated, so there are Schwarts functions $f_v$ on $\mathfrak{g}(\mathbb{Q}_v)$ that vanish on all but one of the nilpotent orbits of a given dimension, hence these too are linearly independent.
\qed

\paragraph{Remark}\setcounter{equation}{0}
Let $v$ be a $p$-adic place. It is conjectured by Jacquet for a general reductive Lie algebra and proven by Waldspurger for $\mathrm{gl}(n)$ that a Schwartz function $f_v$ on $\mathfrak{g}(\mathbb{Q}_v)$ has the property that
\begin{eqnarray}
&&
\exists\lambda\in\mathbb{C}~
\forall X_v\in\mathcal{A}_{\mathrm{G},\mathrm{reg}}(\mathbb{Q}_v)
\\\nonumber\\\nonumber
&&\quad
I_\mathrm{G}^\mathrm{G}(X_v,f_v)=\lambda\cdot\phi_{0,v}(X_v)
\end{eqnarray}
if and only if $I_\mathrm{G}^\mathrm{G}(X_v,f_v)$ and $I_\mathrm{G}^\mathrm{G}(X_v,f_v\hat{~})$ are both supported on the subset of $\mathcal{A}_\mathrm{G}(\mathbb{Q}_v)$ consisting of the characteristic polynomials with coefficients in $\mathbb{Z}_v$.

\paragraph{Definition}\setcounter{equation}{0}
Define the \emph{global Schwartz spaces} $\mathcal{S}_0(\mathcal{A}_\mathrm{G}(\mathbb{A}))$\label{sym152} and $\mathcal{S}_1(\mathcal{A}_\mathrm{G}(\mathbb{A}))$ to be the tensor products of local Schwartz spaces as defined in \eqref{local schwartz space 0} and \eqref{local schwartz space 1}
\begin{eqnarray}
\mathcal{S}_0(\mathcal{A}_\mathrm{G}(\mathbb{A}))
&=&
\bigotimes_v^\mathrm{res} \mathcal{S}_0(\mathcal{A}_\mathrm{G}(\mathbb{Q}_v))
\\\nonumber\\\nonumber
\mathcal{S}_1(\mathcal{A}_\mathrm{G}(\mathbb{A}))
&=&
\bigotimes_v^\mathrm{res} \mathcal{S}_1(\mathcal{A}_\mathrm{G}(\mathbb{Q}_v))
\end{eqnarray}
restricted with respect to the basic functions $\phi_{0,v}$ and $\phi_{1,v}$ as defined in \eqref{local basic function 0} and \eqref{local basic function 1}.

\subsection{The Harish-Chandra transform}

\paragraph{Definition}\setcounter{equation}{0}
Let $v$ be a place of $\mathbb{Q}$. For an element $\varphi_v$ in $\mathcal{S}_0(\mathcal{A}_\mathrm{G}(\mathbb{Q}_v))$, choose a Schwartz function $f_v$ on $\mathfrak{g}(\mathbb{Q}_v)$ such that
\begin{eqnarray}
\forall X_v\in\mathcal{A}_{\mathrm{G},\mathrm{reg}}(\mathbb{Q}_v)
&&
\varphi_v(X_v)=I_\mathrm{G}^\mathrm{G}(X_v,f_v).
\end{eqnarray}
Let $f_v\hat{~}$ be the Fourier transform of $f_v$, denote by $\mathcal{H}_v(\varphi_v)$ the element in $\mathcal{S}_1(\mathcal{A}_\mathrm{G}(\mathbb{Q}_v))$ such that
\begin{eqnarray}
\forall X_v\in\mathcal{A}_{\mathrm{G},\mathrm{reg}}(\mathbb{Q}_v)
&&
\mathcal{H}_v(\varphi_v)(X_v)=\mathcal{I}_\mathrm{max}^\mathrm{G}(X_v,f_v\hat{~}).
\end{eqnarray}
Define the \emph{local Harish-Chandra transform} to be the linear operator\label{sym153} 
\begin{eqnarray}
\mathcal{H}_v:\quad
\mathcal{S}_0(\mathcal{A}_\mathrm{G}(\mathbb{Q}_v)) 
&\longrightarrow&
\mathcal{S}_1(\mathcal{A}_\mathrm{G}(\mathbb{Q}_v))
\end{eqnarray}
defined by
\begin{eqnarray}
\mathcal{H}_v:\quad
\varphi_v
&\mapsto&
\mathcal{H}_v(\varphi_v).
\end{eqnarray}

\paragraph{Definition}\setcounter{equation}{0}
Let $v$ be a place of $\mathbb{Q}$. If $X_v$ is an element of $\mathcal{A}_{\mathrm{G},\mathrm{reg}}(\mathbb{Q})$, then the orbital integral $\mathcal{I}_\mathrm{max}^\mathrm{G}(X_v,~)$ is a vector-valued tempered distribution on $\mathfrak{g}(\mathbb{Q}_v)$. Denote by $\mathcal{I}_\mathrm{max}^\mathrm{G}(X_v,~)\hat{~}$\label{sym154} its componentwise Fourier transform, which is represented by a conjugation invariant function on $\mathfrak{g}(\mathbb{Q}_v)$, hence descends to a function on $\mathcal{A}_\mathrm{G}(\mathbb{Q}_v)$ denoted by
\begin{eqnarray}
Y_v
&\mapsto&
\Big(\mathcal{I}_\mathrm{max}^\mathrm{G}(X_v,~)\hat{~}\Big)(Y_v).
\end{eqnarray}
Define the \emph{local Harish-Chandra kernel function} $\mathcal{K}_v(~,~)$\label{sym155} to be the vector-valued bivariate function on $\mathcal{A}_{\mathrm{G},\mathrm{reg}}(\mathbb{Q}_v)$ such that
\begin{eqnarray}
\label{local harish chandra kernel function}
&&
\forall X_v,Y_v\in\mathcal{A}_{\mathrm{G},\mathrm{reg}}(\mathbb{Q}_v)
\\\nonumber
&&\quad
\mathcal{K}_v(X_v,Y_v)=
|D(Y_v)|_v^{1/2}\Big(\mathcal{I}_\mathrm{max}^\mathrm{G}(X_v,~)\hat{~}\Big)(Y_v).
\end{eqnarray}

\paragraph{Remark}\setcounter{equation}{0}
In the case when $v$ is a $p$-adic place, each component of the vector-valued kernel function $\mathcal{K}_v$ is denoted by $\hat{i}_\mathrm{M}^\mathrm{G}(X_v,Y_v)$\label{sym156} for some Levi subgroup $\mathrm{M}$ by Waldspurger in \cite{Wal}.

\paragraph{Lemma}\setcounter{equation}{0}
\emph{Let $v$ be a place of $\mathbb{Q}$. The local Harish-Chanrda transform $\mathcal{H}_v$ is an integral operator with integral kernel $\mathcal{K}_v$, i.e.
\begin{eqnarray}
&&
\forall \varphi_v\in\mathcal{S}_0(\mathcal{A}_\mathrm{G}(\mathbb{Q}_v))~
\forall X_v\in\mathcal{A}_{\mathrm{G},\mathrm{reg}}(\mathbb{Q}_v)
\\\nonumber\\\nonumber
&&\quad
\mathcal{H}_v(\varphi_v)(X_v)
=
\int_{\mathcal{A}_\mathrm{G}(\mathbb{Q}_v)}
\mathcal{K}_v(X_v,Y_v)\varphi_v(Y_v)~
|D(Y_v)|_v^{1/2}\mathrm{d}Y_v.
\end{eqnarray}
In particular the operator $\mathcal{H}_v$ is well-defined.}

\proof
Every $X_v$ in $\mathcal{A}_{\mathrm{G},\mathrm{reg}}(\mathbb{Q}_v)$ has a $\delta$-sequence in $\mathcal{S}_0(\mathcal{A}_\mathrm{G}(\mathbb{Q}_v))$. More precisely there exists a sequence of functions $\delta_{X_v,1},\delta_{X_v,2},\delta_{X_v,3},\dots$\label{sym157} in $\mathcal{S}_0(\mathcal{A}_\mathrm{G}(\mathbb{Q}_v))$ such that
\begin{eqnarray}
\lim_{i\rightarrow\infty}\delta_{X_v,i}
&=&
\delta_{X_v}
\end{eqnarray}
as distributions on $\mathcal{A}_{\mathrm{G},\mathrm{reg}}(\mathbb{Q}_v)$, where $\delta_{X_v}$ denotes the Dirac distribution at $X_v$. For each natural number $i$ choose a Schwartz function $\widetilde{\delta}_{X_v,i}$ on $\mathfrak{g}(\mathbb{Q}_v)$ such that
\begin{eqnarray}
\forall Y_v\in\mathcal{A}_{\mathrm{G},\mathrm{reg}}(\mathbb{Q}_v)
&&
\delta_{X_v,i}(Y_v)=I_\mathrm{G}^\mathrm{G}(Y_v,\widetilde{\delta}_{X_v,i}).
\end{eqnarray}
Let $f_v$ be a Schwartz function on $\mathfrak{g}(\mathbb{Q}_v)$ such that
\begin{eqnarray}
\forall Y_v\in\mathcal{A}_{\mathrm{G},\mathrm{reg}}(\mathbb{Q}_v)
&&
\varphi_v(Y_v)=I_\mathrm{G}^\mathrm{G}(Y_v,f_v).
\end{eqnarray}

By the local trace formula \eqref{harish chandra local trace formula}
\begin{eqnarray}
&&
\int_{\mathcal{A}_\mathrm{G}(\mathbb{Q}_v)}
\mathcal{I}_\mathrm{max}^\mathrm{G}(Y_v,f_v\hat{~})
I_\mathrm{G}^\mathrm{G}(Y_v,\widetilde{\delta}_{X_v,i})~
|D(Y_v)|_v^{-1/2}\mathrm{d}Y_v
\\\nonumber\\\nonumber
&&\quad
=\int_{\mathcal{A}_\mathrm{G}(\mathbb{Q}_v)}
\mathcal{I}_\mathrm{max}^\mathrm{G}(Y_v,\widetilde{\delta}_{X_v,i}\hat{~})
I_\mathrm{G}^\mathrm{G}(Y_v,f_v)~
|D(Y_v)|_v^{-1/2}\mathrm{d}Y_v,
\end{eqnarray}
hence
\begin{eqnarray}
&&
\mathcal{H}_v(\varphi_v)(X_v)
\\\nonumber\\\nonumber
&=&
\int_{\mathcal{A}_\mathrm{G}(\mathbb{Q}_v)}
\mathcal{H}_v(\varphi_v)(Y_v)\delta_{X_v}(Y_v)~
|D(Y_v)|_v^{-1/2}\mathrm{d}Y_v
\\\nonumber\\\nonumber
&=&
\lim_{i\rightarrow\infty}
\int_{\mathcal{A}_\mathrm{G}(\mathbb{Q}_v)}
\mathcal{H}_v(\varphi_v)(Y_v)\delta_{X_v,i}(Y_v)~
|D(Y_v)|_v^{-1/2}\mathrm{d}Y_v
\\\nonumber\\\nonumber
&=&
\lim_{i\rightarrow\infty}
\int_{\mathcal{A}_\mathrm{G}(\mathbb{Q}_v)}
\mathcal{I}_\mathrm{max}^\mathrm{G}(Y_v,\widetilde{\delta}_{X_v,i}\hat{~})
\varphi_v(Y_v)~
|D(Y_v)|_v^{-1/2}\mathrm{d}Y_v.
\end{eqnarray}
Hence $\mathcal{H}_v$ is an integral operator with integral kernel
\begin{eqnarray}
\lim_{i\rightarrow\infty} 
\mathcal{I}_\mathrm{max}^\mathrm{G}(Y_v,\widetilde{\delta}_{X_v,i}\hat{~})
&=&
|D(Y_v)|_v^{1/2}\left(\mathcal{I}_\mathrm{max}^\mathrm{G}(X_v,~)\hat{~}\right)(Y_v)
\end{eqnarray}
which is independent of the choice of the sequence $\big(\widetilde{\delta}_{X_v,i}\big)_{i=1}^\infty$ and equal to the local Harish-Chandra kernel $\mathcal{K}_v(X_v,Y_v)$.
\qed

\paragraph{Remark}\setcounter{equation}{0}
The local Harish-Chandra transform preserves the basic functions, i.e.
\begin{eqnarray}
\label{local harish chandra transform preserves basic function}
\mathcal{H}_v: \quad
\phi_{0,v}
&\mapsto&
\phi_{1,v}.
\end{eqnarray}

\paragraph{Definition}\setcounter{equation}{0}
For a finite set $S$ of places of $\mathbb{Q}$, define the \emph{$S$-local Harish-Chandra transform} $\mathcal{H}_S$\label{sym158} to be the tensor product
\begin{eqnarray}
\mathcal{H}_S=\bigotimes_{v\in S}\mathcal{H}_v:\quad
\bigotimes_{v\in S}\mathcal{S}_0(\mathcal{A}_\mathrm{G}(\mathbb{Q}_v))
&\longrightarrow&
\bigotimes_{v\in S}\mathcal{S}_1(\mathcal{A}_\mathrm{G}(\mathbb{Q}_v)).
\end{eqnarray}
By \eqref{local harish chandra transform preserves basic function} the limit of $\mathcal{H}_S$ as $S$ approaches infinity defines a linear operator
\begin{eqnarray}
\lim_S\mathcal{H}_S: \quad
\bigotimes_v^\mathrm{res}\mathcal{S}_0(\mathcal{A}_\mathrm{G}(\mathbb{Q}_v))
&\longrightarrow&
\bigotimes_v^\mathrm{res}\mathcal{S}_1(\mathcal{A}_\mathrm{G}(\mathbb{Q}_v)),
\end{eqnarray}
which is defined to be the \emph{global Harish-Chandra transform}
\begin{eqnarray}
\mathcal{H}:  \quad
\mathcal{S}_0(\mathcal{A}_\mathrm{G}(\mathbb{A})) 
&\longrightarrow&
\mathcal{S}_1(\mathcal{A}_\mathrm{G}(\mathbb{A})).
\end{eqnarray}

\subsection{The Poisson summation formula}

\paragraph{Definition}\setcounter{equation}{0}
Let $Z$\label{sym159} be a singular point in the complement of $\mathcal{A}_{\mathrm{G},\mathrm{reg}}(\mathbb{Q})$ in $\mathcal{A}_\mathrm{G}(\mathbb{Q})$. Denote by $\mathfrak{o}_Z$ the $\sim$ equivalence class in $\mathfrak{g}(\mathbb{Q})$ which is the fiber of $Z$. Let $S$ be a finite set of places of $\mathbb{Q}$, let $\varphi_S$ be an element of $\bigotimes_{v\in S}\mathcal{S}_0(\mathcal{A}_\mathrm{G}(\mathbb{Q}_v))$, let $f_S$ be a Schwartz function on $\mathfrak{g}(\mathbb{Q}_S)$ such that
\begin{eqnarray}
\forall X\in\mathcal{A}_{\mathrm{G},\mathrm{reg}}(\mathbb{Q}_S)
&&
\varphi_S(X)
=I_\mathrm{G}^\mathrm{G}(X,f_S)
\end{eqnarray}
By \eqref{general weighted orbital integral}, if $Z'$ is a $\equiv$ equivalence class in $(\mathfrak{m}(\mathbb{Q})\cap\mathfrak{o}_Z)_{\mathrm{M},S}$, then $I_\mathrm{G}^\mathrm{G}(Z',f_S)$ is determined by $\varphi_S(X)$ where $X$ ranges over the elements of $\mathcal{A}_{\mathrm{G},\mathrm{reg}}(\mathbb{Q}_S)$ close to $Z'$. Denote
\begin{eqnarray}
\varphi_S(Z')
&=&
I_\mathrm{G}^\mathrm{G}(Z',f_S).
\end{eqnarray}

\paragraph{Definition}\setcounter{equation}{0}
Let $X$ be an element of $\mathcal{A}_{\mathrm{G},\mathrm{reg}}(\mathbb{Q})$, let $\widetilde{X}$ be a regular semisimple element of $\mathfrak{g}(\mathbb{Q})$ lifting $X$, then by \eqref{tamagawa number} the constant $a^\mathrm{G}(S,\widetilde{X})$ as defined in \eqref{global coefficient} is independent of the finite set of places $S$ and determined by $X$. Denote\label{sym162}
\begin{eqnarray}
a(X)
&=&
a^\mathrm{G}(S,\widetilde{X}).
\end{eqnarray}

\paragraph{Definition}\setcounter{equation}{0}
Let $v$ be a place of $\mathbb{Q}$, let $\varphi_v$ be an element in $\mathcal{S}_0(\mathcal{A}_\mathrm{G}(\mathbb{Q}_v))$, then $\varphi_v$ is said to be \emph{cuspidal} if
\begin{eqnarray}
\varphi_v(X_v)&=&0
\end{eqnarray}
whenever $X_v$ is the image of a point $\widetilde{X}_v$ in $\mathfrak{g}(\mathbb{Q}_v)$ which is regular semisimple and not $\mathbb{Q}_v$-elliptic.

\paragraph{Proposition}\setcounter{equation}{0}
(Poisson summation formula)\\
\emph{Let $\varphi$ be an element of $\mathcal{S}_0(\mathcal{A}_\mathrm{G}(\mathbb{A}))$ which is cuspidal at two distinct places of $\mathbb{Q}$, then}
\begin{eqnarray}
\label{poisson summation formula}
&&
\sum_{X\in\mathcal{A}_{\mathrm{G},\mathrm{reg}}(\mathbb{Q})}
a(X) \varphi(X)+
\\\nonumber\\\nonumber
&&\quad
+\sum_{Z\in\mathcal{A}_\mathrm{G}(\mathbb{Q})
-\mathcal{A}_{\mathrm{G},\mathrm{reg}}(\mathbb{Q})}
\lim_S\Bigg(
\sum_{Z'\in(\mathfrak{g}(\mathbb{Q})\cap\mathfrak{o}_Z)_{\mathrm{G},S}}
a^\mathrm{G}(S,Z')
\varphi_S(Z')\Bigg)
\\\nonumber\\\nonumber
&=&
\sum_{X\in\mathcal{A}_{\mathrm{G},\mathrm{reg}}(\mathbb{Q})}
a(X) \mathcal{H}(\varphi)(X)+
\\\nonumber\\\nonumber
&&\quad
+\sum_{Z\in\mathcal{A}_\mathrm{G}(\mathbb{Q})
-\mathcal{A}_{\mathrm{G},\mathrm{reg}}(\mathbb{Q})}
\lim_S\Bigg(
\sum_{Z'\in(\mathfrak{g}(\mathbb{Q})\cap\mathfrak{o}_Z)_{\mathrm{G},S}}
a^\mathrm{G}(S,Z')
\mathcal{H}_S(\varphi_S)(Z')\Bigg).
\end{eqnarray}

\proof
This follows from the global trace formula \eqref{harish chandra global trace formula} and Remark \ref{cuspidality implies cuspidality}.
\qed

\paragraph{Remark}\setcounter{equation}{0}
The Poisson summation formula \eqref{poisson summation formula} has the general form
\begin{eqnarray}
\label{poisson summation formula general form}
&&
\sum_{X\in\mathcal{A}_{\mathrm{G},\mathrm{reg}}(\mathbb{Q})}a(X)\varphi(X)
+\sum_{Z\in\mathcal{A}_\mathrm{G}(\mathbb{Q})-\mathcal{A}_{\mathrm{G},\mathrm{reg}}(\mathbb{Q})}(\dots.)
\\\nonumber\\\nonumber
&&\quad
=\sum_{X\in\mathcal{A}_{\mathrm{G},\mathrm{reg}}(\mathbb{Q})}a(X)
\mathcal{H}(\varphi)(X)
+\sum_{Z\in\mathcal{A}_\mathrm{G}(\mathbb{Q})-\mathcal{A}_{\mathrm{G},\mathrm{reg}}(\mathbb{Q})}(\dots.).
\end{eqnarray}

\paragraph{Corollary}\setcounter{equation}{0}
\emph{Let $v$ and $w$ be two places of $\mathbb{Q}$. Let $\varphi$ be an element of $\mathcal{S}_0(\mathcal{A}_\mathrm{G}(\mathbb{A}))$ which is cuspidal at two distinct places of $\mathbb{Q}$ such that $\varphi_v$ vanishes on a neighborhood of the complement of $\mathcal{A}_{\mathrm{G},\mathrm{reg}}(\mathbb{Q}_v)$ in $\mathcal{A}_\mathrm{G}(\mathbb{Q}_v)$ and $\mathcal{H}_w(\varphi_w)$ vanishes on a neighborhood of the complement of $\mathcal{A}_{\mathrm{G},\mathrm{reg}}(\mathbb{Q}_w)$ in $\mathcal{A}_\mathrm{G}(\mathbb{Q}_w)$. Then}
\begin{eqnarray}
\sum_{X\in\mathcal{A}_{\mathrm{G},\mathrm{reg}}(\mathbb{Q})}a(X)\varphi(X)
&=&
\sum_{X\in\mathcal{A}_{\mathrm{G},\mathrm{reg}}(\mathbb{Q})}a(X)
\mathcal{H}(\varphi)(X).
\end{eqnarray}

\proof
This follows from \eqref{poisson summation formula general form}.
\qed

\paragraph{Corollary}\setcounter{equation}{0}
\label{poisson summation formula simplified statement}
(Poisson summation formula for $\mathcal{H}_0$)\\
\emph{There exists an endomorphism $\mathcal{H}_0$ of $\mathcal{S}_0(\mathcal{A}_\mathrm{G}(\mathbb{A}))$, defined similarly as $\mathcal{H}$ as a tensor product of local endomorphisms $\mathcal{H}_{0,v}$, such that for each Schwartz function $\varphi$ on $\mathcal{A}_\mathrm{G}(\mathbb{A})$ which is cuspidal at two distinct places of $\mathbb{Q}$}
\begin{eqnarray}
\label{poisson summation formula for h0}
&&
\sum_{X\in\mathcal{A}_{\mathrm{G},\mathrm{reg}}(\mathbb{Q})}
a(X) \varphi(X)+
\\\nonumber\\\nonumber
&&\quad
+\sum_{Z\in\mathcal{A}_\mathrm{G}(\mathbb{Q})
-\mathcal{A}_{\mathrm{G},\mathrm{reg}}(\mathbb{Q})}
\lim_S\Bigg(
\sum_{Z'\in(\mathfrak{g}(\mathbb{Q})\cap\mathfrak{o}_Z)_{\mathrm{G},S}}
a^\mathrm{G}(S,Z')
\varphi_S(Z')\Bigg)
\\\nonumber\\\nonumber
&=&
\sum_{X\in\mathcal{A}_{\mathrm{G},\mathrm{reg}}(\mathbb{Q})}
a(X) \mathcal{H}_0(\varphi)(X)+
\\\nonumber\\\nonumber
&&\quad
+\sum_{Z\in\mathcal{A}_\mathrm{G}(\mathbb{Q})
-\mathcal{A}_{\mathrm{G},\mathrm{reg}}(\mathbb{Q})}
\lim_S\Bigg(
\sum_{Z'\in(\mathfrak{g}(\mathbb{Q})\cap\mathfrak{o}_Z)_{\mathrm{G},S}}
a^\mathrm{G}(S,Z')
\mathcal{H}_{0,S}(\varphi_S)(Z')\Bigg).
\end{eqnarray}

\proof
The local endomorphisms $\mathcal{H}_{0,v}$ are defined by sending the orbital integral $I_\mathrm{G}^\mathrm{G}(X,f_v)$ to the orbital integral of the Fourier transform $I_\mathrm{G}^\mathrm{G}(X,f_v\hat{~})$ for a Schwartz function $f_v$ on $\mathfrak{g}(\mathbb{Q}_v)$. If $v$ is a $p$-adic place, the endomorphism $\mathcal{H}_{0,v}$ is well-defined due to the density of regular semisimple orbital integrals. See \S27 of \cite{Kott}. If $v$ is archimedean, then $\mathcal{H}_{0,v}$ is well-defined due to the second last identity on page 104 of \cite{Var}. The local basic function $\phi_{0,v}$ is preserved by $\mathcal{H}_{0,v}$, hence the global endomorphism $\mathcal{H}_0$ is also well-defined.

The Poisson summation formula \eqref{poisson summation formula for h0}follows from \eqref{poisson summation formula} since the cuspidality condition on $\varphi$ implies that only those $X$ for which the orbital integrals $I_\mathrm{G}^\mathrm{G}(X,f_v)$ and $\mathcal{I}_\mathrm{max}^\mathrm{G}(X,f_v)$ are equal, where $f_v$ denotes a local Schwartz function on $\mathfrak{g}(\mathbb{Q}_v)$ whose orbital integral is equal to $\varphi_v$, contribute to the identity.
\qed

\paragraph{Remark}\setcounter{equation}{0}
The endomorphism $\mathcal{H}_0$ is the Harish-Chandra transform originally considered by Jacquet. If $v$ is a $p$-adic place, then it is conjectured by Jacquet for a general reductive Lie algebra and proven by Waldspurger for $\mathrm{gl}(n)$ that a Schwartz function $\varphi_v$ on $\mathcal{A}_\mathrm{G}(\mathbb{Q}_v)$ is proportional to $\phi_{0,v}$ if and only if $\varphi_v$ and $\mathcal{H}_{0,v}(\varphi_v)$ are both supported on the subset of $\mathcal{A}_\mathrm{G}(\mathbb{Q}_v)$ consisting of the characteristic polynomials with coefficients in $\mathbb{Z}_v$.

\paragraph{Corollary}\setcounter{equation}{0}
\label{local harish chandra transform bijectivity}
\emph{Let $v$ be a place of $\mathbb{Q}$. The local Harish-Chandra transform $\mathcal{H}_v$ is a bijection from $\mathcal{S}_0(\mathcal{A}_\mathrm{G}(\mathbb{Q}_v))$ onto $\mathcal{S}_1(\mathcal{A}_\mathrm{G}(\mathbb{Q}_v))$.}

\proof
The following argument is suggested by Sakellaridis.

The local Harish-Chandra transform $\mathcal{H}_v$ is surjective by definition. For injectivity assume for contradiction that there exists an element $\varphi_v$ in $\mathcal{S}_0(\mathcal{A}_\mathrm{G}(\mathbb{Q}_v))$ and a point $X_v$ in $\mathcal{A}_{\mathrm{G},\mathrm{reg}}(\mathbb{Q}_v)$ such that
\begin{eqnarray}
\forall Y_v\in\mathcal{A}_{\mathrm{G},\mathrm{reg}}(\mathbb{Q}_v)
&&
\mathcal{H}_v(\varphi_v)(Y_v)=0
\end{eqnarray}
and
\begin{eqnarray}
\varphi_v(X_v)
&\neq&
0.
\end{eqnarray}
Since $\varphi_v$ is smooth on $\mathcal{A}_{\mathrm{G},\mathrm{reg}}(\mathbb{Q}_v)$, without loss of generality $X_v$ lies in the dense subset $\mathcal{A}_{\mathrm{G},\mathrm{reg}}(\mathbb{Q})$. By parabolic descent along a suitable parabolic subgroup, the argument is reduced to the special case that $X_v$ is the image of a $\mathbb{Q}_v$-elliptic element of $\mathfrak{g}(\mathbb{Q}_v)$.

Let $\varphi$ be en element of $\mathcal{S}_0(\mathcal{A}_\mathrm{G}(\mathbb{A}))$ which is cuspidal at two other places and whose local component at $v$ is equal to $\varphi_v$, then by the Poisson summation formula \eqref{poisson summation formula}
\begin{eqnarray}
\label{local harish chandra transform bijectivity proof first equality}
\sum_{X\in\mathcal{A}_{\mathrm{G},\mathrm{reg}}(\mathbb{Q})}
a(X) \varphi(X)+
\sum_{Z\in\mathcal{A}_\mathrm{G}(\mathbb{Q})
-\mathcal{A}_{\mathrm{G},\mathrm{reg}}(\mathbb{Q})}
\bigg(\begin{array}{c}\textrm{contribution from}\\\textrm{the singular locus}
\end{array}\bigg)
&=&0
\end{eqnarray}
since $\mathcal{H}(\varphi)$ vanishes identically at the place $v$. There are two cases:
\begin{itemize}
\item If $v$ is $p$-adic, then there exists an integer $N$ such that $X_v$ lies in the image of $\mathfrak{g}_\mathrm{reg.ss}(N^{-1}\mathbb{Z})$ in $\mathcal{A}_{\mathrm{G},\mathrm{reg}}(\mathbb{Q})$ under the natural projection. 

At a finite place $w$ distinct from $v$, denote by $\Lambda_{N,w}$\label{sym163} the lattice
\begin{eqnarray}
\Lambda_{N,w} 
\quad=\quad 
\mathfrak{g}(N^{-1}\mathbb{Z}_w)
&\subset& 
\mathfrak{g}(\mathbb{Q}_w).
\end{eqnarray}
Let $\varphi_w$ be the element of $\mathcal{S}_0(\mathcal{A}_\mathrm{G}(\mathbb{Q}_w))$ such that
\begin{eqnarray}
\forall Y_w\in\mathcal{A}_{\mathrm{G},\mathrm{reg}}(\mathbb{Q}_w)
&&
\varphi_w(Y_w)=I_\mathrm{G}^\mathrm{G}(Y_w,\mathbb{I}_{\Lambda_{N,w}})
\end{eqnarray}
where $\mathbb{I}_{\Lambda_{N,w}}$ denotes the characteristic function of $\Lambda_{N,w}$.

At infinity, the set of rational points where $\bigotimes_{w<\infty}\varphi_w$ is nonzero is contained in the discrete subset
\begin{eqnarray}
\pi_\infty(\mathfrak{g}_\mathrm{reg.ss}(N^{-1}\mathbb{Z}))
&\subset& 
\mathcal{A}_{\mathrm{G},\mathrm{reg}}(\mathbb{R})\end{eqnarray}
where $\pi_\infty$\label{sym164} denotes the natural projection from $\mathfrak{g}(\mathbb{R})$ to $\mathcal{A}_\mathrm{G}(\mathbb{R})$. Choose $\varphi_\infty$ to be a bump function supported away from the complement of $\mathcal{A}_{\mathrm{G},\mathrm{reg}}(\mathbb{R})$ in $\mathcal{A}_\mathrm{G}(\mathbb{R})$ such that $X_v$ is the only point contained in
\begin{eqnarray}
\mathrm{supp}(\varphi_\infty)
&\cap&
\pi_\infty(\mathfrak{g}_\mathrm{reg.ss}(N^{-1}\mathbb{Z})).
\end{eqnarray}

Choose $\varphi$ to be $\bigotimes_w \varphi_w$, then the left hand side of \eqref{local harish chandra transform bijectivity proof first equality} is equal to
\begin{eqnarray}
\qquad\qquad
a(X_v)\varphi(X_v)
&=&
\mathrm{Vol}\left(\mathrm{T}_{X_v}(\mathbb{Q})\backslash\mathrm{T}_{X_v}(\mathbb{A})^1\right)
\lim_S \bigg( \prod_{w\in S}\varphi_w(X_v) \bigg)
\end{eqnarray}
for some torus $\mathrm{T}_{X_v}$\label{sym165}, which is nonzero.

\item If $v$ is archimedean, then choose a regular semisimple element $\widetilde{X}_\infty$ of $\mathfrak{g}(\mathbb{Q})$ lifting $X_\infty$ and choose a Schwartz function $f_\infty$ on $\mathfrak{g}(\mathbb{A})$ such that
\begin{eqnarray}
\forall Y_\infty\in\mathcal{A}_{\mathrm{G},\mathrm{reg}}(\mathbb{R})
&&
\varphi_\infty(Y_\infty)=I_\mathrm{G}^\mathrm{G}(Y_\infty,f_\infty).
\end{eqnarray}
If $N$ is an integer, denote by $\Lambda_{N,\infty}$\label{sym166} the lattice
\begin{eqnarray}
\Lambda_{N,\infty}
\quad=\quad
\mathfrak{g}(N\mathbb{Z})
&\subset&
\mathfrak{g}(\mathbb{R}).
\end{eqnarray}
Since $f_\infty$ is a Schwartz function on $\mathfrak{g}(\mathbb{R})$, the quantity $\varphi_\infty(Y_\infty)$ is rapidly decreasing as $Y_\infty$ approaches infinity in $\mathcal{A}_{\mathrm{G},\mathrm{reg}}(\mathbb{R})$ in such a way that $|D(Y_\infty)|$ is uniformly bounded below, hence for every positive real numbers $\epsilon$ and $r$ there exists a natural number $N_\epsilon$ such that
\begin{eqnarray}
\label{local harish chandra transform bijectivity proof second inequality}
\forall N\geq N_\epsilon
&&
\sum_{\substack{Y\in \pi_\infty(\widetilde{X}_\infty+\Lambda_{N,\infty})\\
Y\neq X_\infty,~|D(Y)|\geq1}}
|\varphi_\infty(Y)|<\epsilon N^{-r}
\end{eqnarray}
where $\pi_\infty$ denotes the natural projection from $\mathfrak{g}(\mathbb{R})$ to $\mathcal{A}_\mathrm{G}(\mathbb{R})$. Choose 
\begin{eqnarray}
\label{local harish chandra transform bijectivity proof third equality}
\epsilon
&=&
\frac{1}{2(n!)^{n\cdot n!}}
|a(X_\infty)\varphi_\infty(X_\infty)|
\end{eqnarray}
where $n$ is the rank of $\mathrm{G}$ and choose $r$ to be $2\mathrm{dim}(\mathfrak{g})$. Let $N$ be a natural number greater than $N_\epsilon$ such that
\begin{eqnarray}
\forall Y\in\pi_\infty(\widetilde{X}_\infty+\Lambda_{N,\infty})
&&
D(Y)\neq0.
\end{eqnarray}

At a finite place $w$, denote by $\Lambda_{N,w}$ the lattice
\begin{eqnarray}
\Lambda_{N,w} 
\quad=\quad 
\mathfrak{g}(N\mathbb{Z}_w)
&\subset& 
\mathfrak{g}(\mathbb{Q}_w).
\end{eqnarray}
Let $\varphi_w$ be the element of $\mathcal{S}_0(\mathcal{A}_\mathrm{G}(\mathbb{Q}_w))$ such that
\begin{eqnarray}
\forall Y_w\in\mathcal{A}_{\mathrm{G},\mathrm{reg}}(\mathbb{Q}_w)
&&
\varphi_w(Y_w)=
I_\mathrm{G}^\mathrm{G}(Y_w,\mathbb{I}_{\widetilde{X}_\infty+\Lambda_{N,w}})
\end{eqnarray}
where $\mathbb{I}_{\widetilde{X}_\infty+\Lambda_{N,w}}$ denotes the characteristic function of the translation of $\Lambda_{N,w}$ by $\widetilde{X}_\infty$.

Choose $\varphi$ to be $\bigotimes_w \varphi_w$, then the left hand side of \eqref{local harish chandra transform bijectivity proof first equality} is equal to
\begin{eqnarray}
\label{local harish chandra transform bijectivity proof fourth equality}
\qquad
&&
a(X_\infty)\prod_{w<\infty}\varphi_w(X_\infty)\cdot\varphi_\infty(X_\infty)+
\\\nonumber\\\nonumber
&&\quad
+\sum_{\substack{Y\in\pi_\infty(\widetilde{X}_\infty+\Lambda_{N,\infty})
\\Y\neq X_\infty,~|D(Y)|\geq1}}
a(Y)\prod_{w<\infty}\varphi_w(Y)\cdot\varphi_\infty(Y)
\\\nonumber\\\nonumber
&=&
C\cdot\Bigg(a(X_\infty)\varphi_\infty(X_\infty)+
\\\nonumber\\\nonumber
&&\quad
+\sum_{\substack{Y\in\pi_\infty(\widetilde{X}_\infty+\Lambda_{N,\infty})
\\Y\neq X_\infty,~|D(Y)|\geq1}}
\bigg(a(Y)\frac{\prod_{w<\infty}\varphi_w(Y)}
{\prod_{w<\infty}\varphi_w(X_\infty)}\bigg)
\varphi_\infty(Y)\Bigg)
\end{eqnarray}
where $C$ is a nonzero constant. Let $Y$ be an element of $\mathcal{A}_{\mathrm{G},\mathrm{reg}}(\mathbb{Q})$. By definition
\begin{eqnarray}
\label{local harish chandra transform bijectivity proof fifth equality}
a(Y)
&=&
\pm |\pi_0(\mathrm{M}_Y)|^{-1}
\mathrm{Vol}(\mathrm{T}_Y(\mathbb{Q})\backslash\mathrm{T}_Y(\mathbb{A})^1)
\end{eqnarray}
where $\mathrm{M}$ is a Levi subgroup in $\mathcal{L}$ and $\mathrm{T}_Y$\label{sym167} is a maximal torus in $\mathrm{G}$. By the Main theorem in \S5 of \cite{Ono},
\begin{eqnarray}
\label{local harish chandra transform bijectivity proof sixth inequality}
\mathrm{Vol}(\mathrm{T}_Y(\mathbb{Q})\backslash\mathrm{T}_Y(\mathbb{A})^1)
&\leq&
|H_\mathrm{Gal}^1(\mathbb{Q},\widehat{\mathrm{T}}_Y)|
\end{eqnarray}
where $H_\mathrm{Gal}^1$ denotes the first Galois cohomology group and $\widehat{\mathrm{T}}_Y$ denotes the Galois module of algebraic characters of $\mathrm{T}_Y$. Let $F$ be the splitting field of $\mathrm{T}_Y$, denote by $\Gamma$ the Galois group of $F$ over $\mathbb{Q}$. Then by the inflation-restriction exact sequence
\begin{eqnarray}
&&
\xymatrix{0\ar[r]& H^1(\Gamma,\widehat{\mathrm{T}}_Y)\ar[r]&
H_\mathrm{Gal}^1(\mathbb{Q},\widehat{\mathrm{T}}_Y)\ar[r]&
H_\mathrm{Gal}^1(F,\widehat{\mathrm{T}}_Y)^\Gamma\ar@{=}[d]\\
&&&\mathrm{Hom}(\mathrm{Gal}(\overline{F}/F),\mathbb{Z}^n)^\Gamma}
\end{eqnarray}
and the compactness of $\mathrm{Gal}(\overline{F}/F)$ which implies that the only continuous homomorphism from $\mathrm{Gal}(\overline{F}/F)$ to $\mathbb{Z}^n$ is the trivial homomorphism, 
\begin{eqnarray}
\label{local harish chandra transform bijectivity proof seventh equality}
|H^1(\Gamma,\widehat{\mathrm{T}}_Y)|
&=&
|H_\mathrm{Gal}^1(\mathbb{Q},\widehat{\mathrm{T}}_Y)|.
\end{eqnarray}
The group $H^1(\Gamma,\widehat{\mathrm{T}}_Y)$ is annihilated by the order of $\Gamma$ which is bounded by the factorial of $n$ since the $\widehat{\mathrm{T}}_Y$ splits over the splitting field of the characteristic polynomial represented by $Y$ which is of degree $n$. By the bar resolution, the group $H^1(\Gamma,\widehat{\mathrm{T}}_Y)$ is a subquotient of $\bigoplus_{g\in\Gamma}\mathbb{Z}^n$ which is generated by at most $n\cdot n!$ elements, hence by \eqref{local harish chandra transform bijectivity proof fifth equality}, \eqref{local harish chandra transform bijectivity proof sixth inequality} and \eqref{local harish chandra transform bijectivity proof seventh equality}
\begin{eqnarray}
\label{local harish chandra transform bijectivity proof eighth inequality}
\forall Y\in\mathcal{A}_{\mathrm{G},\mathrm{reg}}(\mathbb{Q})
&&
|a(Y)|\leq (n!)^{n\cdot n!}.
\end{eqnarray}
As $N$ approaches infinity, either $\prod_{w<\infty}\varphi_w(Y)$ vanishes or
\begin{eqnarray}
N^{-\mathrm{dim}(\mathfrak{g})}\leq\prod_{w<\infty}\varphi_w(Y)\leq1,
\end{eqnarray}
and $\prod_{w<\infty}\varphi_w(X_\infty)$ is nonzero. Hence without loss of generality $N$ is large enough so that
\begin{eqnarray}
\label{local harish chandra transform bijectivity proof ninth inequality}
\forall Y\in\mathcal{A}_{\mathrm{G},\mathrm{reg}}(\mathbb{Q})
&&
\frac{\prod_{w<\infty}\varphi_w(Y)}{\prod_{w<\infty}\varphi_w(X_\infty)} \leq N^r.
\end{eqnarray}
By \eqref{local harish chandra transform bijectivity proof eighth inequality} and \eqref{local harish chandra transform bijectivity proof ninth inequality},
\begin{eqnarray}
&&
\sum_{\substack{Y\in\pi_\infty(\widetilde{X}_\infty+\Lambda_{N,\infty})
\\Y\neq X_\infty,~|D(Y)|\geq1}}
\bigg|a(Y)\frac{\prod_{w<\infty}\varphi_w(Y)}
{\prod_{w<\infty}\varphi_w(X_\infty)}
\varphi_\infty(Y)\bigg|
\\\nonumber\\\nonumber
&\leq&
N^r(n!)^{n\cdot n!}
\sum_{\substack{Y\in\pi_\infty(\widetilde{X}_\infty+\Lambda_{N,\infty})
\\Y\neq X_\infty,~|D(Y)|\geq1}}
|\varphi_\infty(Y)|
\\\nonumber\\
\label{local harish chandra transform bijectivity proof tenth inequality}
&<&
N^r(n!)^{n\cdot n!}\epsilon N^{-r}
\\\nonumber\\
\label{local harish chandra transform bijectivity proof eleventh equality}
&=&
(n!)^{n\cdot n!}\frac{1}{2(n!)^{n\cdot n!}}|a(X_\infty)\varphi_\infty(X_\infty)|
\\\nonumber\\\nonumber
&=&
\frac{1}{2}|a(X_\infty)\varphi_\infty(X_\infty)|
\end{eqnarray}
where the inequality \eqref{local harish chandra transform bijectivity proof tenth inequality} follows from \eqref{local harish chandra transform bijectivity proof second inequality} and the equality \eqref{local harish chandra transform bijectivity proof eleventh equality} follows from \eqref{local harish chandra transform bijectivity proof third equality}. Hence the right hand side of \eqref{local harish chandra transform bijectivity proof fourth equality} is nonzero.
\end{itemize}
In either case the left hand side of \eqref{local harish chandra transform bijectivity proof first equality} is nonzero, which is a contradiction.
\qed

\paragraph{Remark}\setcounter{equation}{0}
Injectivity of $\mathcal{H}_v$ is analogous to the classical result of Harish-Chandra on density of regular semisimple orbital integrals which states that for a $p$-adic reductive Lie algebra $\mathfrak{g}$, if $f$ is a Schwartz function on $\mathfrak{g}$ such that all the regular semisimple orbital integrals of $f$ vanish, then $\mathcal{D}(f)$ vanishes for every invariant distribution $\mathcal{D}$\label{sym168} on $\mathfrak{g}$. See \S27 of \cite{Kott}.

\appendix
\renewcommand{\theparagraph}{(\Alph{section}.\arabic{paragraph})}
\renewcommand{\thesubparagraph}{\Alph{section}.\arabic{paragraph}}

\section{Appendix: Scissors congruence and orbital integrals}\setcounter{paragraph}{0}

In this appendix weighted orbital integrals on $\mathfrak{g}$ are interpreted as scissors congruence classes of polyhedra. Similar results also hold for invariant weighted orbital integrals.

\paragraph{Definition}\setcounter{equation}{0}
Let $\mathbb{E}^n$\label{sym117} denote the $n$-dimensional Euclidean space. A polytope in $\mathbb{E}^n$ that is closed with nonempty interior is said to be \emph{proper}. Let $P$ and $Q$ be proper convex polytopes in $\mathbb{E}^n$. Then $P$ and $Q$ are said to be \emph{translational scissors congruent} if there exist convex polytopes $P_1, P_2,\dots, P_l$ and $Q_1, Q_2,\dots, Q_l$ in $\mathbb{E}^n$ such that
\begin{eqnarray}
P=\bigcup_{i=1}^l P_i
&&
Q=\bigcup_{i=1}^l Q_i
\end{eqnarray}
and $P_i$ is a translation of $Q_i$ for each index $i$ among $1,2,\dots,l$. Define the \emph{scissors group} of $\mathbb{E}^n$, denoted by $\mathbb{S}(\mathbb{E}^n)$, to be the quotient of the free abelian group generated by the proper convex polytopes in $\mathbb{E}^n$ modulo translational scissors congruence.

\paragraph{Definition}\setcounter{equation}{0}
A flag $\Phi$\label{sym118} of linear subspaces in $\mathbb{E}^n$ is said to be \emph{strict of length $l$} if
\begin{eqnarray}
\Phi &=& V_0\supset V_1\supset\dots\supset V_l
\end{eqnarray}
where $V_i$ has codimension $i$. Let $\Phi$ be a strict flag of length $l$. A \emph{rigging} $\mathbf{r}$ of $\Phi$ is a collection
\begin{eqnarray}
\mathbf{r}
&=&
\bigg\{r_1,r_2,\dots,r_l\bigg\}
\end{eqnarray} 
where $r_i$ is a real linear functional on $V_{i-1}$ with kernel $V_i$ for each $i$ among $1,2,\dots,l$. Two riggings $\mathbf{r}$ and $\mathbf{r}'$ are \emph{equivalent} if $r_i$ and $r_i'$ are positive multiples of each other for each $i$ among $1,2,\dots,l$. Denote by $\mathcal{R}ig(\Phi)$ the collection of equivalence classes of riggings of $\Phi$.  A \emph{rigged flag} $\Phi^\mathbf{r}$ is defined to be a strict flag $\Phi$ together with a choice of an element $\mathbf{r}$ in $\mathcal{R}ig(\Phi)$. 

An \emph{orientation} of $\mathbb{E}^n$ is an ordered basis of $\mathbb{E}^n$ defined upto a linear transformation with positive determinant. The product of an orientation of $\mathbb{E}^n$ with a translation invariant measure on $\mathbb{E}^n$ is equal to a volume form on $\mathbb{E}^n$. Fix an ordered basis $\mathcal{B}$\label{sym119}
\begin{eqnarray}
\mathcal{B}&=&\bigg(b_1,b_2,\dots,b_n\bigg)
\end{eqnarray}
of $\mathbb{E}^n$ such that
\begin{eqnarray}
\mathcal{B}'&=&\bigg(b_{l+1},b_{l+2},\dots,b_n\bigg)
\end{eqnarray}
is an ordered basis of $V_l$. Let $\mathbf{r}$ be a rigging of $\Phi$, choose vectors 
\begin{eqnarray}
c_1\in V_0,\quad c_2\in V_1,\quad \dots,\quad c_l\in V_{l-1}
\end{eqnarray}
such that
\begin{eqnarray}
\forall i=1,2,\dots,l&& r_i(c_i)>0.
\end{eqnarray}
Denote by $\mathcal{B}^\mathbf{r}$ the ordered basis
\begin{eqnarray}
\mathcal{B}^\mathbf{r}&=&\bigg(c_1,c_2,\dots,c_l,b_{l+1},b_{l+2},\dots,b_n\bigg)
\end{eqnarray}
of $\mathbb{E}^n$. Define the \emph{sign} of $\mathbf{r}$ by
\begin{eqnarray}
\mathrm{sign}(\mathbf{r})
&=&
\left\{\begin{array}{cl}
1&\textrm{if $\mathcal{B}$ and $\mathcal{B}^\mathbf{r}$ define}\\
&\textrm{the same orientation of $\mathbb{E}^n$,}\\\\
-1 &\textrm{if $\mathcal{B}$ and $\mathcal{B}^\mathbf{r}$ define}\\
&\textrm{opposite orientations of $\mathbb{E}^n$.}
\end{array}\right.
\end{eqnarray}

Let $\Phi^\mathbf{r}$\label{sym120} be a rigged flag of length $l$ in $\mathbb{E}^n$. Let $P$ be a proper convex polytope in $\mathbb{E}^n$. Then the \emph{$\Phi^\mathbf{r}$-boundary} $\partial_{\Phi^\mathbf{r}}$ of $P$ is defined by
\begin{eqnarray}
\partial_{\Phi^\mathbf{r}}P
&=&
r_l^\mathrm{min}\big(
r_{l-1}^\mathrm{min}\big(
\dots\big(
r_2^\mathrm{min}\big(
r_1^\mathrm{min}(P)\big)\dots\big)
\end{eqnarray}
where for each subset $S$ of $V_{i-1}$
\begin{eqnarray}
r_i^\mathrm{min}(S)
&=&
\left\{\begin{array}{cl}
v_i+r_i^{-1}\Big(\min\limits_{s\in S} r_i(s)\Big)
&\textrm{if $v_i$ is a vector in $V_{i-1}$}\\
&\textrm{such that}\\
&\quad v_i+r_i^{-1}\Big(\min\limits_{s\in S} r_i(s)\Big)\\
&\textrm{is a subset of $V_i$ with}\\
&\textrm{nonempty interior,}\\\\
\emptyset
&\textrm{if no such $v_i$ exists,}
\end{array}\right.
\end{eqnarray}
which is a subset of $V_i$ defined upto translation. Fix a translation invariant measure on the Euclidean space $V_l$. This determines the volume of the convex polytope $\partial_{\Phi^\mathbf{r}}P$. Then define the \emph{Hadwiger invariant} $\mathrm{Had}_\Phi$\label{sym121} of $P$ with respect to $\Phi$ by
\begin{eqnarray}
\mathrm{Had}_\Phi(P)
&=&
\sum_{\mathbf{r}\in\mathcal{R}ig(\Phi)}
\mathrm{sign(\mathbf{r})}
\mathrm{Vol}(\partial_{\Phi^\mathbf{r}}P).
\end{eqnarray}

\paragraph{Remark}\setcounter{equation}{0}
Each Hadwiger invariant defines a real-valued additive function on $\mathbb{S}(\mathbb{E}^n)$.

\paragraph{Lemma}\setcounter{equation}{0}
\emph{Let $P$ and $Q$ be proper convex polytopes in $\mathbb{E}^n$. Then $P$ are $Q$ are translational scissors congruent if and only if for every $l$ among $0,1,2,\dots,n$, for each strict flag $\Phi$ of length $l$,}
\begin{eqnarray}
\mathrm{Had}_\Phi(P)&=&\mathrm{Had}_\Phi(Q).
\end{eqnarray}

\proof
See Corollary 2 in \S4 of \cite{MorC}.
\qed

\paragraph{Lemma}\setcounter{equation}{0}
\emph{Let $(H_\Phi)$\label{sym122} be a collection of real numbers indexed by the set of all strict flags $\Phi$ in $\mathbb{E}^n$ which vanishes for all but finitely many $\Phi$. Then there exists an element $[P]$ in the scissors group $\mathbb{S}(\mathbb{E}^n)$ such that
\begin{eqnarray}
\forall \textrm{ strict flag } \Phi \textrm{ in }\mathbb{E}^n
&&
H_\Phi=\mathrm{Had}_\Phi([P])
\end{eqnarray}
if and only if for every $l$ among $0,1,2,\dots,n$, for each strict flag
\begin{eqnarray}
\Phi&=&
V_0\supset V_1\supset\dots\supset V_l
\end{eqnarray}
of length $l$, for every $i$ among $1,2,\dots,l-1$,
\begin{eqnarray}
\label{hadwiger characterization i}
\sum_{\Phi'\in\mathcal{F}_\mathrm{s}(\Phi,i)}H_{\Phi'}&=&0
\end{eqnarray}
where $\mathcal{F}_\mathrm{s}(\Phi,i)$ is the set defined by
\begin{eqnarray}
&&
\mathcal{F}_\mathrm{s}(\Phi,i)
\\\nonumber\\\nonumber
&=&
\bigg\{
\Phi' \textrm{ strict flag in } \mathbb{E}^n:~
\Phi'=V_0\supset\dots\supset V_{i-1}\supset U_i\supset 
V_{i+1}\supset\dots\supset V_l
\bigg\},
\end{eqnarray}
and 
\begin{eqnarray}
\label{hadwiger characterization ii}
\sum_{\Phi'\in\mathcal{F}_\mathrm{s}(\Phi,l)}H_{\Phi'} \wedge \mathbf{d}U_l
&=&0
\end{eqnarray}
where
\begin{eqnarray}
\Phi'&=&V_0\supset V_1\supset\dots\supset V_{l-1}\supset U_l
\end{eqnarray}
and $\mathbf{d}U_l$ is the element of $\bigwedge^{n-l}\mathbb{E}^{n,*}$ defined as the product of the fixed orientation and translation invariant measure on $U_l$.}

\proof
See Corollary 3 in \S4 of \cite{MorC}.
\qed

\paragraph{Remark}\setcounter{equation}{0}
The scissors group $\mathbb{S}(\mathbb{E}^n)$ has the structure of a real vector space.

\paragraph{Definition}\setcounter{equation}{0}
Let $\mathrm{M}$ be a Levi subgroup in $\mathcal{L}$. Let $\mathcal{Y}_\mathrm{M}$ and $\mathcal{Z}_\mathrm{M}$ be positive $(\mathrm{G},\mathrm{M})$-orthogonal sets. Then $\mathcal{Y}_\mathrm{M}$ and $\mathcal{Z}_\mathrm{M}$ are said to be \emph{scissors congruent as $(\mathrm{G},\mathrm{M})$-orthogonal sets} if the convex hull of $\mathcal{Y}_\mathrm{M}$ in $\mathfrak{a}_\mathrm{M}^\mathrm{G}$ is translational scissors congruent to the convex hull of $\mathcal{Z}_\mathrm{M}$ in $\mathfrak{a}_\mathrm{M}^\mathrm{G}$. Define the \emph{scissors group} of $\mathfrak{a}_\mathrm{M}^\mathrm{G}$, denoted by $\mathbb{S}(\mathfrak{a}_\mathrm{M}^\mathrm{G})$\label{sym123}, to be the quotient of the free abelian group generated by the positive $(\mathrm{G},\mathrm{M})$-orthogonal sets modulo translational scissors congruence in $\mathfrak{a}_\mathrm{M}^\mathrm{G}$.

\paragraph{Remark}\setcounter{equation}{0}
The Hadwiger invariants of $\mathcal{Y}_\mathrm{M}$ are supported on the strict flags $\Phi$ in $\mathfrak{a}_\mathrm{M}^\mathrm{G}$ of the form
\begin{eqnarray}
\label{hadwiger remark flag}
\Phi&=&
\mathfrak{a}_\mathrm{M}^{\mathrm{L}^0}\supset
\mathfrak{a}_\mathrm{M}^{\mathrm{L}^1}\supset\dots\supset
\mathfrak{a}_\mathrm{M}^{\mathrm{L}^l}
\end{eqnarray}
where
\begin{eqnarray}
\label{hadwiger remark levi}
\mathrm{L}^0\supset\mathrm{L}^1\supset\dots\supset\mathrm{L}^l
\end{eqnarray}
is a nested chain of Levi subgroups in $\mathcal{L}(\mathrm{M})$. Denote the collection of such $\Phi$ by $\mathcal{F}_\mathrm{s}(\mathfrak{a}_\mathrm{M}^\mathrm{G})$. 

A rigging $\mathbf{r}$ of $\Phi$ is equivalent to a nested chain of parabolic subgroups
\begin{eqnarray}
\label{hadwiger remark rigging}
\mathbf{r}&=&
\mathrm{Q}^0\supset\mathrm{Q}^1\supset\dots\supset\mathrm{Q}^l
\end{eqnarray}
such that $\mathrm{Q}^i$ is a parabolic subgroup in $\mathcal{P}(\mathrm{L}^i)$ for each $i$ among $0,1,2,\dots,l$.

Fix a vector $\xi$ in general position in $\mathfrak{a}_\mathrm{M}^\mathrm{G}$. Then $\xi$ defines a total order on $\Delta_\mathrm{L}^\mathrm{G}$ for each Levi subgroup $\mathrm{L}$ in $\mathcal{L}(\mathrm{M})$, which induces a consistent choice of signs for all rigged flags $\Phi^\mathbf{r}$ with $\Phi$ in $\mathcal{F}_\mathrm{s}(\mathfrak{a}_\mathrm{M}^\mathrm{G})$. More precisely choose an element $s$ of the Weyl group that stablizes $\mathrm{M}$ such that the parabolic subgroup $s\mathrm{Q}^l$ is standard. The nested chain
\begin{eqnarray}
s\mathrm{Q}^0\supset s\mathrm{Q}^1\supset\dots\supset s\mathrm{Q}^l
\end{eqnarray}
determines a sequence of positive roots $\alpha^i$ where $i$ ranges among $1,2,\dots,l$. Let $\sigma$ be the permutation on $l$ letters such that $\alpha^{\sigma(i)}$ is strictly increasing with respect to the total order determined by $\xi$. Then the sign of the rigging $\mathbf{r}$ is equal to
\begin{eqnarray}
\mathrm{sign}(\mathbf{r}) &=&
\mathrm{sign}\big(\mathrm{Det}(s)\big) \mathrm{sign}\big(\sigma\big).
\end{eqnarray}

Each subspace $\mathfrak{a}_\mathrm{M}^\mathrm{L}$ of $\mathfrak{a}_\mathrm{M}^\mathrm{G}$ is equipped with the translation invariant measure determined by the coweight lattice. This choice of orientations and measures determines the numerical values of the Hadwiger invariants $\mathrm{Had}_\Phi(\mathcal{Y}_\mathrm{M})$.

\paragraph{Definition}\setcounter{equation}{0}
Let $\mathrm{M}$ be a Levi subgroup in $\mathcal{L}$, let $X$ be an element of $\mathfrak{m}(\mathbb{Q}_S)$. Define the \emph{scissors-congruence-valued orbital integral}, or \emph{orbital integrohedron} $\mathbb{J}_\mathrm{M}^\mathrm{G}(X,~)$\label{sym124} to be the vector-valued distribution on $\mathfrak{g}(\mathbb{Q}_S)$ taking values in $\bigoplus_{\Phi\in\mathcal{F}_\mathrm{s}(\mathfrak{a}_\mathrm{M}^\mathrm{G})}\mathbb{C}$ such that 
\begin{eqnarray}
\label{orbital integrohedron definition}
&&
\forall f_S\in\mathcal{S}(\mathfrak{g}(\mathbb{Q}_S))
\\\nonumber\\\nonumber
&&
\mathbb{J}_\mathrm{M}^\mathrm{G}(X,f_S)
=
\Bigg(\sum_{\mathbf{r}\in\mathcal{R}ig(\Phi)}\mathrm{sign}(\mathbf{r})
J_\mathrm{M}^{\mathrm{L}^l}(X,f_{S,\mathrm{Q}^l})
\Bigg)_{\Phi\in\mathcal{F}_\mathrm{s}(\mathfrak{a}_\mathrm{M}^\mathrm{G})}
\end{eqnarray}
where $\Phi$, $\mathbf{r}$, $\mathrm{L}^l$ and $\mathrm{Q}^l$ are related as in \eqref{hadwiger remark flag}, \eqref{hadwiger remark levi} and \eqref{hadwiger remark rigging}.

For each Levi subgroup $\mathrm{M}_S$ in $\mathcal{L}^{\mathrm{G}_S}$ define $\mathbb{J}_{\mathrm{M}_S}^{\mathrm{G}_S}(X,~)$ by the analogous formula.

\paragraph{Lemma}\setcounter{equation}{0}
\emph{Let $\mathrm{M}$ be a Levi subgroup in $\mathcal{L}$, let $X$ be an element of $\mathfrak{m}(\mathbb{Q}_S)$, let $f_S$ be a real-valued Schwartz function on $\mathfrak{g}(\mathbb{Q}_S)$. Then $\mathbb{J}_\mathrm{M}^\mathrm{G}(X,f_S)$ defines a unique element of $\mathbb{S}(\mathfrak{a}_\mathrm{M}^\mathrm{G})$.}

\proof
The Schwartz function $f_S$ is real-valued, so $\mathbb{J}_\mathrm{M}^\mathrm{G}(X,f_S)$ is a collection of real numbers indexed by the strict flags $\Phi$ in $\mathcal{F}_\mathrm{s}(\mathfrak{a}_\mathrm{M}^\mathrm{G})$. It suffices to verify \eqref{hadwiger characterization i} and \eqref{hadwiger characterization ii}.

The left hand side of \eqref{hadwiger characterization i} is equal to
\begin{eqnarray}
\label{orbital integrohedron lemma proof first equality}
&&
\sum_{\Phi'\in\mathcal{F}_\mathrm{s}(\Phi,i)}~
\sum_{\mathbf{r}'\in\mathcal{R}ig(\Phi')}
\mathrm{sign}(\mathbf{r}')
J_\mathrm{M}^{\mathrm{L}^l{}'}(X,f_{S,\mathrm{Q}^l{}'})
\\\nonumber\\\nonumber
&=&
\sum_{\substack{\mathrm{Q}^0{}'\supset\mathrm{Q}^1{}'\supset\dots\supset\mathrm{Q}^l{}'\\
\forall j=0,1,2,\dots,l~\mathrm{Q}^j{}'\in\mathcal{P}(\mathrm{L}^j{}')\\
\forall j=0,1,2,\dots,l~j\neq i\Rightarrow\mathrm{L}^j{}'=\mathrm{L}^j}}
\mathrm{sign}(\mathbf{r}')
J_\mathrm{M}^{\mathrm{L}^l{}'}(X,f_{S,\mathrm{Q}^l{}'})
\end{eqnarray}
where $\Phi$, $\mathrm{L}^j$ and $\Phi'$, $\mathbf{r}'$, $\mathrm{L}^j{}'$, $\mathrm{Q}^j{}'$ are related as in \eqref{hadwiger remark flag}, \eqref{hadwiger remark levi} and \eqref{hadwiger remark rigging}. The summands of the right hand side of \eqref{orbital integrohedron lemma proof first equality} with a fixed minimal term $\mathrm{Q}^l{}'$ are in (1,1) correspondence with the sequences of roots $(\alpha^j)$ that are positive with respect to $\mathrm{Q}^l{}'$ such that
\begin{eqnarray}
\alpha^1,\alpha^2,\dots,\alpha^{i-1},\alpha^{i+2},\alpha^{i+3},\dots,\alpha^l
\end{eqnarray}
are determined by $\Phi$. Hence the summand
\begin{eqnarray}
J_\mathrm{M}^{\mathrm{L}^l{}'}(X,f_{S,\mathrm{Q}^l{}'})
\end{eqnarray}
appears twice on the right hand side of \eqref{orbital integrohedron lemma proof first equality} with opposite signs, so \eqref{orbital integrohedron lemma proof first equality} vanishes. Hence $\mathbb{J}_\mathrm{M}^\mathrm{G}(X,f_S)$ satisfies \eqref{hadwiger characterization i}.

The left hand side of \eqref{hadwiger characterization ii} is equal to
\begin{eqnarray}
\label{orbital integrohedron lemma proof second equality}
&&
\sum_{\Phi'\in\mathcal{F}_\mathrm{s}(\Phi,l)}~
\sum_{\mathbf{r}'\in\mathcal{R}ig(\Phi')}
\mathrm{sign}(\mathbf{r}')
J_\mathrm{M}^{\mathrm{L}^l{}'}(X,f_{S,\mathrm{Q}^l{}'})
\wedge\mathbf{d}\mathfrak{a}_\mathrm{M}^{\mathrm{L}^l{}'}
\\\nonumber\\\nonumber
&=&
\sum_{\substack{\mathrm{Q}^0\supset\mathrm{Q}^1\supset\dots\supset\mathrm{Q}^{l-1}\supset\mathrm{Q}^l{}'\\
\mathrm{Q}^l{}'\in\mathcal{P}(\mathrm{L}^l{}')}}
\mathrm{sign}(\mathbf{r}')
J_\mathrm{M}^{\mathrm{L}^l{}'}(X,f_{S,\mathrm{Q}^l{}'})
\wedge\mathbf{d}\mathfrak{a}_\mathrm{M}^{\mathrm{L}^l{}'}
\end{eqnarray}
where $\Phi$, $\mathrm{Q}^j$ and $\Phi'$, $\mathbf{r}'$, $\mathrm{L}^l{}'$, $\mathrm{Q}^l{}'$ are related as in \eqref{hadwiger remark flag}, \eqref{hadwiger remark levi} and \eqref{hadwiger remark rigging}. It suffices to show that \eqref{orbital integrohedron lemma proof second equality} vanishes as a differential form on $\mathfrak{a}_\mathrm{M}^{\mathrm{L}^{l-1}}$. Let $\mathrm{L}$ be a Levi subgroup in $\mathcal{L}^{\mathrm{L}^{l-1}}(\mathrm{M})$ such that
\begin{eqnarray}
\mathrm{dim}(\mathrm{A}_\mathrm{M}/\mathrm{A}_\mathrm{L})
&=&1,
\end{eqnarray}
denote by $\mathbf{d}\mathfrak{a}_\mathrm{L}^{\mathrm{L}^{l-1}}$\label{sym125} the differential form on $\mathfrak{a}_\mathrm{M}^{\mathrm{L}^{l-1}}$ defined by pulling back the volume form on $\mathfrak{a}_\mathrm{L}^{\mathrm{L}^{l-1}}$ along the natural projection
\begin{eqnarray}
\mathfrak{a}_\mathrm{M}^{\mathrm{L}^{l-1}}
&\longrightarrow&
\mathfrak{a}_\mathrm{L}^{\mathrm{L}^{l-1}}.
\end{eqnarray}
Then the orthogonal projection of \eqref{orbital integrohedron lemma proof second equality} onto the one dimensional subspace
\begin{eqnarray}
\mathrm{span}(\mathbf{d}\mathfrak{a}_\mathrm{L}^{\mathrm{L}^{l-1}})
&\subset&
\bigwedge{}^{\mathrm{dim}(\mathrm{A}_\mathrm{M}/\mathrm{A}_\mathrm{G})-l}
~\mathfrak{a}_\mathrm{M}^{\mathrm{L}^{l-1},*}
\end{eqnarray}
is equal to
\begin{eqnarray}
\label{orbital integrohedron lemma proof third expression}
&&
\sum_{\substack{\mathrm{Q}^0\supset\mathrm{Q}^1\supset\dots\supset\mathrm{Q}^{l-1}\supset\mathrm{Q}^l{}'\\
\mathrm{Q}^l{}'\in\mathcal{P}(\mathrm{L}^l{}')}}
\mathrm{sign}(\mathbf{r}')
d_\mathrm{M}^{\mathrm{L}^{l-1}}(\mathrm{L},\mathrm{L}^l{}')
J_\mathrm{M}^{\mathrm{L}^l{}'}(X,f_{S,\mathrm{Q}^l{}'})
\wedge\mathbf{d}\mathfrak{a}_\mathrm{L}^{\mathrm{L}^{l-1}}
\end{eqnarray}
by the definition of the constant $d_\mathrm{M}^{\mathrm{L}^{l-1}}(\mathrm{L},\mathrm{L}^l{}')$ in Remark \ref{angle jacobian}. The summation in \eqref{orbital integrohedron lemma proof third expression} is taken over the set
\begin{eqnarray}
\label{orbital integrohedron lemma proof fourth set}
&&
\bigg\{
(\mathrm{L}^l{}',\mathrm{Q}^l{}'\cap\mathrm{L}^{l-1}):~
\begin{array}{l}
\mathrm{L}^l{}'\in\mathcal{L}^{\mathrm{L}^{l-1}}(\mathrm{M}),~
\mathrm{dim}(\mathfrak{a}_{\mathrm{L}^l{}'}^\mathrm{G})=l,\\
\mathrm{Q}^l{}'\in\mathcal{P}(\mathrm{L}^l{}'),~
\mathrm{Q}^l{}'\subset\mathrm{Q}^{l-1}
\end{array}\bigg\}
\\\nonumber\\\nonumber
&&\quad
\subset\mathcal{L}^{\mathrm{L}^{l-1}}(\mathrm{M})
\times\mathcal{F}^{\mathrm{L}^{l-1}}(\mathrm{M})
\end{eqnarray}
Let $\xi^{l-1}$ be the projection of the vector $\xi$ in $\mathfrak{a}_\mathrm{M}^\mathrm{G}$ used to define the orientation onto $\mathfrak{a}_\mathrm{M}^{\mathrm{L}^{l-1}}$. Then $\xi^{l-1}$ determines a partial map 
\begin{eqnarray}
s^{l-1}:\quad
\mathcal{L}^{\mathrm{L}^{l-1}}(\mathrm{M})\times\mathcal{L}^{\mathrm{L}^{l-1}}(\mathrm{M})
&\longrightarrow&
\mathcal{F}^{\mathrm{L}^{l-1}}(\mathrm{M})\times\mathcal{F}^{\mathrm{L}^{l-1}}(\mathrm{M})
\end{eqnarray}
as in Remark \ref{section s} which is positive in the sense that
\begin{eqnarray}
\mathrm{sign}(\mathbf{r})&=&1
\end{eqnarray}
if the rigging $\mathbf{r}$ corresponds to the element $(\mathrm{L}^l{}',\mathrm{Q}^{l,+})$ in the set \eqref{orbital integrohedron lemma proof fourth set} where $\mathrm{Q}^{l,+}$ denotes the second component of $s^{l-1}(\mathrm{L},\mathrm{L}^l{}')$. Let $\mathrm{Q}^{l,-}$ denote the opposite parabolic of $\mathrm{Q}^{l,+}$, then \eqref{orbital integrohedron lemma proof fourth set} is equal to the disjoint union
\begin{eqnarray}
&&
\bigg\{(\mathrm{L}^l{}',\mathrm{Q}^{l,+}):~
\mathrm{L}^l{}'\in\mathcal{L}^{\mathrm{L}^{l-1}}(\mathrm{M})\bigg\}\coprod
\\\nonumber
&&\quad
\coprod\bigg\{(\mathrm{L}^l{}',\mathrm{Q}^{l,-}):~
\mathrm{L}^l{}'\in\mathcal{L}^{\mathrm{L}^{l-1}}(\mathrm{M})\bigg\},
\end{eqnarray}
hence \eqref{orbital integrohedron lemma proof third expression} is equal to
\begin{eqnarray}
&&
\Bigg(
\sum_{\mathrm{L}^l{}'\in\mathcal{L}^{\mathrm{L}^{l-1}}(\mathrm{M})}
d_\mathrm{M}^{\mathrm{L}^{l-1}}(\mathrm{L},\mathrm{L}^l{}')
J_\mathrm{M}^{\mathrm{L}^l{}'}(X,f_{S,\mathrm{Q}^{l,+}})
\\\nonumber\\\nonumber
&&\quad
-\sum_{\mathrm{L}^l{}'\in\mathcal{L}^{\mathrm{L}^{l-1}}(\mathrm{M})}
d_\mathrm{M}^{\mathrm{L}^{l-1}}(\mathrm{L},\mathrm{L}^l{}')
J_\mathrm{M}^{\mathrm{L}^l{}'}(X,f_{S,\mathrm{Q}^{l,-}})
\Bigg)
\wedge\mathbf{d}\mathfrak{a}_\mathrm{L}^{\mathrm{L}^{l-1}}
\\\nonumber\\
\label{orbital integrohedron lemma proof fifth equality}
&=&
\bigg(J_\mathrm{L}^{\mathrm{L}^{l-1}}(X^\mathrm{L},f_{S,\mathrm{Q}^{l-1}})
-J_\mathrm{L}^{\mathrm{L}^{l-1}}(X^\mathrm{L},f_{S,\mathrm{Q}^{l-1}})\bigg)
\wedge\mathbf{d}\mathfrak{a}_\mathrm{L}^{\mathrm{L}^{l-1}}
\\\nonumber\\\nonumber
&=&0
\end{eqnarray}
where the equality \eqref{orbital integrohedron lemma proof fifth equality} follows from the $S$-local version of \eqref{weighted orbital integral descent} and the fact that $\mathrm{Q}^{l,-}$ is the second component of $s^{l-1,-}(\mathrm{L},\mathrm{L}^l{}')$ where $s^{l-1,-}(~,~)$ is the partial map determined as in Remark \ref{section s} by the vector $-\xi^{l-1}$. Hence \eqref{orbital integrohedron lemma proof third expression} vanishes for every $\mathrm{L}$, so \eqref{orbital integrohedron lemma proof second equality} vanishes as a differential form on $\mathfrak{a}_\mathrm{M}^\mathrm{G}$. Hence $\mathbb{J}_\mathrm{M}^\mathrm{G}(X,f_S)$ satisfies \eqref{hadwiger characterization ii}.
\qed

\paragraph{Remark}\setcounter{equation}{0}
In \S23 of \cite{Kott} Kottwitz defined weight factors and weighted orbital integrals taking values in the complexified $K$-group of the toric variety of the fan of the root hyperplanes in $\mathfrak{a}_\mathrm{M}^{\mathrm{G},*}$. In \S4 of \cite{MorK} Morelli proved that this $K$-group is the additive group of translational scissors congruent classes of positive $(\mathrm{G},\mathrm{M})$-orthogonal sets whose vertices are contained in the coweight lattice in $\mathfrak{a}_\mathrm{M}^\mathrm{G}$.

\paragraph{Definition}\setcounter{equation}{0}
Let $\mathrm{M}$ be a Levi subgroup in $\mathcal{L}$. Define the \emph{total scissors ring} of $\mathfrak{a}_\mathrm{M}^\mathrm{G}$, denoted by $\mathbb{S}(\mathfrak{a}_{\mathcal{L}(\mathrm{M})}^\mathrm{G})$\label{sym126}, to be the direct sum
\begin{eqnarray}
\mathbb{S}(\mathfrak{a}_{\mathcal{L}(\mathrm{M})}^\mathrm{G})
&=&
\bigoplus_{\mathrm{L}\in\mathcal{L}^\mathrm{G}(\mathrm{M})}
\mathbb{S}(\mathfrak{a}_\mathrm{L}^\mathrm{G}).
\end{eqnarray}

Define a bilinear product $\boxtimes$\label{sym127} on $\mathbb{S}(\mathfrak{a}_{\mathcal{L}(\mathrm{M})}^\mathrm{G})$ by
\begin{eqnarray}
&&
\forall \mathrm{L}_1,\mathrm{L}_2\in\mathcal{L}^\mathrm{G}(\mathrm{M})
~\forall [\mathcal{Y}_{\mathrm{L}_i}]\in\mathbb{S}(\mathfrak{a}_{\mathrm{L}_i}^\mathrm{G})
~\textrm{where } i=1,2
\\\nonumber
&&\quad
[\mathcal{Y}_{\mathrm{L}_1}]\boxtimes[\mathcal{Y}_{\mathrm{L}_2}]
\\\nonumber\\\nonumber
&&\qquad
=\left\{\begin{array}{cl}
j^*\Big([\mathcal{Y}_{\mathrm{L}_1}\times\mathcal{Y}_{\mathrm{L}_2}]\Big)
& \textrm{if the natural map}\\
&\quad j:~\mathfrak{a}_{\mathrm{L}_1\cap\mathrm{L}_2}^\mathrm{G}\longrightarrow\mathfrak{a}_{\mathrm{L}_1}^\mathrm{G}\oplus\mathfrak{a}_{\mathrm{L}_2}^\mathrm{G}\\
& \textrm{is an isomorphism,}
\\\\
0 & \textrm{otherwise},
\end{array}\right.
\end{eqnarray}
where the $(\mathrm{G}\times\mathrm{G},\mathrm{L}_1\times\mathrm{L}_2)$-family $\mathcal{Y}_{\mathrm{L}_1}\times\mathcal{Y}_{\mathrm{L}_2}$ is well-defined upto translational scissors congruence in $\mathfrak{a}_{\mathrm{L}_1}^\mathrm{G}\oplus\mathfrak{a}_{\mathrm{L}_2}^\mathrm{G}$, and $j^*$ is the homomorphism
\begin{eqnarray}
j^*:\quad
\mathbb{S}(\mathfrak{a}_{\mathrm{L}_1}^\mathrm{G}\oplus\mathfrak{a}_{\mathrm{L}_2}^\mathrm{G})
&\longrightarrow&
\mathbb{S}(\mathfrak{a}_{\mathrm{L}_1\cap\mathrm{L}_2}^\mathrm{G})
\end{eqnarray}
induced by $j$.

\paragraph{Remark}\setcounter{equation}{0}
The total scissors ring $\mathbb{S}(\mathfrak{a}_{\mathcal{L}(\mathrm{M})}^\mathrm{G})$ is graded by\label{sym128}
\begin{eqnarray}
\label{total scissors ring n grading}
\forall n\in\mathbb{N}
&&
\Big(\mathbb{S}(\mathfrak{a}_{\mathcal{L}(\mathrm{M})}^\mathrm{G})\Big)_n
=\bigoplus_{\substack{\mathrm{L}\in\mathcal{L}^\mathrm{G}(\mathrm{M})\\\mathrm{dim}(\mathfrak{a}_\mathrm{L}^\mathrm{G})=n}}
\mathbb{S}(\mathfrak{a}_\mathrm{L}^\mathrm{G}).
\end{eqnarray}
The grading \eqref{total scissors ring n grading} has a refinement into the $\mathcal{L}^\mathrm{G}(\mathrm{M})$-grading defined by
\begin{eqnarray}
\forall \mathrm{L}\in\mathcal{L}^\mathrm{G}(\mathrm{M})
&&
\Big(\mathbb{S}(\mathfrak{a}_{\mathcal{L}(\mathrm{M})}^\mathrm{G})\Big)_\mathrm{L}=\mathbb{S}(\mathfrak{a}_\mathrm{L}^\mathrm{G})
\end{eqnarray}
where the monoid structure on $\mathcal{L}^\mathrm{G}(\mathrm{M})$ is defined by intersection.

Each Levi subgroup $\mathrm{L}$ in $\mathcal{L}^\mathrm{G}(\mathrm{M})$ defines a homogeneous ideal
\begin{eqnarray}
\label{total scissors ring homogeneous ideal}
\bigoplus_{\substack{\mathrm{L}'\in\mathcal{L}^\mathrm{G}(\mathrm{M})\\
\mathrm{L}'\not\supset\mathrm{L}}}
\mathbb{S}(\mathfrak{a}_{\mathrm{L}'}^\mathrm{G})
&\subset&
\mathbb{S}(\mathfrak{a}_{\mathcal{L}(\mathrm{M})}^\mathrm{G}),
\end{eqnarray}
and taking the quotient of $\mathbb{S}(\mathfrak{a}_{\mathcal{L}(\mathrm{M})}^\mathrm{G})$ modulo \eqref{total scissors ring homogeneous ideal} defines the homomorphism\label{sym129}
\begin{eqnarray}
q_\mathrm{M}^\mathrm{L}:\quad
\mathbb{S}(\mathfrak{a}_{\mathcal{L}(\mathrm{M})}^\mathrm{G})
&\longrightarrow&
\mathbb{S}(\mathfrak{a}_{\mathcal{L}(\mathrm{L})}^\mathrm{G}).
\end{eqnarray}
The map $q_\mathrm{M}^\mathrm{G}$ is the quotient modulo the augmentation ideal onto the graded component $\big(\mathbb{S}(\mathfrak{a}_{\mathcal{L}(\mathrm{M})}^\mathrm{G})\big)_\mathrm{G}$, which is a copy of $\mathbb{R}$.

\paragraph{Definition}\setcounter{equation}{0}
Let $\mathrm{M}$ be a Levi subgroup in $\mathcal{L}$, let $X$ be an element of $\mathfrak{m}(\mathbb{Q}_S)$. Define the \emph{total orbital integrohedron} $\mathcal{J}_\mathrm{M}^\mathrm{G}(X,~)$\label{sym130} to be the $\mathbb{S}(\mathfrak{a}_{\mathcal{L}(\mathrm{M})}^\mathrm{G})\otimes\mathbb{C}$-valued distribution on $\mathfrak{g}(\mathbb{Q}_S)$ such that
\begin{eqnarray}
\label{total orbital integrohedron definition}
&&
\forall f_S\in\mathcal{S}(\mathfrak{g}(\mathbb{Q}_S)
\\\nonumber\\\nonumber
&&\quad
\mathcal{J}_\mathrm{M}^\mathrm{G}(X,f_S)
=\bigg((-1)^{\mathrm{dim}(\mathrm{A}_\mathrm{L}/\mathrm{A}_\mathrm{G})}
|\mathrm{W}_0^\mathrm{L}||\mathrm{W}_0^\mathrm{G}|^{-1}
\mathbb{J}_\mathrm{L}^\mathrm{G}(X^\mathrm{L},f_S)\bigg)_{\mathrm{L}\in\mathcal{L}^\mathrm{G}(\mathrm{M})}
\end{eqnarray}
For each Levi subgroup $\mathrm{M}_S$ in $\mathcal{L}^{\mathrm{G}_S}$ define $\mathcal{J}_{\mathrm{M}_S}^{\mathrm{G}_S}(X,~)$ by the analogous formula.

\paragraph{Lemma}\setcounter{equation}{0}
(Induction and splitting of orbital integrohedra)\\
\emph{Let $\mathrm{M}$ be a Levi subgroup in $\mathcal{L}$, let $X$ be an element of $\mathfrak{m}(\mathbb{Q}_S)$, let $f_S$ be a Schwartz function on $\mathfrak{g}(\mathbb{Q}_S)$.
\begin{itemize}
\item Let $\mathrm{L}$ be a Levi subgroup in $\mathcal{L}(\mathrm{M})$, then
\begin{eqnarray}
\label{orbital integrohedron induction}
\mathcal{J}_\mathrm{L}^\mathrm{G}(X^\mathrm{L},f_S) 
&=& 
q_\mathrm{M}^\mathrm{L}\Big(
\mathcal{J}_\mathrm{M}^\mathrm{G}(X,f_S)\Big).
\end{eqnarray}
\item Let $S$ be the set $\{v_1,v_2\}$, let $f_S$ be of the form $f_{v_1}\otimes f_{v_2}$ where $f_{v_i}$ is a Schwartz function on $\mathfrak{g}(\mathbb{Q}_{v_i})$ where the index $i$ is 1 or 2, then
\begin{eqnarray}
\label{orbital integrohedron splitting}
\qquad\qquad
\mathcal{J}_\mathrm{M}^\mathrm{G}(X,f_{v_1}\otimes f_{v_2}) 
&=& 
\mathcal{J}_\mathrm{M}^\mathrm{G}(X,f_{v_1})
\boxtimes \mathcal{J}_\mathrm{M}^\mathrm{G}(X,f_{v_2}) .
\end{eqnarray}
\end{itemize}
Local identities analogous to \eqref{orbital integrohedron induction} also hold for $\mathrm{G}_v$.}

\proof
The induction identity \eqref{orbital integrohedron induction} follows from the definition \eqref{total orbital integrohedron definition}.

The splitting identity \eqref{orbital integrohedron splitting} follows by comparing the Hadwiger invariants of the two sides of \eqref{orbital integrohedron splitting}. Retain the notations of \eqref{orbital integrohedron definition}. Let $\mathrm{L}$ be a Levi subgroup in $\mathcal{L}(\mathrm{M})$, let $\Phi$ be a strict flag in $\mathcal{F}_\mathrm{s}(\mathfrak{a}_\mathrm{L}^\mathrm{G})$. It suffices to show that the $\Phi$-components of the Hadwiger invariants of the $\mathrm{L}$-components of the two sides of \eqref{orbital integrohedron splitting} are equal.

The $\mathrm{L}$-component of the right hand side of \eqref{orbital integrohedron splitting} is
\begin{eqnarray}
&&
\sum_{\substack{\mathrm{L}_1,\mathrm{L}_2\in\mathcal{L}(\mathrm{L})\\
\mathfrak{a}_\mathrm{L}^\mathrm{G}\rightarrow\mathfrak{a}_{\mathrm{L}_1}^\mathrm{G}\oplus\mathfrak{a}_{\mathrm{L}_2}^\mathrm{G}\\
\mathrm{isomorphism}}}
(-1)^{\mathrm{dim}(\mathrm{A}_\mathrm{L}/\mathrm{A}_\mathrm{G})}
|\mathrm{W}_0^{\mathrm{L}_1}||\mathrm{W}_0^\mathrm{G}|^{-1}
|\mathrm{W}_0^{\mathrm{L}_2}||\mathrm{W}_0^\mathrm{G}|^{-1}\times
\\\nonumber\\\nonumber
&&\quad
\times\mathbb{J}_{\mathrm{L}_1}^\mathrm{G}(X^{\mathrm{L}_1},f_{v_1})
\boxtimes
\mathbb{J}_{\mathrm{L}_2}^\mathrm{G}(X^{\mathrm{L}_2},f_{v_2}),
\end{eqnarray}
whose Hadwiger invariant has $\Phi$-component equal to
\begin{eqnarray}
\label{orbital integrohedron splitting proof first expression}
&&
\sum_{\mathrm{L}_1,\mathrm{L}_2\in\mathcal{L}(\mathrm{L})}
(-1)^{\mathrm{dim}(\mathrm{A}_\mathrm{L}/\mathrm{A}_\mathrm{G})}
|\mathrm{W}_0^{\mathrm{L}_1}||\mathrm{W}_0^\mathrm{G}|^{-1}
|\mathrm{W}_0^{\mathrm{L}_2}||\mathrm{W}_0^\mathrm{G}|^{-1}\times
\\\nonumber\\\nonumber
&&\quad\times
d_\mathrm{L}^{\mathrm{L}^l}(\mathrm{L}_1^l,\mathrm{L}_2^l)
\bigg(\sum_{\mathbf{r}_1\in\mathcal{R}ig(\Phi_1)}\mathrm{sign}(\mathbf{r}_1)
J_{\mathrm{L}_1}^{\mathrm{L}_1^l}(X^{\mathrm{L}_1},f_{v_1,\mathrm{Q}_1^l})\bigg)\times
\\\nonumber\\\nonumber
&&\qquad
\times\bigg(\sum_{\mathbf{r}_2\in\mathcal{R}ig(\Phi_2)}\mathrm{sign}(\mathbf{r}_2)
J_{\mathrm{L}_2}^{\mathrm{L}_2^l}(X^{\mathrm{L}_2},f_{v_2,\mathrm{Q}_2^l})\bigg)
\end{eqnarray}
where $\Phi_i$ denotes the strict flag in $\mathfrak{a}_{\mathrm{L}_i}^\mathrm{G}$ induced by $\Phi$, and $\Phi_i$, $\mathbf{r}_i$, $\mathrm{L}_i^l$ and $\mathrm{Q}_i^l$ are related as in \eqref{hadwiger remark flag}, \eqref{hadwiger remark levi} and \eqref{hadwiger remark rigging} where the index $i$ is 1 or 2. Combine $\mathbf{r}_1$ and $\mathbf{r}_2$ into a rigging $\mathbf{r}$ of $\Phi$, the expression \eqref{orbital integrohedron splitting proof first expression} becomes
\begin{eqnarray}
\label{orbital integrohedron splitting proof second equality}
&&
\sum_{\mathrm{L}_1,\mathrm{L}_2\in\mathcal{L}(\mathrm{L})}
(-1)^{\mathrm{dim}(\mathrm{A}_\mathrm{L}/\mathrm{A}_\mathrm{G})}
|\mathrm{W}_0^{\mathrm{L}_1}||\mathrm{W}_0^\mathrm{G}|^{-1}
|\mathrm{W}_0^{\mathrm{L}_2}||\mathrm{W}_0^\mathrm{G}|^{-1}\times
\\\nonumber\\\nonumber
&&\quad
\times\sum_{\mathrm{Q}_1^l\in\mathcal{P}(\mathrm{L}_1^l)}
\sum_{\mathrm{Q}_2^l\in\mathcal{P}(\mathrm{L}_2^l)}
\mathrm{sign}(\mathbf{r})
d_\mathrm{L}^{\mathrm{L}^l}(\mathrm{L}_1^l,\mathrm{L}_2^l)\times
\\\nonumber\\\nonumber
&&\qquad
\times J_{\mathrm{L}_1}^{\mathrm{L}_1^l}(X^{\mathrm{L}_1},f_{v_1,\mathrm{Q}_1^l})
J_{\mathrm{L}_2}^{\mathrm{L}_2^l}(X^{\mathrm{L}_2},f_{v_2,\mathrm{Q}_2^l})
\\\nonumber\\\nonumber
&=&
\sum_{\mathrm{L}_1,\mathrm{L}_2\in\mathcal{L}(\mathrm{L})}
(-1)^{\mathrm{dim}(\mathrm{A}_\mathrm{L}/\mathrm{A}_\mathrm{G})}
|\mathrm{W}_0^{\mathrm{L}_1}||\mathrm{W}_0^\mathrm{G}|^{-1}
|\mathrm{W}_0^{\mathrm{L}_2}||\mathrm{W}_0^\mathrm{G}|^{-1}\times
\\\nonumber\\\nonumber
&&\quad
\times\sum_{\mathrm{Q}_1^l\in\mathcal{P}(\mathrm{L}_1^l)}
\sum_{\mathrm{Q}_2^l\in\mathcal{P}(\mathrm{L}_2^l)}
\mathrm{sign}(\mathbf{r})
d_\mathrm{L}^{\mathrm{L}^l}(\mathrm{L}_1^l,\mathrm{L}_2^l)\times
\\\nonumber\\\nonumber
&&\qquad
\times\sum_{\mathrm{L}_1^\circ\in\mathcal{L}^{\mathrm{L}_1^l}(\mathrm{L})}
d_\mathrm{L}^{\mathrm{L}_1^l}(\mathrm{L}_1,\mathrm{L}_1^\circ)
J_\mathrm{L}^{\mathrm{L}_1^\circ}(X^\mathrm{L},f_{v_1,\mathrm{Q}_1^{l,\circ}})\times
\\\nonumber\\\nonumber
&&\qquad\quad
\times\sum_{\mathrm{L}_2^\circ\in\mathcal{L}^{\mathrm{L}_2^l}(\mathrm{L})}
d_\mathrm{L}^{\mathrm{L}_2^l}(\mathrm{L}_2,\mathrm{L}_2^\circ)
J_\mathrm{L}^{\mathrm{L}_2^\circ}(X^\mathrm{L},f_{v_2,\mathrm{Q}_2^{l,\circ}})
\end{eqnarray}
by \eqref{weighted orbital integral descent}, where $\mathrm{Q}_i^{l,\circ}$ is the second component of $s_{\mathrm{Q}_i^l}(\mathrm{L}_i,\mathrm{L}_i^\circ)$ for the choice of a collection of partial maps $(s_{\mathrm{Q}_i^l})_{\mathrm{Q}_i^l\in\mathcal{P}(\mathrm{L}_1')}$ such that 
\begin{eqnarray}
\forall \mathrm{Q}_i^l\in\mathcal{P}(\mathrm{L}_i^l)
&&
\mathrm{Q}_i^{l,\circ}\subset\mathrm{Q}_i^l
\end{eqnarray}
where the index $i$ is 1 or 2. The choice of such a collection is equivalent to the choice of a collection of vectors $(\xi_{\mathrm{Q}_i^l})_{\mathrm{Q}_i^l\in\mathcal{P}(\mathrm{L}_1')}$
\begin{eqnarray}
\forall \mathrm{Q}_i^l\in\mathcal{P}(\mathrm{L}_i^l)
&&
\xi_{\mathrm{Q}_i^l}\in\mathfrak{a}_\mathrm{L}^{\mathrm{L}_i^l}
\end{eqnarray}
in general position and positive with respect to $\mathrm{Q}_i^l$. By definition
\begin{eqnarray}
d_\mathrm{L}^{\mathrm{L}^l}(\mathrm{L}_1^l,\mathrm{L}_2^l)
d_\mathrm{L}^{\mathrm{L}_1^l}(\mathrm{L}_1,\mathrm{L}_1^\circ)
d_\mathrm{L}^{\mathrm{L}_2^l}(\mathrm{L}_2,\mathrm{L}_2^\circ)
&=&
d_\mathrm{L}^{\mathrm{L}^l}(\mathrm{L}_1^\circ,\mathrm{L}_2^\circ),
\end{eqnarray}
hence the right hand side of \eqref{orbital integrohedron splitting proof second equality} is equal to
\begin{eqnarray}
\label{orbital integrohedron splitting proof third expression}
&&
(-1)^{\mathrm{dim}(\mathrm{A}_\mathrm{L}/\mathrm{A}_\mathrm{G})}
\sum_{\mathbf{r}\in\mathcal{R}ig(\Phi)}
\mathrm{sign}(\mathbf{r})\times
\\\nonumber\\\nonumber
&&\quad\times
\sum_{\mathrm{L}_1,\mathrm{L}_2\in\mathcal{L}(\mathrm{L})}
\sum_{\substack{\mathrm{L}_1^\circ\in\mathcal{L}^{\mathrm{L}_1^l}(\mathrm{L})\\
\mathrm{L}_2^\circ\in\mathcal{L}^{\mathrm{L}_2^l}(\mathrm{L})}}
|\mathrm{W}_0^{\mathrm{L}_1}||\mathrm{W}_0^\mathrm{G}|^{-1}
|\mathrm{W}_0^{\mathrm{L}_2}||\mathrm{W}_0^\mathrm{G}|^{-1}\times
\\\nonumber\\\nonumber
&&\qquad\times
d_\mathrm{L}^{\mathrm{L}^l}(\mathrm{L}_1^\circ,\mathrm{L}_2^\circ)
J_\mathrm{L}^{\mathrm{L}_1^\circ}(X^\mathrm{L},f_{v_1,\mathrm{Q}_1^{l,\circ}})
J_\mathrm{L}^{\mathrm{L}_2^\circ}(X^\mathrm{L},f_{v_2,\mathrm{Q}_2^{l,\circ}})
\end{eqnarray}
where $(\mathrm{Q}_1^{l,\circ},\mathrm{Q}_2^{l,\circ})$ is the image of $(\mathrm{L}_1^\circ,\mathrm{L}_2^\circ)$ under the partial map $s_{\mathrm{Q}_1^l,\mathrm{Q}_2^l}$ determined by the vector
\begin{eqnarray}
\frac{\xi_{\mathrm{Q}_1^l}}{2}-\frac{\xi_{\mathrm{Q}_2^l}}{2}
&\in&
\mathfrak{a}_\mathrm{L}^{\mathrm{L}^l}=\mathfrak{a}_\mathrm{L}^{\mathrm{L}_1^l}\oplus\mathfrak{a}_\mathrm{L}^{\mathrm{L}_2^l}.
\end{eqnarray}
For each fixed $\mathrm{L}$ and fixed $\Phi$ with minimal term $\mathrm{L}^l$, the map
\begin{eqnarray}
&&
\bigg\{(\mathrm{L}_1,\mathrm{L}_2,\mathrm{L}_1^\circ,\mathrm{L}_2^\circ):~
\begin{array}{l}
\mathrm{L}_1\in\mathcal{L}(\mathrm{L}),~
\mathrm{L}_2\in\mathcal{L}(\mathrm{L}),\\
\mathrm{L}_1^\circ\in\mathcal{L}^{\mathrm{L}_1^l}(\mathrm{L}),~
\mathrm{L}_2^\circ\in\mathcal{L}^{\mathrm{L}_2^l}(\mathrm{L})
\end{array}\bigg\}
\\\nonumber\\\nonumber
&&\quad\longrightarrow\quad
\bigg\{(\mathrm{L}_1^\circ,\mathrm{L}_2^\circ):~
\mathrm{L}_1^\circ\in\mathcal{L}^{\mathrm{L}^l}(\mathrm{L}),~
\mathrm{L}_2^\circ\in\mathcal{L}^{\mathrm{L}^l}(\mathrm{L})\bigg\}
\end{eqnarray}
defined by
\begin{eqnarray}
\label{orbital integrohedron splitting proof fourth map}
(\mathrm{L}_1,\mathrm{L}_2,\mathrm{L}_1^\circ,\mathrm{L}_2^\circ)
&\mapsto&
(\mathrm{L}_1^\circ,\mathrm{L}_2^\circ)
\end{eqnarray}
is surjective. Fix Levi subgroups $\mathrm{L}_1$ and $\mathrm{L}_2$ in $\mathcal{L}(\mathrm{L})$ and parabolic subgroups $\mathrm{P}_1$ in $\mathcal{P}(\mathrm{L}_1)$ and $\mathrm{P}_2$ in $\mathcal{P}(\mathrm{L}_2)$, and let $\mathrm{P}_{1,2}$ be the parabolic subgroup in $\mathcal{P}(\mathrm{L})$ whose unipotent radical $\mathrm{N}_{\mathrm{P}_{1,2}}$ contains $\mathrm{N}_{\mathrm{P}_1}$ and $\mathrm{N}_{\mathrm{P}_2}$, then the map \eqref{orbital integrohedron splitting proof fourth map} restricts to a bijection between Levi subgroups $\mathrm{L}_1^\circ$ and $\mathrm{L}_2^\circ$ that are standard with respect to $\mathrm{P}_1$ and $\mathrm{P}_2$ or $\mathrm{P}_{1,2}$. Hence each fiber of the map \eqref{orbital integrohedron splitting proof fourth map} containing the pair $(\mathrm{L}_1,\mathrm{L}_2)$ has cardinality
\begin{eqnarray}
\frac{|\mathrm{W}_0^\mathrm{G}|\big/|\mathrm{W}_0^\mathrm{L}|}
{\Big(|\mathrm{W}_0^\mathrm{G}|\big/|\mathrm{W}_0^{\mathrm{L}_1}|\Big)
\Big(|\mathrm{W}_0^\mathrm{G}|\big/|\mathrm{W}_0^{\mathrm{L}_2}|\Big)},
\end{eqnarray}
hence \eqref{orbital integrohedron splitting proof third expression} is equal to
\begin{eqnarray}
&&
(-1)^{\mathrm{dim}(\mathrm{A}_\mathrm{L}/\mathrm{A}_\mathrm{G})}
|\mathrm{W}_0^{\mathrm{L}}||\mathrm{W}_0^\mathrm{G}|^{-1}
\sum_{\mathbf{r}\in\mathcal{R}ig(\Phi)}
\mathrm{sign}(\mathbf{r})\times
\\\nonumber\\\nonumber
&&\quad\times
\sum_{\mathrm{L}_1^\circ,\mathrm{L}_2^\circ\in\mathcal{L}^{\mathrm{L}^l}(\mathrm{L})}
d_\mathrm{L}^{\mathrm{L}^l}(\mathrm{L}_1^\circ,\mathrm{L}_2^\circ)
J_\mathrm{L}^{\mathrm{L}_1^\circ}(X^\mathrm{L},f_{v_1,\mathrm{Q}_1^{l,\circ}})
J_\mathrm{L}^{\mathrm{L}_2^\circ}(X^\mathrm{L},f_{v_2,\mathrm{Q}_2^{l,\circ}})
\\\nonumber\\\nonumber
&=&
(-1)^{\mathrm{dim}(\mathrm{A}_\mathrm{L}/\mathrm{A}_\mathrm{G})}
|\mathrm{W}_0^{\mathrm{L}}||\mathrm{W}_0^\mathrm{G}|^{-1}
\sum_{\mathbf{r}\in\mathcal{R}ig(\Phi)}
\mathrm{sign}(\mathbf{r})
J_\mathrm{L}^{\mathrm{L}^l}(X^\mathrm{L},f_{v_1,\mathrm{Q}^l}\otimes f_{v_2,\mathrm{Q}^l})
\end{eqnarray}
by \eqref{weighted orbital integral splitting} where $\Phi$, $\mathbf{r}$, $\mathrm{L}^l$ and $\mathrm{Q}^l$ are related as in \eqref{hadwiger remark flag}, \eqref{hadwiger remark levi} and \eqref{hadwiger remark rigging}, which is equal to the $\Phi$-component of the Hadwiger invariant of the $\mathrm{L}$-component of the left hand side of \eqref{orbital integrohedron splitting}.
\qed

\paragraph{Remark}\setcounter{equation}{0}
The orbital integrohedra $\mathbb{J}_\mathrm{M}^\mathrm{G}(X,~)$ and $\mathcal{J}_\mathrm{M}^\mathrm{G}(X,~)$ could be alternatively defined as distributions on $\mathfrak{g}(\mathbb{Q}_S)$ obtained by integrating against certain polyhedron-valued weight factors with a formula analogous to \eqref{regular weighted orbital integral} if $X$ is regular semisimple.

\section{Appendix: The example of $\mathrm{GL}(2)$}\setcounter{paragraph}{0}

In this appendix some explicit computations for the group $\mathrm{GL}(2)$ are carried out. Similar local computations could be found in \cite{EvEx}. 

\paragraph{Definition}\setcounter{equation}{0}
In this appendix $v$\label{sym169} denotes an odd rational prime $p$, and $|~|$ denotes the $p$-adic absolute value. Let $\mathrm{G}$\label{sym170} denote the group $\mathrm{GL}(2,\mathbb{Q}_p)$, let $\mathrm{M}$ denote the minimal Levi subgroup of $\mathrm{G}$ consisting of the diagonal matrices, then $\mathrm{M}$ is the Levi component of the Borel subgroup $\mathrm{B}$ consisting of the upper triangular matrices. The groups $\mathrm{G}$ and $\mathrm{M}$ are the only standard Levi subgroups of $\mathrm{G}$.

Let $\mathfrak{g}'$\label{sym171} denote the Lie algebra $\mathrm{sl}(2,\mathbb{Q}_p)$, equipped with the adjoint action of $\mathrm{G}$ from the right. Let $\mathfrak{m}'$ denote the intersection of $\mathfrak{m}$ and $\mathfrak{g}'$.  Denote by $\mathcal{A}$\label{sym172} the affine quotient $\mathfrak{g}'/\!\!/\mathrm{G}$ which is identified with the affine line $\mathbb{L}$ via the negative determinant map
\begin{eqnarray}
-\mathrm{det}:\quad
\mathfrak{g}' &\longrightarrow& \mathbb{L}=\mathcal{A}.
\end{eqnarray}

\paragraph{Remark}\setcounter{equation}{0}
Let $\mathfrak{z}$ denote the center of $\mathfrak{g}$ consisting of the diagonal matrices on which $\mathrm{G}$ operates trivially. Since 
\begin{eqnarray}
\mathfrak{g}
&=&
\mathfrak{g}'\oplus \mathfrak{z}
\end{eqnarray}
as representations of $\mathrm{G}$, for computing the local basic functions and the local Harish-Chandra transforms it suffices to consider $\mathfrak{g}'$ instead of $\mathfrak{g}$.

\paragraph{Definition}\setcounter{equation}{0}
The discriminant function $D$ on $\mathcal{A}$ is equal to $4X$ where $X$ denotes the coordinate function on $\mathcal{A}$. Let $\eta$ be an element of $\mathcal{A}_\mathrm{reg}(\mathbb{Q}_p)$ which is identified with the subset of units $\mathbb{Q}_p^\times$ of $\mathbb{Q}_p$. Then $\eta$ is said to be
\begin{itemize}
\item \emph{split} if $\eta$ has a square root in $\mathbb{Q}_p$; 
\item \emph{unramified elliptic} if $\eta$ does not have a square root in $\mathbb{Q}_p$ but the $p$-adic valuation of $\eta$ is even;
\item \emph{ramified elliptic} if the $p$-adic valuation of $\eta$ is odd.
\end{itemize}
The image of $\mathfrak{m}'_\mathrm{reg.ss}(\mathbb{Q}_p)$ in $\mathcal{A}_\mathrm{reg}(\mathbb{Q}_p)$ is equal to the subset of split elements.

\paragraph{Definition}\setcounter{equation}{0}
The algebraic differential form\label{sym173}
\begin{eqnarray}
\mathrm{d}_\eta X
&=&
\frac{\mathrm{d}x\wedge\mathrm{d}y}{y}
=\frac{\mathrm{d}y\wedge\mathrm{d}z}{2x}
=\frac{\mathrm{d}z\wedge\mathrm{d}x}{z}
\end{eqnarray}
on the orbit 
\begin{eqnarray}
\mathfrak{g}'_\eta
&=&
\Bigg\{\left(\begin{array}{cc}x&y\\z&-x\end{array}\right)\in\mathfrak{g}':~
x^2+yz=\eta\Bigg\}
\end{eqnarray}
where $\eta$ lies in $\mathcal{A}_\mathrm{reg}$ is invariant under the action of $\mathrm{G}$.

\paragraph{Lemma}\setcounter{equation}{0}
\emph{Let $f$ be a Schwartz function on $\mathfrak{g}'(\mathbb{Q}_p)$, then}
\begin{eqnarray}
\label{local differential form lemma}
\forall \eta\in\mathcal{A}_\mathrm{reg}(\mathbb{Q}_p)
&&
I_\mathrm{G}^\mathrm{G}(\eta,f)=
\frac{1}{1+p^{-1}}\int_{\mathfrak{g}'_\eta(\mathbb{Q}_p)} f(X)~|\mathrm{d}_\eta X|.
\end{eqnarray}

\proof
There exists a multiplicative constant $\lambda$ such that
\begin{eqnarray}
&&
\forall f\in\mathcal{S}(\mathfrak{g}'(\mathbb{Q}_p))~
\forall \eta\in\mathcal{A}_\mathrm{reg}(\mathbb{Q}_p)
\\\nonumber\\\nonumber
&&\quad
I_\mathrm{G}^\mathrm{G}(\eta,f)=
\lambda\cdot\int_{\mathfrak{g}'_\eta(\mathbb{Q}_p)} f(X)~|\mathrm{d}_\eta X|,
\end{eqnarray}
hence the constant $\lambda$ is determined by
\begin{eqnarray}
\label{local differential form lemma proof first equality}
I_\mathrm{G}^\mathrm{G}(1,\mathbb{I}_{\mathfrak{g}'(\mathbb{Z}_p)})
&=&
\lambda\cdot\int_{\mathfrak{g}'_1(\mathbb{Q}_p)} \mathbb{I}_{\mathfrak{g}'(\mathbb{Z}_p)}(X)~|\mathrm{d}_1 X|.
\end{eqnarray}
By parabolic descent along the Borel subgroup $\mathrm{B}$, the left hand side of \eqref{local differential form lemma proof first equality} is equal to
\begin{eqnarray}
\label{local differential form lemma proof second equality}
I_\mathrm{M}^\mathrm{M}(1,(\mathbb{I}_{\mathfrak{g}'(\mathbb{Z}_p)})_\mathrm{B})
&=&
I_\mathrm{M}^\mathrm{M}(1,\mathbb{I}_{\mathfrak{m}'(\mathbb{Z}_p)})
\end{eqnarray}
which is equal to one. Hence by \eqref{local differential form lemma proof first equality} the reciprocal of $\lambda$ is equal to
\begin{eqnarray}
&&
\int_{\mathfrak{g}'_1(\mathbb{Q}_p)}
\mathbb{I}_{\mathfrak{g}'(\mathbb{Z}_p)}
\bigg(\left(\begin{array}{ll}x&~~y\\z&-x\end{array}\right)\bigg)
~\bigg|\frac{\mathrm{d}x\wedge\mathrm{d}y}{y}\bigg|
\\\nonumber\\\nonumber
&=& \int_{\mathbb{Z}_p}
\bigg(\int_{\big\{y\in\mathbb{Z}_p:~ z=\frac{1-x^2}{y}\in\mathbb{Z}_p\big\}}
\frac{\mathrm{d}y}{|y|}\bigg)~\mathrm{d}x
\\\nonumber\\\nonumber
&=& \int_{|1+x|<1}
\bigg(\int_{|1-x^2|\leq|y|\leq1}\frac{\mathrm{d}y}{|y|}\bigg)~
\mathrm{d}x+
\\\nonumber\\\nonumber
&&\quad
+\int_{|1-x|<1}
\bigg(\int_{|1-x^2|\leq|y|\leq1}\frac{\mathrm{d}y}{|y|}~\bigg)~
\mathrm{d}x+ 
\\\nonumber\\\nonumber
&&\qquad
+\int_{\substack{|1+x|=1\\|1-x|=1}}
\bigg(\int_{|1-x^2|\leq|y|\leq1}\frac{\mathrm{d}y}{|y|}\bigg)~
\mathrm{d}x
\\\nonumber\\\nonumber
&=& 
\sum_{i=1}^\infty
\bigg(\int_{|1+x|=p^{-i}}\mathrm{d}x\bigg)
\bigg(\int_{p^{-i}\leq|y|\leq1}\frac{\mathrm{d}y}{|y|}\bigg)+
\\\nonumber\\\nonumber
&&\quad
+\sum_{j=1}^\infty
\bigg(\int_{|1-x|=p^{-j}}\mathrm{d}x\bigg)
\bigg(\int_{p^{-j}\leq|y|\leq1}\frac{\mathrm{d}y}{|y|}\bigg)+ 
\\\nonumber\\\nonumber
&&\qquad
+\bigg(\int_{\substack{|1+x|=1\\|1-x|=1}}\mathrm{d}x\bigg)
\bigg(\int_{|y|=1}\frac{\mathrm{d}y}{|y|}\bigg)
\\\nonumber\\\nonumber
&=& 
\sum_{i=1}^\infty
\bigg(\frac{1}{p^i}-\frac{1}{p^{i+1}}\bigg)
\bigg(\frac{(i+1)(p-1)}{p}\bigg)+
\\\nonumber\\\nonumber
&&\quad
+\sum_{j=1}^\infty
\bigg(\frac{1}{p^j}-\frac{1}{p^{j+1}}\bigg)
\bigg(\frac{(j+1)(p-1)}{p}\bigg)+
\\\nonumber\\\nonumber
&&\qquad
+\bigg(\frac{p-2}{p}\bigg)\bigg(\frac{p-1}{p}\bigg)
\\\nonumber\\\nonumber
&=& \frac{p-1}{p}\bigg(2\cdot\Big(\frac{1}{p}+\frac{1}{p-1}\Big)+\frac{p-2}{p}\bigg)
\\\nonumber\\\nonumber
&=& 1+p^{-1}.
\end{eqnarray}
\qed

\paragraph{Lemma}\setcounter{equation}{0}
\emph{Let $\eta$ be an element of $\mathcal{A}_\mathrm{reg}(\mathbb{Q}_p)$, then}
\begin{eqnarray}
\phi_{0,v}(\eta)
&=&
\left\{\begin{array}{cl}
1
&\textrm{\emph{if $\eta$ is split}}\\
&~\textrm{ \emph{and contained in $\mathbb{Z}_p-\{0\}$}},\\\\
\displaystyle 1-\frac{2|\eta|^{1/2}}{1+p}
&\textrm{\emph{if $\eta$ is unramified elliptic}}\\
&~\textrm{ \emph{and contained in $\mathbb{Z}_p-\{0\}$}},\\\\
\displaystyle1-p^{-1/2}|\eta|^{1/2}
&\textrm{\emph{if $\eta$ is ramified elliptic}}\\
&~\textrm{ \emph{and contained in $\mathbb{Z}_p-\{0\}$}},\\\\
0&\textrm{\emph{if $\eta$ does not lie in $\mathbb{Z}_p-\{0\}$}}.
\end{array}\right.
\end{eqnarray}

\proof
If $\eta$ does not lie in $\mathbb{Z}_p-\{0\}$, then the orbit $\mathfrak{g}'_\eta(\mathbb{Q}_p)$ is disjoint from the support of the characteristic function $\mathbb{I}_{\mathfrak{g}'(\mathbb{Z}_p)}$, hence $\phi_{0,v}(\eta)$ vanishes.

If $\eta$ is split and contained in $\mathbb{Z}-\{0\}$, then by \eqref{local differential form lemma proof second equality} the value $\phi_{0,v}(\eta)$ is equal to one.

If $\eta$ is unramified elliptic and contained in $\mathbb{Z}_p-\{0\}$, then by \eqref{local differential form lemma}
\begin{eqnarray}
\phi_{0,v}(\eta)
&=&
\frac{1}{1+p^{-1}}
\int_{\mathfrak{g}'_\eta(\mathbb{Q}_p)}
\mathbb{I}_{\mathfrak{g}'(\mathbb{Z}_p)}(X)~
|\mathrm{d}_\eta X|
\\\nonumber\\\nonumber
&=&
\frac{1}{1+p^{-1}}
\int_{\mathfrak{g}'_\eta(\mathbb{Q}_p)}
\mathbb{I}_{\mathfrak{g}'(\mathbb{Z}_p)}
\bigg(\left(\begin{array}{ll}x&~~y\\z&-x\end{array}\right)\bigg)
~\bigg|\frac{\mathrm{d}x\wedge\mathrm{d}y}{y}\bigg|
\\\nonumber\\
\label{local basic function 0 lemma proof first equality}
&=& \frac{1}{1+p^{-1}}
\int_{|x|\leq1}
\bigg(\int_{|\eta-x^2|\leq|y|\leq1}
\frac{\mathrm{d}y}{|y|}\bigg)~\mathrm{d}x
\end{eqnarray}
Decompose the integral on the right hand side of \eqref{local basic function 0 lemma proof first equality} over two disjoint regions:
\begin{itemize}
\item If $|x^2|$ is strictly less than $|\eta|$, then $|\eta-x^2|$ is equal to $|\eta|$, hence
\begin{eqnarray}
&&
\int_{|x^2|<|\eta|}
\bigg(\int_{|\eta-x^2|\leq|y|\leq1}\frac{\mathrm{d}y}{|y|}\bigg)~
\mathrm{d}x
\\\nonumber\\\nonumber
&=& 
\bigg(\int_{|x^2|<|\eta|}\mathrm{d}x\bigg)
\bigg(\int_{|\eta|\leq|y|\leq1}\frac{\mathrm{d}y}{|y|}\bigg)
\\\nonumber\\\nonumber
&=&
\bigg(p^{-\big(\frac{\mathrm{val}(\eta)}{2}+1\big)}\bigg)
\bigg(\frac{(\mathrm{val}(\eta)+1)(p-1)}{p}\bigg)
\end{eqnarray}
where $\mathrm{val}(\eta)$\label{sym174} denotes the $p$-adic valuation of $\eta$.

\item If $|x^2|$ is greater than or equal to $|\eta|$, then $|\eta-x^2|$ is equal to $|x^2|$, hence
\begin{eqnarray}
\qquad
&&
\int_{|\eta|\leq|x^2|\leq1}
\bigg(\int_{|\eta-x^2|\leq|y|\leq1}\frac{\mathrm{d}y}{|y|}\bigg)~
\mathrm{d}x
\\\nonumber\\\nonumber
&=&
\int_{|\eta|\leq|x^2|\leq1}
\bigg(\int_{|x^2|\leq|y|\leq1}\frac{\mathrm{d}y}{|y|}\bigg)~
\mathrm{d}x 
\\\nonumber\\\nonumber
&=&
\int_{|\eta|\leq|x^2|\leq1}(2\mathrm{val}(x)+1)\bigg(\frac{p-1}{p}\bigg)~
\mathrm{d}x
\\\nonumber\\\nonumber
&=& 
\frac{p-1}{p}
\bigg(1-p^{-\big(\frac{\mathrm{val}(\eta)}{2}+1\big)}
+2\int_{|\eta|\leq|x^2|\leq1} \mathrm{val}(x)~\mathrm{d}x \bigg)
\\\nonumber\\\nonumber
&=& 
\frac{p-1}{p}
\Bigg(1-p^{-\big(\frac{\mathrm{val}(\eta)}{2}+1\big)}
+2\sum_{i=0}^\frac{\mathrm{val}(\eta)}{2}
i\bigg(\frac{1}{p^i}-\frac{1}{p^{i+1}}\bigg)\Bigg)
\\\nonumber\\\nonumber
&=& 
\frac{p-1}{p}
\Bigg(1-p^{-\big(\frac{\mathrm{val}(\eta)}{2}+1\big)}+
\\\nonumber\\\nonumber
&&\quad
+2\cdot\bigg(-\Big(\frac{\mathrm{val}( \eta)}{2}+1\Big)
p^{-\big(\frac{\mathrm{val}(\eta)}{2}+1\big)}
+\frac{1-p^{-\big(\frac{\mathrm{val}(\eta)}{2}+1\big)}}{p-1} \bigg)\Bigg).
\end{eqnarray}
\end{itemize}
Hence the right hand side of \eqref{local basic function 0 lemma proof first equality} is equal to
\begin{eqnarray}
&&
\frac{1}{1+p^{-1}}\cdot
\frac{p-1}{p}
\Bigg(\Big(\mathrm{val}(\eta)+1\Big)
p^{-\big(\frac{\mathrm{val}(\eta)}{2}+1\big)}+
\\\nonumber\\\nonumber
&&\quad
-\Big(\mathrm{val}( \eta)+2\Big)
p^{-\big(\frac{\mathrm{val}(\eta)}{2}+1\big)}
+\bigg(1-p^{-\big(\frac{\mathrm{val}(\eta)}{2}+1\big)}\bigg)
\bigg(1+\frac{2}{p-1}\bigg)\Bigg)
\\\nonumber\\\nonumber
&=&
\frac{p-1}{p+1}
\Bigg(-p^{-\big(\frac{\mathrm{val}(\eta)}{2}+1\big)}-
\bigg(\frac{p+1}{p-1}\bigg)p^{-\big(\frac{\mathrm{val}(\eta)}{2}+1\big)}+
\frac{p+1}{p-1}\Bigg)
\\\nonumber\\\nonumber
&=&
1-\frac{2}{p+1}p^{-\frac{\mathrm{val}(\eta)}{2}}.
\end{eqnarray}

If $\eta$ is ramified elliptic and contained in $\mathbb{Z}_p-\{0\}$, then by \eqref{local differential form lemma}
\begin{eqnarray}
\label{local basic function 0 lemma proof second equality}
\phi_{0,v}(\eta)
&=&
\frac{1}{1+p^{-1}}
\int_{|x|\leq1}
\bigg(\int_{|\eta-x^2|\leq|y|\leq1}
\frac{\mathrm{d}y}{|y|}\bigg)~\mathrm{d}x
\end{eqnarray}
Decompose the integral on the right hand side of \eqref{local basic function 0 lemma proof second equality} over two disjoint regions:
\begin{itemize}
\item If $|x^2|$ is less than $|\eta|$, then $|\eta-x^2|$ is equal to $|\eta|$, hence
\begin{eqnarray}
&&
\int_{|x^2|<|\eta|}
\bigg(\int_{|\eta-x^2|\leq|y|\leq1}\frac{\mathrm{d}y}{|y|}\bigg)~
\mathrm{d}x
\\\nonumber\\\nonumber
&=& 
\bigg(\int_{|x^2|<|\eta|}\mathrm{d}x\bigg)
\bigg(\int_{|\eta|\leq|y|\leq1}\frac{\mathrm{d}y}{|y|}\bigg)
\\\nonumber\\\nonumber
&=&
\bigg(p^{-\big(\frac{\mathrm{val}(\eta)+1}{2}\big)}\bigg)
\bigg(\frac{(\mathrm{val}(\eta)+1)(p-1)}{p}\bigg).
\end{eqnarray}

\item If $|x^2|$ is greater than $|\eta|$, then $|\eta-x^2|$ is equal to $|x^2|$, hence
\begin{eqnarray}
\qquad
&&
\int_{|\eta|\leq|x^2|\leq1}
\bigg(\int_{|\eta-x^2|\leq|y|\leq1}\frac{\mathrm{d}y}{|y|}\bigg)~
\mathrm{d}x
\\\nonumber\\\nonumber
&=&
\int_{|\eta|\leq|x^2|\leq1}
\bigg(\int_{|x^2|\leq|y|\leq1}\frac{\mathrm{d}y}{|y|}\bigg)~
\mathrm{d}x 
\\\nonumber\\\nonumber
&=&
\int_{|\eta|\leq|x^2|\leq1}(2\mathrm{val}(x)+1)\bigg(\frac{p-1}{p}\bigg)~
\mathrm{d}x
\\\nonumber\\\nonumber
&=& 
\frac{p-1}{p}
\bigg(1-p^{-\big(\frac{\mathrm{val}(\eta)+1}{2}\big)}
+2\int_{|\eta|\leq|x^2|\leq1} \mathrm{val}(x)~\mathrm{d}x \bigg)
\\\nonumber\\\nonumber
&=& 
\frac{p-1}{p}
\Bigg(1-p^{-\big(\frac{\mathrm{val}(\eta)+1}{2}\big)}
+2\sum_{i=0}^\frac{\mathrm{val}(\eta)-1}{2}
i\bigg(\frac{1}{p^i}-\frac{1}{p^{i+1}}\bigg)\Bigg)
\\\nonumber\\\nonumber
&=& 
\frac{p-1}{p}
\Bigg(1-p^{-\big(\frac{\mathrm{val}(\eta)+1}{2}\big)}+
\\\nonumber\\\nonumber
&&\quad
+2\cdot\bigg(-\Big(\frac{\mathrm{val}( \eta)+1}{2}\Big)
p^{-\big(\frac{\mathrm{val}(\eta)+1}{2}\big)}
+\frac{1-p^{-\big(\frac{\mathrm{val}(\eta)+1}{2}\big)}}{p-1} \bigg)\Bigg).
\end{eqnarray}
Hence the right hand side of \eqref{local basic function 0 lemma proof second equality} is equal to
\begin{eqnarray}
&&
\frac{1}{1+p^{-1}}\cdot
\frac{p-1}{p}
\Bigg(\Big(\mathrm{val}(\eta)+1\Big)
p^{-\big(\frac{\mathrm{val}(\eta)+1}{2}\big)}+
\\\nonumber\\\nonumber
&&\quad
-\Big(\mathrm{val}( \eta)+1\Big)
p^{-\big(\frac{\mathrm{val}(\eta)+1}{2}\big)}
+\bigg(1-p^{-\big(\frac{\mathrm{val}(\eta)+1}{2}\big)}\bigg)
\bigg(1+\frac{2}{p-1}\bigg)\Bigg)
\\\nonumber\\\nonumber
&=&
\frac{p-1}{p+1}
\bigg(1-p^{-\big(\frac{\mathrm{val}(\eta)+1}{2}\big)}\bigg)
\bigg(\frac{p+1}{p-1}\bigg)
\\\nonumber\\\nonumber
&=&
1-p^{-1/2}p^{-\frac{\mathrm{val}(\eta)}{2}}.
\end{eqnarray}
\end{itemize}
\qed

\paragraph{Definition}\setcounter{equation}{0}
The weighted orbital integral $J_\mathrm{M}^\mathrm{G}(\widetilde{\eta},~)$ where $\widetilde{\eta}$ is an element of $\mathfrak{g}'_\mathrm{reg.ss}(\mathbb{Q}_p)$ lifting $\eta$ is defined is $\eta$ is split. For such an $\eta$ choose a square root $\sqrt{\eta}$ of $\eta$ in $\mathbb{Q}_p$ and fix the lift $\widetilde{\eta}$\label{sym175} to be
\begin{eqnarray}
\widetilde{\eta}
&=&
\left(\begin{array}{ll}\sqrt{\eta} & \\ & -\sqrt{\eta}\end{array}\right).
\end{eqnarray} 

\paragraph{Remark}\setcounter{equation}{0}
With respect to the choice of $\mathrm{B}$ as the Borel subgroup and $\mathrm{G}(\mathbb{Z}_p)$ as the maximal compact subgroup, the weight factor $v(~)$ appearing in the definition of $J_\mathrm{M}^\mathrm{G}(\widetilde{\eta},~)$ is the function on $\mathrm{G}(\mathbb{Q}_p)$ such that $v(~)$ is invariant under left translation by $\mathrm{M}(\mathbb{Q}_p)$ and right translation by $\mathrm{G}(\mathbb{Z}_p)$ and
\begin{eqnarray}
&&
\forall u\in\mathbb{Q}_p
\\\nonumber\\\nonumber
&&
v\bigg(\left(\begin{array}{ll}1 & u\\ & 1\end{array}\right)\bigg)=
\left\{\begin{array}{cl} 0 & \textrm{if $u$ is contained in $\mathbb{Z}_p$},
\\\\ -\mathrm{val}(u) & \textrm{otherwise}.\end{array}\right.
\end{eqnarray}
See \S12.1 of \cite{Kott}.

\paragraph{Lemma}\setcounter{equation}{0}
\emph{Let $\eta$ be a split element of $\mathcal{A}_\mathrm{reg}(\mathbb{Q}_p)$, then}
\begin{eqnarray}
\label{local weighted orbital integral lemma}
J_\mathrm{M}^\mathrm{G}(\widetilde{\eta},\mathbb{I}_{\mathfrak{g}'(\mathbb{Z}_p)})
&=&
\left\{\begin{array}{cl}
\displaystyle \frac{\mathrm{val}(\eta)}{2}-\frac{1}{p-1}+\frac{|\eta|^{1/2}}{p-1}
& \textrm{if $\eta$ lies in $\mathbb{Z}_p-\{0\}$},\\\\
0 & \textrm{otherwise}.\end{array}\right.
\end{eqnarray}

\proof
By definition
\begin{eqnarray}
&&
J_\mathrm{M}^\mathrm{G}(\widetilde{\eta},\mathbb{I}_{\mathfrak{g}'(\mathbb{Z}_p)})
\\\nonumber\\\nonumber
&=&
|D(\eta)|^{1/2}
\int_{\mathrm{M}(\mathbb{Q}_p)\backslash \mathrm{G}(\mathbb{Q}_p)}
\mathbb{I}_{\mathfrak{g}'(\mathbb{Z}_p)}(\widetilde{\eta}\cdot\mathrm{ad}(g))~
v(g)\mathrm{d}g 
\\\nonumber\\
\label{local weighted orbital integral lemma proof first equality}
&=&
|\eta|^{1/2}
\int_{\mathbb{Q}_p}
\int_{\mathrm{G}(\mathbb{Z}_p)} 
\mathbb{I}_{\mathfrak{g}'(\mathbb{Z}_p)}
\Bigg(
\left(\begin{array}{ll}\sqrt{\eta}&\\&-\sqrt{\eta}\end{array}\right)
\cdot\mathrm{ad}\bigg(\left(\begin{array}{ll}1 & u\\ & 1\end{array}\right)\bigg)
\mathrm{ad}(k)
\Bigg)\times
\\\nonumber\\\nonumber
&&\quad\times
v\bigg(\left(\begin{array}{ll}1 & u\\ & 1\end{array}\right)\bigg)
\mathrm{d}k\mathrm{d}u
\\\nonumber\\\nonumber
&=&
|\eta|^{1/2}
\int_{\mathbb{Q}_p}
\mathbb{I}_{\mathfrak{g}'(\mathbb{Z}_p)}
\bigg(
\left(\begin{array}{ll}\sqrt{\eta}&2u\sqrt{\eta}\\&-\sqrt{\eta}\end{array}\right)
\bigg)~
v\bigg(\left(\begin{array}{ll}1 & u\\ & 1\end{array}\right)\bigg)\mathrm{d}u
\\\nonumber\\\nonumber
&=&
|\eta|^{1/2}
\int_{1<|u|\leq|2\sqrt{\eta}|^{-1}}
-\mathrm{val}(u)~
\mathrm{d}u
\\\nonumber\\\nonumber
&=&
|\eta|^{1/2}
\sum_{i=1}^{\frac{\mathrm{val}(\eta)}{2}}
i\bigg(\frac{1}{p^{-i}}-\frac{1}{p^{-i+1}}\bigg)
\\\nonumber\\\nonumber
&=&
|\eta|^{1/2}
\bigg(\frac{\mathrm{val}(\eta)}{2}p^{\frac{\mathrm{val}(\eta)}{2}}
-\frac{p^{\frac{\mathrm{val}(\eta)}{2}}-1}{p-1}\bigg)
\end{eqnarray}
where the equality \eqref{local weighted orbital integral lemma proof first equality} follows from the Iwasawa decomposition.
\qed

\paragraph{Lemma}\setcounter{equation}{0}
\emph{Let $\eta$ be a split element of $\mathcal{A}_\mathrm{reg}(\mathbb{Q}_p)$, then}
\begin{eqnarray}
I_\mathrm{M}^\mathrm{G}(\widetilde{\eta},\mathbb{I}_{\mathfrak{g}'(\mathbb{Z}_p)})
&=&
\left\{\begin{array}{cl}
\displaystyle
\frac{\mathrm{val}(\eta)}{2}+\frac{|\eta|^{1/2}}{p-1}-\frac{2p+1}{p^2-1}
& \textrm{if $\eta$ lies in $\mathbb{Z}_p-\{0\}$}, \\\\
\displaystyle
\frac{|\eta|^{-1}}{1-p^{-2}}-\frac{|\eta|^{-{1/2}}}{1-p^{-1}}
& \textrm{otherwise}. \end{array}\right.
\end{eqnarray}

\proof
By definition
\begin{eqnarray}
\label{local invariant weighted orbital integral lemma proof first equality}
&&
I_\mathrm{M}^\mathrm{G}(\widetilde{\eta},\mathbb{I}_{\mathfrak{g}'(\mathbb{Z}_p)})
\\\nonumber\\\nonumber
&=&
J_\mathrm{M}^\mathrm{G}(\widetilde{\eta},\mathbb{I}_{\mathfrak{g}'(\mathbb{Z}_p)})-
\int_{\mathfrak{m}'(\mathbb{Q}_p)}
J_\mathrm{M}^\mathrm{G}(X,\mathbb{I}_{\mathfrak{g}'(\mathbb{Z}_p)})
\psi\big(\mathrm{Tr}(\widetilde{\eta}^\mathrm{T}X)\big)~
\mathrm{d}X
\\\nonumber\\\nonumber
&=&
J_\mathrm{M}^\mathrm{G}(\widetilde{\eta},\mathbb{I}_{\mathfrak{g}'(\mathbb{Z}_p)})-
\int_{\mathbb{Q}_p}
J_\mathrm{M}^\mathrm{G}\bigg(\left(\begin{array}{ll}x\\&-x\end{array}\right),
\mathbb{I}_{\mathfrak{g}'(\mathbb{Z}_p)}\bigg)\psi(2\sqrt{\eta}x)~
\mathrm{d}x
\end{eqnarray}
where $\psi$ denotes the additive character of $\mathbb{Q}_p$ involved in the definition of the Fourier transform on a $p$-adic vector space, the superscript $\mathrm{T}$ denotes the transpose of a $2\times2$ matrix and $\mathrm{Tr}$ denotes the trace of a $2\times2$ matrix.

The Fourier transform of $\mathrm{val}(~)\mathbb{I}_{\mathbb{Z}_p}(~)$ on $\mathbb{Q}_p$ is the function such that
\begin{eqnarray}
\label{local invariant weighted orbital integral lemma proof second equality}
\forall x\in\mathbb{Q}_p
&&
\Big(\mathrm{val}(~)\mathbb{I}_{\mathbb{Z}_p}(~)\Big)\hat{~}(x)
\\\nonumber\\\nonumber
&=&
\bigg(\sum_{i=1}^\infty \mathbb{I}_{p^i\mathbb{Z}_p}(~)\bigg)\hat{~}(x)
\\\nonumber\\\nonumber
&=& 
\sum_{i=1}^\infty \frac{1}{p^i}\mathbb{I}_{p^{-i}\mathbb{Z}_p}(x)
\\\nonumber\\\nonumber
&=&
\sum_{i=\mathrm{max}\{1,-\mathrm{val}(x)\}}^\infty \frac{1}{p^i}
\\\nonumber\\\nonumber
&=&
\frac{p^{-\mathrm{max}\{1,-\mathrm{val}(x)\}}}{1-p^{-1}}
\\\nonumber\\\nonumber
&=&
\mathrm{min}\bigg\{\frac{p^{-1}}{1-p^{-1}},\frac{|x|^{-1}}{1-p^{-1}}\bigg\}.
\end{eqnarray}
The Fourier transform of $|~|\mathbb{I}_{\mathbb{Z}_p}(~)$ on $\mathbb{Q}_p$ is the function such that
\begin{eqnarray}
\label{local invariant weighted orbital integral lemma proof third equality}
\forall x\in\mathbb{Q}_p
&&
\Big(|~|\mathbb{I}_{\mathbb{Z}_p}(~)\Big)\hat{~}(x)
\\\nonumber\\\nonumber
&=&
\bigg(\mathbb{I}_{\mathbb{Z}_p}(~)+
\sum_{i=1}^\infty
\bigg(\frac{1}{p^i}-\frac{1}{p^{i-1}}\bigg)
\mathbb{I}_{p^i\mathbb{Z}_p}(~)\bigg)\hat{~}(x)
\\\nonumber\\\nonumber
&=&
\mathbb{I}_{\mathbb{Z}_p}(x)+
\sum_{i=1}^\infty
\bigg(\frac{1}{p^{2i}}-\frac{1}{p^{2i-1}}\bigg)
\mathbb{I}_{p^{-i}\mathbb{Z}_p}(x)
\\\nonumber\\\nonumber
&=&
\mathbb{I}_{\mathbb{Z}_p}(x)+
\sum_{i=\mathrm{max}\{1,-\mathrm{val}(x)\}}^\infty
\bigg(\frac{1}{(-p)^{2i}}+\frac{1}{(-p)^{2i-1}}\bigg)
\\\nonumber\\\nonumber
&=&
\mathbb{I}_{\mathbb{Z}_p}(x)+
\frac{(-p)^{-\big(2\mathrm{max}\{1,-\mathrm{val}(x)\}-1\big)}}{1+p^{-1}}
\\\nonumber\\\nonumber
&=&
\mathbb{I}_{\mathbb{Z}_p}(x)-
\mathrm{min}\bigg\{\frac{p^{-1}}{1+p^{-1}},\frac{p|x|^{-2}}{1+p^{-1}}\bigg\}.
\end{eqnarray}
Hence by \eqref{local weighted orbital integral lemma} and \eqref{local invariant weighted orbital integral lemma proof first equality} the invariant weighted orbital integral $I_\mathrm{M}^\mathrm{G}(\widetilde{\eta},\mathbb{I}_{\mathfrak{g}'(\mathbb{Z}_p)})$ is equal to
\begin{eqnarray}
&&
\Big(\frac{\mathrm{val}(\eta)}{2}-\frac{1}{p-1}+\frac{|\eta|^{1/2}}{p-1}\Big)
\mathbb{I}_{\mathbb{Z}_p-\{0\}}(\eta)
\\\nonumber\\\nonumber
&&\quad
-\bigg(\Big(\mathrm{val}(~)-\frac{1}{p-1}+\frac{|~|}{p-1}\Big)
\mathbb{I}_{\mathbb{Z}_p}(~)\bigg)\hat{~}(2\sqrt{\eta})
\\\nonumber\\\nonumber
&=&
\Big(\frac{\mathrm{val}(\eta)}{2}-\frac{1}{p-1}+\frac{|\eta|^{1/2}}{p-1}\Big)
\mathbb{I}_{\mathbb{Z}_p-\{0\}}(\eta)+
\frac{1}{p-1}\mathbb{I}_{\mathbb{Z}_p}(\eta)
\\\nonumber\\\nonumber
&&\quad-
\Big(\mathrm{val}(~)\mathbb{I}_{\mathbb{Z}_p}(~)\Big)\hat{~}(2\sqrt{\eta})-
\frac{1}{p-1}\Big(|~|\mathbb{I}_{\mathbb{Z}_p}(~)\Big)\hat{~}(2\sqrt{\eta})
\\\nonumber\\
\label{local invariant weighted orbital integral lemma proof fourth equality}
&=&
\left\{\begin{array}{ll}
\displaystyle
\frac{\mathrm{val}(\eta)}{2}+\frac{|\eta|^{1/2}}{p-1}
\\\\\displaystyle\quad
-\frac{p^{-1}}{1-p^{-1}}-\frac{1}{p-1}\bigg(1-\frac{p^{-1}}{1+p^{-1}}\bigg)
&\textrm{if $\eta$ lies in $\mathbb{Z}_p-\{0\}$},\\\\
\displaystyle
-\frac{|2\sqrt{\eta}|^{-1}}{1-p^{-1}}-\frac{1}{p-1}\bigg(-\frac{p|2\sqrt{\eta}|^{-2}}{1+p^{-1}}\bigg)
&\textrm{otherwise},
\end{array}\right.
\\\nonumber\\\nonumber
&=&
\left\{\begin{array}{ll}
\displaystyle
\frac{\mathrm{val}(\eta)}{2}+\frac{|\eta|^{1/2}}{p-1}
-\frac{1}{p-1}\bigg(\frac{2+p^{-1}}{1+p^{-1}}\bigg)
&\textrm{if $\eta$ lies in $\mathbb{Z}_p-\{0\}$},\\\\
\displaystyle
-\frac{|\eta|^{-1/2}}{1-p^{-1}}+\frac{|\eta|^{-1}}{(1-p^{-1})(1+p^{-1})}
&\textrm{otherwise},
\end{array}\right.
\end{eqnarray}
where the equality \eqref{local invariant weighted orbital integral lemma proof fourth equality} follows from \eqref{local invariant weighted orbital integral lemma proof second equality} and \eqref{local invariant weighted orbital integral lemma proof third equality}.
\qed

\paragraph{Remark}\setcounter{equation}{0}
If $\eta$ is a split element of $\mathcal{A}_\mathrm{reg}(\mathbb{Q}_p)$, then $\phi_{1,v}(\eta)$ has three components indexed by the parabolic subgroups $\mathrm{B}$, $\overline{\mathrm{B}}$ and $\mathrm{G}$, and 
\begin{eqnarray}
\phi_{1,v}(\eta)_\mathrm{B}=\phi_{1,v}(\eta)_{\overline{\mathrm{B}}}=1
&&
\phi_{1,v}(\eta)_\mathrm{G}=-I_\mathrm{M}^\mathrm{G}(\widetilde{\eta},\mathbb{I}_{\mathbb{Z}_p}).
\end{eqnarray}
If $\eta$ is an elliptic element of $\mathcal{A}_\mathrm{reg}(\mathbb{Q}_p)$, then $\phi_{1,v}(\eta)$ has one component corresponding to the parabolic subgroup $\mathrm{G}$ which is equal to $\phi_{0,v}(\eta)$.

\paragraph{Remark: Towards an invariant trace formula}\setcounter{equation}{0}
The simple invariant trace formula \eqref{simple invariant trace formula} holds only for test functions which are cuspidal at two distinct places of $\mathbb{Q}$. An invariant trace formula, once extended to full generality for all Schwartz functions on $\mathfrak{g}(\mathbb{A})$, will likely involve more intricate distributions similar to the invariant weighted orbital integrals $I_\mathrm{M}^\mathrm{G}(X,~)$ which are absent from \eqref{simple invariant trace formula} due to the cuspidality condition. One such identity could be established when $\mathfrak{g}$ is $\mathrm{gl}(2)$ following the arguments of Langlands in \cite{BaseChange}.

Fix a sufficiently large finite set $S$ of places of $\mathbb{Q}$, let $f_S$ be a Schwartz function on $\mathfrak{g}(\mathbb{Q}_S)$. It follows from the non-invariant trace formula \eqref{refined trace formula} that
\begin{eqnarray}
\label{invariant trace formula for gl(2) first equality}
&&\sum_{\mathfrak{o}\in\mathfrak{g}/\sim}~\sum_{X\in\mathfrak{o}_{\mathrm{G},S}}a^\mathrm{G}(S,X)I_\mathrm{G}^\mathrm{G}(X,f_S)+
\\\nonumber\\\nonumber&&\quad
+\frac{1}{2}\sum_{\mathfrak{o}\in\mathfrak{g}/\sim}~\sum_{X\in(\mathfrak{m}\cap\mathfrak{o})_{\mathrm{M},S}}a^\mathrm{M}(S,X)J_\mathrm{M}^\mathrm{G}(X,f_S)
\\\nonumber\\\nonumber
&=&\sum_{\mathfrak{o}\in\mathfrak{g}/\sim}~\sum_{X\in\mathfrak{o}_{\mathrm{G},S}}a^\mathrm{G}(S,X)I_\mathrm{G}^\mathrm{G}(X,f\hat{~}_S)+
\\\nonumber\\\nonumber&&\quad
+\frac{1}{2}\sum_{\mathfrak{o}\in\mathfrak{g}/\sim}~\sum_{X\in(\mathfrak{m}\cap\mathfrak{o})_{\mathrm{M},S}}a^\mathrm{M}(S,X)J_\mathrm{M}^\mathrm{G}(X,f\hat{~}_S).
\end{eqnarray}

The weighted orbital integral $J_\mathrm{M}^\mathrm{G}(X,f_S)$, considered as a function in the variable $X$ on $\mathfrak{m}(\mathbb{Q}_S)$, has at worst logarithmic singularity along the discriminant locus. There is a natural way to subtract from $J_\mathrm{M}^\mathrm{G}(X,f_S)$ an invariant function in $X$ such that the remainder extends to a continuous function on $\mathfrak{m}(\mathbb{Q}_S)$ which is denoted by $J_\mathrm{M}^\mathrm{G,reg}(X,f_S)$\label{sym176}. In general $J_\mathrm{M}^\mathrm{G,reg}(X,f_S)$ is not smooth, as could be seen for instance from \eqref{local weighted orbital integral lemma}. However $J_\mathrm{M}^\mathrm{G,reg}(X,f_S)$, when restricted to a function on the torus $\mathrm{M}(\mathbb{Q}_S)$, has integrable Mellin transform, which follows from the estimates in \S9 of \cite{BaseChange}. By a suitable adelic extension of the arguments of Ferrar in \cite{Ferrar}, it could be shown that the Poisson summation formula applies to the function $J_\mathrm{M}^\mathrm{G,reg}(X,f_S)$ and implies that
\begin{eqnarray}
\label{invariant trace formula for gl(2) second equality}
\sum_{X\in\mathfrak{m}(\mathbb{Q})}J_\mathrm{M}^\mathrm{G,reg}(X,f_S)
&=&
\sum_{X\in\mathfrak{m}(\mathbb{Q})}J_\mathrm{M}^\mathrm{G,reg}(X,f_S)\hat{~}
\end{eqnarray}
where the Fourier transform is taken over $\mathfrak{m}(\mathbb{Q}_S)$.

By definition the coefficients $a^\mathrm{M}(S,X)$ are independent of the point $X$, hence a suitable linear combination of \eqref{invariant trace formula for gl(2) first equality} and \eqref{invariant trace formula for gl(2) second equality} gives the identity
\begin{eqnarray}
\label{invariant trace formula for gl(2) third equality}
&&\sum_{\mathfrak{o}\in\mathfrak{g}/\sim}~\sum_{X\in\mathfrak{o}_{\mathrm{G},S}}a^\mathrm{G}(S,X)I_\mathrm{G}^\mathrm{G}(X,f_S)+
\\\nonumber\\\nonumber&&\quad
+\frac{1}{2}\sum_{\mathfrak{o}\in\mathfrak{g}/\sim}~\sum_{X\in(\mathfrak{m}\cap\mathfrak{o})_{\mathrm{M},S}}a^\mathrm{M}(S,X)\bigg(J_\mathrm{M}^\mathrm{G}(X,f_S)-J_\mathrm{M}^\mathrm{G,reg}(X,f\check{~}_S)\hat{~}\bigg)
\\\nonumber\\\nonumber
&=&\sum_{\mathfrak{o}\in\mathfrak{g}/\sim}~\sum_{X\in\mathfrak{o}_{\mathrm{G},S}}a^\mathrm{G}(S,X)I_\mathrm{G}^\mathrm{G}(X,f\hat{~}_S)+
\\\nonumber\\\nonumber&&\quad
+\frac{1}{2}\sum_{\mathfrak{o}\in\mathfrak{g}/\sim}~\sum_{X\in(\mathfrak{m}\cap\mathfrak{o})_{\mathrm{M},S}}a^\mathrm{M}(S,X)\bigg(J_\mathrm{M}^\mathrm{G}(X,f\hat{~}_S)-J_\mathrm{M}^\mathrm{G,reg}(X,f_S)\hat{~}\bigg).
\end{eqnarray}
Since
\begin{eqnarray}
I_\mathrm{G}^\mathrm{G}(X,f_S)
&\textrm{and}&
J_\mathrm{M}^\mathrm{G}(X,f_S)-J_\mathrm{M}^\mathrm{G,reg}(X,f\check{~}_S)\hat{~}
\end{eqnarray}
are invariant under the adjoint action, the identity \eqref{invariant trace formula for gl(2) third equality} is an invariant trace formula for $\mathrm{gl}(2)$ which is valid for all Schwartz functions on $\mathfrak{g}(\mathbb{Q}_S)$.

In order to extend this construction to a general reductive Lie algebra $\mathfrak{g}$, the singularities of higher rank weighted orbital integrals need to be separated, and the general invariant trace formula for each proper Levi subalgebra $\mathfrak{m}$ needs to be extended to test functions beyond the Schwartz space $\mathcal{S}(\mathfrak{m}(\mathbb{Q}_S))$.

\paragraph{Remark}\setcounter{equation}{0}
Arthur's work \cite{Arthur_Inv1} and \cite{Arthur_Inv2} on the invariant trace formula for reductive groups $\mathrm{G}$ follows an approach which is dual to the approach taken in the previous remark. Instead of working with the geometric side of the trace formula, Arthur devised the cancellation of singularities on the spectral side of the trace formula, utilizing the underlying complex structure and in particular the powerful techniques of residue calculus. However such an alternative is unavailable in the Lie algebra case since both sides of the trace formula are geometric in nature.

\paragraph{Remark}\setcounter{equation}{0}
Distributions similar to $J_\mathrm{M}^\mathrm{G,reg}(X,~)$ has been constructed by Chaudouard in \cite{Chaudb}. He has shown that his version of the weighted orbital integral $J_\mathrm{M}^{\mathrm{G},b}(X,~)$\label{sym177} remains bounded as $X$ traverses $\mathfrak{m}(\mathbb{Q}_S)$ while also enjoying many other nice properties expected of the usual weighted orbital integrals.

\paragraph{Remark}\setcounter{equation}{0}
There has been some important recent work on extending the Arthur-Selberg trace formula beyond the usual space of Schwartz functions. See \cite{FinisLapid} and \cite{Matz}.

\section*{List of symbols}
\addcontentsline{toc}{section}{List of symbols}
\begin{longtable}[l]{p{95pt} p{300pt} p{50pt}}
$\mathrm{G}, \mathfrak{g}$&
a connected reductive group defined over $\mathbb{Q}$ and its Lie algebra&
page \pageref{sym1}\\
$v, S$&
a place of $\mathbb{Q}$ and a finite set of places of $\mathbb{Q}$&
page \pageref{sym2}\\
$\mathbb{Q}_v, \mathbb{Q}_S, \mathbb{A}$&
the completions of $\mathbb{Q}$ at $v$ and $S$ and the ring of adeles of $\mathbb{Q}$&
page \pageref{sym3}\\
$|~|_v, |~|_S, |~|_\mathbb{A}$&
the norms on $\mathbb{Q}_v$, $\mathbb{Q}_S$ and $\mathbb{A}$&
page \pageref{sym4}\\
$\mathrm{G}(\mathbb{Q}_v), \mathrm{G}(\mathbb{Q}_S),$ $\mathrm{G}(\mathbb{A}), \mathrm{G}(\mathbb{Q})$&
the topological groups of $\mathbb{Q}_v$, $\mathbb{Q}_S$, $\mathbb{A}$ and $\mathbb{Q}$-valued points of $\mathrm{G}$&
page \pageref{sym5}\\
$\mathrm{P}_0, \mathrm{M}_0, \mathrm{A}_0$&
a minimal parabolic subgroup, a minimal Levi subgroup and a maximal split torus&
page \pageref{sym6}\\
$\mathrm{P},\mathrm{M}_\mathrm{P},\mathrm{N}_\mathrm{P},\mathrm{A}_\mathrm{P}$&
a parabolic subgroup, its Levi component, its unipotent radical and its split component&
page \pageref{sym7}\\
$\mathrm{P}_i, \mathrm{M}_i, \mathrm{N}_i, \mathrm{A}_i$&
a parabolic subgroup, its Levi component, its unipotent radical and its split component&
page \pageref{sym8}\\
$\overline{\mathrm{P}}, \overline{\mathrm{N}}_\mathrm{P}$&
the parabolic subgroup opposite to $\mathrm{P}$ and its unipotent radical&
page \pageref{sym9}\\
$\mathcal{F}, \mathcal{P}, \mathcal{L}$&
finite sets of parabolic subgroups and Levi subgroups&
page \pageref{sym10}\\
$X(\mathrm{M}_\mathrm{P}), X(\mathrm{M}_\mathrm{P})^*$&
the group of rational characters of $\mathrm{M}_\mathrm{P}$ and its dual&
page \pageref{sym11}\\
$\mathfrak{a}_\mathrm{P}, \mathfrak{a}_\mathrm{P}^*$&
a Euclidean space and its dual space&
page \pageref{sym12}\\
$\Phi_\mathrm{P}, \Delta_\mathrm{P}, \Phi_\mathrm{P}^\vee, \Delta_\mathrm{P}^\vee$&
the sets of roots, simple roots, coroots and simple coroots of $\mathrm{P}$&
page \pageref{sym13}\\
$\Delta_1^2, \Delta_1^{2,\vee}, \hat{\Delta}_1^2, \hat{\Delta}_1^{2,\vee}$&
sets of simple roots and coroots, coweights and weights&
page \pageref{sym14}\\
$\mathfrak{a}_1^2, \mathfrak{a}_1^{2,*}$&
subquotients of $\mathfrak{a}$ and $\mathfrak{a}^*$&
page \pageref{sym15}\\
$\mathrm{W}_0^\mathrm{G}, \mathrm{W}(\mathfrak{a}_1,\mathfrak{a}_2)$&
the Weyl group of $\mathrm{G}$ and the Weyl set from $\mathfrak{a}_1$ to $\mathfrak{a}_2$&
page \pageref{sym16}\\
$\tau_1^2, \hat{\tau}_1^2$&
characteristic functions of positive cones in $\mathfrak{a}$&
page \pageref{sym17}\\
$\mathrm{G}(\mathbb{A})^1$&
the subgroup of $\mathrm{G}(\mathbb{A})$ of elements of norm one&
page \pageref{sym18}\\
$\mathrm{K}, \mathrm{K}_v$&
global and local admissible maximal compact subgroups&
page \pageref{sym19}\\
$H_\mathrm{P}$&
a height function on $\mathrm{G}(\mathbb{A})$&
page \pageref{sym20}\\
$\rho_\mathrm{P}$&
the Weyl vector in $\mathfrak{a}^*$&
page \pageref{sym21}\\
$T',T$&
a fixed negative point and a truncation parameter in $\mathfrak{a}$&
page \pageref{sym22}\\
$\omega$&
a compact subset of $\mathrm{N}_0(\mathbb{A})\mathrm{M}_0(\mathbb{A})^1$&
page \pageref{sym23}\\
$\mathfrak{S}(T',\omega), \mathfrak{S}^T(T',\omega)$&
a Siegel domain and a truncated Siegel domain&
page \pageref{sym24}\\
$F^\mathrm{G}(x,T), F^\mathrm{P}(x,T)$&
characteristic functions of a truncated Siegel domain&
page \pageref{sym25}\\
$\mathcal{S},C_c^\infty$&
spaces of Schwartz functions and compactly supported smooth functions&
page \pageref{sym26}\\
$\mathbb{I}_{\mathfrak{g}(\mathbb{Z}_p)}$&
the characteristic function of the points in $\mathfrak{g}(\mathbb{Q}_p)$ with $p$-adically integral coordinates&
page \pageref{sym27}\\
$\langle~,~\rangle, \psi, \psi_v$&
a nondegenerate invariant bilinear form on $\mathfrak{g}$,  a additive character on $\mathbb{A}$ and its $\mathbb{Q}_v$-component&
page \pageref{sym28}\\
$\wedge,\vee$&
the Fourier transform and the inverse Fourier transform&
page \pageref{sym29}\\
$\sim, \mathfrak{o}$&
an equivalence relation on $\mathfrak{g}(\mathbb{Q})$ and an equivalence class&
page \pageref{sym30}\\
$X_\mathrm{ss}, X_\mathrm{nil}$&
the semisimple and nilpotent components of $X$ under the Jordan decomposition&
page \pageref{sym31}\\
$D, D^\mathrm{M}$&
the discriminant functions on $\mathfrak{g}$ and $\mathfrak{m}$&
page \pageref{sym32}\\
$\mathfrak{g}_\mathrm{reg.ss}$&
the locus of regular semisimple elements of $\mathfrak{g}$&
page \pageref{sym33}\\
$\mathrm{G}_X, \mathrm{G}_X^0, \pi_0(\mathrm{G}_X)$&
the centralizer of $X$ in $\mathrm{G}$, its identity component and its group of connected components&
page \pageref{sym34}\\
$\mathrm{G}_S, \mathrm{G}_v$&
the base changes of $\mathrm{G}$ to $\mathbb{Q}_S$ and $\mathbb{Q}_v$&
page \pageref{sym35}\\
$\mathrm{P}_{S,0}, \mathrm{M}_{S,0}, \mathrm{A}_{S,0},$ $\mathrm{P}_{v,0}, \mathrm{M}_{v,0}, \mathrm{A}_{v,0}$&
a pair of minimal parabolic subgroups, minimal Levi subgroups and maximal split tori of $\mathrm{G}_S$ and $\mathrm{G}_v$&
page \pageref{sym36}\\
$T_S, T_v, \mathfrak{t}_S, \mathfrak{t}_v$&
maximal tori of $\mathrm{G}_S$ and $\mathrm{G}_v$ and the associated Cartan subalgebras&
page \pageref{sym37}\\
$\mathcal{T}_\mathrm{ell}(\mathrm{G}_S)$&
the set of conjugacy classes of elliptic maximal tori of $\mathrm{G}_S$&
page \pageref{sym38}\\
$\mathrm{W}(\mathrm{G}_S,\mathrm{T}_S)$&
the Weyl group of the pair $(\mathrm{G}_S(\mathbb{Q}_S),\mathrm{T}_S(\mathbb{Q}_S))$&
page \pageref{sym39}\\
$I_\mathrm{G}^\mathrm{G}(X,~)$&
the standard normalized orbital integral&
page \pageref{sym40}\\
$K(x,f), K_\mathfrak{o}(x,f)$&
the kernel functions of Chaudouard&
page \pageref{sym41}\\
$I_\mathfrak{o}(f)$&
a variant of the standard orbital integral in the anisotropic case&
page \pageref{sym42}\\
$K_{\mathrm{P},\mathfrak{o}}(x,f), k_\mathfrak{o}^T(~,f)$&
a variant of the kernel function and the truncated kernel function of Chaudouard&
page \pageref{sym43}\\
$J_\mathfrak{o}^T, J^T$&
distributions on $\mathfrak{g}(\mathbb{A})$&
page \pageref{sym44}\\
$\sigma_1^4(~)$&
a combinatorial function on $\mathfrak{a}$&
page \pageref{sym45}\\
$\mathfrak{m}_\mathrm{P}^\mathrm{Q},\mathfrak{m}_\mathrm{P}^\mathrm{Q}(\mathbb{Q})', \mathfrak{n}_2^3$&
the quotient of $\mathfrak{m}_\mathrm{Q}$ by $\mathfrak{m}_\mathrm{P}$, a subset of $\mathfrak{m}_\mathrm{P}^\mathrm{Q}(\mathbb{Q})$ and the quotient of $\mathfrak{n}_2$ by $\mathfrak{n}_3$&
page \pageref{sym46}\\
$\Gamma'_{\mathrm{P}_2}(~,T)$&
the geometric gamma$'$ function&
page \pageref{sym47}\\
$J_\mathfrak{o}^\mathrm{M}$&
a sum of distributions on $\mathfrak{m}(\mathbb{A})$&
page \pageref{sym48}\\
$f_\mathrm{P}$&
the parabolic descent of $f$ along the parabolic subgroup $\mathrm{P}$&
page \pageref{sym49}\\
$\|~\|$&
the Euclidean norm on $\mathfrak{a}$&
page \pageref{sym50}\\
$T_0$&
a distinguished point in $\mathfrak{a}$&
page \pageref{sym51}\\
$w_s$&
a representative of the element $s$ of the Weyl group in $\mathrm{G}(\mathbb{Q})$&
page \pageref{sym52}\\
$J, J_\mathfrak{o}$&
distributions on $\mathfrak{g}(\mathbb{A})$&
page \pageref{sym53}\\
$v_\mathrm{M}(x)$&
the weight factor&
page \pageref{sym54}\\
$f_{\mathrm{P},x}, v'_\mathrm{P}(X)$&
a noninvariant version of the parabolic descent of $f$ along $\mathrm{P}$ and a variant of the weight factor $v_\mathrm{M}(x)$&
page \pageref{sym55}\\
$\theta_\mathrm{P}$&
a polynomial on $i\mathfrak{a}_\mathrm{M}^*$&
page \pageref{sym56}\\
$(c_\mathrm{P}), c_\mathrm{M}$&
a $(\mathrm{G},\mathrm{M})$-family and its associated constant&
page \pageref{sym57}\\
$(c_\mathrm{P}^\mathrm{Q}),(c_\mathrm{P})$&
the $(\mathrm{L},\mathrm{M})$-family and the $(\mathrm{G},\mathrm{L})$-family associated to the $(\mathrm{G},\mathrm{M})$-family $(c_\mathrm{P})$&
page \pageref{sym58}\\
$c'_\mathrm{Q}$&
the function or constant associated to the $(\mathrm{G},\mathrm{M})$-family $(c_\mathrm{P})$&
page \pageref{sym59}\\
$\mathcal{Y}_\mathrm{M}, Y_\mathrm{P}$&
a $(\mathrm{G},\mathrm{M})$-orthogonal set and an element of $\mathcal{Y}_\mathrm{M}$&
page \pageref{sym60}\\
$\big(v_\mathrm{P}(\mathcal{Y}_\mathrm{M})\big)$&
the $(\mathrm{G},\mathrm{M})$-family associated to the $(\mathrm{G},\mathrm{M})$-orthogonal set $\mathcal{Y}_\mathrm{M}$&
page \pageref{sym61}\\
$\mathcal{Y}_\mathrm{M}(x), v_\mathrm{M}(x)$&
a positive $(\mathrm{G},\mathrm{M})$-orthogonal and its associated constant which is the weight factor&
page \pageref{sym62}\\
$J_\mathrm{M}^\mathrm{G}(X,~)$&
the weighted orbital integral&
page \pageref{sym63}\\
$r_\mathrm{P}^\mathrm{L}(x,a)$&
an auxiliary $(\mathrm{L},\mathrm{M})$-family&
page \pageref{sym64}\\
$f_{\mathrm{P},x}, v'_\mathrm{Q}(y)$&
an $S$-local noninvariant parabolic descent of $f$ along $\mathrm{P}$ and its associated weight factor&
page \pageref{sym65}\\
$\mathfrak{g}_\mathrm{nil}, J_\mathrm{nil}, J_\mathrm{nil}^T$&
the nilpotent locus and the distributions supported on $\mathfrak{g}_\mathrm{nil}$&
page \pageref{sym66}\\
$\|~\|, \mathrm{d}(T)$&
a continuous seminorm on $\mathcal{S}(\mathfrak{g}(\mathbb{A}))$ and the distance from $T$ to the root hyperplanes&
page \pageref{sym67}\\
$\Phi_X(x,\overline{Y})$&
a partial Fourier transform of $f$&
page \pageref{sym68}\\
$\Gamma, \mathrm{A}_{0,T'}^{1,T}(\mathbb{R}),\delta_0^2$&
a compact subset of $\mathrm{M}_0(\mathbb{A})^1$, a compact subset of $\mathrm{A}_0^1(\mathbb{R})$ and the modulus function of $\mathrm{P}_2$&
page \pageref{sym69}\\
$\mathcal{D}, \Phi^\mathcal{D}$&
an invariant differential operator on $\mathfrak{g}(\mathbb{R})$ and its action on $\Phi$ via the Fourier transform&
page \pageref{sym70}\\
$N(f)$&
a natural number determined by a Schwartz function $f$&
page \pageref{sym71}\\
$\widehat{\mathbb{Z}}$&
the profinite completion of $\mathbb{Z}$&
page \pageref{sym72}\\
$\|~\|'$&
a seminorm on $\mathcal{S}(\mathfrak{g}(\mathbb{A}))$&
page \pageref{sym73}\\
$Z_\mu, Z_\mathfrak{n}$&
the components of $Z$ on the weight space of the weight $\mu$ and on $\mathfrak{n}$&
page \pageref{sym74}\\
$\phi_\bullet, \Psi(a_0)$&
positive Schwartz functions and the upper bounded determined by the collection of $\phi_\bullet$&
page \pageref{sym75}\\
$\beta_v$&
a bump function on $\mathbb{Q}_v$&
page \pageref{sym76}\\
$\nu, \{p_1,\dots,p_l\}$&
a nilpotent orbit and a collection of polynomials vanishing on $\nu$&
page \pageref{sym77}\\
$f_{\nu,v}^\epsilon$&
a function obtained by truncating $f$ around $\nu$&
page \pageref{sym78}\\
$\|~\|_1$&
a seminorm on $\mathcal{S}(\mathfrak{g}(\mathbb{A}))$&
page \pageref{sym79}\\
$J_\nu^T$&
a distribution supported on $\overline{\nu}$&
page \pageref{sym80}\\
$\|~\|$&
a seminorm on $\mathcal{S}(\mathfrak{g}(\mathbb{A}))$&
page \pageref{sym81}\\
$J_{\overline{\nu}}^T$&
a distribution supported on $\overline{\nu}$&
page \pageref{sym82}\\
$\mathfrak{m}_\mathrm{nil}(\mathbb{Q})_{\mathrm{M},S}$&
the set of $\mathrm{M}(\mathbb{Q}_S)$-conjugacy classes in $\mathfrak{m}_\mathrm{nil}(\mathbb{Q})$&
page \pageref{sym83}\\
$a^\mathrm{M}(S,\nu)$&
a global coefficient that depends on $S$&
page \pageref{sym84}\\
$T^\mathrm{G}$&
a distribution supported on $\mathfrak{g}_\mathrm{nil}(\mathbb{Q}_S)$&
page \pageref{sym85}\\
$\mathfrak{g}_{\mathrm{nil},d}(\mathbb{Q}_S)$&
an open subset of $\mathfrak{g}_\mathrm{nil}(\mathbb{Q}_S)$&
page \pageref{sym86}\\
$T_d^\mathrm{G}, T^{\mathrm{G},d}, T_d^{\mathrm{G},d}$&
distributions constructed from $T^\mathrm{G}$ and $\mathfrak{g}_{\mathrm{nil},d}(\mathbb{Q}_S)$&
page \pageref{sym87}\\
$\Sigma$&
a semisimple element in $\mathfrak{o}$&
page \pageref{sym88}\\
$\mathrm{P}_{1,\Sigma}^0, \mathrm{M}_{1,\Sigma}^0, \mathrm{K}_\Sigma$&
a minimal parabolic subgroup and a minimal Levi subgroup of $\mathrm{G}_\Sigma^0$ and a maximal compact subgroup of $\mathrm{G}_\Sigma^0(\mathbb{A})$&
page \pageref{sym89}\\
$T_{\Sigma,1}$&
a distinguished point in $\mathfrak{a}_1^{\mathrm{G}_\Sigma^0}$&
page \pageref{sym90}\\
$\mathcal{F}^\Sigma, \mathring{\mathcal{F}}_\mathrm{Q}, \bar{\mathcal{F}}_\mathrm{Q}$&
finite sets of parabolic subgroups&
page \pageref{sym91}\\
$\pi_\Sigma$&
a surjection from $\mathcal{F}$ to $\mathcal{F}^\Sigma$&
page \pageref{sym92}\\
$\mathcal{Y}, \mathcal{Y}_\mathrm{Q}$&
variants of a $(\mathrm{G},\mathrm{M})$-orthogonal set&
page \pageref{sym93}\\
$\Gamma'_\mathrm{Q}(,\mathcal{Y}_\mathrm{Q}), \Gamma'_\mathrm{Q}\hat{~}(,\mathcal{Y}_\mathrm{Q})$&
a gamma$'$ function and its Fourier transform&
page \pageref{sym94}\\
$\Phi_{\mathrm{Q},x}^T, v'_\mathrm{Q}, \mathcal{Y}^T(k,x)$&
a variant of the parabolic descent of $f$ along $\mathrm{Q}$ and a weight factor and a $(\mathrm{G},\mathrm{M}$-orthogonal set involved in its definition&
page \pageref{sym95}\\
$T_\Sigma$&
a truncation parameter in $\mathfrak{a}_1$&
page \pageref{sym96}\\
$j_\mathfrak{o}^T(~,f), J_{\mathrm{P},\mathfrak{o}}(x,f)$&
variants of a truncated kernel function and its summand&
page \pageref{sym97}\\
$\Xi$&
a semisimple element of $\mathfrak{m}(\mathbb{Q})$&
page \pageref{sym98}\\
$\big(v_\mathrm{P}(x,T)\big), v'_\mathrm{Q}$&
a $(\mathrm{G},\mathrm{M})$-family and a weight factor&
page \pageref{sym99}\\
$T_{\Xi,\mathrm{M}}$&
a distinguished point in $\mathfrak{a}_{\mathrm{M}_\Xi^0}^{\mathrm{G}_\Xi^0}$&
page \pageref{sym100}\\
$\Phi_{S,\mathrm{Q},x}^T$&
an $S$-local variant of the parabolic descent of $f_S$ along $\mathrm{Q}$&
page \pageref{sym101}\\
$\equiv$&
an equivalence relation on $\mathfrak{m}(\mathbb{Q})\cap\mathfrak{o}$&
page \pageref{sym102}\\
$(\mathfrak{m}(\mathbb{Q})\cap\mathfrak{o})_{\mathrm{M},S}$&
the collection of $\equiv$ equivalence classes in $\mathfrak{m}(\mathbb{Q})\cap\mathfrak{o}$&
page \pageref{sym103}\\
$S_\mathfrak{o}$&
a sufficiently large finite set of places of $\mathbb{Q}$&
page \pageref{sym104}\\
$a^\mathrm{M}(S_\mathfrak{o},X)$&
a global coefficient that depends on $S_\mathfrak{o}$&
page \pageref{sym105}\\
$\mathring{\mathcal{L}}_\Sigma(\mathrm{M}_1)$&
a finite set of Levi subgroups&
page \pageref{sym106}\\
$a^\mathrm{M}(S,X)$&
a global coefficient that depends on $S$&
page \pageref{sym107}\\
$\Gamma, \mathrm{G}_\Gamma, \mathfrak{g}(\mathbb{Q})'$&
a compact subset of $\mathfrak{g}(\mathbb{A})$, a compact subset of $\mathrm{G}(\mathbb{A})^1$ and a subset of $\mathfrak{g}(\mathbb{Q})$&
page \pageref{sym108}\\
$f_{S,\mathrm{P}}$&
the parabolic descent of $f_S$ along $\mathrm{P}$&
page \pageref{sym114}\\
$\mathrm{Ind}_{\mathrm{M},\mathrm{P}}^\mathrm{G}(X),$ $\mathrm{Ind}_\mathrm{M}^\mathrm{G}(X), X^\mathrm{G}$&
the orbit in $\mathfrak{g}(\mathbb{Q}_S)$ parabolically induced from $X$ along $\mathrm{P}$&
page \pageref{sym115}\\
$d_\mathrm{M}^\mathrm{G}(~,~), \xi, s(~,~)$&
a bivariate function on $\mathcal{L}(\mathrm{M})$, a vector in general position in $\mathfrak{a}$ and the partial map determined by $\xi$&
page \pageref{sym116}\\
$\mathrm{G}, \mathrm{M}_0, \mathbb{B}$&
the general linear group $\mathrm{GL}(n,\mathbb{Q})$ and the standard choice of the minimal Levi subgroup and the Borel subgroup&
page \pageref{sym140}\\
$\mathcal{A}_\mathrm{G}, D, \mathcal{A}_{\mathrm{G},\mathrm{reg}}, \pi_\mathrm{M}$&
the space of characteristic polynomials, the discriminant function, the nonvanishing locus of $D$ and the natural map from $\mathcal{A}_\mathrm{M}$ to $\mathcal{A}_\mathrm{G}$&
page \pageref{sym141}\\
$|D(X_v)|_v^{-1/2}\mathrm{d}X_v$&
a measure on $\mathcal{A}_\mathrm{G}(\mathbb{Q}_v)$&
page \pageref{sym142}\\
$\mathcal{A}_{\mathrm{M},\mathrm{reg}}(\mathbb{Q}_v)_\mathrm{ell}$&
the subset of $\mathcal{A}_{\mathrm{M},\mathrm{reg}}(\mathbb{Q}_v)$ consisting of the $\mathbb{Q}_v$-elliptic elements&
page \pageref{sym143}\\
$I_\mathrm{M}^\mathrm{G}(X,~), I_\mathrm{M}^\mathrm{G}(X,~)\hat{~}$&
the invariant weighted orbital integral and its Fourier transform&
page \pageref{sym109}\\
$\mathcal{I}_\mathrm{M}^\mathrm{G}(X,~)$&
the vector-valued orbital integral&
page \pageref{sym144}\\
$\mathcal{I}_\mathrm{max}^\mathrm{G}(X_v,~)$&
the maximal orbital integral&
page \pageref{sym146}\\
$\widetilde{X}_v,$ $\mathrm{M}[\widetilde{X}_v\textrm{-ell}]$&
a lift of $X_v$ in $\mathfrak{g}(\mathbb{Q}_v)$ and a standard Levi subgroup whose Lie algebra contains $\widetilde{X}_v$&
page \pageref{sym147}\\
$\mathcal{S}_0(\mathcal{A}_\mathrm{G}(\mathbb{Q}_v)),$ $\mathcal{S}_1(\mathcal{A}_\mathrm{G}(\mathbb{Q}_v))$&
the local Schwartz spaces&
page \pageref{sym148}\\
$\Lambda_0, \mathbb{I}_{\Lambda_0}, \mathcal{E}, M^\mathrm{T}, \mathrm{Tr}$&
the standard lattice in $\mathfrak{g}(\mathbb{Q}_p)$, its characteristic function, the standard Gaussian on $\mathfrak{g}(\mathbb{R})$, the transpose of the matrix $M$ and the matrix trace&
page \pageref{sym149}\\
$\phi_{0,v}, \phi_{1,v}$&
the local basic functions&
page \pageref{sym150}\\
$Z_v, U_{\varphi_v,Z_v}, \Gamma_\mathrm{G}^\mathrm{G}(~,\nu)$&
a singular point, an open neighborhood of $Z_v$ and the Shalika germ at a nilpotent orbit $\nu$&
page \pageref{sym151}\\
$\mathcal{S}_0(\mathcal{A}_\mathrm{G}(\mathbb{A})),$ $\mathcal{S}_1(\mathcal{A}_\mathrm{G}(\mathbb{A}))$&
the global Schwartz spaces&
page \pageref{sym152}\\
$\mathcal{H}_v$&
the local Harish-Chandra transform&
page \pageref{sym153}\\
$\mathcal{I}_\mathrm{max}^\mathrm{G}(X_v,~)\hat{~}$&
the function representing the Fourier transform of the maximal orbital integral&
page \pageref{sym154}\\
$\mathcal{K}_v(~,~)$&
the local Harish-Chandra kernel&
page \pageref{sym154}\\
$\hat{i}_\mathrm{M}^\mathrm{G}(~,~)$&
Waldspurger's notation for the local Harish-Chandra kernel&
page \pageref{sym156}\\
$\delta_{X_v}, (\delta_{X_v,i}), \widetilde{\delta}_{X_v,i}$&
the Dirac distribution at $X_v$, a $\delta$-sequence converging to $\delta_{X_v}$ and a lift of $\delta_{X_v,i}$ on $\mathfrak{g}(\mathbb{Q}_v)$&
page \pageref{sym157}\\
$\mathcal{H}_S, \mathcal{H}$&
the $S$-local and global Harish-Chandra transforms&
page \pageref{sym158}\\
$Z, \mathfrak{o}_Z, Z', \varphi_S(Z')$&
a singular point in $\mathcal{A}_\mathrm{G}(\mathbb{Q})$, its fiber in $\mathfrak{g}(\mathbb{Q})$, a $\equiv$ equivalence class in $\big(\mathfrak{m}(\mathbb{Q})\cap\mathfrak{o}_Z\big)_{\mathrm{M},S}$ and the value determined by the nearby regular values $\varphi_S(X)$&
page \pageref{sym159}\\
$a(X)$&
the global coefficient for a regular $X$&
page \pageref{sym162}\\
$\Lambda_{N,w}, \mathbb{I}_{\Lambda_{N,w}}$&
a lattice in $\mathfrak{g}(\mathbb{Q}_w)$ for a finite place $w$ and its characteristic function&
page \pageref{sym163}\\
$\pi_\infty$&
the natural projection from $\mathfrak{g}(\mathbb{R})$ to $\mathcal{A}_\mathrm{G}(\mathbb{R})$&
page \pageref{sym164}\\
$\mathrm{T}_{X_v}$&
the centralizer of $X_v$ which is a torus&
page \pageref{sym165}\\
$\Lambda_{N,\infty}, \Lambda_{N,w},$ $\mathbb{I}_{\widetilde{X}_\infty+\Lambda_{N,w}}$&
a lattice in $\mathfrak{g}(\mathbb{R}$, a lattice in $\mathfrak{g}(\mathbb{Q}_w)$ for a finite place $w$, and the characteristic function of the translation of $\Lambda_{N,w}$ by $\widetilde{X}_\infty$&
page \pageref{sym166}\\
$\mathrm{T}_Y, \widehat{\mathrm{T}}_Y, F, \Gamma$&
the centralizer of $Y$ which is a maximal torus, its group of algebraic characters, the splitting field of $\mathrm{T}_Y$ and its Galois group&
page \pageref{sym167}\\
$\mathcal{D}$&
an invariant distribution on a $p$-adic reductive Lie algebra&
page \pageref{sym168}\\
$\mathbb{E}^n, \mathbb{S}(\mathbb{E}^n)$&
the $n$-dimensional Euclidean space and its scissors group&
page \pageref{sym117}\\
$\Phi, V_i, \mathbf{r}, \mathcal{R}ig(\Phi)$&
a strict flag in $\mathbb{E}^n$, the $i$th subspace of $\Phi$, a rigging of $\Phi$ and the collection of equivalence classes of riggings of $\Phi$&
page \pageref{sym118}\\
$\mathcal{B},\mathcal{B}^\mathbf{r}, \mathrm{sign}(\mathbf{r})$&
two ordered bases of $\mathbb{E}^n$ and the sign of $\mathbf{r}$&
page \pageref{sym119}\\
$\Phi^\mathbf{r}, \partial_{\Phi^\mathbf{r}}(P), r_i^\mathrm{min}(~)$&
a rigged flag, the $\Phi^\mathbf{r}$-boundary of the polytope $P$ and an auxiliary function&
page \pageref{sym120}\\
$\mathrm{Had}_\Phi(P)$&
the Hadwiger invariant of the polytope $P$&
page \pageref{sym121}\\
$(H_\Phi), \mathcal{F}_\mathrm{s}(\Phi,i)$&
a collection of real numbers and a collection of strict flags in $\mathbb{E}^n$&
page \pageref{sym122}\\
$\mathbb{S}(\mathfrak{a}_\mathrm{M}^\mathrm{G}), \mathcal{F}_\mathrm{S}(\mathfrak{a}_\mathrm{M}^\mathrm{G})$&
the scissors group of positive $(\mathrm{G},\mathrm{M})$-orthogonal sets and a collection of strict flags in $\mathfrak{a}_\mathrm{M}^\mathrm{G}$&
page \pageref{sym123}\\
$\mathbb{J}_\mathrm{M}^\mathrm{G}(X,~)$&
the scissors-congruence-valued orbital integral&
page \pageref{sym124}\\
$\mathbf{d}\mathfrak{a}_\mathrm{L}^{\mathrm{L}^{l-1}}$&
a differential form on $\mathfrak{a}_\mathrm{M}^{\mathrm{L}^{l-1}}$&
page \pageref{sym125}\\
$\mathbb{S}(\mathfrak{a}_{\mathcal{L}(\mathrm{M})}^\mathrm{G})$&
the total scissors ring of $\mathfrak{a}_\mathrm{M}^\mathrm{G}$&
page \pageref{sym126}\\
$\boxtimes, j, j^*$&
a bilinear product on the total scissors ring and two auxiliary maps&
page \pageref{sym127}\\
$\big(\mathbb{S}(\mathfrak{a}_{\mathcal{L}(\mathrm{M})}^\mathrm{G})\big)_n,$ $\big(\mathbb{S}(\mathfrak{a}_{\mathcal{L}(\mathrm{M})}^\mathrm{G})\big)_\mathrm{L}$&
graded components of the total scissors ring with respect to two different gradings&
page \pageref{sym128}\\
$q_\mathrm{M}^\mathrm{L}$&
a homomorphism from $\mathbb{S}(\mathfrak{a}_{\mathcal{L}(\mathrm{M})}^\mathrm{G})$ to $\mathbb{S}(\mathfrak{a}_{\mathcal{L}(\mathrm{L})}^\mathrm{G})$&
page \pageref{sym129}\\
$\mathcal{J}_\mathrm{M}^\mathrm{G}(X,~)$&
the total orbital integrohedron&
page \pageref{sym130}\\
$v, p, |~|$&
an odd rational prime and the $p$-adic absolute value&
page \pageref{sym169}\\
$\mathrm{G}, \mathrm{M}, \mathrm{B}$&
the group $\mathrm{GL}(2)$ and the standard choice of the minimal Levi subgroup and the Borel subgroup&
page \pageref{sym170}\\
$\mathfrak{g}, \mathfrak{g}', \mathfrak{m}', \mathfrak{z}$&
the Lie algebras $\mathrm{gl}(2)$, $\mathrm{sl}(2)$, $\mathfrak{m}\cap\mathfrak{g}'$ and the center of $\mathfrak{g}$&
page \pageref{sym171}\\
$\mathcal{A}, \mathbb{L}, D, \mathcal{A}_\mathrm{reg}, \mathbb{Q}_p^\times$&
the affine quotient of $\mathfrak{g}'$ by $\mathrm{G}$, the affine line, the discriminant function on $\mathcal{A}$,  the nonvanishing locus of $D$ and the group of units of $\mathbb{Q}_p$&
page \pageref{sym172}\\
$\eta, \mathfrak{g}'_\eta, \mathrm{d}_\eta X$&
a point of $\mathcal{A}_\mathrm{reg}$, the fiber of $\eta$ in $\mathfrak{g}'$ and a differential form on $\mathfrak{g}'_\eta$&
page \pageref{sym173}\\
$\mathrm{val}(\eta)$&
the $p$-adic valuation of $\eta$&
page \pageref{sym174}\\
$\widetilde{\eta}$&
a lift of a split $\eta$ in $\mathfrak{g}'_\mathrm{reg.ss}(\mathbb{Q}_p)$&
page \pageref{sym175}\\
$J_\mathrm{M}^\mathrm{G,reg}(X,~)$&
a regularized version of the weighted orbital integral $J_\mathrm{M}^\mathrm{G}(X,~)$ for the Lie algebra $\mathrm{gl}(2)$&
page \pageref{sym176}\\
$J_\mathrm{M}^{\mathrm{G},b}(X,~)$&
a variant of the weighted orbital integral $J_\mathrm{M}^\mathrm{G}(X,~)$ constructed by Chaudouard&
page \pageref{sym177}
\end{longtable}

\end{document}